\documentclass[11pt]{amsbook}
\usepackage{amsfonts}
\usepackage{graphicx}
\usepackage{enumerate}
\title{The parameter capture map for $V_{3}$}
\author{Mary Rees}
\address{Department of Mathematical Sciences, University of Liverpool, Mathematics Building, Peach St., Liverpool L69 7ZL, U.K.}
\email{maryrees@liv.ac.uk}
\subjclass[2000]{37F10,37B10}

\def\Box{\hbox{$\sqcap \unskip \kern -6.5pt 
\sqcup$}}
\def \Amalg{\mathbin{\raise .5pt
	\hbox{$\scriptstyle \amalg$}}}
\def \cal{\mathcal}
\def \Cbar{\overline{\mathbb C}}	
	    
\numberwithin{equation}{subsection}
\numberwithin{section}{chapter}
\numberwithin{subsection}{section}
\newtheorem{theorem}[section]{Theorem}
\newtheorem{subsectiontheorem}[subsection]{Theorem}

\newtheorem{lemma}[subsection]{Lemma}

\newtheorem{corollary}[subsection]{Corollary}

\begin{document}
\frontmatter
{  

\begin{abstract} This is a study of the Wittner capture construction for critically finite quadratic rational maps for which one critical point is periodic, and the second critical point is in the backward orbit of the first. This construction gives a way of describing rational maps up to topological conjugacy. It is known that representations as Wittner captures are not unique. We show that, in a certain parameter space which we call $V_3$, the set of maps with exactly $2^r$ representations as a Wittner capture is of density bounded from $0$ for each $r\geq 0$, and for each fixed preperiod of the second critical point.\end{abstract}

  \maketitle
  \tableofcontents
}
\mainmatter
{
\chapter{Introduction}

\section{The setting}\label{1.1}

This paper is concerned with the parameter space $V_{3}$ of quadratic rational maps:
$$h_{a}:z\mapsto \frac{(z-a)(z-1)}{z^{2}},\ \ a\in \mathbb C,\ a\neq 
0.$$
It is a sequel to \cite{R5} and addresses some questions at the end of that paper. Roughly speaking, this is  a series of  questions about how many critically finite rational maps in $V_3$ can be represented by a construction called {\em{Wittner capture}}  \cite{W} and in how many ways. More precisely, we shall consider Wittner captures by the {\em{aeroplane polynomial}}. More formally, the questions  concern the {\em{parameter capture map for $V_{3}$ and the aeroplane polynomial}}.  The domain of the parameter capture map is a union of finite sets parametrised by $m$, and for each such restricted finite set domain,  there is a natural restriction of the range to another finite set. We denote the restricted parameter capture map by $\Phi _{m}$. The domain and range of $\Phi_{m}$ both have size of the order  of  $2^{m}$.   We are interested in the image size, and point inverse image sizes of $\Phi _{m}$, for varying $m$.

This is, therefore, a consideration of a problem in a very specific setting. But it is also a particular case of a very basic question in dynamics: how to find simple representations of dynamical systems in a parameter space, and, at the same time, to determine the number of duplicate representations. In complex dynamics, mating, and the more easily analysed Wittner captures, are among the more popular representations. This paper is perhaps the first serious analysis of the number of duplicate representations arising. 

\section{The parameter space, its maps and hyperbolic components}\label{1.2}
The map $h_{a}$   has critical points $0$ and $c_{2}(a)=\frac{2a}{a+1}$, with corresponding critical values $\infty $ and $v_{2}(a)=-\frac{(a-1)^{2}}{4a}$. The point $0$ is of period $3$ under $h_{a}$, with orbit $0\mapsto \infty \mapsto 1\mapsto 0$. Every quadratic rational map with a critical point $c_{1}$ of period $3$ is represented in $V_{3}$ up to conjugation by a M\"obius transformation. The representative is unique if the M\"obius transformation is chosen to map $c_{1}$ to $0$. There are three polynomials in $h_{a}$ up to M\"obius conjugacy: the {\em{rabbit}}, the {\em{antirabbit}} and the {\em{aeroplane}}, corresponding to parameter values $a=a_{0}$ with ${\rm{Im}}(a_{0})>0$ \cite{R5}, $\overline{a_{0}}$ and $a_{1}$, which is real and $<0$. Therefore, with abuse of notation, we shall sometimes refer to $h_{a}$, for $a=a_{0}$, $\overline{a_{0}}$ or $a_{1}$, as ``polynomials''. There are also two {\em{type II}} critically finite maps in $V_{3}$ given by parameter values $a=\pm 1$. We see that $c_{2}(1)=1$ and $c_{2}(-1)=\infty $, so that, for both these parameter values, the critical points lie in the same period three orbit. 

Each of these maps $h_{a}$, for $a=a_{0}$, $\overline{a_{0}}$, $a_{1}$ and $\pm 1$, is a {\em{hyperbolic}} rational map --- an iterate of the map is expanding on the Julia set in the spherical metric --- since the simple equivalent condition for hyperbolicity is trivially satisfied. It is clear that each critical point is in the Fatou set, since each critical point is in a superattracting cycle. Each of these maps $h_{a}$ therefore lies in an open set $H_{a}$ of rational hyperbolic maps  in $V_{3}$, which we call the hyperbolic component, such that each  maps $h_{a'}$ in $h_{a}$ is conjugate to $h_{a}$ on a neighbourhood of the Julia set $J(h_{a'})$ of $h_{a'}$. 

A branched covering of the Riemann sphere is {\em{critically finite}} if the forward orbit of any critical point is finite. Any hyperbolic component contains at most one critically finite map -- always exactly one in $V_{3}$ -- and so the sets $H_{a}$ for $a=a_{0}$, $\overline{a_{0}}$, $a_{1}$ and $\pm 1$ are disjoint.   Identifying $h_{a}$ with $a$, we can regard $H_{a}$ as an open subset of the complex plane. It is simply connected in each case. In fact, the uniformising map in each case is a natural parametrisation of the dynamical variation within the hyperbolic component (\cite{R4}, for example). Also, $H_{a}$ is symmetric about the real axis whenever  $a$ is real. The closures  of the hyperbolic components $H_{\pm 1}$ meet in three points: $0$ (which is in $\mathbb C$, but excluded from $V_{3}$) and points which we shall call $x$ (with ${\rm{Im}}(x)>0$) and $\overline{x}$, which are uniquely determined by the property that $h_{x}$ has a parabolic fixed point with multiplier $e^{2\pi i/3}$. (Of course, it follows immediately that $h_{\overline{x}}$ has a parabolic fixed point with multiplier $e^{-2\pi i/3}$.) 
The closure of $H_{a_{1}}$ also includes $x$ and $\overline{x}$. These two points are accessible from $H_{\pm 1}$ \cite{R4} and $H_{a_{1}}$, and $0$ is accessible from $H_{\pm 1}$ along the real axis. Since $H_{a_{1}}$ is unbounded --- it includes all points in the negative real axis in $(-\infty ,a_{1}]$, for example --- it lies in the unbounded component  of the complement of $H_{1}\cup H_{-1}\cup H_{a_{1}}\cup \{ x,\overline{x}\} $. We call this unbounded complementary component $V_{3}(a_{1},+)$ \cite{R5}, because it meets the positive real axis. Its closure meets $\overline{H_{1}}\setminus \{ x,\overline{x}\} $ and $\overline{H_{a_{1}}}\setminus \{ x,\overline{x}\} $, but not $\overline{H_{-1}}\setminus \{ x,\overline{x}\} $. There are three other complementary components. One is bounded by $H_{-1}\cup H_{a_{1}}\cup \{ x,\overline{x}\} $ and does not meet $\overline{H_{1}}\setminus \{ x,\overline{x}\} $. This complementary component is called $V_{3}(a_{1},-)$ (and meets the negative real axis). The other two complementary components are bounded by $H_{a_{0}}\cup H_{\overline{a_{0}}}\cup \{ x,\overline{x}\} $. One of these, $V_{3}(a_{0})$, is in the upper half-plane, and contains $H_{a_{0}}$. The other, $V_{3}(\overline{a_{0}})$, is in the lower half-plane, and contains $H_{\overline{a_{0}}}$. This is shown in Figure 1.

\begin{figure}
\centering{\includegraphics[width=6cm]{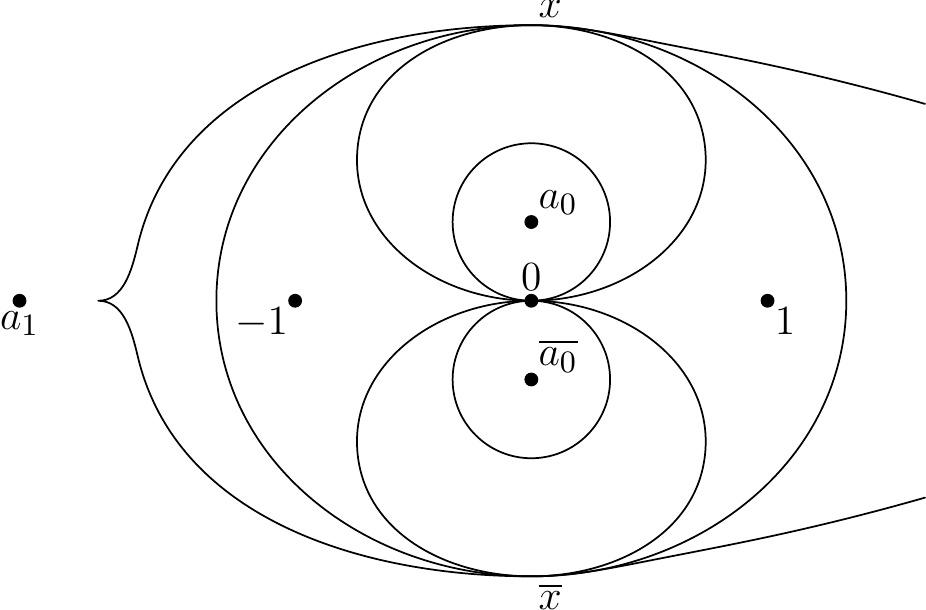}}
\caption{The parameter space $V_{3}$}
\end{figure}

\section{Counting hyperbolic components}\label{1.3}
Each of the sets $V_{3}(a_{0})$, $V_{3}(\overline{a_{0}})$, $V_{3}(a_{1},-)$ and $V_{3}(a_{1},+)$ contains infinitely many hyperbolic components of each of types III and IV, each with a critically finite centre. For maps in a {\em{type III}} hyperbolic component (in the space of quadratic rational maps) there is a single periodic cycle of Fatou components, containing exactly one of the critical values, and the other critical value is in the full orbit of this periodic cycle. For hyperbolic components intersecting $V_{3}$, the cycle is of period $3$, and the second critical value  is in a Fatou component of preperiod $m>0$. For the critically finite centre $h_{a}$, the second critical value $v_{2}(a)$ has preperiod $m$. For {\em{type IV}} hyperbolic components, there are two distinct cycles of periodic Fatou components, each containing a critical point and value. As in \cite{R5}, we are concerned in this paper only with type III hyperbolic components, and, in particular, with their centres. For each fixed $m$, the set of type III centres for preperiod $m$ is finite, and has cardinality
$$\frac{23}{21}2^{m+1}+O(1).$$

Moreover, we can identify the number of type III components of preperiod $m$ in each of the sets $V_{3}(a_{0})$, $V_{3}(\overline{a_{0}})$, $V_{3}(a_{1},-)$ and $V_{3}(a_{1},+)$. We call denote the set of  critically finite centres  of type III components by $P_{3,m}$, and the intersections with each of the above sets of $V_3$ by, respectively, $P_{3,m}(a_0)$, $P_{3,m}(\overline{a_0})$, $P_{3,m}(a_1,+)$ and $P_{3,m}(a_1,-)$. As a consequence of the main result Theorem 2.10 of \cite{R5}, supplemented by the counting in section 3.4 of \cite{R5},  the numbers in each of these sets are, respectively,
$$\frac{3}{7}2^{m+1}+O(1),\ \ \frac{3}{7}2^{m+1}+O(1),\ \ \frac{1}{7}2^{m+1}+O(1),\ \ \frac{2}{21}2^{m+1}+O(1).$$
In fact, we only need the statement of the ``easy'' parts 1 and 2 of Theorem 2.10 of \cite{R5}, together with the counting, to deduce this. Counting results of this type, although not this particular one, are proved in \cite{K-R}.

\section{Thurston equivalence}\label{1.4}
For a branched covering $f$ of $\overline{\mathbb C}$,  we define
$$X(f)=\{ f^n(c)\mid c {\rm{\ critical,\ }}n>0\} .$$
We say that $f$ is {\em{critically finite}} if $X(f)$ is  finite. The convention (not universal) in this paper  is to number the critical values. Adopting this convention, two critically 
finite maps $f_{0}$ and $f_{1}$ with numbered critical values are   \mbox{{\em{(Thurston) equivalent}} }   if there is a homotopy $f_{t}$ from 
$f_{0}$ to $f_{1}$ such that $\# (X(f_{t}))$ is constant in $t$, so that the finite set $X(f_{t})$ varies isotopically with $t$, and the 
isotopy between $X(f_{0})$ and $X(f_{1})$ preserves the numbering of critical values. If $\# 
(X(f_{t}))\geq 3$, this is equivalent to the existence of 
homeomorphisms $\varphi $ and  $\psi :\overline{\mathbb C}\to 
\overline{\mathbb C}$ with $\varphi $ and $\psi $ isotopic via an 
isotopy mapping $X(f_{0})$ to $X(f_{1})$, and preserving the 
numbering of critical values, and such 
that 
$$\varphi \circ f_{0}\circ \psi ^{-1} =f_{1}.$$
We shall write $\simeq $ for the equivalence relation of Thurston equivalence, and will sometimes write simply ``equivalence'' or ``equivalent'' when it is clear that Thurston equivalence is meant. 

Thurston's theorem \cite{D-H2,R3} gives a necessary and sufficient condition for a critically finite rational map $f$ to be equivalent to a rational map $g$, which is usually unique up to M\"obius conjugation. This is certainly true if the forward orbit of every critical point contains a periodic critical point, which is the only case which concerns us here. (In fact, if it is not true, then $\# (X(f))=4$, every critical point  is strictly preperiodic, and maps forward to an orbit of period one or two.)

\section{Capture paths and capture maps}\label{1.5}
A capture path is a particular type of path in the dynamical plane of a critically periodic quadratic polynomial, or any rational map which is M\"obius conjugate to a critically periodic quadratic polynomial, that is, one for which the finite critical point is periodic. In particular, we can define capture paths for the maps $h_a$ for $a=a_0$, $\overline{a_0}$ and $a_1$. We shall also use capture paths for branched coverings which are Thurston equivalent to quadratic polynomials.

So now let $a=a_0$, $\overline{a_0}$ or $a_1$. Recalling that $0$ is of period $3$ under $h_a$ (in fact this is the case for any $h_a\in V_3$), we define
$$Z_m(h_a)=h_a^{-m}(\{ 0,1,\infty \} )=h_a^{-m}(\{ 0,h_a(0),h_a^2(0)\} )$$
and
$$Z(h_a)=\cup _{m\geq 0}Z_m(h_a)=\cup _{i\geq 0}h_a^{-i}(\{ 0,1,\infty \} ).$$
 A {\em{capture path for $h_a$}}  is a path $\beta :[0,1]\to \overline{\mathbb C}$ from the fixed critical point $\beta (0)=c_{2}(a)=\frac{2a}{a+1}$ to a point $\beta (1)=x\in Z(h_a)$. To be a capture path, the path has to cross the Julia set $J(h_{a})$ just once, into the Fatou component containing $x$, and this Fatou component must be adjacent to the ray of entry. Two paths $\beta _{1}$ and $\beta _{1}$ are {\em{equivalent }}if they are homotopic in $\overline{\mathbb C}\setminus \{ h_a^n(x)\mid n\geq 0\} $ via a homotopy fixing endpoints. For $a=a_0$ or $\overline{a_0}$, a capture path is uniquely determined up to equivalence, by its endpoint. This is because the Julia set of $h_a$ in these cases is homeomorphic to the Julia set of the rabbit or antirabbit polynomial respectively, and the forward orbit of $x$ is in just one component of the complement in the filled Julia set, of the closure of the immediate basin of attraction. For $a=a_1$, the forward orbit of $x$ sometimes intersects two components of the complement, in the filled Julia set, of the closure of the  immediate basin  of attraction of $x$. But the set of $x$ of preperiod $m$ for which there are two such components, is of density tending to $0$ as $m$ tends to infinity. We shall be more precise about this below, in \ref{1.7} 

Any capture path is an arc up to equivalence and from now on we assume that capture paths are arcs. If $\beta :[0,1]\to \overline{\mathbb C}$ is an arc then we can define a homeomorphism $\sigma _\beta $ as follows. Take a suitably small disc neighbourhood $U$ of $\beta $. Define $\sigma _\beta $ to be the identity outside $U$, and to map $\beta (0)$ to $\beta (1)$. See Figure 2. It is possible to define $\sigma _{\beta }$ for any continuous path $\beta $, by writing $\beta $ as a union of arcs, up to homotopy. We shall need to employ this later when $\beta $ is a closed loop.

\begin{figure}
\centering{\includegraphics[width=4cm]{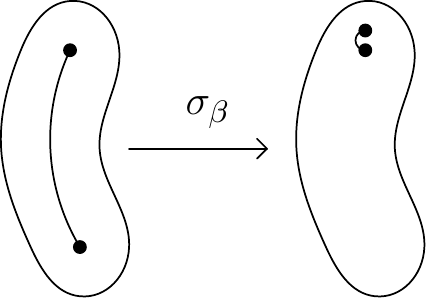}}
\caption{$\sigma _{\beta }$}
\end{figure}

A capture path $\beta $ is {\em{type II}} if the second endpoint is periodic, and {\em{type III}} if the second endpoint is preperiodic. If $\beta $ is a type III capture path for $h_a$, for $a=a_0$, $\overline{a_0}$ or $a_1$, then $\sigma _{\beta }\circ h_a$ is a critically finite branched covering. If $U$ is sufficently small to be disjoint from $\{ h^n(x)\mid n>0\} $ then $\sigma _{\beta }\circ h_a$ is uniquely determined up to Thurston equivalence by the equivalence (homotopy) class of $\beta $. Hence, if $\beta $ is a capture path, $\sigma _{\beta }\circ h_a$ is determined up to equivalence by the endpoint $\beta (1)$, except for a set of endpoints of density tending to $0$ as the preperiod tends to infinity. If $\beta $ is a capture path then we call $\sigma _{\beta }\circ h_a$ a {\em{capture}} or a {\em{Wittner capture}}.  This definition of capture was used in \cite{R1,R2,R3,R5,R6},  following the introduction of captures by Ben Wittner in his thesis \cite{W}, but some authors use ``capture'' more generally for what we call ``type III''. We refer to a hyperbolic component whose centre is a capture in our stricter sense, up to Thurston equivalence, as a {\em{capture component}}. 

If $\beta $ is a type II capture path for $h_a$ then we can define a critically finite branched covering using $\beta $ in a slightly different way (see \cite{R1,R2,R3,R5}). Let $\zeta $ be the uniquely determined path such that $h_a\circ \zeta =\beta $ and such that $\zeta (1)$ is periodic under $h_a$. Then $\sigma _{\zeta }^{-1}\circ \sigma _{\beta }\circ h_a$ is a critically finite branched covering of type II, which is, again, uniquely determined up to Thurston equivalence by the equivalence class of $\beta $.

\section{Laminations and Lamination maps}\label{1.6}

Invariant laminations were introduced by Thurston \cite{T} to 
describe 
the dynamics of polynomials with locally connected Julia sets. The {\em{leaves}} of a lamination $L$  
are straight line segments in $\{ z:\vert z\vert \leq 1\} $. 
Invariance of a lamination means that if 
there is a leaf with endpoints $z_{1}$ and $z_{2}$, then there are also  
 a leaf with endpoints $z_{1}^{2}$ and $z_{2}^{2}$, a leaf with endpoints 
$-z_{1}$ and $-z_{2}$, and a leaf with endpoints $w_{1}$ and $w_{2}$, where 
$w_{1}^{2}=z_{1}$ and $w_{2}^{2}=z_{2}$.  A leaf 
with endpoints $e^{2\pi ia_{1}}$ and $e^{2\pi ia_{2}}$, for $0\leq 
a_{1}<a_{2}<1$,
is then said to have {\em{length}} 
${\rm{min}}(a_{2}-a_{1},a_{1}+1-a_{2})$. {\em{Gaps}} of the 
lamination are components of $\{ z:\vert z\vert <1\} \setminus (\cup 
L)$. 
If the longest leaf of $L$ has length $<{1\over 2}$  then there are 
exactly two with the same image, which is called the {\em{minor 
leaf}}.  A lamination $L$ is {\em{clean}} if finite-sided gaps of $L$ are never 
adjacent. Minor leaves of clean laminations are either equal or have  disjoint interiors. We are only interested, here, in minor leaves which have endpoints which are periodic under $z\mapsto z^2$. If one endpoint of a minor leaf of a clean lamination is periodic, then the other is too, and of the same period. Any two such minor leaves have distinct endpoints. There is a natural partial ordering on minor leaves. We say that $\mu <\mu '$ if $\mu $ separates $\mu '$ from $0$ in the unit disc. For any minor leaf $\mu $ there is a unique minimal minor leaf $\nu =\nu (\mu )$ such that $\nu \leq \mu '$ whenever $\mu '\leq \mu $.

 Each point $e^{2\pi it_{1}}\neq 1$, with $t_{1}$ an odd denominator rational, is an endpoint of the  minor leaf $\mu _{t_1}$ of a unique clean lamination, which we call $L_{t_{1}}$ -- and also $L_{t_2}$, if $e^{2\pi it_2}$ is the other endpoint of $\mu _{t_1}$ (in which case $\mu _{t_1}=\mu _{t_2}$).  We 
can define a {\em{lamination map}} $s_{t_{1}}=s_{t_{2}}$ which maps $L_{t_{1}}$ to $L_{t_{1}}$, 
and such that $s_{t_1}(z)=z^{2}$ for $\vert z\vert \geq 1$, and gaps are mapped to gaps. The gap of $L_{t_{1}}$  containing $0$ has infinitely many sides and is periodic under $s_{t_{1}}$, of period $n$, and is mapped with degree two onto its image by $s$, but the rest of the periodic cycle maps homeomorphically. We can therefore choose $s_{t_{1}}$ so that $0$ is a degree two critical point, of period $n$, and hence is a degree two critically periodic branched covering. Then $s_{t_{1}}$ is Thurston equivalent to a unique critically finite quadratic polynomial, and conversely.

For technical reasons, in \cite{R5}, it turned out to be easier to work with the lamination map which is Thurston equivalent to $h_{a}$, rather than with $h_{a}$ itself. This will be true here, too. The corresponding lamination map is $s_{q}$ for $q=\frac{1}{7}$ or $\frac{6}{7}$ or $\frac{3}{7}$ respectively. Capture paths transfer to this setting. A {\em{capture path for $s_{q}$}} is a path from the fixed critical point $\infty $ of $s_{q}$ to a point $x$ in the full orbit $Z(s_q)=\cup _{n>0}s_{q}^{-n}(\{ 0,s_q(0),s_q^2(0)\} )$ of the critical point $0$, which crosses the unit circle just once, into the infinite-sided gap of $L_{q}$ which contains $x$. The associated type III and type II capture maps are defined in exact analogy to those in \ref{1.5}.

\section{Counting capture maps and components}\label{1.7}

The following information is given in \cite{R5}, with a little more detail in places, but for ease of reference, it is given here. 
A theorem of Tan Lei \cite{TL}, of which a more general version is proved in \cite{R1,R2}, gives a necessary and sufficient criterion for  a capture $\sigma _{\gamma }\circ s_q$ (for any odd denominator rational $q$) to be  Thurston equivalent to a rational map. If $q=\frac{r}{7}$ for $1\leq r\leq 6$, then this rational map is  automatically in $V_3$. Actually, Tan Lei's theorem deals with matings rather than captures, but the result for captures can be proved by the same methods. The criterion is easily described using minor leaves, and therefore  we use the captures $\sigma _{\gamma }\circ s_q$. The critically finite branched covering $\sigma _{\gamma  }\circ s_q$ is Thurston equivalent to a rational map if and only if $\gamma (1)$ is in the larger  region of the unit disc  bounded by the minimal minor leaf $\nu (\mu _q)$ with $\nu (\mu _q)\leq \mu _q$ . The minimal minor leaves $\nu (\mu _q)$ for $q=\frac{1}{7}$, $\frac{6}{7}$ and $\frac{3}{7}$ are $\mu _{1/7}$ (which has endpoints $e^{2\pi i(1/7)}$ and $e^{2\pi i(2/7)}$), $\mu _{6/7}$ (which has endpoints $e^{2\pi i(5/7)}$ and $e^{2\pi i(6/7)}$) and $\mu _{1/3}$ (which has endpoints $e^{\pm 2\pi i(1/3)}$) respectively. 

We write ${\cal{Z}}(q)$ for the set of equivalence classes of capture paths $\gamma $ for $s_q$ such that $\sigma _{\gamma }\circ s_q$ is equivalent to a rational map, and ${\cal{Z}}_m(q)$ for the subset for which the second endpoint of the path is of preperiod $m$. The number $\# ({\cal{Z}}_m(q))$ is closely related to the  the number of points of  $Z_m(s_q)=s_q^{-m}(\{ s_q^i(0)\mid i\geq 0\} )$ in the larger region of the unit disc bounded by $\nu (\mu _q)$. The numbers  of such points for preperiod $m$ are, respectively,
$$\dfrac{3}{7}2^{m+1}+O(1),\ \ \dfrac{3}{7}2^{m+1}+O(1),\ \ \dfrac{1}{3}2^{m+1}+O(1).$$
For $q=\frac{1}{7}$ and $\frac{6}{7}$, it is not hard to show that capture paths $\gamma _1$ and $\gamma _2$ for $s_q$ are equivalent if and only if the  captures $\sigma _{\gamma _1}\circ s_q$ and $\sigma _{\gamma _2}\circ s_q$ are Thurston equivalent. Also, it was shown in \cite{R5} that the equivalent rational maps to captures $\sigma _{\gamma }\circ s_{1/7}$ are precisely the centres of hyperbolic components in $V_3(a_0)$, and similarly for captures $\sigma _{\gamma }\circ s_{6/7}$ and centres of hyperbolic components in $V_3(\overline{a_0})$.  (It was pointed out in \cite{R5} that this can be shown by elementary methods.)

For $q=\frac{3}{7}$, the number of equivalence classes of capture paths is slightly more than the number of endpoints, because , for some choices of endpoint $x$, there are two equivalence classes of capture paths with endpoint $x$. But even for $q=\frac{3}{7}$, there is just one equivalence class of capture path with endpoint $x$, unless the closure of the gap of $L_{3/7}$ which contains $x$ intersects both
$\{ e^{2\pi it}: t\in (\frac{2}{7},\frac{1}{3})\} $ and $\{ e^{2\pi it}: t\in (\frac{2}{3},\frac{6}{7})\} $. In that case, there are exactly two equivalence classes of capture paths for $s_{3/7}$ ending at $x$: those which intersect $\{ e^{2\pi it}\mid t\in (\frac{2}{7},\frac{1}{3}\} $ and those which intersect  $\{ e^{2\pi it}: t\in (\frac{2}{3},\frac{6}{7})\} $. There may then be either one or two Thurston equivalence classes of captures $\sigma _{\gamma }\circ s_{3/7}$ for $\gamma $ with endpoint at $x$. Examples of both cases are given in \cite{R5}. But at any rate, the number of maps in $V_3$ which are equivalent to capture maps  $\sigma _{\gamma }\circ s_{3/7}$ is 
$$\dfrac{1}{3}2^{m+1}(1+o(1)).$$

However, there are Thurston equivalences between some captures $\sigma _{\gamma }\circ s_{3/7}$ and  captures $\sigma _{\zeta }\circ s_q$ for $q=\frac{1}{7}$ or $\frac{6}{7}$. These arise from the so-called {\em{Wittner flip construction}}. 
If $\gamma $ crosses the unit circle at $e^{2\pi it}$ for $t\in (\frac{1}{7},\frac{2}{7})$, then $\sigma _{\gamma }\circ s_{3/7}$ is equivalent to a capture $\sigma _{\zeta }\circ s_{6/7}$, and similarly, if $t\in (\frac{5}{7},\frac{2}{7})$ then $\sigma _{\gamma }\circ s_{3/7}$ is equivalent to a capture $\sigma _{\zeta }\circ s_{1/7}$. This very interesting construction is not used in the results of this paper, so we will not discuss it further.

 Capture paths $\gamma _1$ and $\gamma _2$ which cross the unit circle in $\{ e^{2\pi it}: t\in (-\frac{1}{7},\frac{1}{7})\} $ are equivalent if and only if the   captures $\sigma _{\gamma _1}\circ s_q$ and $\sigma _{\gamma _2}\circ s_p$  are Thurston equivalent. It was shown in \cite{R5} that the rational maps equivalent to these captures are precisely all the centres of the hyperbolic components in $V_3(a_1,-)$. Ideas of Timorin (private communication) indicate that there is a simple direct way of proving this.

Thus, using Thurston equivalence, there are natural one-to-one correspondences between the hyperbolic components in the regions $V_3(a_0)$, $V_3(\overline{a_0})$ and $V_3(a_1,-)$ and equivalence classes of captures by $s_q$ for $q=\frac{1}{7}$, $\frac{6}{7}$ and $\frac{3}{7}$ respectively,  where the capture paths are taken to have endpoints in $\cup _{n>0}s_{q}^{-n}(0)$, and in, respectively:
\begin{itemize}
\item  the larger region of the unit disc bounded by the  leaf $\mu _{1/7}$;  
\item the larger region of the unit disc bounded by the  leaf $\mu _{6/7}$;
\item  the smaller region of the unit disc bounded by the leaf with endpoints $e^{\pm 2\pi i(1/7)}$.
\end{itemize}

But none of this is true for $V_3(a_1,+)$. The only captures which can be equivalent to centres of hyperbolic components in this region are those of the form $\sigma _{\gamma }\circ s_{3/7}$ where $\gamma $ crosses the unit circle at $e^{2\pi it}$ for  $t\in (\frac{2}{7},\frac{1}{3})\cup (\frac{2}{3},\frac{6}{7})$, since as we have seen all the other captures are equivalent to centres in the other regions $v_3(a_0)$ or $V_3(\overline{a_0})$ or $V_3(a_1,-)$. The number of preperiod $m$ endpoints of capture paths in this region is 
$$\dfrac{1}{21}\cdot 2^{m+1}+O(1).$$
It was shown in \cite{R6} that for any integer $n>0$ there are $n$ inequivalent capture paths -- in fact, with distinct endpoints -- all crossing the unit circle at points in $\{ e^{2\pi it}: t\in (\frac{2}{7},\frac{17}{56})\} $.  There are plenty of capture components in $V_3(a_1,+)$, but not every type III hyperbolic component in $V_3(a_1,+)$ is a capture component, since as we have seen in \ref{1.3} there are at least twice as many type III hyperbolic components as can be represented. In fact because of the multiple equivalences it is somewhat less. In fact, the main results of this paper imply that the number of type III hyperbolic components in $V_3(a_1,+)$ which are represented by captures is between $c_12^{m+1}$ and $c_22^{m+1}$ where $0<c_1<c_2<\frac{1}{21}$.  

We now introduce some more notation. We  write ${\cal{Z}}_m(3/7,+)$ for the subset of ${\cal{Z}}_m(3/7)$  such that the intersection of the paths with the unit circle is in 
$$\left \{ e^{2\pi it}: t\in \left(\frac{2}{7},\frac{1}{3}\right)\cup \left(\frac{2}{3},\frac{5}{7}\right)\right \} ,$$
and ${\cal{Z}}_m(3/7,-)$ for the subset of ${\cal{Z}}_m(\frac{3}{7})$ such that the intersection with the unit circle is in 
$$\left \{ e^{2\pi it}: t\in \left(-\frac{1}{7},\frac{1}{7}\right) \right \} .$$

We write ${\cal{Z}}_m(3/7,+,+)$ for the subset of ${\cal{Z}}_m(3/7,+)$ such that the intersection is with  $\{ e^{2\pi it}\mid t\in (\frac{2}{7},\frac{1}{3})\} $. We also write
$${\cal{Z}}(3/7,+)=\cup _{m>0}{\cal{Z}}_m(3/7,+),\ \  {\cal{Z}}(3/7,-)=\cup _{m>0}{\cal{Z}}_m(3/7,-),$$
$${\cal{Z}}(3/7,+,+)=\cup _{m>0}{\cal{Z}}_m(3/7,+,+).$$
Similarly we write ${\cal{Z}}_{m}(3/7,+,p)$ for subset of ${\cal{Z}}_m(3/7,+)$, such that the intersection with the unit circle is in 
$$\left \{ e^{2\pi it}: t\in \left(\frac{1}{3}-2^{-p}\frac{1}{21},\frac{1}{3}-2^{-(p+1)}\frac{1}{21}\right)\cup \left(\frac{2}{3}-2^{-(p+1)}\frac{1}{21},\frac{2}{3}-2^{-p}\frac{1}{21}\right) \right \} .$$
We define ${\cal{Z}}(3/7,+,p)$, and   ${\cal{Z}}(3/7,+,+,p)$, analogously to the above. 

The {\em{parameter capture map}} $\Phi $ (for $V_3$) is defined by
$$\Phi ([\gamma ])=[\sigma _{\gamma }\circ s_q],$$
if $\gamma $ is a capture path for $s_q$, for $q=\frac{1}{7}$, $\frac{6}{7}$ or $\frac{3}{7}$, and where $[\gamma ]$  is the equivalence class of any capture path $\gamma $ for $s_q$ such that $\sigma _{\gamma }\circ s_q$ is equivalent to a rational map, and $[\sigma _{\gamma }\circ s_q]$ denotes the Thurston equivalence class of the critically finite branched covering $\sigma _{\gamma }\circ s_q$. We shall write $\Phi _{m}$ for the parameter capture map restricted to $([\gamma ],s_q)$ for which the second endpoint of $\gamma $ is of preperiod $m$ under $s_q$. As we have seen, the parameter capture map restricted to pairs $([\gamma ],s_{1/7})$ or $([\gamma ],s_{6/7})$ is a bijection onto a naturally defined domain, and similarly restricted to pairs $([\gamma ],s_{3/7})$ with $[\gamma ]\in {\cal{Z}}(\frac{3}{7},-)$. The case of $([\gamma ],s_{3/7})$ for $[\gamma ]\in {\cal{Z}}(3/7,+)$ is quite different. We will usually identify to domain of $\Phi $ in this case with ${\cal{Z}}(3/7,+)$. We shall also write $\Phi _m$ for the capture map $\Phi $ restricted to ${\cal{Z}}_m(3/7,+)$.
 
 We can now state our main theorems.

\begin{theorem}\label{1.8} For any capture path $\gamma $ in ${\cal{Z}}(3/7,+,+,0)$, the capture $\sigma _{\gamma }\circ s_{3/7}$ is equivalent to the centre of a type III hyperbolic component in $P_{3,m}(a_1,+)\subset V_3(a_1,+)$. Thus,
$$\Phi ({\cal{Z}}(3/7,+,+,0))\subset P_{3,m}(a_1,+).$$
\end{theorem}

\begin{theorem}\label{1.9}
There is a constant $c>0$ such that the number of type III hyperbolic components in $V_3(a_1,+)$ with second critical point in an immediate basin of preperiod $m$ and represented by $\sigma _{\gamma _i}\circ s_{3/7}$ for exactly $2^r$ different capture paths $\gamma _i$ in ${\cal{Z}}_m(3/7,+,+,0)$, for $r=1$ or any integer $r\geq 0$, is $\geq \lfloor 2^{-cr+m}\rfloor $. 
Hence, for $c_r=2^{-cr}$,
$$\# (\{ [\gamma ]\in {\cal{Z}}_m(3/7,+,+,0): \#(\Phi _m^{-1}([\sigma _{\gamma }\circ s_{3/7}])=2^r\} ) \geq c_r.2^{m}.$$

\end{theorem}

\section{So what?}\label{1.10}

\par Why does it matter? One of the most basic ideas in dynamics, especially in complex dynamics, is to model parameter space, at least local dynamics, by the parameter plane. The most famous example in complex dynamics is the combinatorial model for the Mandelbrot set for quadratic polynomials, which is based, perhaps surprisingly given the lack of immediate visual resemblance, on the dynamical plane of $z\mapsto z^{2}$. This results from the seminal work of Douady and Hubbard  \cite{D-H1} which was reinterpreted by Thurston \cite{T}. Thurston's framework is used in this article. The information obtained about the relative positions of hyperbolic components, and the dynamics in each hyperbolic component, just from the dynamics of one map, is extraordinarily detailed, and is complete. A conjectural topological model for the parameter space of quadratic polynomials is thus obtained. Important advances in proving that the model is the correct one were initially obtained by Yoccoz (described in \cite{H,Ro3}). Nearly thirty years after the model was first obtained, it remains a conjectural model, although a great deal has been proved about it (for example, \cite{L3}) and important advances have been made recently \cite{L1,L2,Ka,K-L1,K-L2,K-L3}. The topological equivalence of the combinatorial model to the actual parameter space is equivalent to the conjecture that the Mandelbrot set is locally connected (MLC). 

Given that the information about combinatorial structure of quadratic polynomial parameter space is so detailed, it is natural  to attempt a model for parameter space of quadratic rational maps in terms of the one for polynomials. This is the route that has been followed. It would also be natural to try to model parts of the parameter space on the dynamical planes of polynomials within the parameter space by simply using captures. Restricting to the subset $V_3$ of rational parameter space, this means using the dynamical planes of $h_a$ for $a=a_0$, $\overline{a_0}$ and $a_1$, or, equivalently, the dynamical plane $s_q$ for $q=\frac{1}{7}$, $\frac{6}{7}$ and $\frac{3}{7}$. It has been clear for some time that it is possible to do this in some parts of parameter space and essentially impossible in others. If we fix attention on $V_3$, then parts of the Julia sets of $h_a$ for $a=a_0$, $\overline{a_0}$ and $a_1$, are visible in $V_3(a_0)$, $V_3(\overline{a_0})$ and $V_3(a_1,-)$. It can be proved --- following the lines of analogous results of Aspenberger-Yampolsky \cite{Asp-Yam} and of Timorin \cite{Tim} for $V_2$ --- that parts of these Julia sets translate into subsets of  $V_3(a_0)$, $V_3(\overline{a_0})$. There is hope of proving something along these lines for $V_3(a_1,-)$. But the situation for capture paths in ${\cal{Z}}(a_1,+)$ is quite different. Theorem \ref{1.8} shows that the image is contained in $V(a_1,+)$. The map cannot be surjective onto type III centres in $V_3(a_1,+)$, because there are approximately twice as many centres of preperiod $m$ as there are capture paths in ${\cal{Z}}_m(3/7,+)$. It follows from the results of \cite{R6} that the map cannot be injective. But Theorem \ref{1.9} quantifies this. The way in which the proof proceeds also suggests that the capture map from ${\cal{Z}}(3/7,+)$ to $V_3(3/7,+)$ does not extend continuously to a larger set. So a transfer of topological information from dynamical plane to parameter plane is, at least, problematic.  But the proof of \ref{1.9} does  indicates ways in which discontinuities  might  arise.

\section{The parameter capture map is not typical}\label{1.11}
It is interesting to compare  the parameter capture maps $\Phi _m$ on ${\cal{Z}}_m(3/7,+)$  with typical maps between finite sets.  From \ref{1.7} we see that the size of ${\cal{Z}}_m(3/7,+)$ is $\frac{1}{21}2^{m+1}+O(1)$, while from \ref{1.3}, the size of the range $P_{3,m}(a_1,+)$ is $\frac{2}{21}2^{m+1}+O(1)$, that is, twice the size of the domain, to within $O(1)$. So the ratio of the domain size to range size for $\Phi _m$ tends to $\frac{1}{2}$ as $m$ tends to $\infty $. Typical behaviour of gaps between finite sets  is a classical problem in probability, often described in terms of ``balls in boxes'', and involving Stirling numbers of the second kind. Saddle point methods are usually  used in the solution.  It is shown \cite{S-C} that, with probability tending to $1$ as $N$ tends to $\infty $ for maps from a set of  $M$ elements to a set of $N$ elements, and for $M=bN$ with $b>0$, the number of points with $r$ preimages is within $O(\sqrt{N})$ of 
$$\frac{b^{r}}{r!}N,$$
and the image size is within $O(\sqrt{N})$ of $(1-e^{-b})N$. In particular, the maximum inverse image size for a typical map, for $b$ bounded above and below, is $\log N/\log \log N(1+o(1))$. Surprisingly, the maximum inverse image size of $\Phi _m$ is a bit  larger than typical, at least $cm=c_1\log (2^m)$. This is especially intriguing as
 we shall see that 
$$\Phi _{m}\mid {\cal{Z}}_m(s_{3/7},+,+,0)$$
 can be written as a composition of a sequence of maps which converge exponentially fast to the identity, in a reasonable sense.  Therefore one might expect $\Phi _m$ to behave more like a bijection than a typical map. In many respects, this could well be the case. The question of the relation of the image size to domain size is quite open. In principle, this could be achieved by estimation, but the practicality is another matter.
 
\section{Contrast with the subhyperbolic case}\label{1.13}

All the maps considered in this work are hyperbolic critically finite rational maps. Any critically finite rational map is  either hyperbolic or {\em{subhyperbolic}}. Subhyperbolic maps are expanding on their Julia sets, with respect to a suitable adaptation of the spherical metric: Lipschitz equivalent to the spherical metric outside the postcritical set. A critically finite map is hyperbolic if and only if the forward orbit of every critical point contains a periodic critical point. The concept of mating or capture can be extended to critically finite maps which are subhyperbolic and not hyperbolic. In the non-hyperbolic case it might not be useful to distinguish between matings and captures, and the term mating is often used.   Some very interesting studies of these have been made in the non-hyperbolic case, by Daniel Meyer \cite{Mey1,Mey2,Mey3}. It seems that the non-hyperbolic case is very different from the hyperbolic case. In particular, the number of different representations of a map by matings is much larger --- or almost certainly so. No definitive upper bound on the number of equivalent maps captures is given in the hyperbolic case, in the current work.

\section{Organisation of the paper}\label{1.12}

The organisation of this paper is as follows. 

Chapter \ref{2} concerns the Resident's View. The ``Resident's View'' is the term used here, and in earlier papers on the same general subject, for a particular manifestation of  a very basic philosophy in dynamical systems. The idea is that information about variation of dynamics in a parameter space can be obtained by scrutiny of a  fixed dynamical system within the parameter space. It is more appropriate, here, to say ``resident's views'' since the focus of the study is on determining when different models -- or different ``views'' -- represent the same dynamical system.  The representations that we are considering are all Wittner captures. Chapter \ref{2} reinterprets Theorem \ref{1.8} into the context of the Resident's View. The resulting statement, Theorem \ref{2.8}, thus becomes one of the main goals of the paper.  The proofs of both Theorems \ref{1.8} and \ref{1.9} then use Theorem \ref{2.8} and its context.

Theorem \ref{2.8} is statement about  details of a certain group action on paths in the dynamical plane of a quadratic polynomial, a group action which has been the subject of intense study, notably in \cite{R5}, where a fundamental domain was found for the group action in a particular case.  This is the main tool in the current work. Theorem \ref{2.8} gives important information about successive translates of the fundamental domain which are crossed by ``parameter capture paths''. Chapter \ref{3} is devoted to providing the foundations for proving  Theorem \ref{2.8}. The group action  is a twisted group action, where the twisted action comes from certain homemorphisms which can be regarded as conjugating homeomorphisms between different model dynamical spaces. Chapter \ref{3} is essentially devoted to the study of these homeomorphisms.

Chapter \ref{4} contains the proof of Theorem \ref{2.8}. In order to prove this theorem, and in preparation for the proof of Theorem \ref{1.9}, more detailed information about the translates of the fundamental domain are given. Theorem \ref{2.8} implies that  translates of the fundamental domain are paired off, each one with an adjacent translate. The results at the end of Chapter \ref{4}, in particular \ref{4.19}, give good estimates  on the distance between the translates, that is, on the distance between the group elements by which translation is made.  These estimates are vitally important in controlling the effect of the twisted group action, and hence for  proving the theorems \ref{1.8} and \ref{1.9}. 

Chapter \ref{5} presents key examples. The positive density sets which are the subject of Theorem \ref{1.9} are then perturbations of these key examples. The main result is Theorem \ref{5.2}.  This result says, surprisingly, that for  some constant $C>0$, and each integer $k>0$ there is a type III quadratic rational map with one critical point of period $3$ and the other of preperiod $\le C2^k$ which can be represented by exactly $2^k$ different Wittner captures with capture paths in ${\cal{Z}}(3/7,+,+,0)$. The proof is direct and does not need information about the sides of translates of the fundamental domain that are crossed by these capture paths. Nevertheless, considerable information is given about these, in preparation for the positive density results in Chapter \ref{6}, with the proof of theorem \ref{1.9}.

Chapter \ref{6} gives the proof of the positive density results, in particular, \ref{1.9}. It might be interesting to compare the implementation of the technique with the now-classical results  of positive measure sets of interval maps with characteristics of non-uniform hyperbolicity. In all theorems of this type, all estimates depend on some sort of a priori control. It never seems possible to avoid a chicken-and-egg argument, and to prove results without any sort of control. This means that, although results are obtained, concerning the abundance of certain phenomena, no absolute upper bound is achieved, and there are number of open questions.

\chapter{The resident's view of the main theorems }\label{2}

\section{A notation for Thurston equivalence}\label{2.1}

The following notation was used in \cite{R1,R2,R3,R5,R6}. If critically finite branched coverings $f$ and $g$ are Thurston equivalent, then there is an orientation-preserving  homeomorphism $\psi :\overline{\mathbb C}\to \overline{\mathbb C}$ such that $\psi (X(f_0))=X(f_1)$, and $\psi $ preserves the numbering of critical values and  there is a homotopy  $g_t$ from $\psi \circ f\circ \psi ^{-1}=g_0$ and $g_1$ such that $g_t$ is a critically finite branched covering for all $t$ with $X(g_t)=X(g)$ for all $t$. The definition of Thurston equivalence is sometimes stated in this way. If such a $\psi $ exists, then we write $f\simeq _{\psi }g$. The drawback of this notation is that $\simeq _{\psi }$ is not an equivalence relation. However, if $f\simeq _\psi g$ then $g\simeq _{\psi ^{-1}}f$, and if also $g\simeq _{\xi }h$ then $f\simeq _{\psi \circ \xi }h$.

An extension of this notation has also been used, and will be used in this paper. Suppose that $X_1(f)$ is a finite set which contains $X(f)$ and such that $f(X_1(f))\subset X_1(f)$. Suppose that $X_1(g)$ has similar properties relative to $g$. Then $(f,X_1(f))$ is said to be {\em{Thurston equivalent to $(g,X_1(g))$ }}, written $(f,X_1(f))\simeq (g,X_1(g))$ if there is a homotopy through $(f_t,X_1(f_t))$ such that:
\begin{itemize}
\item  $(f_0,X_1(f_0))=(f,X_1(f))$ and $(f_1,X_1(g))=g$;
 \item $f_t$ is a critically finite branched covering with $X(f_t)\subset X_1(f_t)$ and $f_t(X_1(f_t))\subset X_1(f_t)$;
 \item $X_1(f_t)$ moves isotopically with $t$ and numbering of critical points is preserved.
 \end{itemize}
 Similarly to the above, if $(f,X_1(f))\simeq (g,X_1(g))$, then $\psi $ and $g_t$ exist as before, and in addition, $\psi (X_1(f))=X_1(g)$ and $g_t(X_1(g))\subset X_1(g)$ for all $t$. We then write $(f,X_1(f))\simeq _{\psi }(g,X_1(g))$. 
 
 A fact which is used throughout \cite{R1,R2,R3,R5,R6} and which will be used here also, is that if $(f,X_1(f))\simeq _{\psi _0}(g,X_1(g))$, then for each $n\geq 1$ there is $\psi _n$ such that $(f,f^{-n}X_1(f))\simeq _{\psi _n}(g,g^{-n}(X_1(g))$, and such that $\psi _n$ and $\psi _k$ are isotopic via an isotopy which is constant on $f^{-k}(X_1(f))$, for each $0\leq k\leq n$.
 
 \section{The boundary of the universal cover}\label{2.2}

The theory of the Resident's View uses the following facts and notation concerning boundaries of universal covers of Riemann surfaces. Suppose that $V=\overline{\mathbb C}\setminus P$ where $P\subset \Cbar $ is finite and $\#(P)\geq 3$, so that $V$ has negative Euler characteristic. (This is the set-up in all  our examples, but $\Cbar $ could be replaced by any compact Riemann surface $S$, provided that $S\setminus P$ has negative Euler characeristic). The holomorphic universal cover of $V$ is the open unit disc $D$. For fixed $a\in V$, we define $\pi _1(V,P,a)$ to be the set of homotopy classes of paths $\beta :[0,1]\to \Cbar $ such that $\beta ([0,1))\subset V$, $\beta (0)=a$ and $\beta (1)\in P$, where we restrict to homotopies through such paths. Then for any $[\beta ]\in \pi _1(V,P,a)$, and any lift $\tilde \beta $ of $\beta \mid [0,1)$ to $D$, $\lim _{t\to 1}\tilde \beta (t)$ exists as in $D$, and is uniquely determined by $[\beta ]$ and $\tilde \beta (0)$. Therefore, having fixed a lift $\tilde a$ of $a$ to $D$, we can regard $\pi _1(V,a,P)$ as a subset of $\partial D$, by identifying $[\beta ]$ with $\tilde \beta (1)$ for any lift $\tilde \beta $ with $\tilde \beta (0)=\tilde a$.  This subset is dense.

Since $D$ is holomorphic, the action on $D$ of the covering group of $V$  is by M\"obius transformations, and hence this action extends continuously to $\partial D$. The covering group is naturally isomorphic to the fundamental group $\pi _1(V,a)$, which acts naturally on $\pi _1(V,P,a)$ by concatenation, that is, for $[\gamma ]\in \pi _1(V,a)$ and $[\beta ]\in \pi _1(V,P,a)$, we define $[\gamma ].[\beta ]=[\gamma *\beta ]$. This coincides with the action by M\"obius transformations on $\partial D$, using the identification of $\pi _1(V,P,a)$ with a subset of $\partial D$.

\section{Action of mapping class groups}\label{2.3}

If $Y\subset \Cbar$ is finite, then we denote by ${\rm{MG}}(\Cbar, Y)$ the {\em{mapping class group}} of $(\Cbar ,Y)$, that is, the group of isotopy classes of orientation preserving homeomorphisms of $\Cbar $ which leave $Y$ invariant, modulo isotopies which are constant on $Y$. If $Z\subset Y$ then there is a natural homomorphism from ${\rm{MG}}(\Cbar, Y)$ to ${\rm{MG}}(\Cbar, Z)$, given by mapping the isotopy class of $\psi $ in ${\rm{MG}}(\Cbar, Y)$ to the isotopy class of $\psi $ in  ${\rm{MG}}(\Cbar, Z)$.
If $\#(Z)\geq 3$, any lift to the universal cover $D$ of any homeomorphism of $\Cbar $ which preserves $Z$ extends homeomorphically to $D\cup \partial D$. If $y\in Y\setminus Z$ and $\widetilde y$ is a fixed lift of $y$ to $D$,  and $[\psi ]\in {\rm{MG}}(\Cbar ,Y)$ with $\psi (y)=y$ and $\tilde \psi $ is a lift with $\widetilde \psi (\widetilde y)=\widetilde y$, then $\widetilde \psi \mid \partial D$ is uniquely detemined by $\widetilde y$ and  $[\psi ]\in {\rm{MG}}(\Cbar , Y)$.

\section{The Resident's View.}\label{2.4}

The Resident's View  gives information about one punctured sphere in terms of another.  As in \ref{1.3}, we define $P_{3,m}$ to be the set of centres of hyperbolic components of types II, and of type III of preperiod $\leq m$. We also define 
$$V_{3,m}=V_3\setminus P_{3,m}=\Cbar \setminus (\{ 0,\infty \} \cup P_{3,m}).$$
 As we saw in \ref{1.2}, there are just two type II hyperbolic components, with centres $1$ and $-1$. So $V_{3,0}=\mathbb C\setminus \{ 0,\pm 1\} $, and, for every $m\geq 0$, $V_{3,m}$ is the complement in the Riemann sphere of four or more punctures. Hence it is certainly  a surface of negative Euler characteristic, and the holomorphic universal cover identifies with the open unit disc $D$. Therefore, as above, $\pi _1(V_{3,m},P_{3,m},a)$ identifies with a subset of $\partial D$ for a fixed choice $\tilde a\in D$ of lift of $a$, and $\pi _1(V_{3,m},a)$ acts on $D\cup \partial D$ by M\"obius transformations, where this action coincides with the usual action on $\pi _1(V_{3,m},P_{3,m},a)\subset \partial D$ by contenation of paths.

The other punctured sphere, with which we want to compare the first, is a punctured dynamical plane. For $a=a_0$ or $\overline{a_0}$ or $a_1$, we write
$$Z_m(h_a)=h_a^{-m}(\{ 0,\infty ,1\} )=h_a^{-m}(\{ 0,h_a(0),h_a^2(0)\} ),$$
and 
$$Y_m(h_a)=Z_m(h_a)\cup \{ v_2(a)\} ,$$
where $v_2(a)=\frac{2a}{a+1}$ is the second critical point  of $h_a$, which is fixed by $h_a$ for these values of $a$ and hence is also the second critical value. Thus, $Y_0(h_a)$ is the postcritical set $X(h_a)$ and $Y_m(h_a)=h_a^{-m}(Y_0(h_a))$ for all $m\geq 0$, for these values of $a$. Similarly (and as in \ref{1.7}), for $q=\frac{1}{7}$ or $\frac{6}{7}$ or $\frac{3}{7}$, we write
$$Z_{m}(s_q)=s^{-m}(\{ s_q^{n}(0):n\geq 0\} ,$$
and
$$Y_{m}(s_q)=Z_{m}(s_q)\cup \{ \infty\} ,$$
remembering that $\infty  $ is the fixed critical point of $s_q$. Since $\Cbar \setminus Z_m$, for $Z_m=Z_m(h_a)$ or $Z_m(s_q)$ for any of these values of $a$ or $q$, the holomorphic universal cover is $D$ once again, and we can again regard $\pi _1(\Cbar \setminus Z_m(h_a),Z_{m}(h_a),v_2(a))$ (or $\pi _1(\Cbar \setminus Z_m(s_q),Z_m(s_q),\infty )$) as a subset of $\partial D$.  Under our definition of Thurston equivalence, the Thurston equivalence class of $(f,Y)$ for any critically finite branched covering $f$ with $\#(Y)\geq 3$ is simply connected. There is therefore a natural isomorphism between $\pi _1(\Cbar \setminus Z_m(h_a),v_2(a))$ and $\pi _1(\Cbar \setminus Z_m(s_q),\infty )$, and a natural identification between $\pi _1(\Cbar \setminus Z_m(h_a)_,Z_m(h_a),v_2(a))$ and $\pi _1(\Cbar \setminus Z_m(s_q),Z_m(s_q),\infty )$ for $(a,q)=(a_0,\frac{1}{7})$ or $(\overline{a_0},\frac{6}{7})$ or $(a_1,\frac{3}{7})$. 

From now on we fix lifts $\widetilde a$ and $\widetilde \infty $ to the universal cover $D$, of $a\in V_{3,m}$ and $\infty =v_2(s_q)\in \Cbar \setminus Z_m(s_q)$, in order to define the identifications of $\pi _1(V_{3,m},P_{3,m},a)$ and $\pi _1(\Cbar \setminus Z_m(s_q),Z_m(s_q),\infty )$ with subsets of $\partial D$.  
 
Given these identifications, Resident's View gives two maps, both of which are called $\rho $, or $\rho (.,s_q)$ if more precision is needed,  in \cite{R5} (and called $\rho _2$ and $\rho $ in \cite{R3}):
$$\rho :\pi _1(V_{3,m},P_{3,m},a)\to \pi _1(\Cbar \setminus Z_m(s_q),Z_m(s_q),\infty ),$$
$$\rho :\pi _1(V_{3,m},a)\to \pi _1(\Cbar \setminus Z_m(s_q),\infty ).$$
Both these maps $\rho $ are injective. The first map $\rho $ is also monotone from the dense subset $\pi _1(V_{3,m},P_{3,m},a)$ of $\partial D$ into $\pi _1(\Cbar \setminus Z_m(s_q),Z_m(s_q),\infty )\subset \partial D$. The image is not dense, except in the case when $m=0$, and hence the map $\rho $ cannot extend continuously to $\partial D$. However, it does extend continuously except at countably many points, where right and left limits exist.
 We define $D'$ to be the convex hull of the closure of the image of $\rho $ in $\partial D$.  The convex hull can be taken with respect to either the hyperbolic or Euclidean metric on $D$, since the two convex hulls are homeomorphic. 
 There is a slight extension of $\rho $ to the domain $\pi _1(V_{3,m},P_{3,m}\cup \{ 0,\infty \} ,a)$ such that the images of $[\beta ]$ is one of countably many points of $\partial D\cap \partial D'$ if $\beta $ ends at $0$, and is the closure of one of the countably many components of $\partial D'\cap D$ if $\beta $ ends at $\infty $. In the latter case, we identify $\rho ([\beta ])$ with the set of paths from $\widetilde{\infty}=\widetilde {\rho (\beta)} (0)$ to $\overline{\partial D'\cap D}$, modulo homotopy equivalence fixing the first endpoint and keeping the second endpoint  in $\overline{\partial D'\cap D}$. 
 
 The set $\rho (\pi _1(V_{3,m},P_{3,m},a))$ is the set of all paths in $\pi _1(\Cbar \setminus Z_m(s_q),Z_m(s_q),\infty )$ which lift to paths from $\widetilde{\infty }$ to $\partial D\cap \partial D'$. In particular, it includes all paths in ${\cal{Z}}_m(q)$.
 
 The second map $\rho $, with domain $\pi _1(V_{3,m},a)$,  is {\em{not }} a group homomorphism. Nevertheless, the action of $\pi _1(V_{3,m},a)$ on $\partial D$ translates to an action on $\partial D'$, and, indeed on $\partial D$, in the following way.  Suppose that $\gamma \in \pi _1(V_{3,m},a)$  and that $\rho (\gamma )=\alpha $. Then there is $[\psi ]\in {\rm{MG}}(\Cbar,Y_{m+1}(s_q))$ such that
 $$s_q\simeq _{\psi }\sigma _{\alpha }\circ s_q$$
 and for all $\beta \in \pi _1(V_{3,m},P_{3,m},a)$,
 $$\rho ([\gamma  *\beta ])=[\alpha *\psi (\rho (\beta ))]$$
 where $\rho (\beta)$ is any path representing $\rho ([\beta ])$.  Moreover, $[\psi ]$ is uniquely determined and the map 
 $$\Xi:[\gamma ]\mapsto [\psi ]$$
  is a group  isomorphism of $\pi _1(V_{3,m},a)$ onto a subgroup $G$ of
 $$\{ [\psi ]\in {\rm{MG}}(\Cbar,Y_{m+1}(s_q)):\  \psi (\infty )=\infty ,\ \exists \alpha \in \pi _1(\Cbar \setminus Z_m(s_q),\infty ),\ s_q\simeq \sigma _{\alpha }\circ s_q\} .$$
(The notation $\Xi$ has not been used in the earlier papers being quoted, where, for reasons which we will not go into, an anti-isomorphism, rather than an isomorphism, was used.) Here, $[\psi ]$ determines $[\alpha ]$ uniquely, and vice versa. For $[\psi ]\in G$ and $\alpha $ as above, we define $\widetilde \psi $ to be the lift of $\psi $ which maps $\widetilde \infty $ to $\widetilde \alpha (1)$, where the path $\alpha $ has domain $[0,1]$ and the lift $\widetilde \alpha $ of $\alpha $ is such that $\widetilde \alpha (0)=\widetilde \infty $. Then $[\psi ]\to \widetilde \psi \vert \partial D$ is an isomorphism onto its image, and therefore an action of $G$ on $\partial D$ is defined by 
$$[\psi ].x=\widetilde \psi (x)\ \ \ {\rm{(}}x\in \partial D{\rm{).}}$$
If $x=[\zeta ]\in \pi _1(\Cbar \setminus Z_m(s_q),Z_m(s_q),\infty )$, then $\widetilde \psi ([\zeta ])=[\alpha *\psi (\zeta )]$.  The action of $\pi _1(V_{3,m},a)$ on $\partial D$ therefore transfers under $\rho $ to the restriction (to $\partial D'\cap \partial D$) of a natural action of $G$ on $\partial D$.

 \section{Key features of the fundamental domain}\label{2.5}
 
 The fundamental domain for $V_{3,m}$ that was chosen in \cite{R5} is an ideal polygon $F_{3,m}$ in the universal covering space of $V_{3,m}$, that is, a region of $D$ that is bounded by finitely many geodesics (in the hyperbolic metric), with closure intersecting $\partial D$ in just finitely many points, which are all in $\pi _1(V_{3,m},P_{3,m}\cup \{ 0,\infty \} ,a)$. Then $F_{3,m}$ is clearly determined by its vertices. Sides of $F_{3,m}$ are determined by pairs of adjacent vertices. Since $F_{3,m}$ is a fundamental domain, any adjacent pair of vertices has to be identified with some other pair by an element of $\pi _1(V_{3,m},a)$. Any element of $\pi _1(V_{3,m},a)$ identifies at most one pair of vertices of $F_{3,m}$, and the identifying elements freely generate $\pi _1(V_{3,m},a)$. 
 The set of paths chosen in \cite{R5} for the vertices of $F_{3,m}$ is $\Omega _m$, where 
 $$\Omega _m=\gamma _0'*\Omega _m(a_0)\cup \gamma _0*\Omega _m(\overline{a_0})\cup \Omega (a_1,-)\cup \Omega (a_1,+)$$
 where:
 \begin{itemize}
 \item  $\gamma _0'$ and $\gamma _0'$ are paths from $a_1$ to $a_0$ and $\overline{a_0}$;
 \item  the first endpoint   of paths in $\Omega _m(a)$ is at $a$, and the second endpoints  are in $V_3(a)$  for $a=a_0$ and $a=\overline{a_0}$;
\item  the first endpoint of paths in $\Omega _m(a_1,\sigma )$ is at $a_1$ and the second endpoints are in $V_3(a_1,\sigma )$, for $\sigma =\pm $. 
\end{itemize}

 Since $\rho :\pi _1(V_{3,m},P_{3,m},a)\to \pi _1(\Cbar\setminus Z_m(s_q),Z_m(s_q),\infty )$ is monotone and injective, it makes sense to find the subset of $\pi _1(V_{3,m},P_{3,m},a)$ which forms a vertex set for an ideal polygon which is a fundamental domain, by finding the image under $\rho $, because $\rho $ preserves the order of vertices, and adjacency of vertices. This was the technique employed in \cite{R5}. A finite subset of $\pi _1(\Cbar \setminus Z_m(s_{3/7}),Z)_m(s_{3/7}),\infty )$ was found, together with a pairing of pairs of elements of this subset under the action of $G$. It was shown that under very mild extra conditions (which were satisfied) the preimage under $\rho $ is indeed a fundamental domain for $V_{3,m}$, or, more precisely, for the action for the covering group $\pi _1(V_{3,m},a_1)$ on $D$. Actually, $\rho =\rho (.,s_{3/7})$ is not used directly to describe all of $\Omega _m$ in \cite{R5}. Using four particular paths (called $\omega _1$, $\omega _{-1}$, $\omega _{1}'$ and $\omega _{-1}'$) the space $V_{3,m}$ is split into four regions which contain one each of $V_3(a_0)$, $V_3(\overline{a_0})$ and $V_3(a_1,\pm )$. We use $\rho (.,s_{3/7})$ to describe paths of $\Omega _{m}(a_1,\pm )$ --- that is, to describe paths in $\Omega _m(a_1,+)$, which are between $\omega _1$ and $\omega _1'$ (including these) and paths in $\Omega _m(a_1,-)$ which are between $\omega _{-1}$ and $\omega _{-1}'$, including these. We use $\rho (., s_{1/7})$ to describe paths in $\Omega _m(a_0)$, all of which are strictly between $\overline{\gamma _0'}*\omega _{1}'$ and $\overline{\gamma _0'}*\omega _{-1}'$, up to homotopy, and $\rho (.,s_{6/7})$ is used to describe paths in $\Omega _m(\overline{a_0})$, which are all strictly between $\overline{\gamma _0}*\omega _{1}$ and $\overline{\gamma _0}*\omega _{-1}$, up to homotopy. 
 
The description then includes
 $$\rho (\Omega _m(a_0),s_{1/7})={\cal{Z}}_m(1/7),\ \ \rho (\Omega _m(\overline{a_0}),s_{6/7})={\cal{Z}}_m(6/7),$$
 $$\rho (\Omega _m(a_1,-),s_{3/7})={\cal{Z}}_m(a_1,-).$$

 The set  $\Omega _m(a_1,+)$ is defined by its image under $\rho (.,s_{3/7})$. It has some intersection with ${\cal{Z}}_m(a_1,+)$, but  $\rho (\Omega _m(a_1,+))$ neither contains it, nor is contained in it. But part of the result of \ref{1.8} is that each  capture $\sigma _{\gamma  }\circ s_{3/7}$, for $\gamma \in  {\cal{Z}}_m(a_1,+)$, is equivalent to $\sigma _{\beta }\circ s_{3/7}$, for some $\beta \in \Omega _m(a_1,+)$. This is the same as saying that, under the identification of ${\cal{Z}}_m(a_1,+)$ with a subset of $\partial D'\cap \partial D$, the $G$-orbit of $\rho (\Omega _m(a_1,+))$ contains ${\cal{Z}}_m(a_1,+)$.
 
\section{Key features  of  the fundamental domain paths}\label{2.6}
 
 We now summarize some key features  of the paths in $\rho (\Omega _m(a_1,+),s_{3/7})$ from \cite{R5}, and the corresponding generators which identify adjacent pairs. From now on, we write $s$ for $s_{3/7}$ except when this might cause confusion.
 
The paths in $\rho (\Omega _m(a_1,+),s_{3/7})$ are split into two sets, called $R_m$ and $R_m'$. The set $R_m$   is totally ordered using an ordering $<$. One way to describe the ordering is just in terms of the ordering of the endpoints in $\partial D\cap \partial D'$. Apart from  one common element which is maximal in $R_m$, the sets $R_m$ and $R_m'$ are disjoint. Adjacent paths in $R_m$ (or $R_m'$) are adjacent in the ordering. Each adjacent pair in $R_m$ is matched with exactly one adjacent pair in $R_m'$ by exactly one adjacent pair in $R_m'$ by exactly one element of $G$, and vice versa. More precisely, for each adjacent pair $(\beta _1,\beta _2)$ of paths in $R_m$ with $\beta _1<\beta _2'$, there is an adjacent pair $(\beta _1',\beta _2')$ in $R_m'$ and $[\psi ]\in G$ such that $[\psi ].[\beta _j']=[\beta _j]$ for $j=1$, $2$, and $(\beta _1',\beta _2')$. Using the description of the action of $G$ in \ref{2.5} this means that there is $[\alpha ]\in \pi _1(\Cbar \setminus Z_m(s),\infty )$ such that $s\simeq _\psi \sigma _{\alpha }\circ s$ (and $[\alpha ]$ is uniquely determined by this property) and:
$$[\alpha *\psi (\beta _j')]=[\beta _j],\ \ j=1,2.$$
The common maximal path in $R_m$ and $R_m'$ is a path called $\beta _{1/3}$ from $\infty $ to $e^{2\pi i(1/3)}$, and contained in $\{ z:\vert z\vert >1\} $ apart from that endpoint. This is the only path in $R_m\cup R_m'$ whose endpoint is not in a lift of $Z_m(s)$, although it is in $\partial D'$. It is the image under $\rho $ of the path along the negative real axis from $a_1$ to $\infty $ (enlarging the domain of $\rho $ slightly, as described in \ref{2.5}). In general, if $e^{2\pi ix}$ is on the boundary of a gap of $L_{3/7}$, we write $\beta _x$ for the capture path which crosses the unit circle at $e^{2\pi ix}$: this is unambiguous, because $e^{2\pi ix}$ is in the boundary of at most one gap. We write $q_p=\frac{1}{3}-2^{-2p}\frac{1}{21}$. Then $\beta _{q_p}\in {\cal{Z}}_m(a_1,+)$ for each $0\leq 2p\leq m$: the path $\beta _{2/7}$  is a type II capture path, but the others are all of type III. There is a further decomposition of $R_m$ as
$$R_m=\cup _{0\leq 2p\leq m}R_{m,p}\cup \{ \beta _{1/3}\} ,$$
where $\beta _{q_{p+1}}\in R_{m,p}\cap R_{m,p+1}$, and $R_{m,p}$ is the set of all paths $\beta $ in $R_{m}$ with $\beta _{q_p}\leq \beta \leq \beta _{q_{p+1}}$. The minimal path $\beta _{2/7}$ in $R_m$ is matched with $\beta _{5/7}\in R_m'$, and more generally the minimal path $\beta _{q_p}$ in $R_{m,p}$ is matched with $\beta _{1-q_p}$ in $R_{m,p}'$, where we write $R_{m,p}'$ for the set of paths in $R_m'$ which are matched with paths in $R_m'$. 

The paths in $R_{m,p}'$ were never described explicitly in \cite{R5}, and we shall not do so here. They are determined inductively from the paths in $R_{m,0}$, using the ordering on paths of $R_{m}$ and the matching with paths in $R_m'$, starting from the matching of $\beta _{2/7}$ with $\beta _{5/7}$. However, we shall use the following later. Let $\beta _1<\beta _2\in R_{m,p}$, and let $\beta _i$ be matched with $\beta _i'\in R_{m,p}$. Let $(\beta _1,\beta _{1,2})$ and $(\beta _{2,1},\beta _{2,2})$ be adjacent pairs including $\beta _1$ and $\beta _2$, with $\beta _1<\beta _{1,2}$ and $\beta _{2,1}<\beta _2$, matched with $(\beta _1',\beta _{1,2}' )$ and $(\beta _{2,1}',\beta _2')$ in $R_{m,p}'$, and let $(\alpha _1,\psi _1)$ and $(\alpha _2,\psi _2)$ effect the matching, that is,
$$(s,Y_m(s))\simeq _{\psi _j}(\sigma _{\alpha _j}\circ s,Y_m(s))$$
and 
$$[\alpha _j*\psi _j(\beta _j')]=[\beta _j],\ \ [\alpha _1*\psi _1(\beta _{1,2}')]=[\beta _{1,2}],\ \ [\alpha _2*\psi _2(\beta _{2,1}')]=[\beta _{2,1}].$$
Then $\alpha _1*\overline{\alpha _2}$ is trivial in $\pi _1(\Cbar \setminus Z_{p+1}(s),\infty )$. This follows from the fact that, up to  free homotopy in $\Cbar \setminus Z_{2p+2}(s)$, the loop  $\alpha _1*\overline{\alpha _2}$ lies in a set $U^p$ --- which is defined precisely in \cite{R5}, and which is also described in the next chapter  --- and which is disjoint from $Z_{2p+2}(s)$. 

\section{Transferring the fundamental domain under $\rho$}\label{2.7}

We have seen  that 
$$R_m\cup R_m'\subset \rho (\Omega _m).$$
By the monotonicity of $\rho $, and by the form chosen for $F_{3,m}$, the geodesics homotopic to paths $\overline{\rho (\omega _1)}*\rho (\omega _2)$, for all adjacent pairs $(\omega _1,\omega _2)$ in $\Omega $, bound a polygon in $D'$, which by abuse of notation we call $\rho (F)$. The translates of $F$ by $G$ cover $D'$ by the Resident's View and the interiors of the translates of $F$ are disjoint, also by the Resident's View.   We cannot expect the translates to remain disjoint if we project down from $D'$ to $\Cbar \setminus Z(s_q)$, because the covering group of $\Cbar \setminus Z(s_q)$ is not $G$. Nevertheless it is useful to view the translates in the projection. If $(\rho (\omega _1),\rho (\omega _2))=(\beta _1,\beta _2)$ then an adjacent pair $(\beta _1,\beta _2)$ in $R_m$, then 
$$\overline{\rho (\omega _1)}*\rho (\omega _2)=\overline{\beta _1}*\beta _2,$$
and we can identify this with its left to $D'\subset D$, using the lifts of $\beta _1$ and $\beta _2$ which start from $\widetilde{\infty }$ Similarly if $[\gamma ]\in \pi _1(V_{3,m},a_1)$  with $\rho [\gamma ]=\alpha $ and $\Xi ([\gamma ]) =[\psi ]$,then
$$\overline{\rho (\gamma *\omega _1)}*\rho (\gamma *\omega _2)=\psi (\overline{\beta _1}*\beta _2),$$
where, this time, we use the lift of $\alpha *\psi (\beta _1)$ (or $\alpha *\psi (\beta _2)$) starting from $\widetilde{\infty }$ to determine which preimage of $\psi (\overline{\beta _1}*\beta _2)$ we have in $D'$.

Then we have the following theorem. The proof of this will only be completed in Chapter \ref{4}.  Remember (as was recalled in \ref{2.4}) that ${\cal{Z}}_m(3/7)\subset \rho (\pi _1(V_{3,m},P_{3,m},a_1))$, and so for any $\beta \in {\cal{Z}}_m(3/7)$  the point on $\partial D'$ which is identified with $\beta $, must be in the orbit under $\pi _1(V_{3,m},a_1)$ of $\Omega _m$ . This theorem immediately gives Theorem \ref{1.8}, because the branched coverings $\sigma _{\zeta }\circ s$, for $\zeta \in R_m$, are in  $V_m(a_1,+)$, and it gives the structure for proving \ref{1.9}.

\begin{theorem}\label{2.8}

Let $\beta \in {\cal{Z}}_m(3/7,+,+,0)$. Let   
$$((\omega _{2i-1},\omega _{2i}):1\leq i\leq n(\beta ))$$ and 
$$(\gamma _i:1\leq i\leq n(\beta ))$$
 be the unique sequences of, respectively, adjacent pairs in $\Omega _m$ and elements of  $\pi _1(V_{3,m},a_1)$, such that $\beta $ successively crosses the arcs $\overline{\rho (\gamma _i*\omega _{2i-1})}*\rho (\gamma _i*\omega _{2i-1})$, with $\beta =\rho (\gamma _{n(\beta )}*\omega _{2n(\beta )})$. Then the following hold.
\begin{itemize}
\item[1.] For $1\leq i\leq n(\beta )$, there is an adjacent pair $(\beta _{2i-1},\beta _{2i})$ in $R_{m,0}$ and matched with $(\beta _{2i-1}',\beta _{2i}')$ in $R_{m,0}'$ such that
$$(\rho (\omega _{2i-1}),\rho (\omega _{2i-1}))=\begin{cases}(\beta _{2i-1},\beta _{2i}){\rm{\ if\ }}i{\rm{\ is\ odd,}}\\ 
(\beta _{2i-1}',\beta _{2i}'){\rm{\ if\ }}i{\rm{\ is\ even.}}\end{cases}$$
Moreover, $\sigma _{\beta }\circ s$ is Thurston equivalent to $\sigma _{\beta _{2n(\beta )-1}}\circ s$, and $\Phi (\beta )=\beta _{2n(\beta )-1}$.
\item[2.] Let  $[\psi _{2j-1}]$ and $\alpha _{2j-1}$ effect the matching between the adjacent pairs $(\beta _{2j-1}',\beta _{2j}')$ in $R_{m,0}'$ and $(\beta _{2j-1},\beta _{2j})$ in $R_{m,0}$, that is:
$$(s,Y_m)\simeq _{\psi _{2j-1}}(\sigma _{\alpha _{2j-1}}\circ s,Y_m)$$
and
$$\alpha _{2j-1}*\psi _{2j-1}(\beta _{2j-1}')=\beta _{2j-1},\ \ \alpha _{2j-1}*\psi _{2j-1}(\beta _{2j}')=\beta _{2j}.$$
 Write 
$$\psi _{4j-3,4j-1}=\psi _{4j-3}\circ \psi _{4j-1}^{-1},\ \ \alpha _{4j-3,4j-1}=\alpha _{4j-3}*\psi _{4j-3,4j-1}(\overline{\alpha _{4j-1}}),$$
 and, for $j<k$,
$$\psi _{4j-3,4k-1}=\psi _{4j-3,4j-1}\circ \cdots \circ \psi _{4k-3,4k-1}$$
$$\psi _{4j-3,4k+1}=\psi _{4j-3,4k-1}\circ \psi _{4k+1},$$
$$\alpha _{4j-3,4k-1}=\alpha _{4j-3,4j-1}*\psi _{4j-3,4j-1}(\alpha _{4j+1,4k-1})=\alpha _{4j-3,4k-5}*\psi _{4j-3,4k-5}(\alpha _{4k-3,4k-1}),$$
and
$$\alpha _{4j-3,4k+1}=\alpha _{4j-3,4j-1}*\psi _{4j-3,4j-1}(\alpha _{4j+1,4k+1})=\alpha _{4j-3,4k-1}*\psi _{4j-3,4k-1}(\alpha _{4k+1}).$$
Then, for $r=0$ or $1$, taking $\psi _{1,-1}$ to be the identity, and $\alpha _{1,-1}$ to be trivial,
$$\rho (\gamma _i)=\alpha _{1,2i-3},\ \ \Xi (\gamma _i)=[\psi _{1,2i-3}],$$
$$\rho (\gamma _i*\omega _{2i-r})=\begin{cases}\alpha _{1,2i-3}*\psi _{1,2i-3}(\beta _{2i-r}){\rm{\ if\ }}i{\rm{\ is\ odd,}}\\
\alpha _{1,2i-3}*\psi _{1,2i-3}(\beta _{2i-r}')=\alpha _{1,2i-1}*\psi _{1,2i-1}(\beta _{2i-r}){\rm{\ if\ }}i{\rm{\ is\ even.}}\end{cases}$$
In other words
$$\rho (\gamma _i*\omega _{2i-r})=\alpha _{1,4\lfloor i/2\rfloor -1}*\psi _{1,4\lfloor i/2\rfloor -1}(\beta _{2i-r}).$$
Consequently if $i$ is odd then $\rho (\gamma _i*\omega _{2i-1})$ and $\rho (\gamma _i*\omega _{2i})$ are arbitrarily small perturbations of 
$$\beta _1*\psi _1(\overline{\beta _1'}*\beta _3')*\psi _{1,3}(\overline{\beta _3}*\beta _5)\cdots *\psi _{1,2i-5}(\overline{\beta _{2i-5}'}*\beta _{2i-3}')*\psi _{1,2i-3}(\overline {\beta _{2i-3}}*\beta _{2i-1}),$$
and
$$\beta _2*\psi _1(\overline{\beta _2'}*\beta _4')*\psi _{1,3}(\overline{\beta _4}*\beta _6)\cdots *\psi _{1,2i-5}(\overline{\beta _{2i-4}'}*\beta _{2i-2}')*\psi _{1,2i-3}(\overline {\beta _{2i-2}}*\beta _{2i}),$$
and if $i$ is even they are arbitrarily small perturbations of 
$$\beta _1*\psi _1(\overline{\beta _1'}*\beta _3')*\psi _{1,3}(\overline{\beta _3}*\beta _5)\cdots *\psi _{1,2i-5}(\overline{\beta _{2i-5}}*\beta _{2i-3})*\psi _{1,2i-3}(\overline {\beta _{2i-3}'}*\beta _{2i-1}'),$$
and
$$\beta _2*\psi _1(\overline{\beta _2'}*\beta _4')*\psi _{1,3}(\overline{\beta _4}*\beta _6)\cdots *\psi _{1,2i-5}(\overline{\beta _{2i-4}}*\beta _{2i-2})*\psi _{1,2i-3}(\overline {\beta _{2i-2}'}*\beta _{2i}').$$
In particular,
$$\Phi (\beta )=\psi _{1,4\lfloor n(\beta )/2\rfloor -1}^{-1}(\overline{\alpha _{1,4[n(\beta )/2]-1}}*\beta ).$$

 \end{itemize}

 \end{theorem}

\section{Remarks}\label{2.9}

It does not seem to be true that if $\beta \in {\cal{Z}}_m(3/7,+,-,p)$ then $\omega _{2n(\beta )-1}(\beta )\in \Omega _{m,p}$, even for $p=0$. Some calculations in Chapter \ref{6} suggest that it might be true that, for such a $\beta $,  
 $$\omega _{2n(\beta )-1}(\beta )\in \Omega _{m,p}\cup \Omega _{m,p+1}.$$

The main statement in Theorem \ref{2.8} is item 1. Since  an adjacent pair $(\beta _{2i-1},\beta _{2i})$ in $R_{m,0}$  is only matched with an  adjacent pair $(\beta _{2i-1}',\beta _{2i}')$ in $R_{m,0}'$, and indeed any adjacent pair in $\Omega _m$ is only matched with one other, $\gamma _i$ is uniquely determined by $(\omega _{2j-1},\omega _{2j})$ for $j<i$. More precisely $\overline{\gamma _i}*\gamma _{i+1}$
is determined by $(\omega _{2i-1},\omega _{2i})$, and hence the same is true for the images under $\rho $and $\Xi $. This is true whichever pairs  $(\omega _{2i-1},\omega _{2i})$ arise. But item 1 tells us( remembering that $\Xi $ is a group homomorphism and $\rho $ is not one, but is  a twisted group homomorphism), that
 $$\rho (\overline{\gamma _i}*\gamma _{i+1})=\begin{cases}\alpha _{2i-1}{\rm{\ if\ }}i{\rm{\ is\ odd,}}
 \\ \psi _{2i-1}^{-1}(\overline{\alpha _{2i-1}}){\rm{\ if\ }}i{\rm{\ is\ even.}}\end{cases}$$ and 
 $$\Xi  (\overline{\gamma _i}*\gamma _{i+1})=\begin{cases}[\psi _{2i-1}]{\rm{\ if\ }}i{\rm{\  is\ odd,}}\\ 
 [\psi _{2i-1}^{-1}]{\rm{\ if\ }}i{\rm{\  is\ even}}\end{cases}$$
So item 2 follows from item 1.

Part of the proof needs some detailed analysis  of the homeomorphisms $\psi _{2j-1,2k-1}$. The proof of this result is therefore given in two stages, partly in Chapter \ref{3} and partly in Chapter \ref{4}. 

2 of \ref{2.8} also suggests an algorithm for determining when two different captures $\sigma _{\beta }\circ s$ and $\sigma _{\zeta }\circ s$ are equivalent, for $\beta $ and $\zeta \in {\cal{Z}}_m(a_1,+,+,0)$. 

Write 
$$\beta _i(\beta ),\ \ \psi _{2i-1,\beta },\ \ \psi _{2j-1,2k-1,\beta },\ \ \alpha _{2j-1,\beta },\ \ \alpha _{2j-1,2k-1,\beta }$$
 for $\beta _i$, $\psi _{2i-1}$, $\psi _{2j-1,2k-1}$, $\alpha _{2j-1}$ $\alpha _{2j-1,2k-1}$ in Theorem \ref{2.8}. Then $\sigma _{\beta }\circ s$ and $\sigma _{\zeta }\circ s$ are Thurston equivalent if and only if $(\beta _{2n(\beta )-1}(\beta ),\beta _{2n(\beta )}(\beta ))=(\beta _{2n(\zeta )-1}(\zeta ),\beta _{2n(\zeta )}(\zeta ))$. If $\beta \neq \zeta $, then it must be the case that $\beta _i(\beta )\neq \beta _i(\zeta )$ for at least one $i$. Theorem \ref{1.9} will be proved simply by exploring how this is possible.

\chapter{Symbolic dynamics,  paths in $R_{m,0}$, and homeomorphisms}\label{3}

The aim of this chapter is to analyse the homeomorphism which will arise as the homeomorphisms $\psi _{2i-1,2j-1}$ which arise in Theorem \ref{2.8}. This will enable us to do much of the proof of Theorem \ref{2.8} at the end of the chapter.
 
  \section{Symbolic Dynamics}\label{3.1}
  
 Symbolic dynamics with respect to a Markov partition for $h_{a_1}$, or equivalently, for $s_{3/7}$ (using the semiconjugacy from $s_{3/7}$ to $h_{a_1}$), underpin the analysis in this paper. We shall work almost exclusively with $s_{3/7}$, in preference to $h_{a_1}$. 
 
The Markov partition for $s_{3/7} $ is a partition of the unit disc. It is drawn below. This partition was also used in \cite{R5}, and introduced in section 2.9 of that paper. The boundaries of the sets of the partition, within the unit disc, are the vertical lines joining $e^{\pm 2\pi ir/14}$ for $1\leq r\leq 6$. The labels of the closures of the sets of the partition are, as shown in Figure 3, $L_j$  and $R_j$ for $1\leq j\leq 3$, and the central region $C$. We shall write $D(X)$ for the closed set labelled $X$. We also denote by $D(BC)$ and $D(UC)$ the components of the complement in $D(C)$ of the central gap of $L_{3/7}$,   where $D(BC)$ is the component intersecting the lower unit circle and $D(UC)$ is the component intersecting the upper unit circle, and the labels of these sets are $BC$ and $UC$.
    \begin{figure}
\centering{\includegraphics[width=8cm]{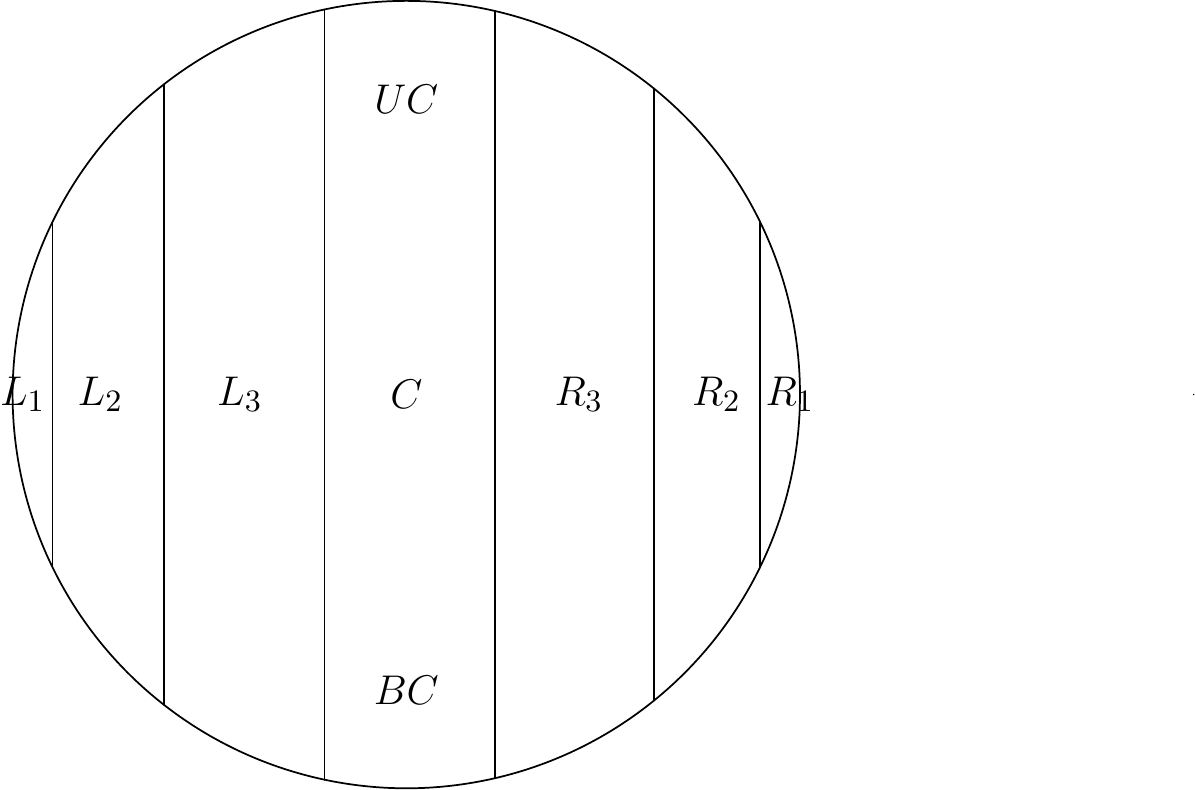}}
\caption{Markov partition for $s_{3/7}$}
\end{figure}

The dynamics under $s_{3/7}$ are as follows:
$$s_{3/7}(D(L_1))=s(D(R_1))=D(R_1)\cup D(R_2),$$
$$s_{3/7}(D(L_2))=s_{3/7}(D(R_2))=D(R_3)\cup D(C),$$
$$s_{3/7}(D(L_3))=s_{3/7}(D(R_3))=D(L_3)\cup D(L_2),\ \ D(C)=D(L_1).$$
Also, $s_{3/7}(D(BC))=s_{3/7}(D(UC))$ is the complement in $D(L_1)$ of the gap containing $s(0)$.  Note that $s_{3/7}$ is a homeomorphism onto image, restricted to $D(L_j)$ or  $D(R_j)$ or $D(BC)$ or $D(UC)$, but a degree two branched cover restricted to $D(C)$. 

As usual in symbolic dynamics, we use words over the alphabet $\{ L_j,R_j: 1\leq j\leq 3\} \cup \{ C,BC,UC\} $ to label subsets in the unit disc, so that for a word $w=X_0X_1\cdots X_{n-1}$ labels the subset $D(w)$ where
$$D(w)=\cap _{j=0}^{n-1}s_{3/7}^{-j}(D(X_j)).$$
Then, as usual, $D(w)\neq \emptyset $ if and only if $w$ is admissible, that is, $D(X_j)\cap s_{3/7}(D(X_{j-1}))\neq \emptyset $ for $0<j\leq n-1$. We will use $BC$ or $UC$ rather than $C$, except for the last letters of words, but do not allow  $BC$ or $UC$ for a last letter, using  $C$ instead. Under this convention, we have, for all admissible words $X_0\cdots X_{n-1}$,
$$s_{3/7}(D(X_0\cdots X_{n-1}))=D(X_1\cdots X_{n-1}).$$

If $w$ is an admissible word of  length $m>1$ (or $m>2$)  ending in $L_2C$ (or  $R_1R_2C$), then $D(w)$ contains a unique point of preperiod $m-1$ (or $m-2$) under $s$. We denote this point by $p(w)$, and in this circumstance we will sometimes talk of $w$ having preperiod $m-1$ (or $m-2$).

 \section{Paths in $R_{m,0}$}\label{3.2}
 
All paths in $R_{m,0}$ have endpoints in $U^0\cap Z_m$, where $U^0$ is the set defined by
$$U^0=\cup _{i=0}^\infty D(L_3^{3i+1}L_2)\cup D(BC)\cup _{i=0}^{\infty }D(L_3^{2+2i}L_2BC)\setminus \cup _{i=0}^{\infty }D(L_3^{1+3i}L_2BC).$$
Not all points in $U^0\cap Z_m$  are endpoints of paths in $R_{m,0}$, and some points are endpoints of two such paths in $R_{m,0}$, but no point in $U^0\cap Z_m$ is the endpoint of more than one path in $R_{m,0}$, and each path in $R_{m,0}$ is uniquely determined up to isotopy preserving $Z_m$ by its first $S^1$-crossing point and its final endpoint. 

The description of the paths uses finite sequences of words 
$$(w_i(w,0):1\leq i\leq n(w)),\ \ (w_i'(w,0):1\leq i\leq n(w)),$$
for each $w$ ending in $C$ with $D(w)\subset U^0$and
$$(w_i'(w,x):1\leq i\leq n(w)).$$
each $w$ ending in $C$ for $D(w)\subset U^x$, for a certain set $U^x$ which will be defined in \ref{3.7}.  
  Since we are only concerned in detail with $R_{m,0}$, and not with $R_{m,p}$ for $p>0$, we write 
$$w_i(w,0)=w_i(w),\ \ w_i'(w,0)=w_i'(w).$$

For $\beta \in R_{m,0}$, we write $w(\beta )$ for the word such that the final endpoint of $\beta $ is $p(w(\beta ))$, and also write 
$$w_i(\beta )=w_i(w(\beta )),\ \ w_1'(\beta )=w_1'(w(\beta ),\ \ w_i'(\beta )=w_i'(w(\beta ),w_1'(\beta )).$$

The $i$'th unit-disc crossing of $\beta $ is determined by $w_i'(\beta )$ as follows. If the last letter of $w_i'(\beta )$ is $L_3$ or $R_3$, then the $i$'th crossing of $\beta $ is along the set $D(w_i'(\beta )(L_2R_3L_3)^\infty )$, which is a leaf of $L_{3/7}$. To determine the direction of the crossing, we use the unique gap $G$ of $L_{3/7}$ which has this leaf in its boundary. The direction of the crossing is in an anticlockwise direction along the boundary of $G$. If $i=n(x)$ then $\beta $ crosses into $G$ and ends at the point of $G\cap Z_m$, starting at what would be the start point of a crossing of the leaf which is part of an anticlockwise circuit round $G$.

Now suppose that the last letter of $w_{i}'(\beta )$ is $BC$. In this case, $i<n(x)$.  Let $u'$ be the word obtained from $w_{i}'$ by replacing the last letter by $UC$. Then $i>1$, and  $u'$ and $w_{i}'$ both have $w_{i-1}'$ as a prefix, while $u'$ is a prefix of $w_{j}'$ for $i<j\leq n(w(\beta ))$.
The words $u'$ and $w_{i}'$ then label disjoint 
subsets $\partial D(u')$ and $\partial D(w_{i}')$ of $\{ z:\vert 
z\vert \leq 1\} $. 
The $i$'th crossing by $\beta 
(w,x,0)$ of $\{ z:\vert z\vert \leq 1\} $ is then from the first 
point of 
$\partial D(w_{i}')$ to the first point of  $\partial D(u')$, 
using anticlockwise order on these intervals.

\section{Properties of the sequences $w_i$, $w_i'$, $w_i'(.,x)$}\label{3.3}
The following properties hold for the sequences. Property 4 allows $w_i(w)$ and $w_i'(w)$ to be defined by induction on $i$, because there is always at least one occurrence of $L_3$ or $BC$ or $UC$ after the last letter of $w_{i-1}(w)$ and at or before the last letter of $w_i(w)$. Similarly, Properties 5 and 6 will allow for an inductive definition of $w_i'(w,x)$ if $D(w)\subset U^x$.

\begin{itemize}
\item[1.]   $w_i(w)$ is a prefix of $w$ of length increasing with $i$, and is also a prefix of $w_{i+1}'(w)$ if $i<n(w)$.
\item[2.] $w_i'(w)$ always ends in $L_3$ or $R_3$ or $BC$,  and is a function only of $i$ and $w_i(w)$.
\item [3.] $w_1'(w,x)=x$ if $D(w)\subset U^x$
\item[4.] Let $w=uy$ and $D(y)\subset U^0$. If $|u|<w_1(w)$, then $w_1(uy)=uw_1(y)$ and $w_1'(uy)=uw_1'(y)$. If $w=uy$ and  $|w_{i-1}(w)|<|u|\leq |w_i(w)|$, then $w_j(w)=uw_{j-i}(y)$ and $w_j'(uy)=uw_{j-i}'(y)$ for $i<j\leq n(w)$. If $w=uUCy_1$ and $w_{i-1}(w)|<|uUC|\leq |uUC|$, then $w_j(w)=uUCy_2$ where $w_{j-i}(BCy)=BCy_2$, for $i<j\leq n(x)$.
\item[5.] If $w=uy$, where $u=u_1u_2$ and $D(v_1u_2)\subset U^x$, with 
$$|w_{i-1}(w)|<|u_1|\leq |w_i(w)|$$ and 
$$|w_1(v_1u_2y)|\leq |v_1|<|w_2(v_1u_2y)|\leq |v_1u_2|,$$
 then
$$w_j'(w,x)=uv_2$$
for $i<j\leq n(x)$, where
$$v_1v_2=w_{j-i+1}'(v_1u_2y,x).$$
\item[6.] If $w_i(w)$ ends in $BC$ or $UC$, then $w_i'(w,x)=w_i'(w)$.
\end{itemize}

Thus $w_i$ and $w_i'$ and $w_i'(.,x)$ can be defined inductively from $w_1$ and $w_1'$ and $w_1'(.,x)$. The definitions are used explicitly in some sections of this article. We therefore include them here. For full details, see Chapter 7 of \cite{R5}.

\section{Definition of $w_1$}\label{3.4}

We define $w_{1}(w)$ for $w$ with $D(w)\subset U^{0}$ and with $w$ ending in $C$. The definition also works for infinite words.  The first letter of $w$ is 
either $BC$ or $L_{3}$. 
If $BC$ is the first letter, then $w_{1}(w)=BC$. 
If $L_{3}$ is the first letter of $w$, and $w$ starts with $L_{3}^{4}$ or $L_{3}^{2}L_{2}BC$, then $w_{1}(w)=L_{3}^{4}$ or $L_{3}^{2}$ respectively. Otherwise, $w$ starts $L_{3}L_{2}$, and 
we look for the first occurrence of one of the following. One of 
these must occur if $w$ is preperiodic, that is, ends in $C$. 

\begin{itemize}
\item[0.] An occurrence of $C$.
\item[1.] An occurrence 
of $xBC$ for any $n\geq 0$, where $x$ is a maximal word in the letters 
$L_{3}$, $L_{2}$, $R_{3}$ with an even number of $L$ letters. 
\item[2.] An occurrence 
of $xUC$, where $x$ is a maximal word in the letters 
$L_{3}$, $L_{2}$, $R_{3}$, with an odd number of $L$ letters. 
\item[3.] An occurrence of $L_{1}R_{2}UC$. 
\item[4.] An occurrence of  $R_{1}R_{2}BC$.  
\item[5.] An occurrence of  a string $v_{1}\cdots v_{n}L_{3}$ which is the end of $w$, or is followed in $w$ by $L_{2}$, where:
\begin{itemize} 
\item $v_{i}=(L_3L_2R_3)^{q_i}L_{3}(L_{2}R_{3})^{m_{i}}L_{3}^{r_{i}-1}$ for $i\geq 2$, where $m_{i}+r_{i}=3$ and $q_i\geq 0$; 
\item either $v_{1}=(L_3L_2R_3)^{q_1}L_{3}(L_{2}R_{3})^{m_{1}}L_{3}^{r_{1}-1}$ for $q_1\geq 0$ and  $m_{1}\geq 1$ and  an odd $m_{1}+r_{1}\geq 3$, or $v_{1}=XL_{1}R_{1}^{p_{1}}R_{2}R_{3}(L_{2}R_{3})^{m_{1}}L_{3}^{r_{1}-1}$ for $3\leq r_{1}\leq 5$, and $m_{1}\geq 0$, and $p_{1}\geq 0$, and $X=BC$ or $UC$;
\item $v_{n}L_{3}$ is either followed by $L_{2}$, or is at the end of $w$; 
\item the string $v_{1}\cdots v_{n}L_{3}$ is maximal with these properties; 
\item $n$ is odd.
\end{itemize}

\end{itemize}

We note that the term $(L_3L_2R_3)^{q_i}$ with $q_i\geq 0$. was missing in the definition of \cite{R5}. This is relevant to the calculations in this paper, because we shall have examples with $q_i=1$ for some $i$.

We define $w_{1}(w)$ to be the prefix 
of $w$ ending in the occurrence listed  -- except in case 1 when $n>0$, and in case 4 when $w$ ends in $L_{1}R_{1}R_{2}BC$, and in case 5.  In case 4, we define $w_{1}(w)$ to be the word up to, and not including, this occurrence of $L_{1}R_{1}R_{2}BC$, and $w_{2}(w,0)$ to be the prefix of $w$ ending in this occurrence of $L_{1}R_{1}R_{2}BC$. In case 5, if the maximal string of words of this type is $v_{1}\cdots v_{n}$, then we define $w_{1}(w,0)$ to be the prefix ending with $v_{n}L_{3}$.

\section{Definition of $w_2(w)$ if $D(w)\subset D(BC)$}\label{3.5}
Now suppose that the first letter of $w$ is $BC$. Then $w$ starts with $(BCL_{1}R_{2})^{k}$ for a maximal $k\geq 0$. 
Then we define 
$w_{2}(w)=w$ if there is no occurrence of cases 1 to 5 above. 
We 
define $w_{2}(w)$ to be the prefix of $w$ defined exactly  as 
$w_{1}(w)$ is 
defined when $w$ starts with $L_{3}$, if there is an occurrence of one of 2 to 5.

\section{Definition of $w_1'$}\label{3.6}
Now we define $w_{1}'(w,0)=w_{1}'(w_{1},0)$. If $w=C$ or if $w_{1}$ starts with $BC$ or $L_{3}^{2}$, then we define $w_{1}'(w,0)=L_{3}$. Now suppose that $w$ starts with $L_{3}L_{2}$, and we assume that one of cases 
0 to 5 holds, which, as pointed out, is always true if $w$ ends in 
$C$. 

\begin{itemize}
\item[0.] $w_{1}'$ is the same length as $w_{1}$, ending in 
    $L_{3}$, if preceded by a maximal subword in the letters $L_{3}$, $L_{2}$ or $R_{3}$ and with an even number of $L$ letters, or ending in $R_{3}$ if preceded by a maximal subword in the letters $L_{3}$, $L_{2}$ or $R_{3}$ and with an even number of $L$ letters, such that $D(w_{1}')$ and $D(w_{1})$ are adjacent.
    \item[1.]  Let $w_{1}''$ be obtained from $w_{1}$ by deleting the last $(UCL_{1}R_{2})^{n}BC$, and adding $C$. Then $w_{1}'$ is the same length as $w_{1}''$, ending in $L_{3}$, such that $D(w_{1}')$ and $D(w_{1}'')$ are adjacent. If $n>0$ then, $w_{2}'$ is the same length as $w_{1}''$, ending in $R_{3}$, such that $D(w_{2}')$ and $D(w_{1}'')$ are adjacent.
    \item[2.] Let $w_{1}''$ be obtained from $w_{1}$ by deleting the last $UC$ and adding $C$. Then  $w_{1}'$ is defined similarly to 1, but ending in 
    $R_{3}$.
    \item[3.] $w_{1}'$ is obtained 
from $w_{1}$ by replacing the last letter by $BC$.
\item[4.] If $w_{1}$ ends in $R_{1}^{2}R_{2}BC$, then $w_{1}'$ is obtained 
from $w_{1}$ by replacing the last letter by $R_{3}L_{2}R_{3}L_{3}$.  If $w_{1}$ is extended in $w$ by $L_{1}R_{1}R_{2}BC$, then $w_{1}'(w,0)$ is the word of the same length as $w_{1}$, with $D(w_{1}')$ to the left of $D(v')$, where $v'$ is obtained from $w_{1}$ by replacing the last letter by $C$. We then also define $w_{2}'(w,0)=v''$, where $v''$ is of the same length as $v'$, with $D(v')$ between $D(w_{1}')$ and $D(v'')=D(w_{2}')$.
\item[5.] The word $w_{1}'$ is obtained 
from $w_{1}$ by deleting all but the first letter of $v_{n}L_{3}$ if $n>1$, or if $n=1$, $m_{1}+r_{1}=3$ and the first letter of $v_{1}$ is $L_{3}$. In the other cases with $n=1$, the word  $w_{1}'$ is obtained from $w_{1}$ by replacing the last $L_{3}^{3}$ by $L_{2}R_{3}L_{3}$ if $r_{1}\geq 3$, by deleting the last $L_{2}R_{3}L_{3}^{2}$ if $r_{1}=2$, and by deleting the last $(L_{2}R_{3})^{2}L_{3}$ if $r_{1}=1$.  

\end{itemize}

\section{Definition of $w_1'(.,x)$}\label{3.7}
This definition generalises the definition given in \cite{R5}, where the only possible value of $w_1'(.,x)$ was $x$, and this was defined only for $w$ with $D(w)\subset U^x$, for a proper subset $U^x$ of $U^0$. 

For $i\geq 1$ we define $w_i'(w,x)=w_i'(w)$ except if $w=uy_1v$ or $uy_2v$, where $w_i(w)=uy_1L_3$ or $uy_2L_3$ and  $D(v)\subset D(L_3L_2)$  is to the right of $D(x(L_2R_3L_3)^\infty )$ and 
$$(y_1,y_2)=(L_3(L_2R_3)^{2j}L_3^4,L_3(L_2R_3)^{2j+2})$$
for some $j\geq 0$, or
$$(y_1,y_2)=(L_3(L_2R_3)^{2j}L_3^2,L_3(L_2R_3)^{2j-1}L_3^2L_2R_3)$$
for some $j\geq 1$. We then define
$$w_i'(uy_1v,x)=w_i'(uy_2v)=w_i'(uy_2L_3L_2),\ \ \ w_i'(uy_2v,x)=w_i'(uy_1v)=w_i'(uy_1L_3L_2).$$

We then define $U^x$ to be the union of all $D(w)$ with $w_1'(w,x)=x$. This agrees with the definition given in \cite{R5}.

The definition of $w_i'(.,x)$ might seem as if it comes from nowhere. It is actually related to certain homeomorphisms which determine the group action on the universal cover of $V_{3,m}$, and hence suggest possible fundamental domains for this action. Homeomorphisms of a similar type also play an important role in the current work, since each capture path is equivalent under the group action to precisely two vertices of the fundamental domain. In both cases, the homeomorphisms act on both paths and on $Z_m$, and, in both cases, {\em{exchangeable pairs}} are important indefining the homeomorphisms.

\section{}\label{3.8}

For the proofs of the main theorems of this paper, we need to analyse the homeomorphisms $\psi  _{4j-3,4j-1,\beta }=\psi _{4j-3,\beta }\circ \psi _{4j-1,\beta }^{-1}$, which are defined in \ref{2.8} and \ref{2.9}. These definitions assume the existence of $\beta _\ell (\beta )$ for $\ell \le 4j$ with the properties claimed in \ref{2.8}. The proof of the existence of the $\beta _\ell (\beta )$ will be inductive, and the existence of $\beta _{4j+t}(\beta )$ for $1\le t\le 4$ will involve an analysis of the homeomorphisms $\psi _{4k-3,4k-1,\beta }$ for $k\le j$.

More generally, we consider homeomorphisms $\psi _{\zeta ,\eta }=\psi _{\zeta }\circ \psi _{\eta }^{-1}$ whenever $(\zeta ,\zeta _{2})$ and $(\eta ,\eta _{2})$ are adjacent pairs in $R_{m,0}$, with $\zeta <\zeta _{2}$ and $\eta <\eta _{2}$, and $\psi _{\zeta }$ matches $(\zeta ,\zeta _{2})$ with the adjacent pair $(\zeta ',\zeta _{2}')$ in $R_{m,0}'$, while $\psi _{\eta }$ matches $(\eta ,\eta _{2})$ with $(\eta ',\eta _{2}')$. This means that there are $[\alpha _\zeta ]$ and $[\alpha _\eta ]\in \pi _1(\Cbar \setminus Z_m,\infty )$ such that 
$$s\simeq _{\psi _\zeta }\sigma _{\alpha _\zeta }\circ s$$
$$[\alpha _\zeta *\psi _\zeta (\zeta ')]=[\zeta ],\ \ [\alpha _\zeta *\psi _\zeta (\zeta _2')]=[\zeta _2]$$
and similarly for $\eta $. 

We start with a very common lemma which has been used many times in the course of preceding work to compute the homeomorphism $\psi $ such that $(s,Y_m)\simeq _\psi (\sigma \alpha \circ s,Y_m)$, and, in particular, to bound the support of $\psi $.
\begin{lemma}\label{3.8.1} Suppose that $\alpha $ is a closed lop based at $v_2=\infty $ which is trivial in $\pi _1(\Cbar \setminus Y_0,v_2)$. and suppose that $\psi $ is a homeomorphism such that $(s,Y_m)\simeq _\psi (\sigma _\alpha \circ s,Y_m)$. Then $\psi $ can be defined inductively by the equations: $\psi _1$ is the identity in ${\rm{Mod}}(\overline{\mathbb C},Y_{0})$, $\psi =\psi _m$, and 
$$s=\sigma _\alpha \circ s\circ \psi _1$$
and 
$$\psi _i=\psi _{i+1}{\rm{\ in\ }}{\rm{Mod}}(\overline{\mathbb C},Y_{i})$$
 and 
$$\psi_i\circ s=\sigma _{\alpha }\circ s\circ \psi _{i+1}$$
Alternatively we can define $\psi _{i+1}$ by
$$\psi _{i+1}=\xi _i\circ \psi _i$$
where 
$$\xi _0=\psi _1$$
 and  
 $$\xi _i={\rm{identity\  in\ }}{\rm{Mod}}(\overline{\mathbb C},Y_{i})$$ and
$$\sigma _\alpha \circ s\circ \xi _{i+1}=\xi _i\circ \sigma _\alpha \circ s.$$

\end{lemma}
\begin{proof}
It is clear that if $\psi _i$ are defined by the first set of equations, and $\psi =\psi _m$, then $\psi $ has the required properties. So it remains to show that the alternative definitions are equivalent to the original ones. If we have 
\begin{equation}\label{3.8.1.1}\psi_i\circ s=\sigma _{\alpha }\circ s\circ \psi _{i+1}\end{equation}
and 
\begin{equation}\label{3.8.1.2}\psi_{i+1}\circ s=\sigma _{\alpha }\circ s\circ \psi _{i+2},\end{equation}
then writing $\psi _{i+1}=\xi _i\circ \psi _i$ and $ \psi _{i+2}=\xi _{i+1}\circ \psi _{i+1}$ in the second equation we have 
$$\xi _i\circ \psi _i\circ s=\sigma _{\alpha }\circ s\circ \xi _{i+1}\circ \psi _{i+1}$$
Combining these, we obtain
$$\xi _i\circ \sigma _{\alpha }\circ s\circ \psi _{i+1}=\sigma _{\alpha }\circ s\circ \xi _{i+1}\circ \psi _{i+1}$$
which gives
\begin{equation}\label{3.8.1.3}\xi _i\circ \sigma _\alpha \circ s=\sigma _{\alpha }\circ s\circ \xi _{i+1}\end{equation}
Conversely, we see that (\ref{3.8.1.3}) implies (\ref{3.8.1.2}) given (\ref{3.8.1.1}), as required. 
\end{proof}

This can be used to bound the support of $\psi _{1,4j-1}$, but, as it turns out, a variant of \ref{3.8.1} is  more useful. We recall from \ref{2.8} that
$$s\simeq _\psi \sigma _\alpha \circ s$$
where $\alpha =\alpha _{1,4j -1}$ where
$$\alpha _{1,4j-1}=\alpha _{1,3}*\psi _{1,3}(\alpha _{5,7})*\cdots *\psi _{1,4j -5})(\alpha _{4j-3,4j-1}).$$
Writing
$$\alpha _{4j-3,4j-1,1}=\psi _{1,4j-5}(\alpha _{4j-3,4j-1}),$$
we have
$$\alpha _{1,4j -1}=\alpha _{1,3}*\alpha _{5,7,1}*\cdots *\alpha _{4j -3,4j-1,1}.$$
We then write $E_{4j-3,4j-1,1}$ for the topological disc enclosed by the closed loop $\alpha _{1,4j-3,4j-1}$. We are only interested in $E_{4j-3,4j-1,1}$ up to homotopy preserving $Z_\infty $. We have
$$\sigma _\alpha =\sigma _{\alpha _{4j-3,4j-1,1}}\circ \cdots \circ \sigma _{\alpha _{1,3}}.$$
Hence if we define $\psi _{1,4j-1,1}$ to be homotopic to the identity relative to $Z_0$, and to satisfy the equation
$$s=\sigma _\alpha \circ s\circ \psi _{1,4j-1,1}$$ then we have
$$\psi _{1,4j-1,1}=\psi _{1,3,1,1}\circ \cdots \circ \psi _{4j-3,4j-1,1,1}$$
where 
$$s=\sigma _{\alpha _{4i-3,4i-1}}\circ s \circ \psi _{4i-3,4i-1,1,1}.$$
the loop $\alpha _{4i-3,4i-1}$ is a perturbation of $\alpha _{4i-3}*\psi _{4i-3,4i-1}(\overline{\alpha _{4i-1}})$.  This is  awkward to work with, simply because the effect of $\psi _{4i-3,4i-1}$ on $\alpha _{4i-1}$ is unclear. We will therefore employ a slightly different approach.

We write $E_{\zeta ,\eta }$ for an arbitrarily  small neighbourhood of the closed disc bounded by a loop homotopic to  $\alpha _{\zeta }*\overline{\alpha _{\eta }}$, where the homotopy preserving $Z_{\infty }$ and the disc does not contain $v_2=c_2=\infty $. Naturally, we wish to apply this in the case $(\zeta ,\eta )=(\beta _{4i-3}(\beta ),\beta _{4i-1}(\beta ))$ for $\beta \in {\cal{Z}}_m(3/7,+,+,0)$, once the existence of $\beta _\ell (\beta )$ as in Theorem \ref{2.8} has been established for $\ell \le 4i$. 
We shall therefore be particularly interested in the following in the case $(\zeta ,\eta )=(\beta _{4i-3}(\beta ),\beta _{4i-1}(\beta ))$,  for varying $i$.

\begin{lemma} \label{3.9} If $k$ is the smallest integer for which $\alpha _{\zeta ,\eta }=\alpha _{\zeta }*\overline{\alpha _{\eta }}$ is nontrivial in $\pi _{1}(\overline{\mathbb C}\setminus Z_{k},v_{2})$, then $\psi _\zeta \circ \psi _\eta ^{-1}=\psi _{\zeta ,\eta }=\psi _{\zeta ,\eta ,m+1}$ can be defined inductively as an element of ${\rm{Mod}}(\overline{\mathbb C},Y_{m+1})$ by the equations
\begin{equation}\label{3.9.1}\sigma _{\alpha _{\zeta ,\eta }}\circ s\circ \psi _{\zeta ,\eta ,k+1}=s,\end{equation}
and for $k<\ell \leq m$
\begin{equation}\label{3.9.2}\psi _{\zeta ,\eta ,\ell +1}=\xi _{\zeta ,\eta ,\ell }\circ \psi _{\zeta ,\eta ,\ell }\end{equation}
\begin{equation}\label{3.9.3}\sigma _{\alpha _{\zeta }}\circ s\circ \xi _{\zeta ,\eta ,\ell }=\xi _{\zeta ,\eta ,\ell -1}\circ \sigma _{\alpha _{\zeta} }\circ s.\end{equation}
It follows that the support of $\psi _{\zeta ,\eta ,\ell }$, up to isotopy preserving $Z_{m}$, is 
$$\cup _{1\leq i\leq \ell -r}(\sigma _{\alpha _\zeta }\circ s)^{-i}(E_{\zeta ,\eta }).$$
 
 \end{lemma}

\noindent{\em{Proof.}} We have 
\begin{equation}\label{3.9.4}\sigma _{\alpha _{\zeta }}\circ s\circ \psi _{\zeta }=\psi _{\zeta }\circ s{\rm{\ rel\ }}Y_{m+1},\end{equation}
\begin{equation}\label{3.9.5}\sigma _{\alpha _{\eta }}\circ s\circ \psi _{\eta }=\psi _{\eta }\circ s{\rm{\ rel\ }}Y_{m+1}.\end{equation}
Here, $[\psi _\zeta ]$ and $[\psi _\eta ]$ are the only elements of ${\rm{Mod}}(\Cbar ,Y_{m+1})$ satisfying (\ref{3.9.4}) and (\ref{3.9.5}) respectively. Now (\ref{3.9.5}) gives
\begin{equation}\label{3.9.6}\psi _{\eta }^{-1}\circ \sigma _{\alpha _{\eta }}\circ s=s\circ \psi _{\eta }^{-1}{\rm{\ rel\ }}Y_{m+1}.\end{equation}
Then combining (\ref{3.9.4}) and (\ref{3.9.6}), we have
\begin{equation}\label{3.9.7}\sigma _{\alpha _{\zeta} }\circ s\circ \psi _{\zeta ,\eta }=\psi _{\zeta }\circ s\circ \psi _{\eta }^{-1}=\psi _{\zeta ,\eta }\circ \sigma _{\alpha _{\eta} }\circ s{\rm{\ rel\ }}Y_{m+1}.\end{equation}
If this equation is satisfied by $[\chi ]$ replacing $[\psi _{\zeta ,\eta }]$, then (\ref{3.9.4}) is satisfied with  $[\chi \circ \psi _{\eta }]$ replacing $[\psi _{\zeta }]$, which implies that $[\chi ]=[\psi _{\zeta ,\eta }]$. So (\ref{3.9.6}) uniquely determines $[\psi _{\zeta ,\eta }]\in {\rm{MG}}(\Cbar ,Y_{m+1})$. 

Now we define $\psi _{\zeta ,\eta ,k+1}$ by (\ref{3.9.1}) and $[\psi _{k+1}]=[{\rm{identity}}]$ in ${\rm{MG}}(\Cbar ,Y_k))$, and  define $\psi _{\eta ,\zeta ,k}$ to be the identity, we have 
$$\xi _{\zeta ,\eta ,k}=\psi _{\zeta ,\eta ,k+1}.$$
Then we clearly want to define $\psi _{\zeta ,\eta ,\ell +1}$ in terms of $\psi _{\zeta ,\eta ,\ell }$ by $[\psi _{\zeta ,\eta ,\ell +1}]=[\psi _{\zeta ,\eta ,\ell }]$ in ${\rm{MG}}(\Cbar ,Y_\ell )$ and:
\begin{equation}\label{3.9.8}\sigma _{\psi_{\zeta ,\eta ,\ell }(\alpha _\eta )}^{-1}\circ \sigma _{\alpha _{\zeta }}\circ s\circ \psi _{\zeta ,\eta ,\ell +1}=\psi _{\zeta ,\eta ,\ell }\circ s,\end{equation}
If we do this for $\ell \leq m+1$, then we will have $[\psi _{\zeta ,\eta ,m +2}]=[\psi _{\zeta ,\eta ,m+1}]$ in ${\rm{MG}}(\Cbar ,Y_{m+1})$ and (\ref{3.9.6}) is satisfied by $\psi _{\zeta ,\eta ,m+1}$ replacing $\psi _{\zeta ,\eta }$, which, as we have seen, implies that $[\psi _{\zeta ,\eta ,m+1}]=[\psi _{\zeta ,\eta }]$
Now (\ref{3.9.8}) is equivalent to:
\begin{equation}\label{3.9.9}\sigma _{\psi_{\zeta ,\eta ,\ell }(\alpha _{\eta })}^{-1}\circ \sigma _{\alpha _{\zeta }}\circ s\circ \xi _{\zeta ,\eta ,\ell }\circ \psi _{\zeta ,\eta ,\ell }=\xi _{\zeta ,\eta ,\ell -1}\circ \psi _{\zeta ,\eta ,\ell -1}\circ s.\end{equation}
But then this gives, by (\ref{3.9.8}) with $\ell -1$ replacing $\ell $, 
$$\sigma _{\psi_{\zeta ,\eta ,\ell }(\alpha _{\eta })}^{-1}\circ \sigma _{\alpha _{\zeta }}\circ s\circ \xi _{\zeta ,\eta ,\ell }\circ \psi _{\zeta ,\eta ,\ell}=\xi _{\zeta ,\eta ,\ell -1}\circ \sigma _{\psi_{\zeta ,\eta ,\ell -1}(\alpha _{\eta })}^{-1}\circ \sigma _{\alpha _{\zeta }}\circ s\circ \psi _{\zeta ,\eta ,\ell }$$
Then this gives
$$\sigma _{\psi_{\zeta ,\eta ,\ell }(\alpha _{\eta })}^{-1}\circ \sigma _{\alpha _{\zeta }}\circ s\circ \xi _{\zeta ,\eta ,\ell }\circ \psi _{\zeta ,\eta ,\ell}=\sigma _{\psi_{\zeta ,\eta ,\ell}(\alpha _{\eta })}^{-1}\circ \xi _{\zeta ,\eta ,\ell -1}\circ \sigma _{\alpha _{\zeta }}\circ s\circ \psi _{\zeta ,\eta ,\ell },$$
which gives
$$\sigma _{\alpha _{\zeta }}\circ s\circ \xi _{\zeta ,\eta ,\ell }=\xi _{\zeta ,\eta ,\ell -1}\circ \sigma _{\alpha _{\zeta }}\circ s$$
which gives (\ref{3.9.3}), from which we see that the support of $\xi _{\zeta ,\eta ,\ell }$ is 
$$(\sigma _{\alpha _{\zeta }}\circ s)^{k-\ell}(E_{\zeta ,\eta }),$$
 and hence we also have the required expression for the support of $\psi _{\zeta ,\eta ,\ell}$.

\Box

\section{Basic exchangeable  pairs}\label{3.10}

From \ref{3.9}, we see that $\psi _{\zeta ,\eta }$ is a composition of exchanges between sets. Restricted to $Z_m$, each of these exchanges is just an involution. If $\alpha _{\zeta }*\overline{\alpha _{\eta }}$ is a simple closed loop up to homotopy, then $\psi _{\zeta ,\eta }$ is a composition of disc exchanges, first between the two halves of $E_{\zeta ,\eta }$ and then between two discs in each component of $(\sigma _{\alpha _{\zeta }}\circ s)^{-i}(E_{\zeta ,\eta })$ for increasing $i$.

An {\em{exchangeable pair (for $\zeta $)}} is a notation for keeping track of inverse images of sets $E=E_{\zeta ,\eta }$ under  $(\sigma _{\alpha_\zeta }\circ s)^{n}$, and hence for describing the support of $\psi _{\zeta ,\eta }$. We shall see later that a component of
$(\sigma _{\alpha _{\zeta }}\circ s)^{-n}(E)$ coincides with a component of 
$$s^{-n}(E)\cup (\sigma _\zeta \circ s)^{1-n}(s^{-1}\zeta ).$$
In all cases that we shall consider, we shall have $E\subset  D(L_3L_2)$ or $E\subset D(BC)$, and the intersection of $\zeta $ with the unit disc will be with $D(L_3L_2)$, or the union of $D(BC)$ and the single leaf of $L_{3/7}$ joining $e^{\pm 2\pi i(2/7)}$ respectively.

It follows any component  $C$  is coded by an {\em{ exchangeable pair for $\zeta $}} which is denoted by  $v_1\leftrightarrow v_2$, possibly with the words ``top'' or ``bottom'' attached to $v_1$ and $v_2$, where $v_1E$ and $v_2E$ are the components of $s^{-n}E$ within $C$, and ``top'' or ``bottom'' indicates that the component of $(\sigma _{\zeta }\circ s)^{1-n}(s^{-1}(\zeta ))$ attached to ($v_1E$ or $v_2E$) at a point on the upper or lower unit circle.  We shall write $C=C(v_1\leftrightarrow v_2,E)$ or $C(v_1\leftrightarrow v_2)$ if it is clear what $E$ is. Strictly speaking, we should include the path $\zeta $ in the notation, but in the cases which we consider,  $\zeta $ is usually a path in $R_{m,0}$ for which $w_1'(\zeta )$ is known, and the words $w_j'(\zeta )$, which encode the unit-disc crossings of $\zeta $, coincide with $w_j'(v_1)=w_j'(v_2)$, and hence $\zeta $ can be computed from $v_1$ (or $v_2$). At any rate, $v_1$ and $v_2$ are admissible of words of the same length $n$ in the letters 
$\{ L_j,R_j,BC,UC:1\leq j\leq 3\} $. Of course $C(v_1\leftrightarrow v_2)=C(v_2\leftrightarrow v_1)$, and so $v_1\leftrightarrow v_2$ and $v_2\leftrightarrow v_1$ are naturally equivalent.  If $C=C(v_1\leftrightarrow v_2)$ does not intersect $\zeta $, then $(\sigma _\zeta \circ s)^{-1}(C)=s^{-1}(C)$. If in addition the left-most letters of $v_1$ and $v_2$ are the same, then the components of $(\sigma _\zeta \circ s)^{-1}(C)=s^{-1}(C)$ are $C(u_1\leftrightarrow u_1v_2)$ and $C(u_2v_1\leftrightarrow u_2v_2)$, where $u_1$ and $u_2$ are the letters such that $u_1v_j$ and $u_2v_j$ are admissible for one, and hence both, $j$. 

The {\em{length}} of an exchangeable pair $v_1\leftrightarrow v_2$ is the length $|v_1|$ of $v_1$ --- which is equal to the length $|v_2|$ of $v_2$. 

{\em{Basic exchangeable pairs}} are defined inductively. The only basic exchangeable pairs of length $1$ are $L_3\leftrightarrow R_3$ and $L_2\leftrightarrow R_2$. We define $C=C(L_i\leftrightarrow R_i)$ to be the set obtained by joining $D(L_i)$ and $D(R_i)$ by an arc along the unit circle  which joins a point on $\partial D(L_i)$ in the lower unit circle to a point in $\partial D(R_i)$ in the upper unit circle. Now let $n\ge 2$. Then
$v_1\leftrightarrow v_2$ is a basic exchange of length $n$ if and only if  a component of $s^{-n}(C)$ is the union of sets $D(v_1)$ and $D(v_2)$ and an arc joining them. Thus if $v_1\leftrightarrow v_2$ is a basic exchange of length $n$, and $v_1'$ and $v_2'$ are the suffixes of $v_1$ and $v_2$ of length $m$ for some $m<n$, then $v_1'\leftrightarrow v_2'$ is a basic exchange of length $m$. In particular the last letters of $v_1$ and $v_2$ are $L_i$ and $R_i$ for $i=2$ or $3$. If $i=3$, we say that $v_1\leftrightarrow v_2$ is {\em{ a basic exchange for $L_3$}} --- because $s^{-1}(D(L_3))=D(L_3^2)\cup D(R_3L_3)$ --- and if $i=2$ we say that $v_1'\leftrightarrow v_2'$ is {\em{a basic exchange for $BC$}} --- because $s^{-1}(D(BC))=D(L_2BC)\cup D(R_2BC)$. 

It follows that a basic exchange $v_1\leftrightarrow v_2$ is an exchange for $\zeta $ if, for any $m$ with $1<m<n$, the arc in $C(v_1'\leftrightarrow v_2')$ does not intersect $\zeta $. It turns out that a basic exchange  $v_1\leftrightarrow v_2$ is automatically an exchange for $\zeta \in R_{m,0}$ if all the unit disc crossings of $\zeta \in R_{m,0}$ are in $D(L_3)$, and either $v_1$ and $v_2$ have no common letters in the the same positions, or there are common letters, and, among  the common letters, if there is an occurrence of $L_3$, it is the first letter of each word.  This becomes apparent from computing all the basic exchanges $v_1\leftrightarrow v_2$. We now do this.
This calculation was also carried out in \cite{R5}. If $v_1\leftrightarrow v_2$ is a basic exchange, and the first letters of $v_1$ and $v_2$ are the same, or can admissibly be preceded by the same letter, then $uv_1\leftrightarrow uv_2$ is a basic exchange for any word $u$ such that $uv_1$ is admissible. Therefore we only need to compute those basic exchanges $v_1\leftrightarrow v_2$ such that all the letters in the same positions of $v_1$ and $v_2$ are different. Also, of course, we only have to compute up to equivalence, that it, up to interchanging $v_1$ and $v_2$. 

We start with basic exchangeable pairs for $L_3$ . The unique exchangeable pair of length one is
\begin{equation}\label{3.10.1}L_{3}{\rm{(bot)}}\leftrightarrow {\rm{(top)}}R_{3}\end{equation}
Since $L_3E\subset D(L_3^2L_2)$ and $R_3E\subset D(R_3L_3L_2)$, for any $\zeta \in R_{m,0}$, the only possible intersection between $C(L_3\leftrightarrow R_3)$ and $\zeta $ is with the initial segment of $\zeta $, before the first intersection with $S^1$. Up to homotopy preserving $S^1$ and $Z_\infty $, we can assume that $C(L_3\leftrightarrow R_3)\cap \{ z:\vert z\vert >1\} $ is close to the clockwise arc of the unit circle between $R_3E$ and $L_3E$, mostly close to the lower unit circle, and not intersecting $\zeta $. The only exchanges  of length $2$ are 
\begin{equation}\label{3.10.2}L_{3}^{2}{\rm{(top)}}\leftrightarrow {\rm{(bot)}}L_{2}R_{3}\end{equation}
and
\begin{equation}\label{3.10.3} R_{3}L_{3}{\rm{(bot)}}\leftrightarrow {\rm{(top)}}R_{2}R_{3}\end{equation}
These are both basic. The only exchanges of length $3$ are
\begin{equation}\label{3.10.4}R_3L_3^2{\rm{(top)}}\leftrightarrow {\rm{(top)}}L_3L_2R_3,\end{equation}
\begin{equation}\label{3.10.5}L_3^3{\rm{(bot)}}\leftrightarrow {\rm{(bot)}}R_3L_2R_3,\end{equation}
\begin{equation}\label{3.10.6}L_2R_{3}L_{3}{\rm{(top)}}\leftrightarrow {\rm{(bot)}}L_1R_{2}R_{3},\end{equation}
\begin{equation}\label{3.10.7}R_2R_{3}L_{3}{\rm{(bot)}}\leftrightarrow {\rm{(top)}}R_1R_{2}R_{3}.\end{equation}
Of these, only the component represented by (\ref{3.9.4}) might intersect $\zeta $. If we assume, as will always be true, that $L_3L_2R_3E$ is disjoint from $\zeta $ -- in fact it is always disjoint from $E\cup \zeta $-- then this too is basic.

The following was computed in \cite{R5}, and is easily checked. The exchanges are grouped together so that those  in (\ref{3.11.1}) are derived from  (\ref{3.10.4}), those in (\ref{3.11.2})  are derived from (\ref{3.10.5}) and those in (\ref{3.11.3}) from (\ref{3.10.6}), with (\ref{3.10.7}) being the last one in the group (\ref{3.11.3}).

\begin{lemma}\label{3.11} The complete list of basic exchangeable pairs $a\leftrightarrow b$ for $\zeta $,  for which the first letters of $a$ and $b$ are not the same, but for which the letters preceding $a$ and $b$ in any inverse image has to be the same, is as follows, where $n$ is an integer $\geq 0$ . In all cases, if one of $a$ and $b$ has all letters in the set $\{L_3,L_2,R_3\} $, then $b$ has been chosen to have this property, and if both $a$ and $b$ have all letters in $\{ L_3,L_2,R_3\} $ then $D(b)$ is chosen to be to the left of $D(a)$, if the exchange is along the bottom of the unit circle, and to the right if the exchange is along the top of the unit circle. Any other basic exchanges are represented by suffixes of these. In (\ref{3.11.1}) to (\ref{3.11.3}), the connection in $C$ between $D(a)$ and $D(b)$ is between the nearest component of $\partial D(a)\cap S^1$ and $\partial D(b)\cap S^1$, while in (\ref{3.11.4}) the connection is between the furthest component of $\partial D(a)\cap S^1$ and $\partial D(b)\cap S^1$.

\begin{equation}\label{3.11.1}\begin{array}{l}
L_{2}(UCL_{1}R_{2})^{n}R_{3}L_3^2{\rm{(bot)}}\leftrightarrow {\rm{(bot)}}L_{3}(L_{3}L_{2}R_{3})^{n+1}\\
BCL_1R_2(UCL_{1}R_{2})^{n}R_{3}L_3L_{3}{\rm{(bot)}}\leftrightarrow {\rm{(bot)}}R_3L_2R_3(L_{3}L_{2}R_{3})^{n+1}\\
R_1R_2(UCL_{1}R_{2})^{n}R_{3}L_3L_{3}{\rm{(top)}}\leftrightarrow {\rm{(top)}}R_2R_3(L_{3}L_{2}R_{3})^{n+1}\\
\end{array}
\end{equation} 

\begin{equation}\label{3.11.2}\begin{array}{l}L_{2}(BCL_{1}R_{2})^{n}R_{3}L_{2}R_{3}{\rm{(top)}}\leftrightarrow {\rm{(top)}}L_{3}(L_{3}L_{2}R_{3})^{n}L_{3}^{3}, \\ 
UCL_{1}R_{2}(BCL_1R_2)^nR_{3}L_{2}R_{3} {\rm{(top)}}\leftrightarrow  {\rm{(top)}}R_{3}L_{2}R_{3}(L_3L_2R_3)^nL_{3}^{3},\\
R_1R_2(BCL_1R_2)^nR_{3}L_{2}R_{3} {\rm{(top)}}\leftrightarrow  {\rm{(top)}}R_2R_3(L_3L_2R_3)^nL_{3}^{3},\\
,\end{array}\end{equation}

\begin{equation}\label{3.11.3}\begin{array}{l}L_2(BCL_1R_{2})^{n+1}R_{3}{\rm{(top)}}\leftrightarrow {\rm{(top)}}L_3(L_3L_2R_{3})^{n+1}L_{3}\\
UCL_1R_2(BCL_1R_{2})^nR_{3}{\rm{(top)}}\leftrightarrow {\rm{(top)}}R_3L_2R_3(L_3L_2R_{3})^nL_{3}\\
R_1R_2(BCL_1R_{2})^{n+1}R_{3}{\rm{(bot)}}\leftrightarrow {\rm{(bot)}}R_2R_3(L_3L_2R_{3})^{n+1}L_{3},\\
R_2R_{3}L_{3}{\rm{(bot)}}\leftrightarrow {\rm{(top)}}R_1R_{2}R_{3}\end{array}
\end{equation}

For a word $u$ ending in $L_3$ or $R_3$, we define $i=i(u)$ to be the largest integer such that $w_i(uL_3L_2)$ is a proper prefix of $uL_3L_2$ --- defining $i(u)=0$ if $w_1(uL_3L_2)=uL_3L_2$. 
The exchanges $a\leftrightarrow b$ for which $v$ can be chosen so that $i(vbL_3L_2)=i(vaL_3L_2)=1$ and $w_1'(vbL_3L_2)\neq w_1'(vaL_3L_2)$, and both $D(va)$ and $D(vb)$ intersect the upper unit circle  in $D(L_3L_2)$, are precisely those for which $v=v_0v_1$, where  all letters of $v_1$ and $a$ and $b$ are in $\{ L_3,L_2,R_3\} $ and $v_0v_1$ starts $L_3L_2$ and $v_1$ starts $L_3$ and all subsequent letters are in $\{ L_2,R_3\} $, and $i(v_0L_3L_2)=v_0L_3L_2$.  These are derived from the first lines of (\ref{3.11.1}) and (\ref{3.11.2}) and (\ref{3.11.3}). The exchanges $v_1a\leftrightarrow v_1b$ are all of the form
\begin{equation}\label{3.11.4}\begin{array}{l}
L_3(L_2R_3)^{2t+1}L_{3}^2L_{2}R_{3}{\rm{(bot)}}\leftrightarrow {\rm{(bot)}}L_3(L_2R_3)^{2t+2}L_{3}^2,\\
L_3(L_2R_3)^{2t+2}{\rm{(bot)}}\leftrightarrow {\rm{(bot)}}L_3(L_2R_3)^{2t}L_3^4.\end{array}\end{equation}
For these exchanges we have
$$i(v_1b)=i(v_1a)= 1=i(va)=i(vb),$$
$$w_1'(vaL_3L_2)=v_0w_1'(v_1aL_3L_2),\ \ w_1'(vbL_3L_2)=v_0w_1'(v_1bL_3L_2),$$
and 
$$w_1'(v_1a)\neq w_1'(v_1b).$$
In fact 
$$w_1'(L_3(L_2R_3)^{2t+2}L_3L_2)=L_3(L_2R_3)^{2t},$$
$$w_1'(L_3(L_2R_3)^{2t}L_3^5L_2)=L_3(L_2R_3)^{2t}L_3^2,$$
$$w_1'(L_3(L_2R_3)^{2t+1}L_{3}^2L_{2}R_{3}L_3^2L_2)=L_3(L_2R_3)^{2t},$$
$$w_1'(L_3(L_2R_3)^{2t+2}L_{3}^3L_2)=L_3(L_2R_3)^{2t+3}L_3.$$

The exchanges $a\leftrightarrow b$ for which $v$ can be chosen so that $vb$ starts with $L_3L_2$ with 
$$i(vbL_3L_2)=i(vaL_3L_2)=1,\ \ w_1'(vbL_3L_2)=w_1'(vaL_3L_2),\ \ w_1(vbL_3L_2)\neq w_1(vaL_3L_2),$$
 and such that $D(vb)$ intersects the upper unit circle, and $D(va)$ does not, are obtained from the second lines of (\ref{3.11.1}) and (\ref{3.11.2}) and the first two lines of (\ref{3.11.3}). They are as follows, where $n$ is any integer $\geq 0$
 \begin{equation}\label{3.11.5}\begin{array}{l}
L_3(L_2R_3)^{2t}L_2BCL_1R_2(UCL_{1}R_{2})^{n}R_{3}L_3^2{\rm{(bot)}}\leftrightarrow {\rm{(bot)}}L_3(L_2R_3)^{2t+2}(L_{3}L_{2}R_{3})^{n+1}\\
 (BCL_{1}R_{2})^{n+1}R_{3}L_{2}R_{3}{\rm{(bot)}}\leftrightarrow {\rm{(bot)}}(L_{3}L_{2}R_{3})^{n+1}L_{3}^{3}, \\ 
(BCL_1R_{2})^{n+1}R_{3}{\rm{(top)}}\leftrightarrow {\rm{(bot)}}(L_3L_2R_{3})^{n+1}L_{3}\end{array}\end{equation}
The first of these can be prefixed by any word $v$ such that $w_1(vL_3L_2)=vL_3L_2$. The last two can be prefixed by $(v_1L_2,v_1L_3)$ or by $(v_2UCL_1R_2,v_2R_3L_2R_3)$ for any $v_1$ and $v_2$ such that $w_1(v_1L_3^{2}L_2)=v_1L_3^{2}L_2$ and $w_1(v_2R_3L_2R_3L_3L_2)=v_2R_3L_2R_3L_3L_2$.

\end{lemma}

The pairs $(a,b)$ as in (\ref{3.11.4}) are precisely those for which $w_1'(vbL_3L_2)=y_2$ and $w_1'(vaL_3L_2)=y_1$ with $y_1\neq y_2$, and such that  the set $U^{y_2}$ contains $D(vbv_2)$ if and only if  $w_1'(v_2)>w_1'(y_2)$,  and contains $D(vav_2)$ if and only if  $w_1'(v_1)\leq w_1'(y_2)$, while the set $U^{y_1}$ contains $D(vav_2)$ if and only if $w_1'(v_2)>y_1$, and contains $D(vbv_2)$ if and only if $w_1'(v_1)\leq y_1$. It follows that $D(vav_1)\subset U^{y_1}\cap U^{y_2}$ if $y_1<w_1'(v_1)\leq y_2$ and that $D(vbv_1)$ is contained in neither $U^{y_1}$ nor $U^{y_2}$ if $y_1<w_1'(v_1)\leq y_2$. In fact, $D(vbv_2)$ is not contained in any $U^y$ if $y_1<w_1'(v_2)\leq y_2$. 

\section{Exchangeable pairs for $BC$}\label{3.12}
For convenience, we list the basic exchangeable pairs for $BC$. These are obtained by changes to the later letters in the pairs above and are as follows:

\begin{equation}\label{3.12.1}\begin{array}{l}BCL_{1}R_{2}{\rm{(bot)}}\leftrightarrow {\rm{(bot)}}L_3^2L_{2},\\ 
L_{3}^{3}L_{2}{\rm{(top)}}\leftrightarrow {\rm{(top)}}L_{2}BCL_{1}R_{2}.
\end{array}
\end{equation} 

\begin{equation}\label{3.12.2}L_{3}(L_{3}L_{2}R_{3})^{n}L_{3}^{2}L_{2}{\rm{(top)}}\leftrightarrow {\rm{(top)}}L_{2}(BCL_{1}R_{2})^{n+1} \end{equation}
for some $n\geq 1$.

\begin{equation}\label{3.12.3}\begin{array}{l}
R_{2}R_{3}L_{3}^{2}L_{2}\leftrightarrow R_{1}R_{2}BCL_{1}R_{2},\\ 
R_{3}L_{2}R_{3}L_{3}^{2}L_{2}{\rm{(top)}}\leftrightarrow  {\rm{(top)}}UCL_{1}R_{2}BCL_1R_2.\\
\end{array}\end{equation}

\begin{equation}\label{3.12.4}\begin{array}{l}R_{2}R_{3}L_{2}\leftrightarrow R_{1}^{2}R_{2}),\\ 
 R_{3}L_{2}R_{3}L_{3}\leftrightarrow UCL_{1}R_{2}R_{3}),\\ 
 L_{3}(L_{3}L_{2}R_{3})^{n}L_{2}{\rm{(top)}}\leftrightarrow {\rm{(top))}}L_{2}(BCL_{1}R_{2})^{n-1}BCL_{1}R_{1}R_{2},\\
 R_{2}R_{3}(L_{3}L_{2}R_{3})^{n}L_{2}\leftrightarrow R_{1}R_{2}(BCL_{1}R_{2})^{n-1}BCL_{1}R_{1}R_{2}),\\
 R_{3}L_{2}R_{3}(L_{3}L_{2}R_{3})^{n}L_{2}{\rm{(top)}}\leftrightarrow {\rm{(top)}}UCL_{1}R_{2}(BCL_{1}R_{2})^{n}.\end{array}\end{equation}

 \section{Exchangeable pairs for $(\zeta ,\eta )$}\label{3.13}
 
 Now we return to the question of when a  component $C$ of $s^{2-i}(\sigma _\zeta \circ s)^{-2}(E)$ coincides  with  a component of $(\sigma _{\alpha _\zeta }\circ s)^{-i}(E)$, for $E=E_{\zeta ,\eta }$, and can therefore be coded by basic exchangeable pairs. Coincidence fails when $C$ is in $(\sigma _{\alpha _\zeta }\circ s)^{j-i}(C')$ for a component $C'$ of $s^{2-j}(\sigma _\zeta \circ s)^{-2}(E))$ for some $j<i$ and $C'\cap \alpha _\zeta \neq \emptyset $. Now $C'=C_1'\cup C_2'\cup \ell _r$ where $C_1'$ and $C_2'$ are components of $s^{-j}(E)$ and $r=1$ or $2$, and $\ell _1$ and $\ell _2$ are the  component of $s^{2-j}((\sigma _{\zeta }\circ s)^{-1}(s^{-1}(\zeta ))$. In our cases it will always be true that
 $$(C_1'\cup C_2')\cap \alpha _\zeta =\emptyset $$
 and
 $$\ell \cap \alpha _\zeta =\ell \cap \zeta .$$
This will be proved in \ref{4.16}. This means that components of $(\sigma _{\alpha _\zeta }\circ s)^{-i}(E)$ will be components of $(\sigma _{\zeta _0}\circ s)^{2-i}(\sigma _\zeta \circ s)^{-2}(E))$.
We shall call these {\em{exchangeable pairs for $(\zeta ,\eta )$}}.
   
    This coding, as before, is inductive. Exchangeable pairs for $(\zeta ,\eta )$ (and $L_{3}$ or $BC$) are the same as basic exchangeable pairs for length $\leq 3$, and indeed all the exchangeable pairs $a\leftrightarrow b$ for $(\zeta ,\eta $ for which corresponding letters in $a$ and $b$ are all different are the same as the basic ones. In general the exchangeable pairs $a\leftrightarrow b$ of length $n$ are of the form $xa_{1}\leftrightarrow xb_{1}$ for exchangeable pairs $a_{1}\leftrightarrow b_{1}$ and letters $x$ such that $xa_{1}$ is an admissible word. (In that case $xb_{1}$ is admissible also as we are assuming that $a_{1}$ and $b_{1}$ can be preceded by the same letter.) The exception is when the component of  $s^{2-n}(\ell _{r}))$ corresponding to $a_{1}{\rm{(top)}}\leftrightarrow {\rm{(top)}}b_{1}$ intersects $\zeta $. The inverse image exchangeable pairs (for $L_{3}$)are then
  $$L_{3}a_{1}{\rm{(bot)}}\leftrightarrow R_{3}b_{1},\ \ L_{3}b_{1}{\rm{(bot)}}\leftrightarrow R_{3}a_{1}$$
 Then the coding proceeds inductively. 
 
 Relatively few letters can be different in an exchangeable pair $a\leftrightarrow b$. For a basic exchangeable pair the only pairs for which more than the last $5$ letters differ are
 \begin{itemize}
 \item $UCL_{1}R_{2}R_{3}L_{2}R_{3}\leftrightarrow R_{3}L_{2}R_{3}L_{3}^{3}$,
 \item pairs  $u_{1}a_{1}\leftrightarrow u_{2}b_{1}$ for which $a_{1}$ (and $b_{1}$) has length $\leq 4$  and the first two letters of $(u_{1},u_{2})$ are $(L_{2}BC,L_{3}^{2})$ or $(UCL_{1},R_{3},L_{2})$ or $(R_{1}R_{2},R_{2}R_{3})$, the last letters are $(R_{2},R_{3})$ and the rest is a suffix of $((BCL_{1}R_{2})^{n},(L_{3}L_{2}R_{3})^{n})$.
 \end{itemize}

The exchangeable pairs for $(\zeta ,\eta )$ are the same except that some more letters in $a$ and $b$ can be different. If $w_{k}'(\beta _{1})=u^{k}$ then any exchangeable pair for $(\zeta ,\eta )$ is of the form
 $$u_{r}a_{r}u_{r-1}a_{r-1} \cdots u_{1}a_{1}\leftrightarrow u_{r}b_{r}u_{r-1}b_{r-1} \cdots u_{1}b_{1}$$
 where
 \begin{itemize}
 \item $a_{i}\leftrightarrow b_{i}$ is a basic exchangeable pair
 \item $D(u^{k_i})$ is between $D(u_{i}a_{i})$ and $D(u_{i}b_{i})$ for some $k_i$ and $2\leq i\leq r$.
 \end{itemize}
 
 We have the following extension to \ref{3.11}, which follows from the definition of the $w_j(u)$ and $w_j'(u)$.
 \begin{lemma}\label{3.14.1} Let $v\leftrightarrow v'$, be an exchangeable pair for $(\zeta ,\eta )$. For any prefixes $u$ and $u'$ of $vL_{3}L_{2}$ and $v'L_{3}L_{2}$ with $\vert u\vert =\vert u'\vert $, we have the following. The largest integer $i=i(u)$ such that $w_{i}(u)$ is a proper prefix of $u$ is the same as the largest integer $k$ such that $w_{i}(u')$ is a proper prefix of $u'$. Write
 $$v=u_{1}a_{1}\cdots u_{t}a_{t},\ \ v'=u_{1}b_{1}\cdots u_{t}b_{t}.$$
 If $k(u_{1}L_{3}L_{2})\geq 1$ then 
 $$w_{1}'(u_{1}a_{1})=w_{1}'(u_{1}b_{1})=w_{1}'(v)=w_{1}'(v'),$$
 and this is also true if $i(u_{1}a_{1})\geq 1$, except for certain choices of $(a_{1},b_{1})$ as specified in (\ref{3.11.5}).\end{lemma}

\section{The homeomorphisms $\psi _{m,q}$.}\label{3.15}
As in \cite{R6}, we let  $\beta _q:[0,1]\to \overline{\mathbb C}$ be the capture path which crosses the unit circle at $e^{2\pi iq}$ and ends at the point $\beta _q(1)\in Z(s)$ in the unique gap of $L_{3/7}$ with $e^{2\pi iq}$ in its boundary. Thus, $\beta _q$ is defined for all $q$ such that $e^{2\pi iq}$ is in  the boundary of a gap of $L_{3/7}$. If $\beta _q(1)$ is strictly preperiodic, then $\beta _q$ is a type III capture map and $\sigma _{\beta _q}\circ s$ is   type III branched covering. If $\beta _q(1)$, is periodic then we define $\zeta _q:[0,1]\to \Cbar$ to be the unique path such that $s\circ \zeta _q=\beta _q$ and $\zeta _q(1)$ is periodic under $s$. In this case $\sigma _{\zeta _q}^{-1}\circ \sigma _{\beta _q}\circ s$ is a type II branched covering and $\beta _q$ is a type II capture path (by definition). An example is given by $q=\frac{2}{7}$. This is the most important example for us.

In 3.3 of \cite{R6} we gave a proof (not new in any essential way, of course) that
$$\sigma _{\zeta _{5/7}}^{-1}\circ \sigma _{\beta _{5/7}}\circ s \simeq _{\psi _{0,2/7}}\sigma _{\zeta _{2/7}}^{-1}\circ \sigma _{\beta _{2/7}}\circ s $$
We then used $\psi _{0,2/7}$ to define a sequence of homeomorphisms $\psi _{m,2/7}$ and a sequence of closed loops $\alpha _{m,2/7}$, for all $m\ge 0$, and satisfying the following properties.
 
\begin{itemize}
\item[1.] $\alpha _{m,2/7}$ is an arbitrarily small perturbation of $\beta _{2/7}*\psi _{m,2/7}(\overline{\beta _{5/7}})$ which bounds a disc disjoint from the common endpoint $0$ of $\beta _{2/7}$ and $\beta _{5/7}$, for all $m\ge 0$.
\item[2.] $\psi _{m+1,2/7}=\xi _{m,2/7}\circ \psi _{m,2/7}$ for all $m\ge 0$.
\item[3.] $\psi _{m+1,2/7}$ and $\psi _{m,2/7}$ are isotopic relative to $Y_m$ for all $m\ge 0$ and $\xi _m$ is isotopic to the identity relative to $Y_m$ for all $m\ge 0$.
\item[4.] $\sigma _{\alpha _{m,2/7}}\circ s\circ \psi _{m+1,2/7}=\psi _{m,2/7}\circ s$ and  and $\sigma _{\beta _{2/7}}\circ s\circ \xi _{m+1,2/7}=\xi _{m,2/7}\circ \sigma _{\beta _{2/7}}$ for all $m\ge 1$.
\end{itemize}

 In the definition of $\psi _{m,2/7}$ and $\alpha _{m,2/7}$ in \cite{R5}, we started with the definition of $\psi _{0,2/7}$. But we might as well start from $\psi _{2,2/7}$ and $\xi _{1,2/7}$, because then 1 to 4 can be used to defined  $\psi _{m,2/7}$ and $\xi _{m,2/7}$ and $\alpha _{m,2/7}$ inductively for all $m\ge 2$. (Of course, this will also give the definition of $\psi _{1,2/7}$, as $\psi _{2,2/7}=\xi _{1,2/7}\circ \psi _{1,2/7}$.) The homeomorphism $\psi _{2,2/7}$ is a composition of anticlockwise Dehn twists round the boundaries of four topological discs, two larger disjoint discs and two smaller discs, one contained in each of the others. The two larger discs are called $D_+$ and $D_-$, both contained in the unit disc, symmetrically placed about the origin, with $D_+$ bounded by the leaf of $L_{3/7}$ --- just a vertical line up to isotopy preserving $Z$ --- joining $e^{2\pi i(5/14)}$ and $e^{2\pi i(9/14)}$ -- and $D_-$ is its symmetric copy bounded by the leaf joining $e^{2\pi i(1/7)}$ and $e^{2\pi i(6/7)}$. The two smaller discs are the components of $s^{-1}(D_-)$.  The homeomorphism $\xi _{1,2/7}$ is the composition of anticlockwise Dehn twist round the boundaries of the components of $s^{-1}(D_{-,2/7})$,  and of a clockwise twist in an annulus $A_{1,2/7}$. This annulus bounds a disc disjoint from $\infty $ and containing $D_{+}$. Its intersection with the unit disc is the set bounded by one vertical leaf  joining $e^{2\pi i(2/7)}$ and $e^{2\pi (5/7)}$, and one joining $e^{2\pi i(5/14)}$ and $e^{2\pi i(9/14)}$. A clockwise twist in the annulus means, up to isotopy, a clockwise twist round the outer boundary composed with an anticlockwise twist round the inner boundary.

From the definitions, it will be seen that $\psi _{2,2/7}$ is the identity on $\beta _{5/7}$ and hence $\alpha _{2/2/7}$ is an arbitrarily small perturbation of $\beta _{2/7}*\overline{\beta _{5/7}}$. For any path $\gamma $ based at $\infty $ we also have 
$$\alpha _{m,2/7}*\psi _{m,2/7}(\gamma )=\beta _{2/7}*\psi _{m,2/7}(\overline{\beta _{5/7}}*\gamma ).$$
For $\gamma \in {\cal{Z}}(3/7,+,0)$, $\psi _{2,2/7}$ is the identity on $\overline{\beta _{5/7}}*\gamma $ and hence for all $m\ge 3$,
$$\alpha _{m,2/7}*\psi _{m,2/7}(\gamma )=\beta _{2/7}*\xi _{m-1,2/7}\circ \cdots \circ \xi _{2,2/7}(\overline{\beta _{5/7}}*\gamma ).$$

In order to proceed further we need to give some more information about the homeomorphisms $\xi _{i,2/7}$ for $i\ge 2$. The support of $\xi _{2,7}$ is a union of annuli $A_{i,2/7}\cup C_{i,2/7}$ where 
$$A_{i,2/7}=(\sigma _{\beta _{2/7}}\circ s)^{-1}(A_{i-1,2/7})$$
for all $i\geq 2$ and 
$$C_{i,2/7}=(\sigma _{\beta _{2/7}}\circ s)^{-1}(C_{i-1,2/7})$$
for all $i\ge 1$. 

Previously in this chapter, we have used homeomorphisms which are compositions of disc exchanges. The homeomorphisms $\xi _{i,2/7}$ are a bit different, because each component of the support is a topological annulus rather than a topological disc. Within each annulus of $A_{i,2/7}$, components of intersection with the unit disc are cyclically permuted by $\xi _{2,2/7}$, up to isotopy preserving $Z_m$, for any fixed $m\ge i$, with each component of intersection sent to the next one in the clockwise direction. Within each annulus of $C_{i,2/7}$ the components of 
$$s^{1-i}(\{ z:|z|\le 1\} \setminus \partial C_1\cap \{ z:|z|<1\} )$$
are cyclically permuted, with each one sent to the next one round in an anticlockwise direction. It is therefore possible to code the action of these homeomorphisms using words. This was done in section 7.2 of \cite{R6}. 

We have already defined $A_{1,2/7}$. The set  $C_{0,2/7}$ is simply the complement in $D_{-,2/7}$ of a sufficiently small disc neighbourhood of $s(v_1)$ that the intersection of the neighbourhood with $Z$ is simply $s(v_1)$. The sets $A_{i,2/7}$ and $C_{i,2/7}$ therefore have empty intersection with $Z_i\cup \{ \infty \} =Y_i$, but do intersect $Z_{i+1}$. The branched covering  $\sigma _{\beta _{2/7}}\circ s$ is of degree one or two on each component of $A_{i,2/7}$ and $C_{i,2/7}$, and both are possible. 

We now give some pictures of these sets. 
    \begin{figure}\includegraphics[width=8cm]{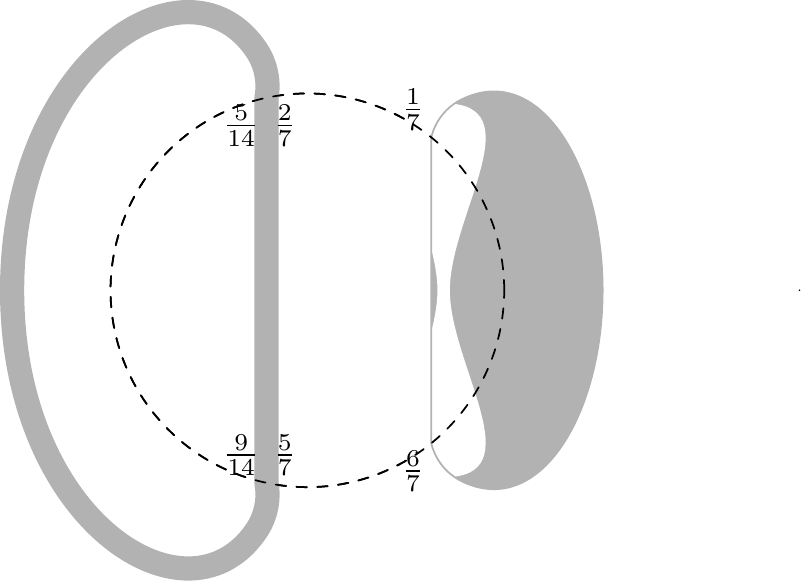}
\caption{$A_{1,2/7}$ and $C_{0,2/7}$}
\end{figure}

\begin{figure}
\includegraphics[width=4cm]{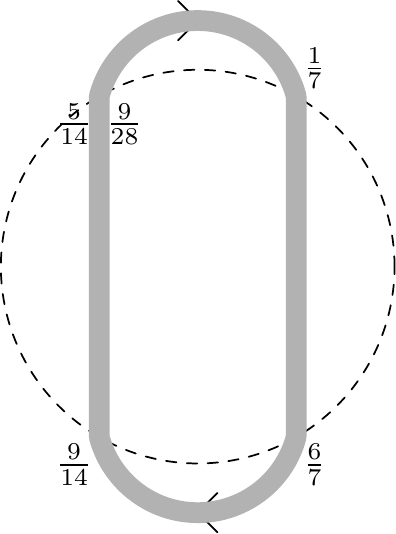}
\caption{$A_{2,2/7}$}\end{figure}

\begin{figure}
\includegraphics[width=6cm]{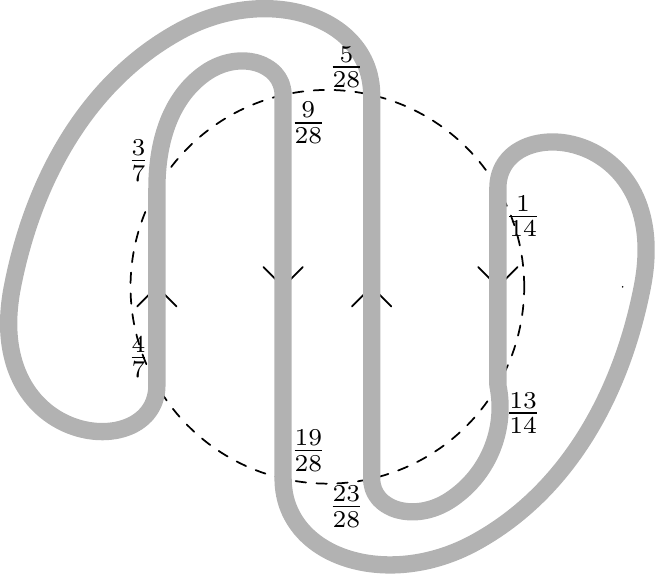}
\caption{$A_{3,2/7}$}\end{figure}

\begin{figure}
\includegraphics[width=6cm]{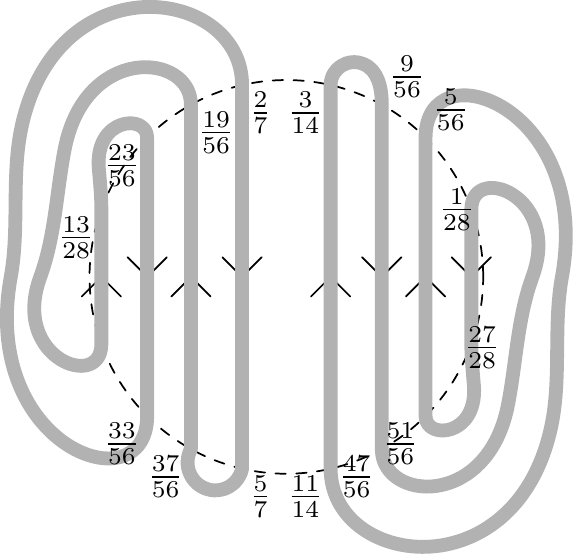}
\caption{$A_{4,2/7}$}\end{figure}

\begin{figure}
\includegraphics[width=6cm]{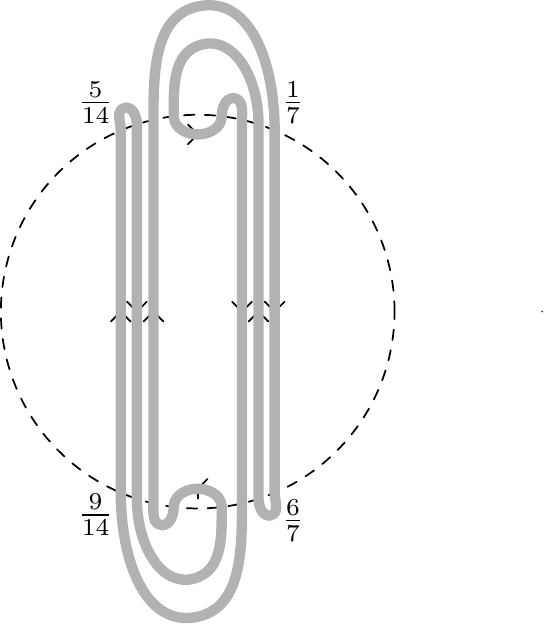}
\caption{$A_{5,1,2/7}$}\end{figure}

\begin{figure}
\includegraphics[width=6cm]{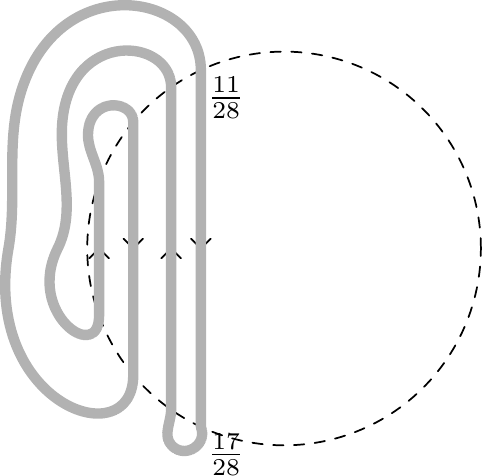}
\caption{$A_{5,2,2/7}$}\end{figure}

\begin{figure}
\includegraphics[width=4cm]{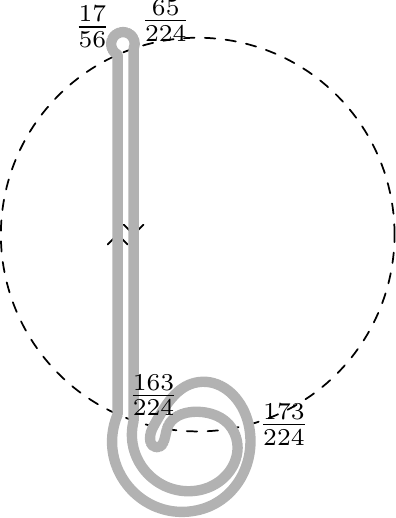}
\caption{$A_{6,2,2/7}$}\end{figure}

The permutations of sets $D(w)$ in these annuli include the following. The clockwise twist in $A_{1,2/7}$ preserves the single component of intersection $D(L_3)$ with the unit disc, up to isotopy preserving $Z_m$, for any fixed $m$. In $A_{2,/7}$ the permutation of components of intersection is indicated by:
$$L_{3}^{2}\to {\rm{ (top)}}R_{3}L_{3}\to {\rm{(bot)}}L_{3}^{2}.$$
Here, ``top'' immediately preceding the word indicates a component of $A_{2,2/7}\cap \{ z:|z|>1$ ending on the top of the unit circle and ``bot'' immediately preceding the word indicates component of $A_{2,2/7}\cap \{ z:|z|>1$ ending on the bottom of the unit circle. In fact, in this example, one component of $A_{2,2/7}\cap \{ z:|z|>1$ is above the unit circle and one is below. Similarly,  the action of $\xi _{3,2/7}$ in $A_{3,2/7}$ is indicated by 
$$L_{2}R_{3}L_{3}{\rm{(top)}}\to  {\rm{(top)}}L_{3}^{3}{\rm{bot}}\to {\rm{(top)}}R_{2}R_{3}L_{3}{\rm{(bot)}}$$
$$\to  {\rm{(bot)}}R_{3}L_{3}^{2}{\rm{(top)}}\to {\rm{bot}}L_{2}R_{3}L_{3}.$$
For the two components of $A_{4,2/7}$, which we call $A_{4,1,2/7}$ and $A_{4,2,1/7}$, we have 
$$L_{1}R_{2}R_{3}L_{3}{\rm{(top)}}\to  {\rm{(top)}}L_{2}R_{3}L_{3}^{2}{\rm{(bot)}}\to  
{\rm{top)}}L_{3}L_{2}R_{3}L_{3}{\rm{(bot)}}$$
$$\to {\rm{(bot)}}L_{3}^{4}{\rm{(top)}}\to {\rm{(bot)}}L_{1}R_{2}R_{3}L_{3},$$
$$R_{1}R_{2}R_{3}L_{3}{\rm{(bot)}}\to  {\rm{(bot)}}R_{2}R_{3}L_{3}^{2}{\rm{(top)}}\to  
{\rm{(bot)}}R_{3}L_{2}R_{3}L_{3}{\rm{(top)}}$$
$$\to {\rm{(top)}}R_{3}L_{3}^{3}{\rm{(bot)}}\to {\rm{(top)}}R_{1}R_{2}R_{3}L_{3}.$$

For the single component $A_{5,1,2/7}$ of $A_{4,1,2/7}$, the action of $\xi _{5,2/7}$ is indicated by 
$${\rm{(bot)}}BCL_{1}R_{2}R_{3}L_{3}\to  
{\rm{(bot)}}L_{3}L_{2}R_{3}L_{3}^{2}{\rm{(top)}}\to  
{\rm{(top)}}
R_{3}L_{3}L_{2}R_{3}L_{3}{\rm{(bot)}}\to $$
$${\rm{(bot)}}R_{3}L_{3}^{4}\to {\rm{(top)}}UCL_{1}R_{2}R_{3}L_{3}\to 
{\rm{(top)}}R_{3}L_{2}R_{3}L_{3}^{2}\to  
{\rm{(bot)}}
$$
$$L_{3}^{2}L_{2}R_{3}L_{3}\to {\rm{(top)}}L_{3}^{5}\to {\rm{(bot)}}BCL_{1}R_{2}R_{3}L_{3}.$$
Since $\beta _{2/7}$ does not intersect any of the components of $A_{4,2,2/7}\cap \{ z:|z|\ge 1\} $,  the components of $(\sigma _{\beta _{2/7}}\circ s)^{-1}(A_{4,2,2/7})$ are actually components of $s^{-1}(A_{4,2,2/7})$. In fact since none of the points $e^{2\pi ir/7}$ for $r\in \{1,2,4\} $ are in any interval of the unit circle bounded by a component of $A_{4,2,2/7}\cap \{ z:|z|\ge 1\} $, all components of $A_{i,2/7}$ in $(\sigma _{\beta _{2/7}}\circ s)^{4-i}(A_{4,2,2/7})$, for all $i\ge 5$, are actually components of $s^{4-i}(A_{4,2,2/7})$. In fact all components of $s^{-p}(A_{i,2/7}\cap \{ z:|z|\ge 1\} $ are disjoint from $\beta _{2/7}$ for all $p\ge 0$ except for those components represented by:
$$(L_3L_2R_3)^nL_{3}^{2}{\rm{(top)}}\to {\rm{ (top)}}(R_3L_3L_2)^nR_{3}L_{3},$$
if $i=3n+2$ for any $n\ge 0$;
$$(R_3L_3L_2)^nR_3L_3^2{\rm{(top)}}\to {\rm{(bot)}}(L_2R_3L_3)^{n+1}$$
if $i=3n+3$;
$$(L_2R_3L_3)^nL_3{\rm{(bot)}}\to {\rm{(top)}}L_3(L_2R_3L_3)^n$$
if $i=3n+4$. For any other component represented by $Y\to Y'$, any component of any $A_{j,2/7}\cap \{ z:|z|\ge 1\} $ in the backward orbit is represented by 
$ZXY\to ZX'Y$, where $(X,X')$ is one of the following, or obtained from one of these by adjoining $BC,L_3$ or $(BCL_1,L_3L_2)$, or a suffix of one of these 
 if $Z$ is empty:
$$(R_1R_2,R_2R_3),\ (UCL_1R_2(BCL_1R_2)^p,R_3L_2R_3(BCL_1R_2)^p)$$
for any $p\ge 0$.  So, altogether, these give the same pairs as the basic exchanges in \ref{3.11}.

As for components of $C_{i,2/7}$, we have a similar pattern. Here are some pictures of some of the sets $C_{i,2/7}$.
\begin{figure}
\includegraphics[width=4cm]{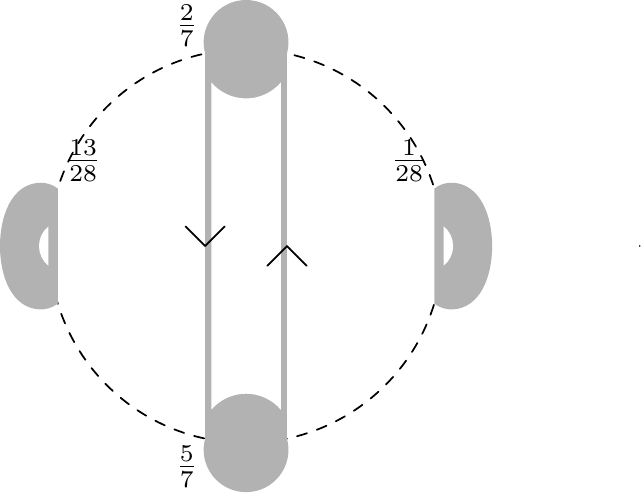}
\caption{$C_{2,2/7}$}\end{figure}

\begin{figure}
\includegraphics[width=4cm]{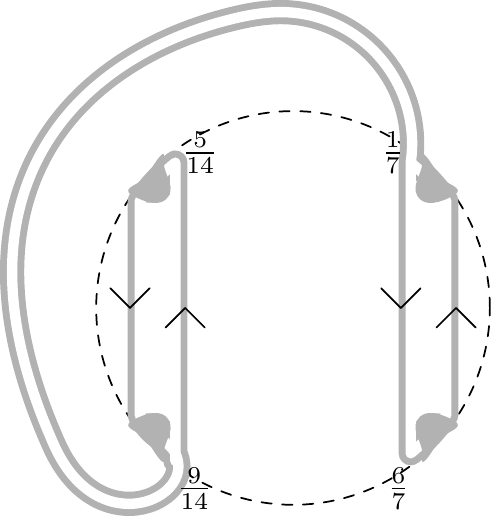}
\caption{$C_{3,2/7}$}\end{figure}

\begin{figure}
\includegraphics[width=4cm]{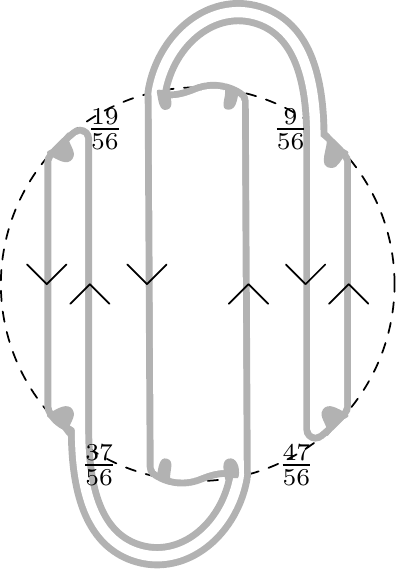}
\caption{$C_{5,2/7}$}\end{figure}

The homemorphism $\xi _{2,2/7}$ restricted to   $C_{2,2/7}$ is represented by:
$$UC\to BC\to UC.$$
The preimage $C_{3,2/7}$ of $C_{2,2/7}$ is connected and $\xi _{3,2/7}$  restricted to $C_{3,2/7}$ is represented by:
$$R_{2}BC\to R_{2}UC\to L_{2}BC\to L_{2}UC\to R_{2}BC.$$
This, and the following representations, should be read in conjunction with the pictures. In $C_{3,2/7}$, for example, the piece of $C_{3,2/7}$ between $D(R_2BC)$ and $D(R_2UC)$ is inside the unit circle, along the right boundary of $D(R_2C)$, while the piece of $C_{3,2/7}$ between $D(L_2UC)$ and $D(R_2UC)$ traces an approximate diagonal path from the lower unit circle (adjacent to $D(L_2UC)$, then enters the unit disc again on the upper boundary of $D(R_2C)$ and travels down the left vertical boundary of $D(R_2C)$ to $D(R_2BC)$. The set $C_{4,2/7}$ has two components  $C_{4,1,2/7}$  and $C_{4,2,2/7}$. The homeomorphism $\xi _{4,2/7}$ restricted to $C_{4,1,2/7}$ and $C_{4,2,2/7}$ is represented by, respectively:
$$L_{1}R_{2}BC\to L_{1}R_{2}UC\to 
L_{3}L_{2}BC\to 
L_{3}L_{2}UC\to L_{1}R_{2}BC,$$
$$R_{1}R_{2}BC\to R_{1}R_{2}UC\to 
R_{3}L_{2}BC\to 
R_{3}L_{2}UC\to R_{1}R_{2}BC.$$
The preimage of $C_{4,1,2/7}$ in $C_{5,2/7}$ has a single component $C_{5,1,2/7}$ and the homeomorphism $\xi _{5,2/7}$ restricted to $C_{5,1,2/7}$ is 
 represented by: 
$$BCL_{1}R_{2}BC\to BCL_{1}R_{2}UC\to 
L_{3}^{2}L_{2}BC\to L_{3}^{2}L_{2}UC\to $$
$$UCL_{1}R_{2}BC\to UCL_{1}R_{2}UC\to 
R_{3}L_{3}L_{2}BC\to 
R_{3}L_{3}L_{2}UC\to BCL_{1}R_{2}BC.$$
None of the components of $C_{4,2,2/7}\cap \{ z:|z|\ge 1\} $ intersect $\beta _{2/7}$ and neither do any of the components of $s^{-n}(C_{4,2,2/7}\cap \{ z:|z|\ge 1\} )$ for any $n\ge 0$. Therefore the components of $(\sigma _{\beta _{2/7}}\circ s)^{-n}(C_{4,2,2/7}\cap \{ z:|z|\ge 1\} )$ are components of $s^{-n}(C_{4,2,2/7}\cap \{ z:|z|\ge 1\} )$. The same is true for all but two of the components of $C_{4,1,2/7}\cap \{ z:|z|\ge 1\} $. The exceptional pair are the components represented by 
$$L_1R_2UC\to L_3L_2BC$$
and
$$L_3L_2UC\to L_1R_2BC,$$
 where both the components intersect the unit circle in $D(L_1R_2UC)$ and $D(L_3L_2UC)$, and there are adjacent components of $C_{4,1,2/7}$ in the boundaries of $D(L_1R_2C)$ and $D(L_3L_2C)$ connecting $D(L_1R_2UC)$ to $D(L_1R_2BC)$, and $D(L_3L_2UC)$ to $D(L_3L_2BC)$. Of the preimages of these in $C_{5,1,2/7}$, only those represented by
 $$UCL_1R_2UC\to R_3L_3L_2BC$$
 and
 $$R_3L_3L_2UC\to UCL_1R_2BC$$
  have backward orbits which intersect $\beta _{2/7}$. Similarly in $C_{i,,2/7}\cap \{ z:|z|\ge 1\} $, the  only components whose backward orbits intersect $\beta _{2/7}$ represented by
 $$(R_2UCL_1)^nR_2UC\to (L_2R_3L_3)^nL_2BC$$
 and
 $$(L_2R_3L_3)^nL_2BC\to (R_2UCL_1)^nR_2UC$$
 if $i=3n+3$, or 
 $$(L_1R_2UC)^{n+1}\to (L_3L_2R_3)^nL_3L_2BC$$
 and
 $$(L_3L_2R_3)^nL_3L_2BC\to (L_1R_2UC)^{n+1}$$
 if $i=3n+4$, or 
 $$(UCL_1R_2)^nUC\to (R_3L_3L_2)^nBC$$
 and
 $$(R_3L_3L_2)^nBC\to (UCL_1R_2)^nUC$$ 
 if $i=3n+5$.

The pre-images in $C_{i+n,2/7}$ of any component of $C_{i,2/7}\cap \{ z:|z|\ge 1\} $ represented by $Y\to Y'$ whose backward orbit does not intersect $\beta _{2/7}$, that is, any but the cycle given above, are obtained by adjoining $Z$ to $Y\to Y'$ if the first letters of $Y$ and $Y'$ are the same, or are in the set $\{ R_3,BC,UC\} $, and adjoining  $ZX$ and $ZX'$ to $Y\to Y'$, giving $ZXY\to ZX'Y'$, where $(X,X')$ is of one of the following forms, or a  suffix of one of these if $Z$ is empty:
$$(R_1,R_2),\ (UC(L_1R_2BC)^pL_1,R_3(L_2R_3L_3)^pL_2),$$
$$(L_2(BCL_1R_2)^pBCL_1R_1,L_3(L_3L_2R_3)^{p+1})$$
for any $p\ge 0$. These  are the same forms as for the basic exchangeable pairs for $BC$ in \ref{3.12}.

\section{$\psi _{m,q}$}\label{3.16}

We saw in 3.3 of \cite{R5} that the equivalence between the type III captures  $\sigma _{\beta _{1-q}}$ and $\sigma _{\beta _q}\circ s$ can be realised similarly for any $q\in (\frac{2}{7},\frac{1}{3})$ such that $\beta _q$ is a type III capture path, that is, $e^{2\pi iq}$ is on the boundary of a gap of $L_{3/7}$. The homeomorphisms $\psi _{m,q}$ and $\xi _{m,q}$ and loops $\alpha _{m,q}$ are defined in the same inductive fashion as the homeomorphisms $\psi _{m,2/7}$ and $\xi _{m,2/7}$ and loops $\alpha _{m,2/7}$ in \ref{3.15}. The support of $\xi _{m,q}$ is again a union of annuli 
$$A_{m,q}\cup C_{m,q}=(\sigma _{\beta _{q}}\circ s)^{2-m}(A_{2,q}\cup C_{2,q})$$ 
and $\psi _{1,q}$ is again a composition of anticlockwise Dehn twists round the boundaries of discs $D_{\pm ,q}$ and $s^{-1}(D_{-,q})\subset D_{+,q}\cup D_{-,q}$.

\section {First stage in the proof of Theorem \ref{2.8}}\label{3.17}

We are now ready to prove parts of Theorem \ref{2.8}. We prove the base of the induction and part of the inductive step. We cannot prove the whole of the inductive step at this stage, however. That will have to wait until Chapter 4. Nevertheless, the part of the inductive step that can be proved holds for ${\cal{Z}}_m(3/7,+,p)$ for any $p\ge 0$. We have the following.

\begin{lemma}\label{3.17.1} Let $\beta \in {\cal{Z}}_m(3/7,+,p)$. Let   
$$((\omega _{2i-1},\omega _{2i}):1\leq i\leq n(\beta ))$$ and 
$$(\gamma _i:1\leq i\leq n(\beta ))$$
 be the unique sequences of, respectively, adjacent pairs in $\Omega _m$ and elements of  $\pi _1(V_{3,m},a_1)$, such that $\beta $ successively crosses the arcs $\overline{\rho (\gamma _i*\omega _{2i-1})}*\rho (\gamma _i*\omega _{2i-1})$, with $\beta =\rho (\gamma _{n(\beta )}*\omega _{2n(\beta )})$. Then the following hold.
 \begin{itemize}
 \item[1.] $(\rho(\omega _1),\rho (\omega _2))=(\beta _1,\beta _2)\in R_{m,p}$ if $\beta \in {\cal{Z}}_m(+,+,p)$ and $(\rho (\omega _1),\rho (\omega _2)=(\beta _1',\beta _2')\in R_{m,p}'$ if $\beta \in {\cal{Z}}_m(+,-,p)$. 
 \item[2.] Let $\beta \in {\cal{Z}}_m(3/7,+,+,p)$. Suppose that it is true that, for $j\le k$, there are adjacent pairs $(\beta _{2j-1},\beta _{2j})$  in $R_{m,p}$ matched with $(\beta _{2j-1}',\beta _{2j}')$ in $R_{m,p}'$, and 
 $$(\rho (\omega _{2j-1}),\rho (\omega _{2j-1}))=\begin{cases}(\beta _{2j-1},\beta _{2j}){\rm{\ if\ }}j{\rm{\ is\ even,}}\\ 
(\beta _{2j-1}',\beta _{2j}'){\rm{\ if\ }}j{\rm{\ is\ odd.}}\end{cases}$$
Then $(\omega _{2k+1},\omega _{2k+2})$ is an adjacent pair in $R_m$ if $k$ is even and in $R_m'$ if $k$ is odd. 
\item[3.] An analogous statement holds for $\beta \in {\cal{Z}}_m(3/7,+,-,p)$, with odd and even interchanged. 
\end{itemize}
\end{lemma}
\subsection{Remarks}\label{3.17.2} The inductive step in 2 (and 3) is not quite strong enough to deduce the full statement of Theorem \ref{2.8}. The proof would be complete if $R_m$ and $R_m'$ could be replaced by $R_{m,p}$ and $R_{m,p}'$. It has only proved to be possible to do that for $p=0$, because of lack of detailed knowledge about $R_{m,p}$ for $p>0$. In any case, the completion of the proof of Theorem \ref{2.8} §will have to wait until Chapter \ref{4}.
\subsection{}\label{3.17.3}
\begin{proof}  The paths in $R_{m,p}$ have never been described in detail for $p>0$, although the set of endpoints lie in a set $U^p$ which is described in Theorem 2.10 of \cite{R5}.  We now present the only facts which we shall need. All statements are up to homotopy preserving $Y_\infty $.
 
 \begin{itemize}
 \item All unit disc crossings of paths in $R_m$ and the last component of intersection with the unit disc,  are either along leaves of $L_{3/7}$ or in gaps of $L_{3,7}$. The same is true for paths in $R_m'$.
 \item The first unit-disc-crossings of paths in $R_{m,p}$ (apart from the path $\beta _{q_{p+1}}$) start at a point $e^{2\pi it}$ for $t\in [q_p,q_{p+1})$. 
\item The sets $U^p$  are not all disjoint, not all points in $U^p\cap Z_m(s)$ are second endpoints of paths in $R_{m,p}$, and it is possible for at most two different paths in $R_{m,p}$ to have the same endpoint. However, each path in $R_m$ is uniquely determined by its first unit-disc crossing and its final endpoint.
\end{itemize}

Suppose first that $\beta \in {\cal{Z}}_m(3/7,+,+,p)$, that is, the $S^1$ crossing point is at $e^{2\pi ix}$ for $x\in (q_{p},q_{p+1})$ for some $p\geq 0$. We do not need to consider $x=q_p$ for any $p$ because $\beta _{q_p}\in R_{m,p}$. 

First, we prove 1. By definition, a capture path has only one intersection with the unit circle before entering the gap of $L_{3/7}$ containing the second endpoint of $\beta $. So the first essential intersections with $S^1$ of the paths  $\rho (\gamma _i*\omega _{2i-1})$ and  $\rho (\gamma _i*\omega _{2i})$ must be distinct, and opposite sides of the intersection with $S^1$ of $\beta$, although one of them could coincide with this. If  $\omega _i\in\Omega _m(a_0)$ (or  $\omega _i\in\Omega _m(\overline{a_0})$) for at least one of $i=1$ or $2$,  then both $\rho (\omega _1)$ and $\rho (\omega _2)$ have first essential $S^1$-crossing at $e^{2\pi i(6/7)}$ (or $e^{2\pi i(1/7)}$). This is is incompatible with $\beta $ crossing $\overline{\rho (\omega _1)}*\rho (\omega _2)$.  The paths in $\Omega _m(a_1,-)$ are ruled out more simply as candidates for $\rho (\omega _1)$ or $\rho (\omega _2)$, because these are simply the paths ${\cal{Z}}_m(3/7,-)$, which have first and only crossings of $S^1$ at points $e^{2\pi ix}$ for $x\in (-\frac{1}{7},\frac{1}{7})$. Similarly we can rule out $\Omega _m(a_1,-)$.   So we must have $\omega _i\in \Omega _m(a_1,+)$, and $\rho (\beta _i)\in R_{m,p}$ for $i=1$ and $2$, in order to have $\beta $ crossing $\overline{\beta _1}*\beta _2$ before the first $S^1$-crossing points of those paths. 

Now we consider the inductive step. So we assume the inductive hypothesis of the lemma, which is assuming that  suppose that  1 of \ref{2.8} is true for $\beta _{i}$ for $i\leq 2k$, for some $k\geq 1$. Then we consider what $(\omega _{2k+1},\omega _{2k+2})$ can be.  First, let $k$ be even. Then the inductive hypothesis gives $\rho (\gamma _{k+1})=\alpha _{1,2k-1}$ in the notation of \ref{2.8}, that is, $\alpha _{1,2k-1}$ is an arbitrarily small perturbation of 
$$\alpha _{1,2k-1}=\alpha _{1,3}*\psi _{1,3}(\alpha _{5,7}*\cdots *\psi _{1,2k-5}(\alpha _{2k-3,2k-1})$$
where $\alpha _{4j-3,4j-1}$ is an arbitrarily small perturbation of 
$$\beta _{4j-3}*\psi _{4j-3}(\overline {\beta _{4j-3}'}*\beta _{4j-1}')*\psi _{4j-3,4j-1}(\overline{\beta _{4j-1}}).$$
This means that the number of vertical intersections of $\rho (\gamma _{k+1})=\alpha _{1,2k-1}$ with the unit disc is even. So $\psi _{1,2k-1}(\rho (\omega _{2k+1}))$ and $\psi _{1,2k-1}(\rho (\omega _{2k+2}))$ must first  cross the unit disc in the same direction as $\beta $, that is, across the upper unit circle first. Also, we now know that $\psi _{1,2k-1}$ is a composition of disc exchanges, and we have a complete list of disc exchanges, in particular of those which intersect the unit circle within the set $\{ e^{2\pi it}:t\in [\frac{2}{7},\frac{1}{3})\} $.   None of these supports can cancel first  $S^1$-crossing of $\rho (\omega _{2k+1})$ or $\rho (\omega _{2k+2})$ at $e^{2\pi ix}$ for $x\notin [\frac{2}{7},\frac{1}{3})$. So we must have $(\rho (\gamma _{2k+1},\gamma _{2k+2})=(\beta _{2k+1},\beta _{2k+2})$ where $(\beta _{2k+2},\beta _{2k+2})$ is an adjacent pair in $R_m$. 
 
 If $k$ is odd then the proof is similar. But this time, 
 $$\alpha _{1,2k-1}=\alpha _{1,3}*\psi _{1,3}(\alpha _{5,7}*\cdots *\psi _{1,2k-1}(\alpha _{2k-1})$$
 and the number of intersections with the unit circle is odd. 

  Now we consider the case of $\beta \in {\cal{Z}}_m(3/7,+,-,p)$, that is, with $S^1$-intersection at $e^{2\pi ix}$ for $x\in (-q_{p+1},-q_p)$.  For the base of the induction, we consider $\alpha _{q_p}*\psi _{q_p}(\beta )$. Then the first arc of the form $\rho(\overline{\eta _1})*\rho(\eta _2)$ crossed by this, for an adjacent pair $(\eta _1,\eta _2)\in \Omega _m$, is for an adjacent  pair $(\rho (\eta _1),\rho (\eta _2))=(\beta _1,\beta _2)$ in $R_{m,p}$. Then let $\psi _1,\alpha _1)$ effect the matching between $(\beta _1,\beta _2)$ and $(\beta _1'\beta _2')$ where $(\beta _1',\beta _2')$ is the matched adjacent pair in $R_{m,p}$. Then $\overline{\beta _1}*\beta _2=\psi _1(\overline{\beta _1'}*\beta _2')$ is also the first arc crossed by $\alpha _1*\psi _1(\beta )$, which means that $\overline{\beta _1'}*\beta _2'$ is the first arc crossed by $\beta $. Then we proceed inductively as before.  \end{proof}

\chapter{Quadruples of paths}\label{4}

\section{Notation}\label{4.1}

Our goal in this work is to prove our main theorem, Theorem \ref{1.9}, concerning capture paths $\beta \in {\cal{Z}}_m(3/7,+,+,0)$. This means analysing the sequence of translates of the fundamental domain crossed by $\beta $, as explained in Chapter 1. A particular structure for this sequence is claimed in theorem \ref{2.8}. In particular the sequence of pairs $(\omega _{2i-1},\omega _{2i})$ is related to a sequence $(\beta _{2i-1}.\beta _{2i})$ of adjacent pairs $R_{m,0}$, and a matching sequence of adjacent pairs $(\beta _{2i-1}',\beta _{2i}')$ in $R_{m,0}'$. Theorem \ref{2.8} has not yet been completely proved. But the base case of the induction has been completed, and  much of the inductive step has been established. This is enough to work with. We will now assume the existence of $(\beta _{2j-1},\beta _{2j})$ for $j\le i$, with the properties claimed in \ref{2.8}, and will work to extend these. In the process, the proof of Theorem \ref{2.8} will be completed.

From now on we consider capture paths $\beta \in {\cal{Z}}_m(3/7,+,+,0)$, and the associated sequences $\beta _j(\beta )$ for $j\le 2i$ in $R_{m,0}$, and the homeomorphisms $\psi _{j,\beta }$ and $\psi _{j-2k,\beta }$ for any $k\ge 1$, as defined in \ref{2.9}. We write $Y_m$ and $Z_m$ for $Y_m(s)$ and $Z_m(s)$. If $\beta $ ends in the point $p(x)$ for a word $x$ (see \ref{2.1}) then we will also write $\beta _j(x)$ and $\psi _{j,x}$ and $\psi _{j-2k,j,x}$ for these quantities. 

A natural way to proceed is to consider what happens if $\psi _{j-2k,j,\beta }$ is always the identity. If that were the case then the sequence of arcs $\overline{\beta _{2j-1}}*\beta _{2j}=\psi _{2i-1}(\overline{\beta _{2j-1}'}*\beta _{2j}')$ which is crossed by $\beta $ can be computed just by studying $R_{m,0}$. Much of this chapter is devoted to studying this sequence, which is denoted by $\zeta _j$. The precise definition is given below. We call such sequences  {\em{$\zeta $-sequences}}. One aim is to relate $\beta _j(\beta )$ to various $\zeta $-sequences. In the process, we also obtain information about quadruples $(\beta _{4j+t}(\beta ):-3\le t\le 0)$. In the process we are led to study {\em{$\zeta $-quadruples}} $(\zeta _{4j+t}:-3\le t\le 0)$.

\section{Principal arcs and levels}\label{4.2}

Let $(\zeta _{1},\zeta _{2})$ be any pair of paths starting at $\infty $ and with endpoints in $Z_{m}$, such that unit-disc-crossings are either along leaves of $L_{3/7}$ or in gaps of $L_{3/7}$. Thus, $\zeta _{1}$ and $\zeta _{2}$ could be paths in $R_{m,0}$ or $R_{m,0}'$, and  be possibly adjacent in one of these, but not necessarily so. We define $k=k(\zeta _{1},\zeta _{2})$ to be the least integer  such that the $k$'th unit-disc-crossings of $\zeta _{1}$ and $\zeta _{2}$ are distinct, up to isotopy preserving $Z_{\infty }\cup \{ \infty \} $. We call $k(\zeta _{1},\zeta _{2})$ the {\em{level}} of $(\zeta _{1},\zeta _{2})$. 
 
 If $\zeta _{1}$ and $\zeta _{2}$ are in $R_{m,0}$, this means that $w_{i}'(\zeta _{1})=w_{i}'(\zeta _{2})$ for $i<k$, but $w_{k}'(\zeta _{1})\neq w_{k}'(\zeta _{2})$. Note that, if $\zeta _{j}\in R_{m,0}$, the $i$'th unit disc crossing of $\zeta _{j}$ is the same as the $i-1$'th crossing of $\psi _{\zeta _{j}}(\zeta _{j}')$ for $i\geq 2$, unless it is the last half-crossing of both $\zeta _{j}$ and $\psi _{\zeta _{j}}(\zeta _{j}')$, in which case they meet in the middle. 
 
Now we specialise to $(\zeta _{1},\zeta _{2})$ being an adjacent pair in $R_{m,0}$. Then the first unit-disc-crossing of $\zeta _{1}$ starts at a point $e^{2\pi it}$ for $t\in [\frac{2}{7},\frac{9}{28}]$ and similarly for $\zeta _{2}$.  If $k$ is odd, then the $k$'th crossing of $\zeta _{1}$ starts to the right of that of $\zeta _{2}$ and if $k$ is even it starts to the left. If $k>1$ and all letters of $w_{k-1}'(\zeta _1)=w_{k-1}'(\zeta _2)$ are in the set $\{ L_3,L_2,R_3\} $, then the $k$'th crossings of $\zeta _1$ and $\zeta _2$ start on the upper unit circle if and only if $k$ is odd. If the word $w_{k-1}'(\zeta _1)$ contains some other letter, then it must contain $BC$ or $UC$, and the prefix ending in this first occurrence of $BC $ or $UC$ is a prefix of $w_j'(\zeta _1)=w_j'(\zeta _2)$ for all $i\leq j\leq  k-1$ for some $i$. If the first occurrence is of $BC$ is preceded by an odd number of letters in $\{ L_3,L_2\} $, or the first occurrence is $UC$, preceded by an even number of letters in $\{ L_3,L_2\} $, then $i$ is odd, and all crossings from the $i$'th onwards have both endpoints on the upper unit circle. In the opposite cases, $i$ is even and all crossings from the $i$'th onwards have both endpoints on the upper unit circle.

We define $D(\zeta _1,\zeta _2)$ to be the region between the $k$'th  unit-disc intersections of $\zeta _1$ and $\zeta _2$, where $k=k(\zeta _1,\zeta _2)$, replacing the $k$'th unit-disc intersection by the leaf of $L_{3/7}$ containing the point on $S^1$ at the start of the intersection, if is is not a complete crossing. 

The {\em{principal arc of $(\zeta _1,\zeta _2)$}} is the arc of $\overline{\zeta _{1}}*\zeta _{2}$ between the  $k$'th  unit-disc-crossings of $\zeta _1$ and $\zeta _2$.We denote it by $I(\zeta _{1},\zeta _{2})$. It is the unique arc bounded by endpoints of essential unit-disc crossings of $\zeta _{1}$ and $\zeta _{2}$ which are also essential unit-disc crossings  on $\overline{\zeta _{1}}*\zeta _{2}$, up to homotopy. The principal arc can have both endpoints on either the upper or lower unit circle. The arc of the unit circle between the endpoints of $I(\zeta _1,\zeta _2)$ is in $D(\zeta _1,\zeta _2)$.

There are two types of principal arcs of adjacent pairs in $R_{m,0}$:
\begin{itemize}
\item[Type 1.] $w_k'(\zeta _1)=u_1$ and $w_k'(\zeta _2)=u_2$ are adjacent values of $w_k'$, subject to the  fixed value of $w_{k-1}'(\zeta _1)=w_{k-1}'(\zeta _2)$ if $k>1$, and $w_{k-1}'(\zeta _1)$ is not between $w_k'(\zeta _1)$ and $w_k'(\zeta _2)$.
\item[Type 2.] $k>1$, and   for the fixed value of $w_{k-1}'(\zeta _1)=w_{k-1}'(\zeta _2)=u$,  with the $k$'th  crossing of $\zeta _1$ starting furthest left and that of $\zeta _2$ starting furthest right if $k$ is odd, and that of $\zeta _1$ starting furthest right and that of $\zeta _2$ starting furthest left if $k$ is even.

\end{itemize}

For both types,  the path  $\zeta _1$ is maximal among among paths in $R_{m,0}$ with $w_k'(\zeta _1)=u_1$, while $\zeta _2$ is minimal among paths $\zeta \in R_{m,0}$  with $w_k'(\zeta _2)=u_2$. We shall refer to $(\zeta _1,\zeta _2)$ being {\em{of type 1}} (or 2) if $I(\zeta _1,\zeta _2)$ is of type 1 (or 2). The figures show examples of a principal arcs $I(\zeta _1,\zeta _2)$. On the first figure, the set $D(\zeta _1,\zeta _2)$ below $I(\zeta _1,\zeta _2)$ is also marked. The principal arc in the first figure could be of type 1 or type 2, depending on whether $\zeta _1<\zeta _2$ or $\zeta _2<\zeta _1$. If $\zeta _2<\zeta _1$ then there could be at least one path $\zeta _3$ with $k(\zeta _1,\zeta _2)=k(\zeta _1,\zeta _3)=k(\zeta _2,\zeta _3)$ and $w_{k-1}'(\zeta _3)=w_{k-1}'(\zeta _1)=w_{k-1}'(\zeta _2)$, and such that $\zeta _1<\zeta _3$, and the $k$'th unit-disc crossing of $\zeta _3$ is strictly between the $k$'th unit-disc crossings of $\zeta _1$ and $\zeta _2$. The second figure must be a type two arc, because, for such a configuration, we must have $\zeta _1<\zeta _2$, and a path $\zeta _3$ with $\zeta _2<\zeta _3$ and $w_{k-1}'(\zeta _1)=w_{k-1}'(\zeta _2)=w_{k-1}'(\zeta _3)$, the the $k$'th unit-disc crossing of $\zeta _3$ between those of $\zeta _1$ and $\zeta _2$, must exist.
  \begin{figure}
\centering{\includegraphics[width=4cm]{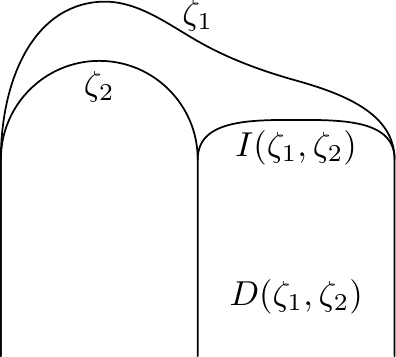}}
\caption{$D(\zeta _1,\zeta _2)$ and the principal arc $I(\zeta _1,\zeta _2)$}
\end{figure}

 \begin{figure}
\centering{\includegraphics[width=4cm]{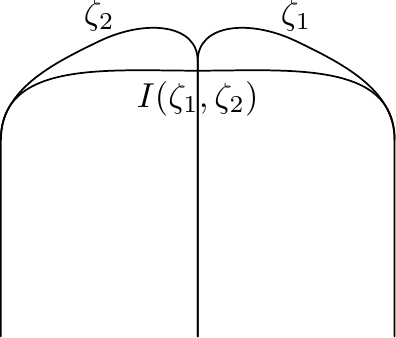}}
\caption{A type 2 principal arc}
\end{figure}

\section{Principal arcs for $[\psi ]$ and $R_{m,0}$}\label{4.3}

For the results to be proved, we need principal arcs of $\alpha *\psi (R_{m,0})$, for $[\psi ]\in G$ and  $[\alpha ]\in \pi _1(\Cbar \setminus Z_m(s),\infty )$ with $s\simeq _\psi \sigma _\alpha \circ s$.  In fact we define {\em{principal arcs for  $[\psi ]$ and $R_{m,0}$ over $\beta $,}} for $\beta \in {\cal{Z}}_m(3/7,+,+)$ as follows, for any adjacent pair $(\zeta _1,\zeta _2)$ in $R_{m,0}$ such that $\psi (\overline{\zeta _1}*\zeta _2)$ crosses $\beta  $ at least once. We define exactly one principal arc for $[\psi ]$ over $\beta $ for each such pair. This is simply the first arc of $\psi (\overline{\zeta _1}*\zeta _2)$ which is crossed by $\beta $, up to isotopy preserving $S^1\cup Z_m(s)$, where $\psi (\overline{\zeta _1}*\zeta _2)$ is chosen up to isotopy preserving $Z_m(s)$ to have only essential intersections with $S^1$.  If $\psi =\psi _{1,4\lfloor i/2\rfloor -1,\beta }$, and the existence of the sequence $\beta _j(\beta )$ as in Theorem \ref{2.8} has been established for $j\le 2i$, then the principal arc for  $(\beta _{2i-1}(\beta ),\beta _{2i}(\beta ),[\psi ])$ is indeed the arc of $\psi (\overline{\beta _{2i-1}(\beta )}*\beta _{2i}(\beta ))$ such that the lift determined by $\alpha _{1,4\lfloor i/2\rfloor -1}*\psi _{1,4\lfloor i/2\rfloor -1,\beta }(\beta _{2i-1}(\beta ))$ is crossed by  the lift of $\beta $. This is the first arc of $\psi (\overline{\zeta _1}*\zeta _2)$, for any adjacent pair $(\zeta _1,\zeta _2)$ in $R_{m,0}$, which is crossed by $\beta $ after crossing the principal arc for $(\beta _{2i-3}(\beta ),\beta _{2i-2}(\beta ),[\psi _{1,4\lfloor (i-1)/2\rfloor -1}])$.  This is useful to know if $i\geq 3$, because, as we shall see in \ref{4.21}, for $i$ even and $\geq 4$, $\psi _{2i-3,2i-1,\beta }$ is the identity on the principal arc for $[\psi _{1,2i-1 ,\beta }]$ over $\beta $ which is  immediately beneath the principal arc for $(\beta _{2i-3}(\beta ),\beta _{2i-2}(\beta ),[\psi _{1,2i-5 ,\beta }])$ over $\beta $.

\section{The $\zeta $-sequence $\zeta _i(x)$}\label{4.4}

Fix any infinite admissible word $x$ ending in $(L_2R_3L_3)^\infty $ and such that $D(x)$ intersects $\{ e^{2\pi it}:t\in (\frac{2}{7},\frac{9}{28})$. We are thinking of $p(x)$ as being the last $S^1$ crossing-point of a path in ${\cal{Z}}(3/7,+,+,0)$. We define the {\em{$\zeta $-sequence for $x$}}, $\zeta _j(x)$, for $\leq 2n(x)$ by the following properties:
\begin{itemize}
\item[1.] $(\zeta _{2i-1}(x),\zeta _{2i}(x))$ is an adjacent pair in $R_{m,0}$ with $\zeta _{2i-1}(x)<\zeta _{2i}(x)$;
\item[2.] $D(x)\subset D(\zeta _{2i-1}(x),\zeta _{2i}(x))$ and $I(\zeta _{2i-1}(x),I(\zeta _{2i}(x))$ is over $D(x)$;
\item[3.] $I(\zeta _1(x),\zeta _2(x))$ is the outermost arc as in 2, and for $i\geq 1$, the principal arc $I(\zeta _{2i+1},\zeta _{2i+1}(x))$, if it exists, is not separated from $I(\zeta _{2i-1}(x),\zeta _{2i}(x))$ by any other principal arc $I(\omega _1,\omega _2)$.
\end{itemize}
There is one more property, which we now explain.The pair $(\zeta _1(x),\zeta _2(x))$ is always defined, but it can happen that $(\zeta _{2i+1}(x),\zeta _{2i+2}(x))$ is not defined for any $i\geq 1$, satisfying conditions 1 to 3. This happens whenever $x$ is not the endpoint of any path in $R_{m,0}$. In the language of \cite{R5}, this happens precisely when  $D(x)$ is not contained in $U^y$ for any value $y$ of $w_1'(.)$. The $x$ for which $\zeta _i(x)$ is defined only for $i=1$, $2$, and for which $x\neq w(\zeta _1(x))$, is precisely the set of $x$ such that $D(x)\subset D(vbv_2)$, for $v_2$ such that $y_1<w_1'(v_2)\leq y_2$, and $y_1$ and $y_2$ as in \ref{4.11}.
 Also, more simply, $(\zeta _{2i+1}(x),\zeta _{2i+2}(x))$ is not defined if $x=w(\zeta _{2i-1}(x))$. In that case, we define $i=n(x)$. If $(\zeta _3(x),\zeta _4(x))$ is defined, then $x=w(\zeta _{2n(x)-1}(x))$.  If there is more than one adjacent pair $(\eta _1,\eta _2)$  in $R_{m,0}$ such that $I(\eta _1,\eta _2)$ is bounded by $I(\zeta _{2i-1}(x),\zeta _{2i}(x))$ and is not separated from $I(\zeta _{2i-1}(x),\zeta _{2i}(x))$ by any other principal arc, and $D(x)$ is bounded by $I(\eta _1,\eta _2)$,  then there are exactly two such. If $i>1$, then it is always true that  $w_1'(\zeta _{2i-1}(x))=w_1'(\zeta _{2i}(x)$ and there is exactly one such pair $(\eta _1,\eta _2)$ with $w_1'(\eta _1)=w_1'(\eta _2)=w_1'(\zeta _{2i-1}(x))$.  and in this case we choose $(\eta _1,\eta _2)$ to be this $(\zeta _{2i+1},\zeta _{2i+2}(x))$. If $i=1$ then the two adjacent pairs $(\eta _1,\eta _2)$ and $(\eta _3,\eta _4)$ in $R_{m,0}$ are  such that  $I(\eta _1,\eta _2)=I(\eta _3,\eta _4)$ and $\eta _1<\zeta _1(x)<\zeta _2(x)<\eta _3(x)$. Condition 4 is then as follows:
\begin{itemize}
\item[4.]  Write $w_2'(x)=w_1(x)u_2$. Then $(\zeta _3(x),\zeta _4(x))$ is chosen so that $\zeta _3(x)<\zeta _1(x)$ if  $L_3u_2 >w_1'(\zeta _1(x))$ and $\zeta _3(x)>\zeta _1(x)$ if $L_3u_2\leq w_1'(x)$. If $i>1$ then $w_1'(\zeta _{2i+1}(x))=w_1'(\zeta _{2i+2}(x))$.  \end{itemize}

The reason for the rather complicated definition of $(\zeta _3,\zeta _4)$ is that $w_2'(.)$ is not injective on all words, as shown by the definitions in section \ref{3}. But there are never more than two adjacent pairs $(\eta _1,\eta _2)$ and $(\eta _3,\eta _4)$ such that $I(\eta _1,\eta _2)$ and $I(\eta _3,\eta _4)$ both bound $D(x)$ and are bounded by $I(\zeta _{2i-1}(x),\zeta _{2i}(x))$, and are not separated from $I(\zeta _{2i-1}(x),\zeta _{2i}(x))$ by any other principal arc.  In this case there are $y_1\neq y_2$ such that 
$$w_1'(\eta _1)=w_1'(\eta _2)=y_1,\ \ w_1'(\eta _3)=w_1'(\eta _4)=y_2,$$
and 
$$I(\eta _1,\eta _2)=I(\eta _3,\eta _4)\subset (U^{y_1}\cap U^{y_2}).$$

We can also define a {\em{lower $\zeta $ sequence }}for $x$ such that  $D(x)$ intersects $\{ e^{2\pi it}:t\in (\frac{19}{28},\frac{5}{7})\} $. The definition is similar but the first pair $(\zeta _1,\zeta _2)$ in a lower $\zeta $ sequence satisfies $k(\zeta _1,\zeta _2)=2$.

We can also define the sequence $\zeta _i(x)$ for $x\in D(BC)$. and $x=p(\gamma )$ with  $\gamma \in {\cal{Z}}(3/7,+,-)$.  For the moment , we just do it for $x=p(\beta )$ with $\beta \in {\cal{Z}}(3/7,+,-)$ and for  $x\in D(BC)\setminus D(BCL_1R_2UC)$. The outermost principal arc $I(\zeta _{3}(x),\zeta _{4}(x))$ is the same for all such $x$, and is the unique principal arc $I(\zeta _{3},\zeta _4)$ such that $w_1'(\zeta _3)=w_1'(\zeta _4)$ and the third intersections of $\zeta _3$ and $\zeta _4$ with the unit circle are the furthest right and left respectively. We then define the sequence $\zeta _{2i-1}(x),\zeta _{2i}(x)$  for $i\geq 3$ as before  satisfying conditions 1 to 3 above. The complications of condition 4 do not arise, as there is always exactly one choice for $(\zeta _{2i+1}(x),\zeta _{2i+2}(x))$ so that 1 to 3 are satisfied, once $(\zeta _{2i-1}(x),\zeta _{2i}(x))$ is given. 

Finally for $x\in D(BCL_1R_2)^kUC\setminus D((BCL_1R_2)^{k+1}UC)$, any pair $(\eta _1,\eta _2)$ such that $I(\eta _1,\eta _2)$ bounds $D(x)$ is such that $\eta _1$ and $\eta _2$ have at least $k+1$ common unit-disc crossings.  We therefore only define $(\zeta _{2i-1}(x),\zeta _{2i}(x))$ for $i\geq k$.

\begin{lemma} \begin{itemize}\label{4.5}
\item[1.] $k(\zeta _1,\zeta _2)=1$
\item[2.] The principal arc $I(\zeta _{2i-1},\zeta _{2i})$ is of Type 1 as in  \ref{4.2}  if $i$ is odd, and of type 2 if $i$ is even.
\item[3.] $k(\zeta _{2i-1},\zeta _{2i})$ is always odd, $1<k(\zeta _{2i-1},\zeta _{2i})$ if $i>1$, and for all $i<n(x)$, 
$$k(\zeta _{2i-1},\zeta _{2i})\leq k(\zeta _{2i+1},\zeta _{2i+2})\leq k(\zeta _{2i-1},\zeta _{2i})+2.$$
\item[4.] If $i$ is even, and $i<n(x)$, then $k(\zeta _{2i+1},\zeta _{2i+2})=k(\zeta _{2i-1},\zeta _{2i})$.

\end{itemize}
\end{lemma}

\begin{proof} Clearly $I(\zeta _1,\zeta _2)$ has to be of type 1, and $k(\zeta _1,\zeta _2)=1$. Also, it is clear that if $i>1$ then $k(\zeta _{2i-1},\zeta _{2i})>1$.

Now  let $i>1$ and $k=k(\zeta _{2i-1},\zeta _{2i})$ and assume inductively that $k$ is odd.  If the longest prefix of $w_k'(\zeta _{2i-1})$ which ends in $BC$ or $UC$ does not coincide with the longest such prefix of $w_k'(\zeta _{2i})$, then there is no principal arc strictly inside $I(\zeta _{2i-1},\zeta _{2i}(x))$ and $i$ is maximal. Now we assume that $i$ is not maximal, because otherwise there is nothing to prove. 

Now suppose that $\zeta \in R_{m,0}$ and $j$ is minimal with $w_j'(\zeta )$ strictly between $w_k'(\zeta _{2i-1})$ and $w_k'(\zeta _{2i})$, and such that the $j$'th unit-disc-crossing of $\zeta $ intersects the unit-circle arc bounded by $I(\zeta _{2i-1},\zeta _{2i})$. Then the $j$'th unit-disc-crossing of $\zeta $ must intersect both the circle arcs bounded by the $k$'th unit-disc crossings of $\zeta _{2i-1}$ and $\zeta _{2i}$. Write $u=w_j'(\zeta )L_2$. We claim that $j\ge k-1$. Since $\zeta $ intersects the circle arcs bounded by the $k$'th unit-disc crossings of $\zeta _{2i-1}$ and $\zeta _{2i}$, and letter $BC$ or $UC$ in $w_j'(\zeta )$ must be in the longest common prefix of $w_k'(\zeta _{2i-1})$ and $w_k'(\zeta _{2i})$. So $w_j'(\zeta )$ must end in $L_3$ or $R_3$. Then $u=w_j'(\zeta )L_2$ is an admissible word. Now either $w_{k-1}'(\zeta _{2i-1})=w_{k-1}'(u,w_1'(x))$, or $w_{k-1}(\zeta _{2i-1})$ is a proper prefix of $w_{k-1}(u,w_1'(x))$ and $w_{k-1}'(\zeta _{2i-1})$ is a proper prefix of $w_{k-1}'(u,w_1'(x))$. The latter only happens  if, for $n\geq 3$ and odd,
$$w_{k-1}(u)=u_0u_1u_2\cdots u_n,\ \ w_{k-1}(\zeta _{2i-1})=w_{k-1}(\zeta _{2i})=u_0u_1,$$
$$w_{k-1}'(u)=u_0\cdots u_{n-1},\ \ w_{k-1}'(\zeta _{2i-1})=w_{k-1}'(\zeta _{2i})=u_0,$$ 
where $u_\ell =(L_3L_2R_3)^{r_\ell }u_\ell '$ for some $r_\ell \geq 0$ and 
$$u_\ell '\in \{ L_3L_2R_3L_3,\ L_3(L_2R_3)^2,L_3^3\} $$
 for $1\leq \ell \leq n-1$.  Now $u_0u_1u_2$ cannot be a prefix of $w_k'(\zeta _{2i-1})$ or $w_k'(\zeta _{2i})$ --- which both extend $w_{k-1}(\zeta _{2i-1})=u_0u_1$ --- because otherwise $w_{k-1}'(\zeta _{2i-1})$ would be defined differently. It follows that in this case $D(u_0u_1u_2)$ is strictly between $w_k'(\zeta _{2i-1})$ and $w_k'(\zeta _{2i})$. Also $w_{k-2}(u)$ is a prefix of $u_0$. Hence $j\ge k-1$.

 Also, $j$ must have the opposite parity to $k$, because the $j$'th crossing of $\zeta $ must start  on the circle arc bounded by $w_k'(\zeta _{2i-1})$ and $w_k'(\zeta _{2i})$ which is not bounded by $I(\zeta _{2i-1},\zeta _{2i})$. But then $w_{k-1}'(\zeta )$ is equal to, or extends, $w_{k-1}'(\zeta _{2i-1})$. We have seen that if  $w_{k-1}'(\zeta _{2i-1})$ is a proper prefix of $w_{k-1}'(\zeta )$, then $D(w_{k-1}'(\zeta _{2i-1})u_1L_3L_2)$ is strictly between $w_k'(\zeta _{2i-1})$ and $w_k'(\zeta _{2i})$ and hence $w_{k-1}'(\zeta )$ is strictly between  $w_k'(\zeta _{2i-1})$ and $w_k'(\zeta _{2i})$ and hence $j=k-1$. If $w_{k-1}'(\zeta _{2i-1})=w_{k-1}'(\zeta )$ then we claim that $j=k+1$. For suppose not so.  Then $D(w_{k+1}(\zeta ))$ intersects the region between $w_k'(\zeta _{2i-1})$ and $w_k'(\zeta _{2i})$ and $D(w_{k+1}'(\zeta ))$ intersects the complement. Since $w_k(\zeta )$ is a prefix of both $w_{k+1}(\zeta ))$ and $D(w_{k+1}'(\zeta )$, it must be the case that $w_k(\zeta )$ is a prefix of one of $w_k'(\zeta _{2i-1})$ and $w_k'(\zeta _{2i})$. But then $w_k(\zeta )=w_k(\zeta _{2i-1})$ is a prefix of $w_k'(\zeta _{2i-1})$, or similarly for $\zeta _{2i}$. Both of these are impossible. 

It follows that the $k-1$'th or the $k+1$'th unit disc crossing of $\zeta _{2i+1}$ and $\zeta _{2i+2}$ is strictly between the $k$'th unit-disc crossings of $\zeta _{2i-1}$ and $\zeta _{2i}$. If the $k-1$'th crossing is between and   $I(\zeta _{2i-1},\zeta _{2i})$ is of type 1, then $w_{k-1}'(\zeta _{2i+1})\neq w_{k-1}'(\zeta _{2i-1})$ and $I(\zeta _{2i-1},\zeta _{2i+1})$ must be of type 2. If the $k+1$'th unit disc crossing of $\zeta _{2i+1}$ and $\zeta _{2i+2}$ is  the first to be strictly between, then again $I(\zeta _{2i+1},\zeta _{2i+2})$ must be of type 2. Now suppose that $I(\zeta _{2i-1},\zeta _{2i})$ is of type 2. Then if $i<n(x)$,  there is at least one other path $\zeta $ between $\zeta _{2i-1}$ and $\zeta _{2i}$   with $w_{k-1}'(\zeta )=w_{k-1}'(\zeta _{2i-1})$. It follows that $w_{k-1}'(\zeta _{2i+1})=w_{k-1}'(\zeta _{2i+2})$ and $k(\zeta _{2i+1},\zeta _{2i+1})=k(\zeta _{2i-1},\zeta _{2i})$ and $I(\zeta _{2i+1},\zeta _{2i+2})$ is of type 1. This completes the proof of 2, 3 and 4. 
\end{proof}

\section{Quadruples}\label{4.6}
We have the following immediate corollary to Lemma \ref{4.5}
\begin{corollary}\label{4.6.1} For any quadruple $(\zeta _{4i+s}(x):-3\leq s\leq 0)$, for any $i\geq 1$, $(\zeta _{4i-3}(x),\zeta _{4i-2}(x))$ is of type 1 and $(\zeta _{4i-1}(x),\zeta _{4i}(x))$ is of type 2. In addition, there are two possibilities for the quadruple $(\zeta _{4i+s}:-3\leq s\leq 0)$ for $i\geq 1$.
\begin{itemize}
\item[A.] $k(\zeta _{4i-1}(x),\zeta _{4i}(x))=k(\zeta _{4i-3}(x),\zeta _{4i-2}(x))+2$;
\item[B.] $k(\zeta _{4i-1}(x),\zeta _{4i}(x))=k(\zeta _{4i-3}(x),\zeta _{4i-2}(x))$.
\end{itemize}
\end{corollary}
We shall refer to these as {\em{quadruples}} of {\em{type A}} and {\em{type B}}. Any quadruple  $(\zeta _{4i+t}(x):-3\le t\le 0)$ in a sequence $\zeta _j(x)$ will be referred to as a {\em{$\zeta $-quadruple}}. If $\zeta _j(x)$ is an upper $\zeta $ sequence then we will refer to $(\zeta _{4i+s}(x):s\le 0$ as an {\em{upper $\zeta $ quadruple}}. Lower $\zeta $ quadruples are similarly defined. We recall the notation $E_{\zeta ,\eta }$ of  \ref{3.8} to denote an arbitrarily small disc neighbourhood of the closed loop $\alpha _\zeta *\overline{\alpha _\eta }$, for paths $\zeta $ and $\eta \in R_{m,0}$. We continue to use this notation for quadruples of various types. So if $(\zeta _i:1\le i\le 4)$ is a quadruple, we use the notation $E_{\zeta _1,\zeta _3}$ to denote an arbitrarily small neighbourhood of $\alpha _{\zeta _1}*\overline{\alpha _{\zeta _3}}$. The figures show examples of type A and type B quadruples, together with the shaded region $E_{\zeta _1,\zeta _3}$ in each case. In the examples drawn we have $\zeta _3<\zeta _4<\zeta _1<\zeta _2$ for the type A quadruple, and $\zeta _1<\zeta _2<\zeta _3<\zeta _4$ for the type B quadruple. Examples like this will arise later, but these are not the only possibilities.

 \begin{figure}
\centering{\includegraphics[width=6cm]{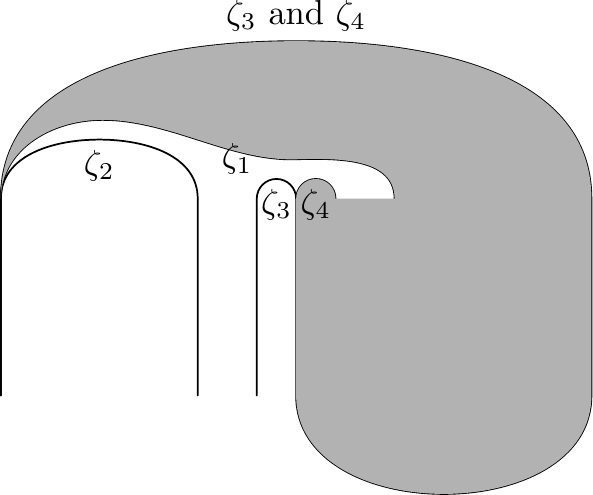}}
\caption{Example of a type A quadruple}
\end{figure}

 \begin{figure}
\centering{\includegraphics[width=6cm]{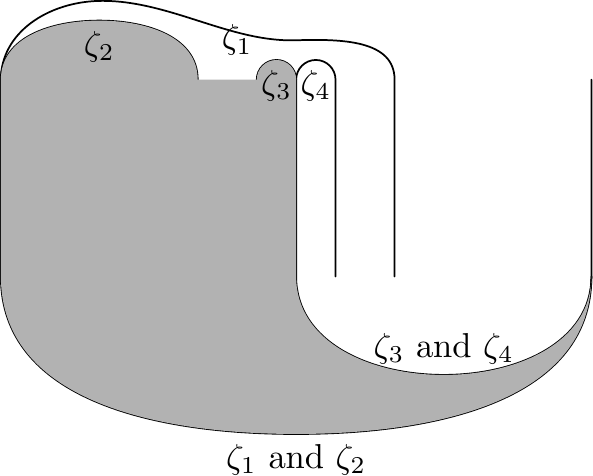}}
\caption{Example of a type B quadruple}
\end{figure}

\section{Quadruples of types C and AC }\label{4.12}
 
We shall see later that the sets of the loops  $(\beta _{4i+j}(\beta ):1\leq j\leq 4)$ are always closely related to quadruples or lower quadruples, of  types A or  B or {\em{type C}} or {\em{type AC}}, which we now describe.
 We say that $(\eta _i:1\leq i\leq 4)$ is of {\em{type C}} if  $(\eta _1,\,\eta _2)$ and $(\eta _3,\eta _4)$ are two adjacent pairs with $\eta _1<\eta _2$ and $\eta _3<\eta _4$ and  one of the following holds.
 
 \begin{itemize}
 \item[C-1] There is a basic exchangeable pair $(a,b)$ which is one of the pairs of  \ref{3.11.4},  such that   the following hold.
 $$k(\eta _1,\eta _2)=1,\ \ k(\eta _3,\eta _4)=3,$$
 $$w_1(\eta _3)=v_1aL_3,\ \ w_1'(\eta _3)=x=w_1'(v_1aL_3L_2),\ \ w_2'(\eta _3,x)=v_1av',$$
 $$w_1'(v_1aL_3L_2)(L_2R_3L_3)^\infty <v'(L_2R_3L_3)^\infty <w_1'(v_1bL_3L_2)(L_2R_3L_3)^\infty ,$$
$$w_1'(\eta _1)(L_2R_3L_3)^\infty <v_1bL_3L_2$$
and $\eta _1$ is maximal with respect to this property, and $(\eta _1,\eta _2)$ of type 1 and $(\eta _3,\eta _4)$ of type 2. Moreover, 
$$I(\eta _3,\eta _4)\subset uE_{\eta _1,\eta _3},$$
where $u\leftrightarrow u'$ is an exchange for $\eta _1$ which is prefixed by $v_1a\leftrightarrow v_1b$ such that $w_2(uL_3L_2)=uL_3L_2$ and $w_2(u'L_3L_2)=u'L_3L_2$.
Such a quadruple is said to be {\em{type C of level 1}}

 \item[C-2.] There is a quadruple $(\zeta _i:1\leq i\leq 4)$ which is of type C of level 1 such that
  $$w_1'(\eta _i)=w_1'(\zeta _3)=w_1'(\zeta _4),$$
$(\eta _1,\eta _2)$ is of type 1 and $(\eta _3,\eta _4)$ is of type 2 and for some odd $k_1\geq 3$, 
   $$k(\eta _1,\eta _2)=k_1,\  \ k(\eta _3,\eta _4)=k_1+2,$$ 
    and, writing $v=v_1aL_3L_2v_2$ for a word $v_2$,
    $$I(\eta _j,\eta _{j+1})=vI(\zeta _j,\zeta _{j+1}),\ \ E_{\eta _1,\eta _3}=vE_{\zeta _1,\zeta _3}.$$
 
 \end{itemize}
 
The definition of type AC is exactly similar to type C except that the basic exchange $va\leftrightarrow vb$ is one of the types listed in \ref{3.11.4}. This makes for some important differences, the chief of which are that  $I(\eta _3,\eta _4)$ spans an arc on the lower unit circle, and that now $w_1'(vaL_3L_2)=w_1'(vbL_3L_2)$.

\section{Intersections between sets $E_{\eta _1,\eta _3}$ and $s^{-n}E_{\omega _1,\omega _3}$}\label{4.9}

If $E$ any set with  $E\subset D(L_3)$, then a component of $s^{-n}(E)$ is of the form $uE$ for some word $u$ such that $uL_3$ is admissible if $E\subset D(L_3)$. This means simply that $uE=\{ p(ux):p(x)\in E\} $. 
The definition and properties of $w_j(.)$ and $w_j'(.)$ and $w_j'(.,y)$, for any value $y$ of $w_1'$, ensure that if $(\eta _i:1\leq i\leq 4)$ is an upper or lower quadruple  and  $w_1'(uL_3L_2)=w_1'(\eta _3)=y$, then there is an upper or lower $\zeta $ quadruple $(\eta  _i:5\leq i\leq 8)$ such that $w_1'(\eta  _i)=y$ for all $5\leq i\leq 8$ and 
$$uE_{\eta _1,\eta _3}=E_{\eta _5,\eta   _7}\subset  U^y$$
 The condition 4 in the definition of the zeta sequence ensures that $w_{k_1}'(.,y)$ is constant on
$$\{ \eta \in R_{m,0}:(\eta ,\eta ^2){\rm{\ is\ an\ adjacent\ pair\ in\ }}R_{m,0},\ \ I(\eta ,\eta ^2)\subset E_{\eta _5,\eta _7}\} ,$$
 where $k_1=k(u)+k(\eta _1,\eta _2)$ where $k(u)$ is the largest integer $k$ such that $w_k(uL_3L_2)$ is a proper prefix of $uL_3L_2$. (This condition is only needed if $k(\eta _1,\eta _2)=1$.) A necessary and sufficient condition for $(\eta _i:5\leq i\leq 8)$ to be an upper $\zeta $ quadruple with $k(\eta _5,\eta _6)\geq 3$ is that $k(u)$ be even and $k(u)+k(\eta _1,\eta _2)\geq 3$. 

\begin{lemma}\label{4.10} Let $(\eta _{i}:1\leq i\leq 4 )$ and $(\omega _{j}:1\leq j\leq 4)$ be two upper or lower $\zeta $ quadruples, such that
$$w_1'(\eta _3)=w_1'(\omega _3),$$
and $(\eta _1,\eta _2)=(\omega _1,\omega _2)$ if $k(\eta _1,\eta _2)=k(\omega _1,\omega _2)=1$. Write $E=E_{\eta _1,\eta _3}$. Let $E'$ be a component of   $s^{-n}(E_{\omega _1,\omega _3})$ for some $n>0$.   Then  $E\cap E'\neq \emptyset $ if and only if one of $E'$ and $E$ is contained in the other.

\end{lemma}

{\textbf{Remark}} 
\begin{itemize}
\item[1.] For this lemma to be useful, we still have to show that for $\alpha =\alpha _{\eta _1,\eta _3 }$, a component of $(\sigma _{\alpha }\circ s)^{-n}(E_{\eta _1,\eta _3})$ is, up to homotopy preserving $Z_\infty $, the union of two components of $s^{-n}(E_{\eta _1,\eta _3})$ and a path joining the two. This will be shown in later lemmas. 
\item[2.] The conditions of a $\zeta $ quadruple ensure that the second unit-disc-crossing of $\eta _3$ is not contained in $\cup _{n\geq 0}s^{-n}(E_{\eta _1,\eta _3})$.
\end{itemize}

\begin{proof}  Let $u$ be the word  such that $E'=uE_{\omega _1,\omega _3}$. Let $(\omega _i:5\leq i\leq 8)$ be the quadruple such that $uE_{\omega _1,\omega _3}=E_{\omega _5,\omega _7}$. If $k(u)\geq 1$, then $w_1'(\omega _i)=w_1'(u)$ for all $5\leq i\leq 8$ and $k(\omega _7,\omega _8)\geq k(\omega _5,\omega _6)>1$. 

If $k(u)=0$, then $w_1'(\omega _i)=uw_1'(\omega _{i-4})$  for $5\leq i\leq 8$. If $k(\omega _1,\omega _2)>1$ then $w_1'(\omega _i)=w_1'(\omega _1)$ for all $1\leq i\leq 4$ and 
$$w_1'(\omega _i)=w_1'(\omega _5)=uw_1'(\omega _3)=uw_1'(\eta _3)$$
 for $5\leq i\leq 7$. It follows that $E'$ is contained in or disjoint from $E$ in this case, except possibly when $k(\eta _1,\eta _2)=1$. But then the only case to consider is when $w_1'(\omega _i)=w_1'(\eta _1)$ for all $5\leq i\leq 8$ or $w_1'(\omega _i)=w_1'(\eta _2)$ for all $5\leq i\leq 8$. But then $\eta _1$ (or $\eta _2$) is either the maximal (or minimal) path in $R_{m,0}$ with this value of $w_1'$, and depending on whether $\eta _4\leq \eta _1$ or $\eta _2\leq \eta _3$, as well as whether $w_1'(\omega _i)=w_1'(\eta _1)$ or $w_1'(eta _2)$, $E'$ is contained in or disjoint from $E$.
 
  So now suppose that $k(u)>0$. Write 
 $$k_1=k(\eta _1,\eta _2),\ \ k_2=k(\eta _3,\eta _4),\ \ k_3=k(\omega _5,\omega _6),\ \ k_4=k(\omega _7,\omega _8).$$
 
 If $I(\omega _5,\omega _6)=I(\eta _1,\eta _2)$, then we have disjointness of $E$ and $E'$  if $\omega _7$ and $\eta _3$ are on opposite sides of $\eta _1$, and one of $E$ and $E'$ is contained in the other if they are on the same side. So now suppose that $I(\omega _5,\omega _6)\neq I(\eta _1,\eta _2)$. Now we consider the boundary of $E$, up to $Z_\infty $-preserving homotopy. Either the $k_1-1$'th or the $k_1-2$'th unit-disc-crossings of the $\eta _i$ are in common (up to $Z_\infty $-preserving homotopy) depending on whether the quadruple is of type A or B. Lock these common crossings together. Lock $\eta _1$ and $\eta _2$ together up to the last common unit-disc crossing, and similarly for $\eta _3$ and $\eta _4$.   Finally, lock together $I(\eta _1,\eta _2)$ and $I(\eta _3,\eta _4)$. In the figures shown earlier of the sets $E_{\zeta _1,\zeta _3}$, for both type A and type B quadruples,   the horizontal arcs of the boundary of $E_{\zeta _1,\zeta _3}$ denote the locking of the principal arcs $I(\zeta _1,\zeta _2)$ and $I(\zeta _3,\zeta _4)$. Returning to the quadruple $(\eta _i:1\le i\le 4)$, we obtain a set bounded by the locked-together principal arcs and by common segments of the paths $\eta _i$. This set is $E$ up to homotopy, and so we denote it by $E$. It has directed boundary with direction given by the paths $\eta _i$. No path in $\zeta \in R_{m,0}$ with $w_1'(\zeta )=w_1'(\eta _3)$ travels through any of the ``locks'' except for the lock between the common crossings of all the $\eta _i$, as $(\eta _1,\eta _2)$ is of type 1 and $(\eta _3,\eta _4)$ is of type 2. So the only way for a path $\zeta \in R_{m,0}$ to enter $E$ is between $\eta _4$ and $\eta _1$,  or between $\eta _2$ and $\eta _3$ (depending on whether $\eta _3<\eta _1$ or $\eta _1<\eta _3$). We choose locks similarly for $E'$. Then either one of $E$ and $E'$ is contained in the other or they are disjoint. For suppose none of these hold. Then the boundaries of $E$ and $E'$ must intersect transversally, and since both $E$ and $E'$ are topological discs, there must be two transversal intersections. But the only place for a  transversal intersection is between the common lock of all the $\eta _i$ and the common lock of all the $\omega _i$.  So this is impossible, and one of $E$ and $E'$ must be contained in the other and they are disjoint.
 
 We have to use a slight variant on this if $k_1=1$, because in that case there is no common crossing of all the $\eta _i$. But there is still only one way in: through the top.

\end{proof}

\begin{corollary}\label{4.11} Let $(\eta _i:1\le i\le 4)$ be an upper or lower $\zeta $ quadruple,  let $E=E_{\eta _1,\eta _3}$ and let $E'$ be a component of $s^{-n}E$ for some $n>0$. Then $E'\cap I(\eta _j,\eta _{j+1})=\emptyset $ for $j=1$, $3$ and either $E'\subset E$ or $E'\cap E=\emptyset $.\end{corollary}

\begin{proof} We have $E'=E_{\eta _5,\eta _7}$ for a quadruple $(\eta _{4+j}:1\leq j\leq 4)$ and with $k(\eta _5,\eta _6)\geq k(\eta _1,\eta _2)$. If $E'$ contains $I(\eta _1,\eta _2)$, then $k(\eta _5,\eta _6)=k(\eta _1,\eta _2)$. We clearly cannot have $E\subset E'$, so we must have $E'\subset E$ and $E'\neq E$. Since $I(\eta _1,\eta _2)$ is the furthest left (or right ) of the principal arcs in $E$ of level $k_1$, then if $I(\eta _1,\eta _2)$ is in $E'$ it must be equal to $I(\eta _5,\eta _6)$ which implies that $I(\eta _1,\eta _2)$ is fixed by $s^i$. This is impossible. So $I(\eta _1,\eta _2)$ is not contained in $E'$. Now suppose that $I(\eta _3,\eta _4)\subset E'$. Since $I(\eta _3,\eta _4)$ intersects $\partial E$, and $E'\subset E$, this is only possible if $I(\eta _3,\eta _4)\subset \partial E'$. But this, in turn is only possible if $I(\eta _3,\eta _4)\subset s^{-i}(\partial E)$. But $\partial E$ is a union of principal arcs of level $k_1$, with $I(\eta _1,\eta _2)$ at one end, and the $k_1$'th unit-disc crossing of  $\eta _3$ at the other end, and $I(\eta _3,\eta _4)$. We cannot have $I(\eta _3,\eta _4)\subset s^{-i}(I(\eta _3,\eta _4))$, and anything else in $\partial E$ mapping to $I(\eta _3,\eta _4)$ under $s^{-i}$ leads to $E'$ not contained in $E$. So $I(\eta _3,\eta _4)$ is not contained in $E'$ either.
\end{proof}

We need  a variant on \ref{4.10} for quadruples of type C or AC.

\begin{lemma}\label{4.14}Let $(\eta _i:1\leq i\leq 4)$ be of type C or AC. Write $k_1=k(\eta _1,\eta _2)$ and $E=E_{\eta _1,\eta _3}$.
\begin{itemize}  \item[1.] If $k_1>1$ then $I(\eta _1,\eta _2)$ and $I(\eta _3,\eta _4)$ are not contained in $s^{-n}(E$ for any $n>0$.
 \item[2.] . Let $u_1$ and $u_2$ be admissible words. Then the only way that $u_1E$ and $u_2E$  can  intersect without one being contained with the other is if $u_2=u_1va$ or $u_1=u_2va$.

 \item[3.] If  $k_1=1$, then $I(\eta _1,\eta _2)$ is disjoint from $s^{-n}(E)$ and $I(\eta _3,\eta _4)$ is contained in a component of $s^{-n}(E)$ for a single value of $n>0$.
 
\end{itemize}

\end{lemma}
\begin{proof}

\noindent 1. Suppose that $I(\eta _3,\eta _4)\subset uE$ for some nontrivial $u$. Since $k(\eta _1,\eta _2)>1$ by the definition of type C or AC we have $w_1'(\eta _1)=w_1'(\eta _3)=x$, say.There is at most one $y$with $U^x\cap U^y\neq \emptyset $ and they have a common prefix apart form the last at most four letters. If $I(\eta _1,\eta _2)$ or $I(\eta _3.\eta _4)$ is in $uE$ then  we need $ux=x$ or $uy=x$, both of which are impossible.  

\noindent 2. We can assume without loss of generality that $u_2$ is trivial. The question is then how it is possible to have $u_1E\cap E\neq \emptyset $ . We have seen in \ref{4.11} that it is impossible to have the first segment of $\eta _1$ in the unit disc in $u_1E_{\eta _1,\eta _3}$. It is similarly impossible to have the first segment of $\eta _3$ in the unit disc in $u_1E$ because then $w_1'(\eta _3)$ would have to be between $u_1w_1'(\eta _3)$ and $u_1'(w_1'(\eta _1)$. The second unit-disc-crossing of $\eta _3$ is in $vaE$ and . So $va$ must be a prefix of $u_1$, and since $I(\eta _3,\eta _4)\subset \partial vaE$, we must have $va=u_1$.

\noindent 3.  We now have $k(\eta _1,\eta _2)=1$.  From the definition of type C,  or AC, there is a basic exchange $a\leftrightarrow b$ and  an admissible word $v_1$ such that $I(\eta _3,\eta _4)\subset v_1aE$.  WE have $v_1aE\setminus E\ne \emptyset $, and by 2, $v_1a$ is the unique nontrivial word with these properties.

\end{proof}

\section{A basic result about the quadruple homeomorphisms}\label{4.13.0}

One of our aims is to obtain a good bound on the supports of the homeomorphisms $\psi _{4i-3,4i-1,\beta }$, for any capture path $\beta \in {\cal{Z}}_\infty (3/7,+,+,0)$. Any bound obviously has to be obtained inductively, because the support of $\psi _{4i-3,4i-1,\beta }$ is obtained from $E_{\beta _{4i-3},\beta _{4i-1},\beta }$, which is, in turn, obtained from $(\beta _{4i+t}:-3\leq t\leq 0)$, which in turn depends on $\psi _{4j-3,4j-1,\beta }$ for $j<i$. It is clearly important to know how the supports of  $\psi _{4j-3,4j-1,\beta }$, for varying $j$, intersect each other, and also how they intersect the boundaries of sets $E_{\beta _{4i-3},\beta _{4i-1},x}'$ for varying $i$. 

We start with a preliminary lemma.

\begin{lemma}\label{4.8}Let $(\zeta _1,\zeta _2)$ be an adjacent pair in $R_{m,0}$ and let $\zeta _{1,0}$ be the initial segment of $\zeta _1$. No component of $R_{m,0}\cap \{ z:|z|\geq 1\} $ or of $(\sigma _{\zeta _{1,0}}\circ s)^{-n}(s^{-1}(\zeta _{1,0}))$ comes between $I(\zeta _1,\zeta _2)$ and $I(\zeta _3,\zeta _4)$ for any quadruple $(\zeta _i:1\leq i\leq 4)$ of type A or B, unless it is a component of $\sigma _{\zeta _{1,0}}\circ s)^{-n}(s^{-1}(\zeta _{1,0}))$ which is isotopic to $I(\zeta _3,\zeta _4)$ up to isotopy preserving $Z_\infty \cup S^1$.\end{lemma}
\begin{proof} Any component of $(\sigma _{\zeta _{1,0}}\circ s)^{-n}(s^{-1}(\zeta _{1,0}))$ is isotopic, up to isotopy preserving $Z_\infty \cup S^1$, to $I(\zeta _5,\zeta _6)$, for some adjacent pair $(\zeta _5,\zeta _6)$ in $R_{m,0}$. By definition, there are no principal arcs separating $I(\zeta _1,\zeta _2)$ and $I(\zeta _3,\zeta _4)$. The result follows.\end{proof}

We have the following.

\begin{subsectiontheorem}\label{4.13} Let $(\zeta ^i:1\leq i\leq 4)$ be an upper or lower quadruple or type C or AC.   Write 
$$\zeta ^1=\zeta ,\ \ \psi =\psi _{\zeta ,\zeta ^3},\ \ E=E_{\zeta ,\zeta ^3},\ \ k_1=k(\zeta ,\zeta^2).$$ 

Write $\zeta _{n}$ for the portion of $\zeta  $ up to  the $n$'th intersection with $S^1$. 

Then each component $C$ of $(\sigma _{\alpha _\zeta }\circ s)^{-i}(E)$ is a component of $(\sigma _{\zeta _{2k_1-1}}\circ s)^{-i}(E)$, and takes the form 
$$E_1\cup E_2\cup \omega _{1,m(C),n(C,1)}\cup \omega _{2,m(C),n(C,2)}\cup \ell ,$$
where:
\begin{itemize}
\item $E_1$ and $E_2$ are components of $s^{-i}(E)$;
\item $C$ is either contained in $U^0$, or disjoint from it, or $\ell (C)$ is an arc with endpoints on the unit circle on either side of $e^{2\pi i(2/7)}$, and for one of $r=1$ or $2$, $E_r\cup \omega _{r,m(C),n(C,r)}$ is contained in $U^0$ apart from the first segment on $\omega _{r,m(C),n(C,r)}$;
\item If $E_r\cup \omega _{r,n(C),m(C)}$ intersects $U^0$, then  $\omega _{r,n(C),m(C)}$ is the segment of $\omega _r\in R_{m,0}$ from the $m(C)$'th intersection with the unit circle to the $n(C,r)$'th for $r=1$ and $2$, and $m(C)\leq  2k_1$;
\item $\ell $ is a component of $s^{-j}(\zeta _0)$ for some $j\leq i$;
\item the only intersections between $C$ and $\alpha _\zeta $ are between $\ell $ and $\zeta _{2k_1-1}$. 
\end{itemize}
In particular,  the support of $\psi $ is contained in $\cup _{i\geq 1}(\sigma _{\zeta _{2k_1-1}}\circ s)^{-i}(E)$. \end{subsectiontheorem}

\begin{proof}

We start by showing that any intersections between $(\sigma _{\alpha _\zeta }\circ s)^{-i}(E)$ and $\alpha _\zeta $ must be intersections between $(\sigma _{\alpha _\zeta }\circ s)^{-i}(E)$ and $\zeta $, and in fact must be intersections with $\zeta _{2k_1}$. Recall that $\alpha _\zeta $ is an arbitrarily small perturbation of $\zeta *\psi _\zeta (\zeta ')$. We need to cut down the intersections between $(\sigma _{\alpha _\zeta }\circ s)^{-i}(E)$ and $\alpha _\zeta $. We first concentrate on cutting down the intersection between  $(\sigma _{\alpha _\zeta }\circ s)^{-i}(E)$ and $\psi _\zeta (\zeta ')$. Now $E\subset U^0$, where $U^0$ is defined in \ref{3.2}, and apart from the component of $\alpha _\zeta \cap \{ z:\vert z\vert >1\} $ which contains $\infty $, the loop $\alpha _\zeta $ is  also contained in $U^0$. It follows by induction on $i$ that each component of $(\sigma _{\alpha _\zeta }\circ s)^{-i}(E)$ is contained in a component of  
$$(\sigma _{\zeta _1}\circ s)^{-i}(U^{0}),$$
up to isotopy preserving $Z_m$ and $S^1$, where $\zeta _{1}$ is the initial segment of $\zeta $, up to the first intersection with $S^{1}$. These are the same as the components of 
$$(\sigma _{\beta _{2/7}})^{-i}(U^{0}).$$
Next we show that there are no intersections between this set and the piece of $\psi _\zeta (\zeta ')$ which is  common with $\psi _{m,2/7}(\beta _{5/7})$ up to $S^1$ and $Z_m$-preserving isotopy, where $\zeta '$ is the path in $R_{m,0}'$ which is matched with $\zeta $.  To do this, we simply show that there are no intersections at all with $\psi _{m,2/7}(\beta _{5/7})$ .
To do this, since $U^0\cap Z_2=\emptyset $, it suffices to show that that  $\psi _{m,2/7}(\beta _{5/7})$ is disjoint from $T^{0}$ where 
$$T^{0}=\cup _{i\leq m-3}(\sigma _{\beta _{2/7}}\circ s)^{-i}(U^{0}).$$
It suffices to prove this by induction on $m$. Let $D_1$ be the disc bounded by $\overline{\alpha _{2/7}}*\alpha _{9/28}$, that is, bounded by an arbitrarily small perturbation of $\overline{\beta _{2/7}}*\beta _{5/7}*\overline{\beta _{19/28}}*\beta _{9/28}$. We denote by $\partial 'D_1$ the part of the boundary of $D_1$ which is an arbitrarily small perturbation of $\beta _{5/7}*\overline{\beta _{19/28}}$.  The support of $\psi _{2,2/7}$, illustrated in the figures in \ref{3.15} is disjoint from $D_1$ and so $\partial 'D_1=\psi _{2,2/7}\partial 'D_1$. We note that $\psi _{2,2/7}(\beta _{5/7})$  is disjoint from $D_{1}$, and $\psi _{3,2/7}(\beta _{5/7})$ is disjoint from $D_{1}\cup s^{-1}(D_{1})$. Now $U^{0}$ is obtained by adjoining to $D_{1}$ the region bounded by $\partial 'D_{1}$ and $\psi _{m,2/7}(\partial 'D_{1})$. It follows that $\psi _{m,2/7}(\beta _{5/7})$ is disjoint from $\psi _{m,2/7}(\partial 'D_{1}\cup s^{-1}(\partial 'D_{1})$. But the intersection of $\psi _{m,2/7}(\beta _{5/7})$ with the unit disc can be homotoped  to $\partial  'U^{0}=\psi _{m,2/7}(\partial 'D_1)$ via a homotopy preserving $Z$. Each component of $s^{-i}U^{0}$ is either contained in, or disjoint from, $U^{0}$. It follows that $\psi _{m,2/7}(\beta _{5/7})$ is disjoint from $T_{0}$. Hence, any intersection between $\psi _{\zeta }(\zeta ')$ and $(\sigma _{\alpha _{\zeta }}\circ s)^{-i}(E)$ must be contained in $U^0$ and coincides, up to isotopy preserving $Z\cup S^1$,  with an intersection between $\psi _{m,2/7}(\overline{\beta _{5/7}})*\psi _{\zeta }(\zeta ')$ and  $(\sigma _{\alpha _{\zeta }}\circ s)^{-i}(E)$. Following \cite{R5}, we denote by $\partial 'U^0$ the subset of $\partial U^0$ which is isotopic, up to isotopy preserving $Z\cup S^1$, to $\psi _{m,2/7}(\overline{\beta _{5,2/7}})*\psi _{m,9/28}(\beta _{19/28})$. Then $\psi _{m,2/7}(\overline{\beta _{5/7}})*\psi _\zeta (\zeta ')$ can be isotoped, up to isotopy preserving $Z\cup S^1$, into $\partial 'U^0\cup \zeta $. Hence the same is true for $\alpha _\zeta $. So any intersection between $\alpha _\zeta $ and  $(\sigma _{\alpha _{\zeta }}\circ s)^{-i}(E)$ must be contained in $\partial 'U^0\cup \zeta $. But $\partial 'U$ has no transversal intersections with $(\sigma _{\beta _{2/7}}\circ s)^{-i}(U^{0})$, and hence none with  $(\sigma _{\alpha _\zeta }\circ s)^{-i}(E)$ either. (In fact there is at most one component of $(\sigma _{\beta _{2/7}}\circ s)^{-i}(U^{0})$ which intersects $U^0$ without being contained in it, for which the corresponding component of $\sigma _{\beta _{2/7}}\circ s)^{-i}(\beta _{2/7})$ intersects the unit circle at $e^{2\pi i(2/7)}$. This is possible if and only if $i$ is divisible by $3$.) It follows that any intersection between $\alpha _\zeta $ and $(\sigma _{\alpha _\zeta }\circ s)^{-i}(E)$ is an intersection between $\zeta $ and $(\sigma _{\alpha _\zeta }\circ s)^{-i}(E)$, up to isotopy preserving $Z\cup S^1$. Furthermore, since $\zeta \setminus \zeta _{2k_1}$ coincides with $\psi _{\zeta }(\zeta ')$ , all intersections between $(\sigma _{\alpha _\zeta }\circ s)^{-i}(E)$ and $\alpha _{\zeta }$ must be with $\zeta _{2k_1}$.

Suppose that the theorem is true for $i$, and we consider pre-images of a component $C$ of $(\sigma _{\alpha _{\zeta }}\circ s)^{-i}(E)$ under $\sigma _{\alpha _\zeta }\circ s$. So  by the inductive hypothesis $C$ is a component of $(\sigma _{\zeta _{2k_1-1}}\circ s)^{1-i}(s^{-1}(E)\cup s^{-1}(\zeta _{2k_1-1}))$ and 
$$C=E_1\cup E_2\cup \omega _{1,n(C),m(C,1)}\cup \omega _{2,n(C),m(C,2)}\cup \ell _C$$
where $E_1$ and $E_2$ are components of $s^{-i}E$ and $\omega _r\in R_{m,0}$ if $\omega _{r,n(C),m(C)}\cap U^0\ne \emptyset $.  and $\omega _{r,n(C),m(C)}$ denotes the segment of $\omega _r$ from the $m(C)$'th intersection with $S^1$ to the $n(C,r)$'th. It is possible that $m(C)=n(C,r)$, in which case by convention  these segments are absent.  If $C$ intersects $U^0$ then it is contained in $U^0$, unless $\ell (C)$ passes over $e^{2\pi i(2/7)}$. In that case, exactly one of the sets $E_r\cup \omega  _{r,n(C),m(C,r)}$ intersect $U^0$. We assume without loss of generality that this is true for $r=1$. If this happens then $w_1'(\omega _1)\ne w_1'(\zeta )$, and $E_1\cap \zeta =\emptyset $. In fact we have $s^{-i}E\cap \zeta _1=\emptyset $ for all$ i>0$. Let the first unit disc crossings of $\zeta ^1$ and $\zeta ^3$ be represented by 
$$D(v_1)=D((L_3L_2R_3)^{k_1}u_1(L_2R_3L_3)^\infty )$$ and 
$$D(v_2)=D((L_3L_2R_3)^{k_2}u_2(L_2R_3L_3)^\infty )$$ respectively, where $k_1$ and $k_2$ are maximal. Then  $k_1=k_2$ or $k_2=k_1\pm 1$ and $D(v_2)$ is strictly between $D(v_1)$ and $D((L_3L_2R_3)^{\pm 1}v_1)$. This is all that is needed.

Since any intersection between $\zeta _{2k_1}$ and $\ell (C)$ must occur outside the unit disc, it cannot occur on $\zeta _{2k_1}\setminus \zeta _{2k_1-1}$, and must occur on $\zeta _{2k_1-1}$. There are no transversal  intersections between $\zeta _{2k_1}$ and $\omega _{r,n(C),m(C,r)}$, because if $\omega _{r,n(C),m(C,r)}$ intersects $U^0$ then $\omega _r$ is a path in $R_{m,0}$ and cannot have transversal intersections with $\zeta \in R_{m,0}$, because paths in $R_{m,0}$ have no transversal intersections \cite{R5}.
It remains to consider how $\zeta _{2k_1}$ can intersect $E_1\cup E_2$.  By \ref{4.6}, there are no intersections between $\partial E$ and $\partial E_1\cup \partial E_2$ if $(\zeta ^j:1\leq j\leq 4)$ is an upper  or lower quadruple. and $\zeta _{2k_1}\setminus \zeta _{2k_1-1}$ is in $\partial E$   Since $E_1$ and $E_2$ are contained in the unit disc, and $\zeta _{2i+1}\setminus \zeta _{2i}$ is in the exterior for all $i\geq 1$, the only possible intersections between $\zeta _{2k_1}$ and $E_1\cup E_2$ are with $\zeta _{2k_1-2}$. We now show that no intersection is possible in these cases. The case of $k_1=1$ is dealt with since then $2k_1-2=0$ and $\zeta _0$ is empty. So we can assume that $k_1>1$. The subsets $E_1$ and $E_2$ of $C$ are bounded by paths $\omega _r$ and $\omega  _{r+2}$ in $R_{m,0}$ (for $r=1$ and $2$), where $\omega _1$ and $\omega _2 $ are as above. By the invariance properties, $k(\omega _r,\omega  _{r+2})\geq k(\zeta ^1,\zeta ^3)\ge  k_1-1$ or $\geq k_1$, depending on whether the quadruple $(\zeta ^i:1\leq i\leq 4)$ is of  type A or of type B. We can extend the notion of type A or B to lower quadruples, and if $k_1=2$ then $(\zeta ^i:1\le i\le 4$ is a type A quadruple. Suppose, for contradiction, that  $\zeta _{2k_1}$  intersects the component of $E_r$ between $\omega _r$ and $\omega  _{r+2}$. This is certainly impossible if  $w_1'(\zeta )=w_1'(\omega _r)$ because then $\zeta _{2k_1-2}$ or $\zeta _{2k_1-4}$ has to coincide with $\omega _{r,2k_1-2}$ or $\omega _{r,2k_1-4}$ (or $\omega _{r+2,2k_1-2}$ or $\omega _{r+2,2k_1-4}$). WE necessarily have $\zeta _{2k_1-2}$ coinciding with  $\omega _{r,2k_1-2}$  or $\omega _{r+2,2k_1-2}$ if $k_1=2$, because then $(\zeta ^j:1\le j\le 4)$ is of type A.  It is also impossible if $w_1'(\zeta )\neq w_1'(\omega _r)$ because although paths $\zeta $ and $\omega $ in $R_{m,0}$ can occasionally intersect if $w_1'(\zeta )\neq w_1'(\omega )$, intersections can only be between $r$'th components of intersection with the exterior of the unit disc for the same $r$. So $\zeta _{2k_1}$cannot intersect $E_1\cup E_2$ This completes the proof for  upper and lower quadruples. 

Now we consider the case of $(\zeta ^i:1\leq i\leq 4)$ being of  type C or AC. There is no intersection if $k(\omega _r,\omega _{r+2})>k_1$. Since $k(\zeta ,\zeta ^3)=k_1$ the only possibility is that $k(\omega _r,\omega _{r+2})=k_1$. If $k_1>1$ then there can only be an intersection if $w_1'(\omega _j)=y$ and $w_1'(\zeta )=x$ and $U^x\cap U^y\neq \emptyset $. But $y$ also has to be in $s^{-i}(x)$, and this is never true in the cases when $U^x\cap U^y\neq \emptyset $, except in the case $i=0$. So the only possibility is that $k_1=1$. The proof for type AC is then exactly as for type A, so now we consider type C only, and we write $w_1'(\zeta ^3)= uL_3$ and $w_2'(\zeta ^3)=uau'$ for a basic exchangeable pair $(a,b)$. Then  $D(w_1'(\zeta )(L_2R_3L_3)^\infty )$ is between $D(uL_3(L_2R_3L_3)^\infty )$ and $D(ubL_3L_2)$ and is maximal to the right of $D(ubL_3L_2)$. If $D(w_1'(\zeta )(L_2R_3L_3)^\infty )$ is also between $D(vuL_3(L_2R_3L_3)^\infty )$ and $D(vubL_3L_2)$ then $u$ must be a prefix of $vu$ and so $u=v^nv'$ for some $n\geq 1$ or $v=uu''$. But $D(ubL_3L_2)$ also has to be in $D(vu)$. So $v^{n+1}v'$ has to be a prefix of $v^nv'bL_3L_2$. So $vv'$ has to be a prefix of $v'bL_3L_2$. This is impossible because $w_1(vv')=vv'$ and $w_1(v'bL_3L_2)$ is a proper prefix of $vbL_3L_2$.

\end{proof}

\section{The first four elements of the $\beta _i $ sequence}

In this section we prove Theorem \ref{2.8} for the first four elements of the sequence $\omega _i=\omega _i(\beta )$, given a capture path $\beta \in {\cal{Z}}_m(3/7,+,+,0)$.  As in \ref{2.8}, we let $\gamma _i$ be the sequence of elements of $\pi _1(V_{3,m},a_1)$ and $(\omega _{2i-1},\omega _{2i})$ the sequence of adjacent pairs in $\Omega _m$ such that $\rho(\beta )$ successively crosses $\overline{\rho (\gamma _i*\omega _{2i-1})}*\rho (\gamma _i*\omega _{2i})$.

\begin{lemma}\label{4.15} Then there is a quadruple $(\beta _i:1\le i\le 4)$ of type A, AC or C such that
$$(\rho (\omega _{1}),\rho (\omega _{2}))=(\beta _{1},\beta _{2}),$$
$$(\rho (\omega _{3}),\rho (\omega _{4}))=(\beta _{3}',\beta _{4}'),$$
where $(\beta _3',\beta _4')$ is the adjacent pair in $R_{m,0}'$ which is matched with $(\beta _3,\beta _4)$. Moreover the following holds. Either $\psi _{1,3,\beta}=\psi _{\beta _1,\beta _3}$ preserves $D(\beta _1,\beta _2)$ or there is a single basic exchange $a\leftrightarrow b$ and prefix $v$ such that $C(va\leftrightarrow vb ,D(X)$ intersects both $D(\beta _1,\beta _2)$ and its complement, where $X=L_3$ or $BC$, depending on whether $a\leftrightarrow b$ is a basic exchange for $L_3$ or $BC$.
If $(\beta _i:1\le i\le 4)$ is of type C or AC then $a\leftrightarrow b$ is the basic exchange occurring in the definition of types C and AC.
  
 \end{lemma}

\begin{proof} We already know that $(\beta _1,\beta _2)=(\zeta _1(x),\zeta _2(x))$, where $(\zeta _i(x))$ is the $\zeta $-sequence for $x$. Write $\zeta _i(x)=\zeta _i$ and let $(\zeta _{2i-1}',\zeta _{2i}')$ be the adjacent pair in $R)_{m,0}$ which is matched with $(\zeta _{2i-1},\zeta _{2i})$.  Write $E=E_{\zeta _1,\zeta _3}$. Clearly,  if $\zeta _i(x)$ is defined for $i=3$, then $(\rho (\omega _3),\rho (\omega _4))=(\zeta _3',\zeta _4')$ and $(\zeta _3,\zeta _4)=(\beta _3,\beta _4)$ has the required properties, unless $I(\zeta _3,\zeta _4)$ intersects the support of $\psi _{\zeta _1,\zeta _3}$. The only way that this can happen is if $C(u\leftrightarrow u',E)$ intersects $I(\zeta _3,\zeta _4)$ for some exchange $u\leftrightarrow u'$ for $\zeta _1$ which is prefixed by some basic exchange $va\leftrightarrow vb$. The only way in which $C(va\leftrightarrow vb,D(L_3L_2))$ (or $C(va\leftrightarrow vb,D(BC))$) cannot be contained in $D(\zeta _1,\zeta _2)=D(\beta _1,\beta _2)$ is if $a\leftrightarrow b$ is one of the exchanges of (\ref{3.11.4}) (or (\ref{3.11.5})). 
If one of these happens then we have
$$D(u'X)\cap D(\beta _1,\beta _2)\ne \emptyset $$
where $X=L_3L_2$ or $X=BC$, as appropriate. 
  Let $\xi _1$ be the exchange supported by $C(u\leftrightarrow u',E_{\beta _1,\zeta _3})$ if $\zeta _i(x)$ is defined for $i\ge 3$ and supported by $C(u\leftrightarrow u',D(L_3L_2))$ if $\zeta _3(x)$ is not defined fo $i\ge 3$. In both cases there is an adjacent pair $(\beta _3,\beta _4)$ in $R_{m,0}$ such that  $I(\beta _3,\beta _4)\subset D(uJL_3L_2)$ and $\xi _1(I(\beta _3,\beta _4))$ is under $I(\beta _1,\beta _2)$ and not separated from it by anything else in the image of $\xi _1$. If $\zeta _i(x)$ is defined for $i\ge 3$ then $\xi _1$ interchanges $I(\zeta _3(x),\zeta _4(x))=I(\zeta _3,\zeta _3)$ and $I(\beta _3,\beta _4)$, and also  $E_{\beta _1,\zeta _3}\subset E_{\beta _1,\beta _3}$. In any case if  $\xi _1'$ is the exchange supported by $C(u\leftrightarrow u',E_{\beta _1,\beta _3})$, then 
  $$\xi _1'(I(\beta _3,\beta _4))=\xi _1(I(\beta _1,\beta _3).$$
  Also, since $I(\beta_3,\beta _4)$ is in the boundary of $E_{\beta _1,\beta _3}$ it is clear that $I(\beta _3,\beta _4)$ is not in any component of $s^{-n}(E_{\beta _1,\beta _3})$ for any $n>0$, apart from $uE_{\beta _1,\beta _3}$. It follows that  
  $$\rho (\omega _i)=\beta _i'{\rm{\ for\ }}i=3,\ 4,$$
  where $(\beta _3',\beta _4')$ is the adjacent pair of $R_{m,0}'$ which is matched with $(\beta _3,\beta _4)$. So in all cases we have a quadruple $(\beta _i:1\le i\le 4)$ of type A, C or AC with the required properties.   It is also true that
 \begin{equation}\label{4.15.1}\begin{array}{l}D(\beta _3,\beta _4)\subset uD(X)),\\ u'D(\beta _3,\beta _4)\subset D(\beta _1,\beta _2)\end{array}\end{equation}
 where $X=L_3L_2$ or $BC$, as appropriate. 
 
We claim that there is at most one basic exchange $va\leftrightarrow vb$ such that  $C=C(va\leftrightarrow vb,D(X))$ intersects both  $D(\zeta _1,\zeta _2)=D(\beta _1,\beta _2)$ and its complement. We see this as follows. The set $C$ never intersects $I(\zeta _1,\zeta _2)$ because there are no values of $w_1'$ between $w_1'(\zeta _1)$ and $w_1'(\zeta _2)$, and $w_1'(\zeta )$ is never between  two values of $w_i'(.,y)$ for any value $y$ of $w_1'$ and $i\geq 1$. We only need to consider exchanges $va\leftrightarrow vb$ as in (\ref{3.11.4}) for which $w_1'(vaL_3L_2)\neq w_1'(vbL_3L_2)$, or as in (\ref{3.11.5}) for which $w_1(vaL_3L_2)\neq w_1(vbL_3L_2)$. For all of these exchanges, since there are no intersections with $I(\zeta _1,\zeta _2)$, we only need to consider those with $w_1(vaL_3L_2)=vaL_3$ and $w_1(vbL_3L_2)=vbL_3L_2$. This condition ensures that for  any two such, say $v_1a_1\leftrightarrow v_1b_1$ and $v_2a_2\leftrightarrow v_2b_2$, neither of $v_1a_1L_3L_2$ or $v_1b_1L_3L_2$ can be a proper prefix of $v_2a_2L_3L_2$ or $v_2b_2L_3L_2$. It follows that if these are distinct basic pairs then
$$(D(v_1a_1L_3L_2)\cup D(v_1b_1L_3L_2))\cap (D(v_2a_2L_3L_2)\cup D(v_2b_2L_3L_2))=\emptyset .$$
Now we  show that $D(v_1b_1L_3L_2)$ and $D(v_2b_2L_3L_2)$ are not in $D(\zeta _1,\zeta _2)$ for the same pair $I(\zeta _1,\zeta _2)$. But if  $D(v_1b_1L_3L_2)\subset D(\zeta _1,\zeta _2)$ then $w_1'(\zeta _1)$ must be the maximal value of $w_1'$ to the right of  $D(v_1b_1L_3L_2)$ which means that $w_1'(\zeta _1)$ has a common prefix  with $v_1b_1$ of length $\vert v_1b_1\vert -1$. It is not possible for this to be true also for $v_2b_2$ if $v_2b_2\neq v_1b_1$. So the claim is proved.

If $va\leftrightarrow vb$ exists as in (\ref{3.11.4}) then, using the same argument, the only basic exchange $v_2a_2\leftrightarrow v_2b_2)$ for which  the homeomorphism $\xi  $ supported by $C(v_2a_2\leftrightarrow v_2b_2)$ does not preserve $vaD(\zeta _1,\zeta _2)$ is $vava\leftrightarrow vavb$. By induction it follows that $\psi _{\beta _1,\beta _3}$ preserves $\cup _{n\geq 0}(va)^nD(v)$

\end{proof}

\section{A general property of $\beta $ quadruples}

In order to prove Theorem \ref{2.8} in its entirety, we need some  inductive information about the form that a quadruple $(\beta _{4i+t}(\beta ):-3\le t\le 0)$ can take, assuming the existence of $\beta _j(\beta )$ as in \ref{2.8} for $j\le 4i$.

\begin{lemma}\label{4.20} Corollary \ref{4.11} and Lemma \ref{4.13} remain true  if we generalise to include quadruples  $(\zeta ^t:1\leq t\leq4)$ of the following form. There are two quadruples $(\zeta ^i:-1\le i\le 2)$ and $(\zeta ^{-1},\zeta ^0,\zeta ^3,\zeta ^4)$ which are either both $\zeta $quadruples or both type C, and $E=E_{\zeta ^1,\zeta ^3}$ is defined to be $E_{\zeta ^{-1},\zeta ^1}\setminus E_{\zeta ^{1},\zeta ^3}$. Also any component of $(\sigma _{\alpha _\zeta }\circ s)^{-n}(E)$ coincides with a component of $(\sigma _{\zeta _k}\circ s)^{-n}(E)$, where $\zeta _k$ is the longest common segment of $\zeta ^1$ and $\zeta ^3$. \end{lemma}
\begin{proof} First we consider \ref{4.11}. Without loss of generality we have $E_{\zeta ^1,\zeta ^3}\subset E_{\zeta ^{-1},\zeta ^3}$, and hence any component of $s^{-n}(E_{\zeta ^1,\zeta ^3}$ is contained in the corresponding component of $s^{-n}(E_{\zeta ^{-1},\zeta ^3})$. By \ref{4.10}, any component of $s^{-n}(E_{\zeta ^{-1},\zeta ^3}$ is either contained in or disjoint from each of $E_{\zeta ^{-1},\zeta ^3}$ and $E_{\zeta ^{-1},\zeta ^1}$. In the latter case  we use the fact that  since it is not possible for a component of $s^{-n}(E_{\zeta ^{-1},\zeta ^3})$ to contain $E_{\zeta ^{-1},\zeta ^1}$, for $n>0$. So any component of $s^{-n}(E_{\zeta ^{-1},\zeta ^3}$ is either contained in or disjoint from $E_{\zeta ^1,\zeta ^3}=E_{\zeta ^{-1},\zeta ^3}\setminus E_{\zeta ^1,\zeta ^3}$. So the same is true for any component of $s^{-n}(E_{\zeta ^{1},\zeta ^3}$

The proof of the analogue of \ref{4.13} is essentially the same as the proof of \ref{4.13}. The proof that any component of $\sigma _{\alpha _\zeta }\circ s)^{-n}(E)$ coincides with a component of $\sigma _{\zeta _k}\circ s)^{-n}(E)$ is the same as before. For the remainder of the statement of \ref{4.13}, we only need to check that $\zeta _k$ does not cut through the interior of components of $s^{-n}(E)$ Since $E\subset E_{\zeta ^{-1},\zeta ^3}$, this is the same as before. 
 \end{proof}

\begin{subsectiontheorem}\label{4.19} Let $\beta \in {\cal{Z}}_\infty (3/7,+,+,0)$ with endpoint $p(x)$. Suppose that $\beta _j(\beta )$ exists as in \ref{2.8} for any $j\le 4i$. For any $i\geq 1$, let $\psi _{1,4i-1,\beta }=\xi _n\circ \cdots \circ \xi _1$, where each $\xi _\ell $ is a disc exchange with support $C(v_\ell \leftrightarrow v_\ell ',E_\ell )$, where $E_\ell =E_{4j-3,4j-1,\beta }$ for some $j<i$ depending on $\ell $. For any $1\leq \ell $, let $\psi =\xi _\ell \circ \cdots \xi _1$.    Each component $J_1$ of $\psi (\overline{\tau _1}*\tau _2)\cap \{ z:|z|\ge 1\} $ over $x$, where $(\tau _1,\tau _2)$ runs over all adjacent pairs $(\tau _1,\tau  _2)$in  $R_{m,0}$, is paired with some other arc $J_2$ which is similarly defined with respect to some adjacent pair $(\tau _3,\tau _4)$, with the following properties.
\begin{itemize}
\item[1.] Paired arcs $J_1$ and $J_2$ both cross $\beta $, with one inside the other.
\item[2.] Choosing the numbering so that arcs of $\psi (\tau _1)$ and $\psi (\tau _3)$ adjoin $J_1$ and $J_2$ at one end, there is no path in $R_{m,0}$ which cuts between $\psi (\tau _1)$ and $\psi (\tau _3)$, including cutting between $J_1$ and $J_2$.
\item[3.] Every pair of arcs $J_1$ and $J_2$ is of one of the following forms:
\begin{enumerate}[a)]
\item     subarcs of $\psi (\overline{\tau _1}*\tau _2)$ and $\psi (\overline{\tau _3}*\tau _4)$, in either order, where $(\tau _i:1\leq i\leq 4)$ is a quadruple or lower quadruple or type C or AC;
\item  subarcs of $\psi (\overline{\tau _3}*\tau _4)$ and $\psi (\overline{\tau _1}*\tau _2)$ where, for some adjacent pair $(\tau _{-1},\tau _0)$, both $(\tau _i,\:-1\le i\le 2)$ and $(\tau _{-1},\tau _0,\tau _3,\tau _4)$ are upper quadruples, or both type C or AC;
\end{enumerate}
\end{itemize}
\end{subsectiontheorem}

{\textbf{Remark}} 
We shall refer to $J_1$ and $J_2$ as being {\em{paired over $x$}}. Of course a  local $S$ inverse of $s^n$ often maps an upper or lower $\zeta $ quadruple $(\zeta _{4j+t}(x_1):-3\leq t\leq 0)$ to another upper or lower $\zeta $  quadruple $(\zeta _{4\ell +t}(x_2):-3\leq t\leq 0)$. Since $S$ might map the upper unit circle into the lower unit circle,  upper quadruples can be mapped to lower quadruples, and vice versa. 
\begin{proof} 

The proof is by induction on $i$ and $\ell $. We assume inductively that $i=1$ or that $i>1$ and that the statement is already proved for $i-1$ replacing $i$. By \ref{4.15} to \ref{4.18}, we already know that $(\beta _{t}(\beta ):1\leq t\leq 4)$ is a quadruple of type A or C or AC. 

So now we need to prove the lemma for $i$, and for each $\ell \geq 1$.  We use induction on $\ell $. So let $\psi '$ be the identity if $\ell =1$ and $\psi '=\xi _{\ell -1}\circ \cdots \circ \xi _1$ if $\ell >1$. Writing $\xi =\xi _\ell $, we assume inductively that the lemma is true for $\psi '$ replacing $\psi $ --- which is trivial if $\psi '$ is the identity --- and prove it for $\psi =\xi \circ \psi '$. Let the support of $\xi $ be $C=C(v\leftrightarrow v',E)$, where $vE=E_{\delta _1,\delta _3}$ and $v'E=E_{\delta _1',\delta _3'}$, where $(\delta _i:1\le i\le 4)$ and $(\delta _i':1\le i\le 4$ are both quadruples of the inductive type. This means that either $(\delta _i:1\le i\le 4)$ is a $\zeta $ quadruple, or of type C or AC, or there are two quadruples $(\delta _i:-1\le i\le 2)$ and $\{ \delta _{-1},\delta _0,\delta _3,\delta _4)$ which are both upper $\zeta $ quadruples or both type C, and $E_{\delta _1,\delta _3}=E_{\delta _{-1},\delta _1}\setminus E_{\delta _{-1},\delta _3}$, and similar properties hold for $(\delta _i':1\le i\le 4)$.  We know from \ref{4.13} that the connecting arc of $C$ between $vE$ and $v'E$, apart from the central arc, consists of common arcs from each of the quadruples $(\delta _i:1\le i\le 4)$ and $(\delta _i':1\le i\le 4)$.   Let $J_1'$ and $J_2'$ be paired sub-arcs of $\psi '(\overline{\tau _1'}*\tau _2')$ and $\psi '(\overline{\tau _3'}*\tau _4')$ where $(\tau _i':1\le i\le 4)$ satisfies the inductive hypothesis. 

If $C$ does not intersect $J_1'\cup J_2'$, then there is nothing to prove. If neither $vE$ nor $v'E$ intersects $J_1'\cup J_2'$, but $C$ does intersect   $J_1'\cup J_2'$, then the intersection has to be of the arc of $C$ between $vE$ and $v'E$, with $J_1'\cup J_2'$.  But, using the inductive hypothesis, by property 2 of the lemma for $J_1'\cup _2'$, no part of the connecting arc of $C$ can separate $J_1'$ and $J_2'$. 

If $vE$ or $v'E$ contains $J_1'\cup J_2'$, then $\xi (J_1'\cup J_2')$ is a single pair of arcs $J_1\cup J_2$ with the required properties. In fact $(J_1,J_2)$ is a paired arc of exactly the same type as $(J_1',J_2')$ because $\xi $ maps quadruples in $vE$ to quadruples in $v'E$. 

So it remains to consider what happens when one of $vE$ or $v'E$ intersects $J_1'\cup J_2'$ but does not contain it.   Suppose without loss of generality that $vE$ intersects, but does not contain, $J_1'\cup J_2'$. Then, by the inductive hypothesis,  there is a $\zeta $ quadruple $(\eta _i':1\le i\le 4)$ such that $J_1'$ and $J_2'$ both intersect $E_{\eta _1',\eta _3'}$ and $vE=E_{\delta _1,\delta _3}$ also intersects $E_{\eta _1',\eta _3'}$. Also by the inductive hypothesis on $J_1'$ and $J_2'$, we have  $J_2'=I(\eta _3',\eta _4')$, and either $J_1'=I(\eta _1',\eta _2')$, or $J_1'$ intersects the region  between  $J_2'$ and the nearest segment of one of $\eta _1'$ and $\eta _2'$ to $I(\eta _1',\eta _2')$. By the inductive hypothesis on $vE$ we have $(\delta _1,\delta _2)=(\eta _1',\eta _2')$ or $(\delta _{-1},\delta _0)$ depending on whether or not $(\delta _i:1\le i\le 4)$ is itself a $\zeta $ quadruple or type C quadruple.  So then there are various configurations to consider.  The figure shows the configuration when $(\delta _i:1\le i\le 4)$ is a $\zeta $ quadruple and $\psi '$ is the identity on $vE=E_{\delta _1,\delta _3}$ --- which means that $\psi '$ is also the identity on $vE$ --- and  the new pairing is of $(\delta _3,\delta _4)$ with $(\eta _3',\eta _4')=(\tau _3',\tau _4')$ -- and also of $(\delta _3,\delta _4)$ with $(\eta _1',\eta _2')=(\tau _1',\tau _2')$, but that is not shown in the picture. Now suppose that $\psi '$ is not necessarily the identity on $vE$. If there is an intersection between $vE$ and $J_1'\cup J_2'$ then necessarily at least one of  $I(\delta _3,\delta _4)$ and $I(\delta _1,\delta _2)$ is under $J_1'$, say $I(\delta _3,\delta _4)$. If $I(\delta _3,\delta _4)\subset \psi '(E_{\tau _1',\tau _3'})$, then $I(\delta _3',\delta _4')=\psi '(I(\tau _5'\tau _6')$ for some $I(\tau '_5,\tau _6')\subset E'$, where $E'$ is in the support of $\psi '$ and $E_{\tau _1',\tau _3'}\subset E'$, where $(\tau _5',\tau _6',\tau _3',\tau _4')$ is another quadruple of the same type, and the pairing for $\xi \circ \psi '=\psi $ is of $\xi (J_2')$ with an arc of $\xi \circ \psi '(I(\tau _5',\tau _6'))$. If $I(\delta _3,\delta _4)$ is not in $\psi '(E')$ then it is separated from $\psi (E')$ by $\psi '(I(\tau _5',\tau _6'))$ for some such $(\tau _5',\tau _6')$, and the new pairing is still of $(\tau _5',\tau _6')$ with each of $(\tau _3',\tau _4')$ and $(\tau _1',\tau _2')$. Then $(\tau '_i:3\le i\le 6)$ and $(\tau _1',\tau _2',\tau _5',\tau _6')$ are quadruples of the required types, and this completes the proof of the inductive hypothesis. 

 \begin{figure}
\centering{\includegraphics[width=6cm]{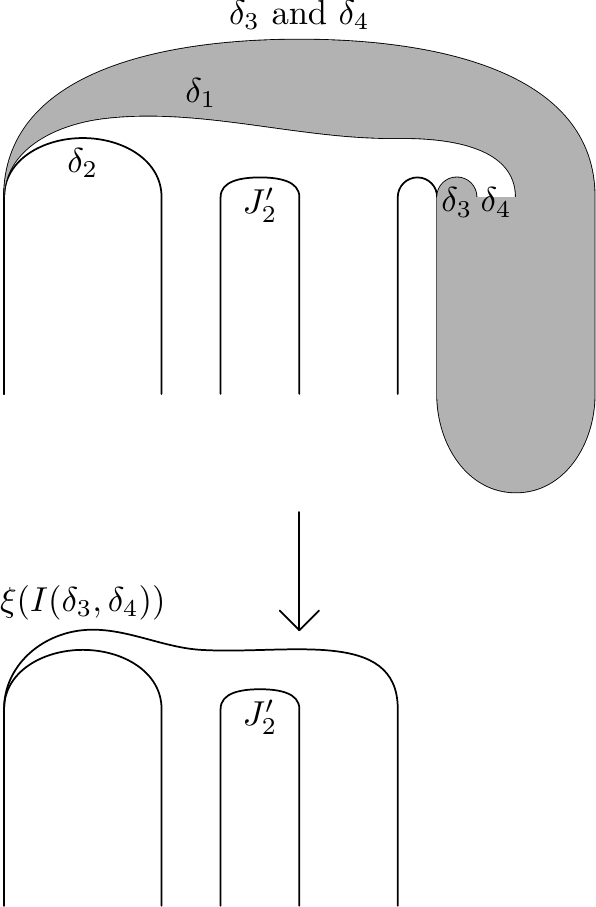}}
\caption{Pairing of type 3b)}
\end{figure}

\end{proof} 

\subsection{Satellite arcs}

 We have seen in \ref{4.19} that there are a number of ways in which paired principal arcs for $[\psi ]$ and $R_{m,0}$ can differ from  a $\zeta $ quadruple, up to isotopy preserving $S^1\cup Z_\infty $.  We have seen that if $\psi =\psi _{1,4i-1,\beta }$ for any capture path $\beta \in {\cal{Z}}(3/7,+,+,0)$, the paired principal arcs for $\psi (R_{m,0})$ have no transversal intersection with principal arcs for $R_{m,0}$, and and pairs of principal arcs for $\psi (R_{m,0})$ are not separated by any other principal arcs. If both endpoints coincide, then we say that this is a pair of {\em{satellite arcs (for $\psi $)}}, or a {\em{satellite pair}}.Pairs of arcs which are not a satellite pair are called  a {\em{non-satellite pair}}. If $I(\psi _{1,4i-5,\beta }(\beta _{4i-3}),\psi _{1,4i-5,\beta }(\beta _{4i-2}(\beta ))$ and  $I(\psi _{1,4i-1,\beta }(\beta _{4i-1}),\psi _{1,4i-1,\beta }(\beta _{4i}(\beta ))$ are a satellite (or non-satellite) pair over $p(w(\beta ))$, then we call $(\beta _{4i+t}(\beta ):-3\leq t\leq 0)$ a {\em{satellite quadruple}} (or {\em{non-satellite quadruple}}.

\begin{corollary}\label{4.21} Assume the existence of $\beta _j(\beta )$ with the properties of Theorem \ref{2.8} for $j\le 4i$. Then  $\psi _{4i-3,4i-1,\beta }$ is the identity on $I(\beta _{4i-1}(\beta ),\beta _{4i}(\beta ))$ for all $i\geq 2$, and $(\beta _{4i-3}(\beta ),\beta _{4i-2}(\beta ))$ and $(\beta _{4i-1}(\beta ),\beta _{4i}(\beta ))$ are the adjacent pairs in $R_{m,0}$ whose principal arcs for $[\psi _{1,4i-5,\beta }]$ are the next two inside the principal arc of $(\beta _{4i-5}(\beta ),\beta _{4i-4}(\beta ))$ for $[\psi _{1,4i-5,\beta }]$ which is crossed by $\beta $.\end{corollary}

\begin{proof} By \ref{4.11}, \ref{4.20} and \ref{4.19},  then any component of \\ $s^{-n}(E_{4i-3,4i-1,\beta })$, for any $n>0$, is either strictly contained in $E_{4i-3,4i-1,\beta }$ or is disjoint from it. Hence it cannot intersect $I(\beta _{4i-1}(\beta ),\beta _{4i}(\beta ))$, which is in the boundary. So if $C(v\leftrightarrow v',E)$ is the support of any disc exchange in the composition for $\psi _{4i-3,4i-1,\beta }$ then $vE$ and $v'E$ are disjoint from $I(\beta _{4i-1}(\beta ),\beta _{4i}(\beta ))$. There are also no transversal intersections between  $I(\beta _{4i-1},\beta _{4i})$ and  the two segments of paths in$R_{m,0}$ in $C$ and no intersections between $I(\beta _{4i-1},\beta _{4i})$ and the central arc $\ell (C)$ of $C$ .\end{proof}

\section{Satellite arcs}\label{4.22}

It will be clear from \ref{4.19} that satellite arcs are formed by the action of disc exchanges. From the form given for the support of a disc exchange $\xi $ in \ref{4.13}, and \ref{4.20}, we see that satellite arcs for $\psi _{1,4i-1,\beta }$  are always parallel to components of 
$$(\omega _{1,n(C),m(C)}\cup \omega _{2,n(C),m(C)}\cup \ell (C))\setminus S^1$$
 for varying $C$, where $C$ is the support of some disc exchange in the composition for $\psi _{1,4i-1,\beta }$. To see this , we need to use the fact that $w_1'(\beta _i)$ is constant for $i\geq 3$. So then there are never transverse intersections between $\ell (C)$ and $\ell (C')$, or between $\omega _{i,n(C),m(C)}$ and $\omega _{j,n(C'),m(C')}$, for different $C$ and $C'$, and for $i$, $j\in \{ 1,2\} $. 
 So satellite arcs are always parallel either to arcs $\ell (C)$ or to segments in the exterior of the unit disc of paths in $R_{m,0}$, not including the first segments of such paths. 

\section{Final proof of Theorem \ref{2.8}}\label{4.23}

We are now finally ready to prove Theorem \ref{2.8} in full generality.We do this for the different cases of the quadruple $(\beta _i(\beta ):1\le i\le 4)$. Let $\beta \in {|cal{Z}}_m(3/7,+,+,0)$ In what follows, we use the notation $E_{4j-3,4j-1,\beta }$ for the set $E_{\beta _{4j-3},\beta _{4j-1},\beta }$ if $2j\le i$ and the existence of the sequence $\beta _\ell (\beta )$ for $\ell \le 2i$ as in Theorem \ref{2.8} has been established. We write $E_{4j-3,4j-1}$ for $E_{4j-3,4j-1,\beta }$ if no confusion can arise.

\begin{subsectiontheorem}\label{4.16} Let $(\beta _i:1\leq i\leq 4)$ be of type A. Then $(\beta _1,\beta _2),(\beta _3,\beta _4))$ extends to a sequence  $(\beta _{2i-1},\beta _{2i})$ of adjacent pairs in $R_{m,0}$ such that, if  $\gamma _i$ is the sequence of elements of $\pi _1(V_{3,m},a_1)$ and $(\omega _{2i-1},\omega _{2i})$ the sequence of adjacent pairs in $\Omega _m$ such that $\rho(\beta )$ successively crosses $\overline{\rho (\gamma _i*\omega _{2i-1})}*\rho (\gamma _i*\omega _{2i})$, then
$$(\rho (\omega _{2i-1}),\rho (\omega _{2i}))=\begin{cases}(\beta _{2i-1},\beta _{2i}){\rm{\ if\ }}i{\rm{\ is\ odd,}}\\(\beta _{2i-1}',\beta _{2i}'){\rm{\ if\ }}i{\rm{\ is\ even,}}\end{cases}$$
where $(\beta _{2i-1}',\beta _{2i}')$ is the adjacent pair in $R_{m,0}'$ matched with $(\beta _{2i-1},\beta _{2i})$.  Then $I(\beta _{2j-1},\beta _{2j})$ is contained in $D(\beta _1,\beta _2)$ for all $j\geq 2$.\end{subsectiontheorem}

\begin{proof} Naturally, this is proved by induction on $j$. It is true for $j=2$ by the assumption that $(\beta _j:1\le j\le 4)$ is of type A.  Using \ref{4.13}, this means that $I(\beta _3,\beta _4)$ is disjoint from $(\sigma {\beta _{1}}\circ s)^{-n}(E_{1,3})$ for all $n>0$, and $I(\beta _3,\beta _4)\subset D(\beta _1,\beta _2)$. Suppose that we have proved the Theorem for $\beta _j$ for $j\le 4i$. Then we consider $\psi _{1,4i-1,\beta }$. We claim that $\psi _{1,4i-1,\beta }$ does not map any  principal arc under $I(\psi _{1,4i-1,\beta }(\beta _{4i-1}),\psi _{1,4i-1,\beta }(\beta _{4i})$ out of $D(\beta _1,\beta _2)$. To see this, as usual, we write 
$$\psi _{1,4i-1,\beta }=\psi _{1,3,\beta }\circ \cdots \circ \psi _{4i-3,4i-1,\beta }$$
 for $1\le j\le i$, and we write  $\psi _{4j-3,4j-1,\beta }$ as a composition of disc exchanges with supports $C(u_{j,\ell}\leftrightarrow u_{j,\ell }',E_{4j-3,4j-1,\beta }$ where $E_{4j-3,4j-1,\beta }=E_{\beta _{4j-3}(\beta ),\beta _{4j-1}(\beta )}$ and $u_{j,\ell }\leftrightarrow u_{j,\ell }'$ is an exchange for $\beta _{4j-3}$. This means $u_{j,\ell }\leftrightarrow u_{j,\ell '}$ is prefixed by a basic exchange $v_{j,\ell }a_{j,\ell }\leftrightarrow v_{j,\ell }b_{j,\ell }$, where the first letters of $a_{j,\ell }$ and $b_{j,\ell }$ are different, that is, each $a_{j,\ell}\leftrightarrow b_{j,\ell }$ is one of the basic exchanges listed in \ref{3.11}. From what we now know about the quadruple $(\beta _{4j+t}(\beta ):-3\le t\le 0)$ from \ref{4.19}, and from the inductive hypothesis, we have 
 $$E_{4j_1-3,4j_1-1}\cap (L_3L_2R_3)^nE_{4j_2-3,4j_2-1,\beta }=\emptyset  $$
 for all $n>0$ and all $j_1$ and $j_2$. It follows that we can only have a nontrivial intersection between $u_{j_1,\ell _1}E_{4j_1-3,4j_1-1,\beta }$ and $u_{j_2,\ell _2}E_{4j_2-3,4j_2-1,\beta }$ if one of the following holds.
 \begin{itemize}
 \item $v_{j_1,\ell _1}=v_{j_2,\ell _2}$ and $\{ a_{j_1,\ell _1},b_{j_1,\ell _1}\} =\{ a_{j_2,\ell _2},b_{j_2,\ell _2}\} $;
 \item One of $v_{j_1,\ell _1}a_{j_1,\ell _1}$ or $v_{j_1,\ell _1}b_{j_1,\ell _1}$ is a prefix of $v_{j_2,\ell _2}$ or vice versa.
\end{itemize}
We also know from \ref{4.15} that there is at most one basic exchange $a\leftrightarrow b$ and prefix $v$ such that $C(a\leftrightarrow vb, D(X))$ which intersects both $D(\beta _1,\beta _2)$ and its complement. If $a\leftrightarrow b$  and $v$ do not exist, then $\psi _{1,4i-1,\beta}$ preserves $D(\beta _1,\beta _2)$. So suppose  that  $a\leftrightarrow b$  and $v$ do exist. Then there is an exchange $u\leftrightarrow u'$ for $\beta _1$ which is prefixed by $va\leftrightarrow vb$ and $u'E_{1,3}\cap D(\beta _1,\beta _2)\ne \emptyset $ (choosing $u'$ without loss of generality). Then by hypothesis, since we are considering type A quadruples  $I(\beta _3,\beta _4)\cap u'E_{1,3}=\emptyset $, and then we also have 
$$D(\beta _3,\beta _4)\cap u'(E_{1,3}\cup D(\beta _1,\beta _2)=\emptyset .$$
Hence if $(\eta _1,\eta _2)$ is an adjacent pair such that $I(\eta _1,\eta _2)\subset D(\beta _3,\beta _4)$ and $\psi _{1,4j-1}^{-1}(I(\eta _1,\eta _2))$ is not contained in $D(\beta _3,\beta _4)$,  then the  first disc exchange in the composition for $\psi _{1,4i-1,\beta }^{-1}$  which moves $I(\eta _1,\eta _2)$ out of $D(\beta _3,\beta _4)$  must have support $C(u_{j,\ell }\leftrightarrow u_{j,\ell }'E_{4j-3,4j-1})$ where $v_{j,\ell }a_{j,\ell }$ and $v_{,\ell }b_{j,\ell }$ do not have $va$ or $vb$ as a prefix, and similarly with the roles reversed. So then further disc exchanges cannot move out of $D(\beta _1,\beta _2)$.  It follows that for any $j\le i$ we have 
$$\psi _{1,4i-1}^{-1}(I(\eta _1,\eta _2)\subset D(\beta _3,\beta _4)\cup D(v_{j,\ell }a_{j,\ell })\cup D(v_{j,\ell }b_{j,\ell })$$
for some $(j,\ell )$ and $\{ a_{j,\ell },b_{j,\ell } \} \ne \{ a,b\} $. 
If
$$I(\eta _1,\eta _2)\subset D(v_{j,\ell }a_{j,\ell })\cup D(v_{j,\ell }b_{j,\ell })$$
for $\{ a_{j,\ell }, b_{j,\ell } \} \ne \{ a,b\} $ then it cannot be moved into the domain of the disc exchange with support  $C(va\leftrightarrow vb D(X))$  by $\psi _{1,4i-1,\beta }^{-1}$ It follows that we have a quadruple $(\beta _{4i+t}:1\le t\le 4)$ such that
$$(\rho (\omega _{4i+1}),\rho (\omega _{4i+2}))=(\beta _{4i+1},\beta _{4i+2}),$$
$$(\rho (\omega _{4i+3}),\rho (\omega _{4i+4}))=(\beta _{4i+3}',\beta _{4i+4}')$$
and $(\beta _{4i+t}:1\le t\le 4)$ has all the required properties. 

\end{proof}

Now we deal with the cases of $(\beta _i:1\le i\le 4)$ being a quadruple of type C or AC

\begin{subsectiontheorem}\label{4.17} Let $\beta \in {\cal{Z}}_m(3/7,+,+,0)$ and let $p(\beta )$ denote the endpoint of $\beta $. Let $(\beta _i(\beta ):1\le i\le 4)=(\beta _i:1\leq i\leq 4)$ be of  AC. Let $va\leftrightarrow vb$ be the basic exchange as in \ref{4.15}, prefixing the exchange $u\leftrightarrow u'$ for $\beta 1$, also as in \ref{4.15}.  Then $\beta j$ exists as in Theorem \ref{2.8}
$$D(\beta _{2j-1},\beta _{2j})\subset (D(\beta _1,\beta _2)\cap E_{1,3})\cup _{1\le i }u^i(E_{1,3}$$
for all $j$. \end{subsectiontheorem}
\noindent {\textbf{Remark}} We automatically have $n\ge 1$, as $(\beta _i:1\le i\le 4)$ is of type AC.
\begin{proof} The proof follows the same lines as the proof for type A quadruples. By the definition of type AC we have 
$$D(\beta _3,\beta _4)\subset u(E_{1,3}\cap D(L_3L_2)).$$
We also have 
$$u'E_{1,3}\subset D(\beta _1,\beta _2)\cap E_{1,3}$$
but 
$$uE_{1,3}\cap E_{1,3}=\emptyset $$
This might seem strange because the first common unit disc crossing of $\beta _3$ and $\beta _4$ is on one side of $uE_{1,3}$, and the second common crossing is on the other side of $uE_{1,3}$. But the paths of $R_{m,0}$ with first two unit disc crossings in common with $\beta _3$ and $\beta _4$, and the third crossing strictly between this and the second and the first, are not in $E_{1,3}$. The paths in $uE_{1,3}$ are of this type. Related to this, provided $v$ is nontrivial, or $u\ne va$ we have
$$L_3L_2R_3E_{1,3}\cap E_{1,3}=\emptyset $$
and , further,
$$L_3L_2R_3(\cup _{i\ge 0}u^iE_{1,3})\cap (\cup _{i\ge 0}u^iE_{1,3})=\emptyset $$
If $v$ is trivial and $u=a$ then we still have
$$L_3L_2R_3(E_{1,3}\cap D(BC))\cap (E_{1,3}\cap D(BC))=\emptyset $$
and 
$$L_3L_2R_3(\cup _{i\ge 1}u^iE_{1,3})\cap (\cup _{i\ge 1}u^iE_{1,3})=\emptyset .$$
In all cases, we have
$$D(\beta _1,\beta _2)\cap L_3L_2R_3(E_{1,3}\cup D(\beta _1,\beta _2))=\emptyset .$$
In view of this we define
$$E=\cup _{i\ge 1}u^iE_{1,3}\cup D(\beta _1,\beta _2).$$
Then in all cases we have
$$E\cap L_3L_2R_3E=\emptyset $$
The only exchange $u_2\leftrightarrow  u_2'$ for $\beta _3$ such that $C(u_2\leftrightarrow u_2',E)$  intersects $uE_{1,3}\cup D(\beta _1,\beta _2)$ without being contained in it is $u\leftrightarrow u'$. Similarly the only exchange $u_2\leftrightarrow u_2'$ such that  $C(u_2\leftrightarrow u_2',E)$  intersects $\cup _{1\le i\le n}u^iE_{1,3}\cup D(\beta _1,\beta _2)$ without being contained in it is $u^{n+1}\leftrightarrow u^nu'$. It follows by induction on $i$, using \ref{4.15}, \ref{4.13}, \ref{4.20} and \ref{4.19}  as in \ref{4.16}, that $\beta _\ell $ exists as in the statement of \ref{2.8} for all $\ell \le 2i$, for all $i$, with $I(\beta _{2j-1},\beta _{2j})\subset E$ for all $j\le i$. These properties imply that $w_1'(\beta _j)-=w_1'(\beta _3)$ for all $j\ge 3$ and hence the exchanges for $\beta _j$ are closely related to the exchanges for $\beta _3$. In particular it is still true that $C(u_2\leftrightarrow u_2',E)$ is contained in $E$ whenever it intersects it, where $u_2\leftrightarrow u_2'$ is an exchange for $\beta _j$, for any $j\ge 3$, since $\overline{\beta _j}*\beta _3\subset E$.
\end{proof}

\begin{subsectiontheorem}\label{4.18} Let $\beta \in {\cal{Z}}_m(3/7,+,+,0)$ and let $p(\beta )$ denote the endpoint of $\beta $. Let $(\beta  _i(x):1\leq i\leq 4)$ be of type C, with associated basic exchange $va\leftrightarrow vb$ prefixing  the exchange $u\leftrightarrow  u'$ for $\beta _3$. Then $\beta _i$ exists as in Theorem \ref{2.8} with 
\begin{equation}\label{4.18.1}\psi _{1,3}(I(\zeta _{2i-1},\zeta _{2i}(x))\subset  D(\beta _1,\beta _2)\end{equation} for all $i\geq 2$ and 
\begin{equation}\label{4.18.2}D(\beta _{2i-1},\beta _{2i})\subset (D(\beta _1,\beta _2)\cap E_{1,3})\cup _{n\geq 1}(u)^nE_{\beta _1,\beta _3}.\end{equation}
for all $i\geq 3$.\end{subsectiontheorem}

\begin{proof}  By hypothesis, $(\beta _i:1\le i\le 4)$ is of type C. We can write $E_{1,3}=E_1\cup E_2$ where 
$$E_1=\{ \zeta \in E_{1,3}:w_1'(\zeta _3)<w_1'(\zeta )<\leq w_1'(\zeta _1)\} $$
$$E_2=\{ \zeta \in E: w_1'(\zeta )=w_1'(\zeta _3)\} $$
Then the conditions of type C are such that 
$$uE_{1,3}\cap E_{1,3}=uE_1\cap E_2.$$
It is also the case that $u'D(L_3L_2)\cap E_{1,3}=\emptyset $ -- which was not the case in \ref{4.16}. As before we have
$$E_{1,3}\cap L_3L_2R_3E_{1,3}=\emptyset $$
and the only exchanges $u_2\leftrightarrow u_2'$ for which $C(u_2\leftrightarrow u_2')$ intersects $E_{1,3}$ without being contained in it are $u\leftrightarrow u'$ and $uu\leftrightarrow uu'$. in fact, we now have
$$uuE_{1,3}\cup uu'E_{1,3})\cap E_{1,3}=\emptyset ,$$
but $C(uu\leftrightarrow uu',E_{1,3}$, which is the support of a disc exchange in the composition for $\psi _{1,3}$, has transversal intersection with $I(\beta _3,\beta _4)$. It follows that $(\beta _5,\beta _6)$ exists and either $I(\beta _5,\beta _6)\subset D(\beta _3,\beta _4)\cap uE_{1,3}$ or $I(\beta _5,\beta _6)\subset u^2E_{1,3}$. In both cases $(\beta _7,\beta _8)$ exists with all the required properties and $I(\beta _7,\beta _8)\subset D(\beta _3,\beta _4)\cap uE_{1,3}$ or $I(\beta _7,\beta _8)\subset u^2E_{1,3}$ respectively. Then similarly to \ref{4.17}, we define
$$E=D(\beta _1,\beta _2)\cup \cup _{n\ge 1}u^iE_{1,3}$$
Once again  if $u_2\leftrightarrow u_2'$ is any exchange for $\beta _3$, then $C(u_2\leftrightarrow u_2')$ is contained in $E$ whenever it intersects it. So by induction $\beta _\ell $exists  as i Theorem \ref{2.8} for $\ell \le 2i $ for all $i\ge 1$, with $I(\beta _{2i-1},\beta _{2i})\subset E$,  and $C(u_2\leftrightarrow u_2')$ is contained in $E$ whenever it intersects it for any exchange $u_2\leftrightarrow  $  for $\beta _j$, for any $j\ge 3$.

\end{proof}

\chapter{Specific examples}\label{5}

In this chapter we expand on the examples in \cite{R6}. We show, somewhat surprisingly, that the Thurston equivalence classes found in \cite{R6} are considerably larger than was demonstrated there. We will use this class of examples in Theorem \ref{6.2}. 

\section{Some old and new notation}\label{5.1}
 For any word $x$ which ends in $C$, we write $p(x)$ for the point in $Z\cap D(x)$ of lowest possible preperiod under $s$. 

It is convenient to modify and extend the notation of \cite{R6}. As in \cite{R6} we define
 $$a=L_{3}(L_{2}R_{3})^{2},\ \ b=L_{3}^{5},$$
 $$c=L_{3}^{3}L_{2}C,$$
$$d=L_{3}^{2}a,$$
$$v_0=L_3L_2R_3a,\ \ w_0=L_3L_2R_3b,$$
$$u_{0}=L_{3}L_{2}R_{3}L_3L_2C,$$
$$t_{0}=L_{3}L_{2}R_{3}d,$$
and inductively, for $k\geq 0$, we define
$$v_{k+1}=v_kt_ka,\ \ t_{k+1}=v_kt_kd,\ \ u_{k+1}=v_kt_kc.$$
Now we generalise this construction.  Let $\alpha =(\alpha (k))\in \{ a,\b\} ^{\mathbb N}$. We define words $v_{k,\alpha }$, $t_{k,\alpha }$ and $u_{k,\alpha }$ for all $k\geq 0$. To start the induction we define 
$$v_{0,\alpha }=\begin{cases}v_0{\rm{\ if\ }}\alpha (0)=a,\\ w_0{\rm{\ if\ }}\alpha (0)=b,\end{cases}\ \ t_{0,\alpha }=t_0.$$
Then we define, for $k\geq 0$,
$$v_{k+1,\alpha }=v_{k,\alpha }t_{k,\alpha }\alpha (k),\ \ t_{k+1,\alpha }=v_{k,\alpha }t_{k,\alpha }d,\ \ u_{k+1,\alpha }=v_{k,\alpha }t_{k,\alpha }c.$$
Thus, the sequences $v_k$, $t_k$ and $u_k$ are $v_{k,\alpha}$, $t_{k,\alpha }$, and $u_{k,\alpha }$ in the special case where $\alpha =\alpha ^{\infty }$, where $\alpha ^{\infty }(k)=a$ for all $a$. We also define $\alpha ^{\ell }$ for $\ell \geq 0$ to be the sequence such that 
$$\alpha ^\ell (k)=\begin{cases}a{\rm{\ if\ }}k<\ell ,\\ b{\rm{\ if\ }}\ell \geq k.\end{cases}$$
Then the word which was called $w_{\ell ,r}$ in \cite{R6} for $0\leq \ell \leq r$ is $v_{r,\alpha ^\ell}$ in our current notation. It will sometimes be convenient to write
$$v_{k,\alpha ^\ell }=v_{k,\ell },\ \ t_{k,\alpha ^\ell }=t_{k,\ell },\ \ u_{k,\alpha ^\ell }=u_{k,\ell }$$
which is the opposite order of indices to that used in \cite{R6}, but since different letters are used, hopefully this will not cause confusion We define, for $r\geq 0$
$$x_r=v_ru_r$$
and, more generally, 
$$x_{r,\alpha }=v_{r,\alpha }u_{r,\alpha },\ \ x_{r,\alpha ^\ell }=x_{r,\ell }.$$
We define $\gamma ^{r,\alpha }$ to be the path in  ${\cal{Z}}(a_1,+,+,0)$ with endpoint $p(x_{r,\alpha })$ and 
$$\gamma ^{r,\ell }=\gamma ^{r,\alpha ^\ell },\ \ 0\leq \ell \leq r+1.$$
 Since $x_{r,\alpha }$ depends on $\alpha (\ell )$ only for $\ell \leq r$, there are only finitely many choices for $x_{r,\alpha }$: in fact, $2^{r+1}$ choices.  In \cite{R6}, it was proved that the captures $\sigma _{\gamma ^{r,\ell }}\circ s$, for $0\leq \ell\leq r+1$  are all Thurston equivalent.  (The notation here has been changed from that in \cite{R6}, as $\beta $with various indices is rather  heavily used in this paper).  In \cite{R6} the proof was direct. We shall now give a more indirect proof in keeping with the methods of this paper.  Note that the preperiod of the points $x_{r,\alpha }$ is the same for all $\alpha $. We fix $m$ so that $p(x_{r,\alpha })$ is of pre-period $\leq m$ under $s$. Write $n_{r,\ell }=n(\gamma ^{r,\ell })$, and $n(r,\alpha )=n(\gamma ^{r,\alpha })$.  We will show the following.
 
 \begin{theorem}\label{5.2} Fix $r\ge 0$. Let $\alpha _1$ and $\alpha _2\in \{ a,b\} ^{r+1}$ be such that, for some $p<r$, 
 $$\alpha _1(i)=\alpha _2(i){\rm{\ for\ }} i\le p$$
 $$\alpha _1(i)a,\alpha _2(i)=b{\rm{\ for\ all\ }}p<i.$$
 Then for $\gamma _i=\gamma ^{r,\alpha _i}$, the captures $\sigma _{\gamma _1}\circ s$ and $\sigma _{\gamma _2}\circ s$ are Thurston equivalent. Hence the captures $\sigma _{\gamma ^{r,\alpha }}\circ s$ are Thurston equivalent for all $\alpha \in \{ a,b\} ^{r+1}$. \end{theorem}

The final claim of Theorem \ref{5.2} follows easily from the first part, because we see inductively that $\sigma _{\gamma ^{r,\alpha }}$ is Thurston equivalent to $\sigma _{\gamma ^{r,\infty }}$ for all $\alpha \in \{ a,b\} ^{\mathbb N}$, by moving from $\alpha ^\infty $ to $\alpha $ in a series of steps, $\alpha _i$ for $0\le i\le n$. To do this, we define $\alpha _0=\alpha ^\infty $, and let $p_i$ be the increasing sequence of non-negative integers such that $\alpha (j)$ is constant for $p_i<j\le p_{i+1}$ and $\alpha (i)=a$ for $i<p_1$. Then we define $\alpha _i$ by induction on $i$ by  $$\alpha _{i+1}(j)=\alpha _i(j){\rm{\  if\ }}j\le p_i,$$
$$\alpha _{i+1}(j)=b\ne a=\alpha _i(j){\rm{\ if\ }}p_i<j{\rm{\ and\ if\ }}i{\rm{\ is\ even,}}$$
$$\alpha _{i+1}(j)=a\ne b=\alpha _i(j){\rm{\ if\ }}p_i<j{\rm{\ and\ if\ }}i{\rm{\ is\ odd.}}$$
 Then $\alpha _n=\alpha $.  We shall use Theorem \ref{5.2} in Chapter \ref{6} to show that the set of capture paths $\gamma $ in ${\cal{Z}}_m$ such that $\sigma _\gamma \circ s$ is equivalent to $\sigma _{\gamma _2}\circ s$ for precisely $2^{r+1}$ capture paths $\gamma _2$ in ${\cal{Z}}_m$, is of positive density, as $m\to \infty $, for each $r\ge 0$ (in fact, also for $r=-1$, but that case does not use this construction).  The idea of that proof will be to show that, for $\gamma =\gamma ^{r,\alpha }$, for each $\alpha \in \{ a,b\} ^{\mathbb N}$, the density of paths $\omega \in R_{m,0}$ in a  sufficiently small neighbourhood $U$ of $\beta $ such that $\sigma _{\omega }\circ s$ is equivalent to $\sigma _{\gamma _2}\circ s$ for exactly one path $\gamma _2$ in a given neighbourhood $U'$ of $\gamma $, approaches $1$ as $U$ is made arbitrarily small, and that any capture path $\gamma _3\in {\cal{Z}}$ such that $\sigma _\omega \circ s$ is equivalent to $\sigma _{\gamma _3}\in {\cal{Z}}$, has to be in a given neighbourhood $U'_\ell $ of $\gamma ^{r,\ell }$ or $U'_{\ell ,-}$ of $\gamma ^{r,\ell ,-}$, for $U$ sufficiently small. 

We now prove \ref{5.2}. This proof follows the lines of the proof of the corresponding result in \cite{R6} pretty exactly, with some minor differences of notation, to accommodate the notation of this paper. I had actually hoped to provide a proof which gave an explicit calculation of the path $\beta \in R_{m,0}$ such that $\sigma _\beta \circ s$ is equivalent to $\sigma _\gamma \circ s$ for $\gamma =\gamma ^{r,\alpha }$, for all $\alpha $. But the calculation is simply too complicated.

\begin{proof} Write $x=x_{r,\alpha _2}$ and $y=x_{r,\alpha _1}$. Exactly as in \cite{R6}, $D(y)$ is to the right of $D(x)$.  Now we  imitate the proof of Lemma 2.6 of \cite{R6}. We claim that the suffixes $u$ of $y$  such that $D(u)$ is between $D(x)$ and $D(y)$ are precisely those which start with $v_{\alpha _1,p+1}$ or $t_{\alpha _1,p+1}=v_{\alpha _2,p+1}$ or $u_{\alpha _1,p+1}=u_{\alpha _2,p+1}$. We see this as follows. The largest common prefix of $x$ and $y$ is $v_{\alpha _1,p}t_{\alpha _1,p}L_3=v_{\alpha _2,p}t_{\alpha _2,p}L_3$. So this must also be a prefix of $u$, if $D(u)$ is between $D(x)$ and $D(y)$. Since $u$ is also a suffix of $y$, it must start with either $v_{\alpha _1,p+1}$ or $t_{\alpha _1,p+1}$ or $u_{\alpha _1,p+1}$. It remains to show that this is a sufficient condition for $D(u)$ to be between $D(x)$ and $D(y)$. If the prefix of $u$ is $t_{\alpha _1,p+1}$ or $u_{\alpha _1,p+1}$, then this is clear. So now suppose that the prefix of $u$ is $v_{\alpha _1,p+1}$. It is clear that $D(u)$ is to the right of $D(x)$. So we need to show that $D(u)$ is to the left of $D(y)$. For some $k$ with $p+1<k<r$, it must be the case that $t_{k,\alpha _1}$ is a prefix of $u$. Since $D(t_{k,\alpha _1})$ is to the left of $D(v_{k,\alpha _1})$, which contains $D(x)$, the claim is proved. 

Now the proof continues exactly as in \cite{R6}. We choose $z$ with $p(z)\in Z_\infty $ such that  all the letters of $z$ apart from the last letter $C$  are in $\{ L_3,L_2,R_3\} $,  such that $D(z)$ is to the right of $D(y)$, and such that the forward orbit of $p(z)$ does not intersect the subset  $D'$ of the unit disc between $D(x)$ and $D(y)$, that is
$$p(z)\notin\cup _{n\ge 0}s^{-n}(D(x)\cup D'\cup D(y)).$$
 We let $\gamma (z)$ denote the capture path in ${\cal{Z}}(3/7,+,+,0)$ which crosses the upper unit circle into $D(z)$. Now we define paths $\eta _1$   from $p(z)$ to $p(y)$, and a closed path $\omega $ based at $p(z)$. The path $\eta _1$ first  crosses the unit circle  out of the gap of $L_{3/7}$ containing $p(z)$, and then crosses the second time into the gap containing $p(y)$. The closed loop $\omega $ also crosses the upper unit circle out of the gap containing $p(z)$. It then crosses the upper unit circle again at the top edge of the right-most leaf of $L_{3/7}$ in the boundary of the gap containing $p(x)$. It traces an anticlockwise path round the boundary of the subset of the unit disc bounded by $D(x)$ on the left and $D(y)$ on the right, as far as the top of the left-most leaf in the boundary of the gap containing $p(y)$. It then re-crosses the upper unit circle into the gap of $L_{3/7}$ containing $p(z)$, and ends at $p(z)$. 
 
 \begin{figure}
\centering{\includegraphics[width=4cm]{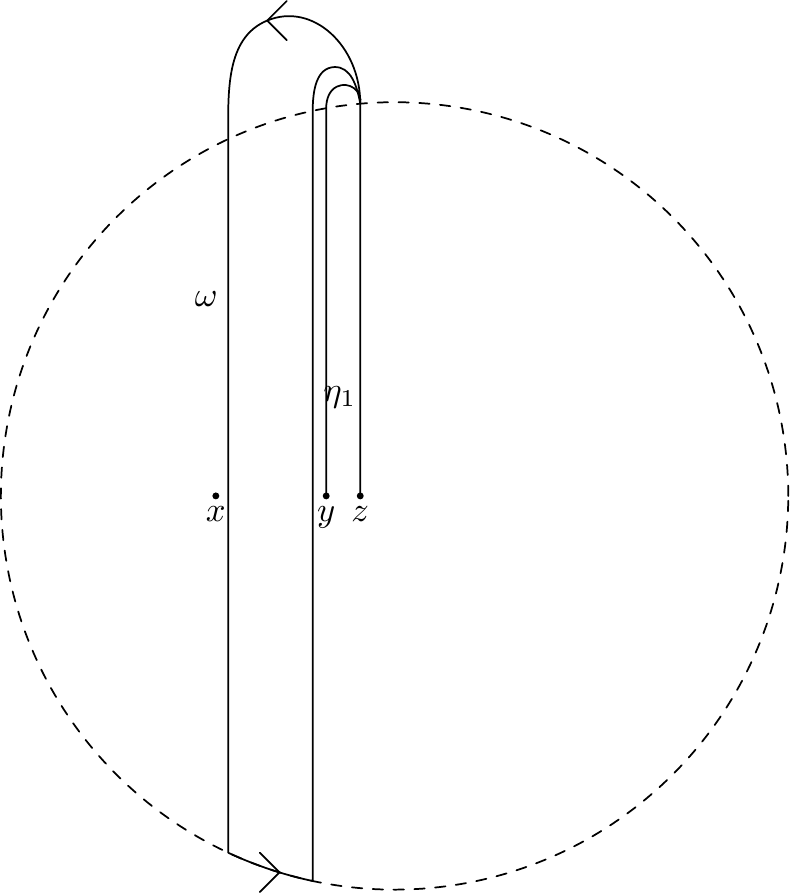}}
\caption{$\omega $}
\end{figure}
It is then immediate that
$$\sigma _{\gamma _1}\circ s\simeq \sigma _{\eta _1}\circ \sigma _{\gamma _z}\circ s.$$
It should be noted that the notation here is different from that in \cite{R6}, where the paths which are here called $\eta _1$ and $\omega $ were called $\gamma _1$ and $\alpha $. For most of this paper, paths $\beta $, with various indexes, denote paths in $R_{m,0}$, and, from now on, paths $\gamma $, with various indexes, denote capture paths, usually in ${\cal{Z}}(3/7,+,+,0)$. In particular, at present, the paths $\gamma _1$ and $\gamma _2$ denote the capture paths $\gamma ^{r,\alpha _1}$ and $\gamma ^{r,\alpha _2}$. 

Since the disc enclosed by $\omega $ is disjoint from the forward orbit of $z$, we have
$$(\sigma _{\gamma (z)}\circ s,Y_m)\simeq _{\psi _m}(\sigma _{\omega }\circ \sigma _{\gamma (z)} \circ s,Y_m)$$
for some homeomorphism $\psi _m$, for any $m$ such that $p(x)$ and $p(y)\in Z_m$. Hence, if 
$$\eta _2=\omega *\psi _m(\eta _1),$$
 we have
$$(\sigma _{\eta _1}\circ \sigma _{\gamma (z)}\circ s,Y_m)\simeq _{\psi _m}(\sigma _{\eta _2} \circ \sigma _{\gamma (z)}\circ s,Y_m).$$
The aim then, as in \cite{R6}, is to compute $\psi _m$ and to show that 
\begin{equation}\label{5.2.1}(\gamma (z)*\eta _2, O^+(p(x)))\simeq (\gamma _2,O^+(p(x))),\end{equation}
where, here, $\simeq $ denotes homotopy constant on $O^+(p(x))$, the forward orbit of $p(x)$. We then have
$$\sigma _{\eta _2} \circ \sigma _{\gamma (z)}\circ s\simeq \sigma _{\gamma _2}\circ s$$
and hence
$$\sigma _{\gamma _2}\circ s\simeq \sigma _{\gamma _1}\circ s,$$
as required.  As usual, $\psi _m$ is a composition of disc exchanges with supports $(\sigma _{\gamma (z)}\circ s)^{-n}(D')$ for varying $n>0$. We only need to consider those disc exchanges which have an effect on $p(y)$. These are the disc exchanges with supports $C(v'a\leftrightarrow v'b,D')$ where $v'b$ is a prefix of $x_{r,\alpha _2}$ of the same length as a prefix $v$ of $x_{r,\alpha _1}$ such that $vau=x_{r,\alpha _1}$. It follows that $\psi _m(\eta _1)$ has endpoint $p(x)$, as required, because the successive changes from $a$ to $b$ are made preceding suffixes $u$ of $x_{r,\alpha _1}$ with $D(u)\subset D'$, that is, precisely the changes that are needed to change $y=x_{r,\alpha _1}$ to $x=x_{r,\alpha _2}$. These disc exchanges also introduce hooks of $\psi _m(\eta _1)$ round the sets $D(v'a)$. But these do not matter, because such a $v'a$ never occurs as a subword of $x$, that is, none of these sets intersects the forward orbit $O^+(p(x))$ of $p(x)$. We see this as foilows. Each  $v'b$ is a prefix of $x_{r,\alpha _2}$ which ends in $x_{r,\alpha _2}$ before a subword $v_{p+1,\alpha _2}$ or $t_{p+1,\alpha _2}$. So $v'b$ must end in a subword $v_{k,\alpha _2}$ for some $k\ge p+1$. It follows that the only occurrences of $v'$ in $x_{r,\alpha _2}$ must be followed by $b$, since $\alpha _2(i)=b$ for all $i\ge p+1$. So (\ref{5.2.1}) follows, as required.

\end{proof}

For the remainder of this section, we obtain  information as possible about the paths $\beta _i(\gamma ^{r,\alpha })$, in preparation for studying perturbations of the capture paths $\gamma ^{r,\alpha }$ in Chapter \ref{6}.. We denote by $n(\alpha )$ the largest integer such that $(\beta _{2n(\alpha )}(\gamma ^{r,\alpha })$ is defined. Even though, by Theorem \ref{5.2},. we have 
$$\beta _{2n(\alpha _1)-1}(\gamma ^{r,\alpha _1})=\beta _{2n(\alpha _2)-1}(\gamma ^{r,\alpha _2}),$$
for all $\alpha _1,\alpha _2\in \{ a,b\} ^{r+1}$, it does not follow, and is certainly not true, that $\beta _i(\gamma ^{r,\alpha _1})=\beta _i(\gamma ^{r,\alpha _2})$ for all $i$.

\section{Lemmas about the $\zeta $ sequence for $x_{r,\alpha }$}\label{5.3}

\begin{lemma}\label{5.3.1} 
\begin{itemize}
\item[1.] $$w_1'(v_0u_0)=L_3,\ \ w_1'(w_0u_0)=L_3L_2R_3L_3^3,$$
$$w_2'(v_0u_0)=v_0L_3L_2R_3L_3^3,\ \ w_2'(w_0u_0)=w_0L_3L_2R_3L_3^3.$$
More generally, if $\alpha (0)=a$ then, for all $r>0$,
$$w_1'(x_{r,\alpha })=L_3,\ \ w_2'(x_{r,\alpha })=v_0t_0aL_3,$$
and if $\alpha (0)=b$ then for all $r>0$,
$$w_1'(x_{r,\alpha })=L_3L_2R_3L_3^3,\ \ w_2'(x_{r,\alpha })=w_0L_3L_2R_3L_3^3L_2y$$
for some $y$.

\item[2.]  For any $r>0$ and any $\alpha $, and any prefixes $u_1$ and $u_1u_2$ of $x_{r,\alpha }$ such that each of $u_1$ and $u_2$ ends in $a$ or $b$ (or $u_1$ can be empty),   the number of $j$ such that 
$$\vert u_1\vert <\vert w_j'(x_{r,\alpha })\vert \leq \vert u_1\vert   +\vert u_2\vert $$
is equal to the number of occurrences of $da$ or $db$ or $dd$ in $u_2$. If $u_2$ ends in $c$, that is, if $x_{r,\alpha }=u_1c$ then we add two for the last occurrence of $dc$.
\item[3.] In particular, for any $\ell >0$, the number of  $j$ for which $\vert w_j'(x_{r,\alpha })\vert \leq \vert v_\ell \vert +1 $ is
$2^\ell $ for all $0<\ell \leq r$, the number for which $\vert w_j'(x_{r,\alpha })\vert \leq \vert v_\ell \vert $ is $2^\ell -1$ if $\alpha (\ell )=a$ and $2^\ell $ if $\alpha (\ell )=b$, and the number of $j$ for which $w_j'(x_{r,\alpha })$ is defined is $2^{r+1}$.

\end{itemize}
\end{lemma}

\begin{proof} 
\noindent 1. This is a direct calculation. 

\noindent 2. Apart from the prefix $v_0$ or $w_0$ of $x_{r,\alpha }$ (depending on whether $\alpha (0)=a$ or $b$) every other occurrence of $v_0$ (or $w_0$) occurs in a subword $dav_0t_0$ or $dbv_0t_0$  (or $daw_0t_0$ or $bw_0t_0$). Also $x_{r,\alpha }$ starts with $v_0t_0$ or $w_0t_0$, and then is made up of subwords $u$ with some overlapping: 
$$dd,\ \ dav_0t_0,\ \ dbv_0t_0,\ \ daw_0t_0,\ \ dbw_0t_0,\ \ \ dc.$$
For each such occurrence  $u$,  if $z$ is the prefix of $x_{r,\alpha }$ which precedes the occurrence of $u$ in $x_{r,\alpha }$, then, apart from $dc$, by the rules described in \ref{3.4} to \ref{3.7},  the corresponding word $w_j'(x_{r,\alpha })$ is:
$$zdL_2R_3L_3,\ \ zdaL_3,\ \ zdL_2R_3L_3,\ \ zdaL_3L_2R_3L_3^3,\ \ zdbL_3L_2R_3L_3^3.$$
  As for $dc$, if $x_{r,\alpha }=zdc$ then the two  words $w_j'(x_{r,\alpha })$ corresponding to $dc$ are:
$$zdL_2R_3L_3,\ \ zdL_3^5.$$
So, as claimed, we simply need to count occurrences of $dd$, $da$ and $db$, and  two for the last occurrence of $dc$. Since the occurrences of $db$ in $x_{r,\alpha }$ are simply obtained by replacing occurrences of $da$ in $x_r$, it suffices to count occurrences of $da$ and $dd$ in $x_r$.

3. A simple calculation shows that the number $N(\ell )$ of occurrences of $dd$ or $da$ or $db$ in $t_k$ is the same as the number in $v_k$, for all $k\geq 1$, and $N(1)=1$. Also, $N(\ell +1)=2N(\ell )+1$. So, by induction, $N(\ell )=2^\ell -1$, and the number of $j$ for which $|w_j'(x_r)|\leq |v_\ell |$ is $2^\ell -1$, adding in $w_1'(x_r)=L_3$, but discounting $w_{2^\ell }'(x_r)=v_\ell L_3$, since this corresponds to the occurrence of $da$ at the end of $v_\ell $. The number of $j$ for which $|w_j'(x_r)|\leq |v_\ell |+1$ is $2^\ell $. The number of occurrences of $dd$ or $dc$ in $u_r$ is the same as the number of $da$ or $dd$  in $v_r$. So the number of words $w_j'(x_r)$ is $2\cdot (2^r-1)+1=2^{r+1}-1$. 

Now we consider $x_{r,\alpha }$ in general. The number of occurrences of $dd$, $da$ and $db$ in $v_{\ell ,\alpha }$ is the same as the number of occurrences of $dd$ or $da$ in $v_\ell $, and the number of occurrences of $dd$, $da$, $db$ or $dc$ in $x_{r,\alpha }$ is the same as the number of $dd$, $da$ or $dc$ in $x_r$. The only difference is that $w_{2^\ell }'(x_{r,\alpha })|<|v_{\ell ,\alpha }|$ if and only if $\alpha (\ell )=b$. So the number of $j$ such that $|w_j'(x_{r,\alpha })|< |v_\ell|$ is the same as the number with $|w_j'(x_{r,\alpha })|\le |v_\ell|+1$ if $\alpha (\ell )=b$.

\end{proof}
 
 Now we consider the sequence $\zeta _i(x)$ of \ref{4.4} for $x=x_{r,\alpha }$. From \ref{5.3} and the definitions in \ref{4.6} and \ref{4.12}, it follows that  if $\alpha (0)=a$, then  $\zeta _i(x)$ is defined for all $i$ with $2i\leq n(x)$  and $(\zeta _i(x_{r,\alpha }):1\leq i\leq 4)$ is of type A. If $\alpha (1)=b$ then $\zeta _i(x)$ is only defined for $i=1$, $2$, but $w_2'(x)=w_1(x)u$ where $L_3u=t_0u'$ for some $u$ and $u'$, and since $v_0<t_0<w_0$, the conditions of \ref{4.14} are satisfied for the existence of a type C quadruple $(\zeta _j:1\leq i\leq 4)$ such that $\zeta _j(x)=\zeta _j$ for $j=1$, $2$, and $\zeta _j=\zeta _j(x_{\alpha ^\ell })$ for $j=3$, $4$, for any $\ell \geq 1$. For most of the time, we will consider $\zeta _i(x_{r,\alpha })$ for $\alpha $ with $\alpha (0)=a$. We will see later that, if $\alpha (0)=b$, then the sequence to consider is $\zeta _i(x_{r,\alpha '})$, where $\alpha '(0)=a$ and $\alpha '(\ell )=\alpha (\ell )$ for all $\ell >0$.
 
For the moment, suppose that $x=x_{r,\alpha }$ with $\alpha (0)=a$. It is reasonable to expect that all but the last few disc crossings by $\zeta _{2i-1}(x)$ and $\zeta _{2i}(x)$ are in the sequence of crossings determined by $w_j'(x)$. Recall that $k=k(\zeta _{2i-1},\zeta _{2i})$ is the first integer for which the $k$'th unit-disc-crossings of $\zeta _{2i-1}$ and $\zeta _{2i}$ do not coincide. For the moment, we write $k(i,x)$ for this integer, so that $w_j'(\zeta _{2i-1}(x))=w_j'(\zeta _{2i}(x))$ for $j<k(i,x)$. The integer $k(i,x)$ is always odd,  and $k(i,x)\leq k(i+1,x)\leq k(i,x)+2$ for all $i$.
 
  \begin{lemma}\label{5.5} 
  Let $x=x_{r,\alpha }$ for any $r\geq 0$ and any $\alpha $ with $\alpha (0)=a$. We write $k(i,x)=k(\zeta _{2i-1}(x),\zeta _{2i}(x))$.
   \begin{itemize}
   \item[1.] Apart from $(\zeta _1,\zeta _2)$, there are two pairs $(\zeta _{4i-1}(x),\zeta _{4i}(x))$ and \\ 
   $(\zeta _{4i+1}(x),\zeta _{4i+2}(x))$, with $i$ even with 
   $$w_{k-1}'(\zeta _{2i+s}(x))=w_{k-1}'(x),\ \ -1\leq s\leq 2$$
   for each odd $k\geq 5$. There are two extra pairs for each prefix $zv_0L_3L_2R_3L_3^3$ of $x_{r,\alpha }$ such that $z$ has an even number of letters $L_3$ and $L_2$.
   \item[2.] $(\zeta _{4i+t}(x):-3\leq t\leq 0)$ is of type A if $w_{k(2i-1,x)-1}'(\zeta _{2i-1}(x))=w_{k(2i-1,x)-1}'(x)$ and of type B otherwise, in which case 
   $$w_{k(2i-1,x)-1}'(\zeta _{2i-1}(x))=zv_0L_3L_2R_3L_3^3$$
    is a prefix of $x_r$ such that $z$ has an even number of letters $L_3$ and $L_2$. 
   \item[3.]    For any $\alpha $, and any $r\geq 1$, the pair  $(\zeta _{2i-1}(x),\zeta _{2i}(x))$ is defined for  
  $$i\leq 2^{r+1}+2^{r-1}-1.$$
 Moreover,  we have 
 $$\vert w_{k(i,x)-1}'(\zeta _{2i-1}(x))\vert <  \vert v_1 \vert $$
  for $i\leq 3$, and for any $2\leq \ell \leq r$, we have 
  $$\vert w_{k(i,x)-1}'(\zeta _{2i-1}(x))\vert <  \vert v_\ell \vert $$
   for $i\leq 2^\ell +2^{\ell -1}-1$ if $\alpha (\ell )=a$, and for $i\leq 2^\ell +2^{\ell -1}$ if  $\alpha (\ell )=b$.
 If $\alpha (\ell )=a$ then for  the next pair $(\zeta _{2i-1}(x),\zeta _{2i}(x))$, we have 
    $$|w_{k(i,x)-1}'(\zeta _{2i-1}(x))|=|w_{k-1}'(\zeta _{2i}(x)|=|v_\ell |+1.$$
      \item[4.] 
 For $1\leq \ell \leq r$, we have 
 $$\vert w_{k(i,x)-1}'(\zeta _{2i-1}(x))\vert \leq \vert v_\ell \vert +1\ \ \Leftrightarrow\ \ k(i,x)\leq 2^\ell ,$$
  and  
 $$\vert w_{k(i,x)-1}'(\zeta _{2i-1}(x))\vert =2^{r+1}-1\ \ \Leftrightarrow\ \ k(i,x)=2^{r+1}-1.$$
  \end{itemize}
 \end{lemma}

\begin{proof} First we consider the case of $x=x_r$. We write $\zeta _i(v_ru_r)=\zeta _i$ and $k_i=k(i,v_ru_r)$. Recall that, apart from the occurrence at the beginning of the word, each occurrence of $v_0t_0$ in $v_ru_r$ is preceded by $da=L_3^2a^2$. Now we consider the possible words  $w_j'(v)$ for any word $v$ for prefixes $v$ of $v_ru_r$ ending in $L_3L_2$, for any even value of $j$. These are the possible words for $w_{k_i-1}'(\zeta _{2i-1})=w_{k_i-1}'(\zeta _{2i})$, for any $i$. Of course the words  $w_j'(x)$ are possibilities, for even $j<2^{r+1}$. We have seen in \ref{5.3.1} that each word $w_j'(x_r)$ --- without the restriction to even $j$ --- is determined by an occurrence  of $da$ or of $dd$, apart from $j=1$ (which corresponds to the first occurrence of $v_0$) and the three largest possible values of $j$. More precisely, if $z_1=zdd$ is a prefix of $v_ru_r$, then the largest $j$ with $\vert w_j'(v_ru_r)\vert < \vert z_1\vert $ gives $w_j'(v_ru_r)=zdL_2R_3L_3$. If $z_2=zdav_0t_0$, then the largest $j$ with $\vert w_j'(v_ru_r)\vert <\vert z_2\vert $ gives $w_j'(v_ru_r)=zdaL_3$.   In both cases, if $z=z'd$, then $zL_2R_3L_3=w_{j-1}'(v_ru_r)$, but this is associated to the previous occurrence of $dd$ in $v_ru_r$.  The other alternative is that $z=z'v_0L_3L_2R_3$, and in this case $w_{j-1}'(v_ru_r)\vert <\vert z\vert $. By induction, there are $2^r-1$ occurrences of  $da$ in $x_r$, and $2^r-2$ occurrences of $dd$. Of these, $2^{\ell -1}$ and $2^{\ell -1}-1$ respectively occur in the prefix $v_\ell $ of $x_r$ if $\ell \geq 1$.  We have seen that $j=2^\ell $ for the occurrence of $dav_0t_0$ preceded  in $v_ru_r$ by $z$, for which the $zda=v_\ell $ for any $\ell \geq 2$. In particular, this  $j$ is even. It follows that if $zv_0t_0$ is a prefix of $v_\ell $, then the integers $j_1$ and $j_2$ such that $w_{j_1}'(v_ru_r)=zL_3$ and $w_{j_2}'(v_ru_r)=v_\ell zL_3$ have opposite parities. Hence, by induction, for each $r\geq \ell \geq 2$, the number of prefixes $zda$ of $v_\ell  $ with $z$ nontrivial and for which $w_j'(zdav_0t_0)=zdaL_3$ and for which $j$ is even, is the same as the number of such prefixes with $j$ odd. So both these numbers are $2^{\ell -2}$.  Similarly, the numbers for such prefixes $zv_0$ of $v_ru_r$,   are $2^{r-1}$ for both $j$ even and $j$ odd, for all $r\geq 1$. 

We have $k(1,x)=1$ and $k(i,x)=3$ for $2\leq i\leq 5$. We have
$$w_2'(\zeta _j)=v_0L_3L_2R_3L_3^3{\rm{\ for\ }}3\leq j\leq 6,$$
and
$$w_2'(\zeta _j)=v_0t_0aL_3{\rm{\ for\ }}7\leq j\leq 10.$$
The reason for the extra pairs is that 
$$w_3'(v_0L_3L_2R_3L_3L_2BC)=w_3'(v_0L_3L_2R_3L_3^6L_2)=v_2L_3L_2R_3L_3^3,$$ 
and 
$v_0L_3L_3R_3L_3L_2BCL_1R_2R_3$ and $v_0L_3L_2R_3L_3^6L_2$  are prefixes of $w_3'(\zeta _3)$ and $w_3'(\zeta _4)$ respectively. It is therefore clear that $I(\zeta _3,\zeta _4)$ is over $D(x_r)$. Similarly we have extra pairs for each prefix $zv_0L_3L_2R_3L_3^3$ of $x_r$ such that $z$ has an even number of letters $L_3$ and $L_2$. We do not have any other extra pairs, because although there are other prefixes $z_1L_3^3$ of $x_r$, these are prefixes of the form $z_2dL_3^3=z_2L_3^3(L_2R_3)^2L_3^3$ rather than $z_2v_0L_3L_2R_3L_3^3$.  It follows from the detail of the rules in \ref{3.4} to \ref{3.7} that there are no extra pairs in this case.

So  the pairs $(\zeta _{2i-1}(x),\zeta _{2i}(x)$ correspond two-to-one with $w_k'(x_r)$ for odd values of $k\geq 5$ and with $k\leq 2^{r+1}-1$, and two extra pairs for each prefix $zv_0L_3L_2R_3L_3^3$ of $x_r$ such that $z$ has an even number of letters $L_3$ and $L_2$. We then have 

$$w_{k-2}'(zv_0L_3L_2R_3L_3^3L_2)=zL_3=w_{k-2}'(x_r),$$
but $w_{k-1}'(x_r)\neq w_{k-1}'(zv_0L_3L_2R_3L_3^3L_2)$. If $k$ is odd, and $k-1$ is even, there are pairs $(\zeta _{4i-1},\zeta _{4i})$ and $(\zeta _{4i+1},\zeta _{4i+2})$, for some $i$, such that $k(2i,x)=k(2i+1,x)$ and $w_{k-1}'(\zeta _{4i+s})=zL_3^3$ for $-1\leq s\leq 2$, and $w_q'(\zeta _{2i+s})=w_q'(v_ru_r)$ for $q\leq k-2$. The exceptional feature of these pairs is that $(\zeta _{4i+s}:1\leq s\leq 4)$ is of type B, with 
$$k=k(\zeta _{4i+s-1},\zeta _{4i+s})=k(\zeta _{4i-1},\zeta _{4i})$$
for $0\leq s\leq 3$, and
$$w_{k-1}'(\zeta _{4i+t})\neq w_{k-1}'(\zeta _{4i+u})=w_{k-1}'(x),\ \ -1\leq t\leq 2,\ \ 3\leq y\leq 6.$$
So there are four pairs corresponding to each such value of $w_{k-1}'(x)$ with $k-1$ even. These are pairs $(\zeta _{4i+2s-1}(x),\zeta _{4i+2s}(x))$ for $0\leq s\leq 3$ for some $i$. 

The number of occurrences of $v_0t_0$ in $x_r=v_ru_r$ is $2^{r+1}$ and the number of occurrences in $v_\ell =v_{\ell -1}t_{\ell -1}a$ is $2^{\ell }$, for all $\ell \geq 1$. For all $r\geq 1$, if $uv_0t_0$ is a prefix of $v_r$, then $v_ruv_0t_0$ is a prefix of $x_r$, and vice versa. Moreover, exactly one of $v_r$ and $u_r$ has an even number of letters $L_3$ and $L_2$, and hence the number of prefixes  $zv_0t_0$ of $x_r$ for which $z$ has an even number of letters $L_3$ and $L_2$ is $2^r$. Similarly for all $\ell \geq 2$ --- which we need so that $v_0t_0$ is a prefix of $v_{\ell -1}$ and $t_{\ell -1}$ --- the number of prefixes  $zv_0t_0$ of $v_{\ell }$ for which $z$ has an even number of letters $L_3$ and $L_2$ is $2^{\ell -1}$. 

 It follows that the number of pairs $(\zeta _{2i-1}(x),\zeta _{2i}(x))$ is 
$$1+4\cdot 2^{r}=1+2^{r+1}+2^r,$$
and the number for which $\vert w_{k-1}'(\zeta _{2i-1})\vert <\vert v_1 \vert $ is $3$, and, for $\ell \geq 2$, the number for which $\vert w_{k-1}'(\zeta _{2i-1})\vert <\vert v_\ell \vert $  is 1 ---  plus twice times the number of odd $k$ (even $k-1$) such that $|w_{k-1}'(x)|<|v_\ell |$, which is $2^{\ell }-2$ -- plus twice  for the even number $k-1$ such that $|w_{k-1}'(x)|=|v_\ell |+1$ -- which is $2^{\ell -1}-1$ --- which gives $2^{\ell +1}-1$, as required. The final pair for which $k(i,x)\leq 2^{\ell }+1$ is $i=2^{\ell }+2^{\ell -2}$.
 $$i=2^{\ell }+2^{\ell -1}.$$ 
The possible values of $k(i,x)$ are precisely the odd $k$ for which $w_k'(x_{r,\alpha })$ is defined. Thus, 4 for $x_r$ follows from \ref{5.3} for $x_r$, and 1 to 4 are proved for $x_r$. The proof for $x_{r,\alpha }$ is exactly the same, apart from minor differences  because some occurrences of $dav_0t_0$ are replaced by $dbv_0t_0$.

\end{proof}

Let $x=x_{r,\alpha }$ for any $r\geq 0$ and any $\alpha $ and let $\gamma (x)$ be the associated capture path in ${\cal{Z}}_m(3/7,+,+,0)$ and write $\beta _i(x)$ for the associated sequence $\beta _i(\gamma (x))$ as in \ref{3.3}. We write $E_{4j-3,4j-1,x}$ for the sequence $E_{4j-3,4j-1,\beta (x)}=E_{\beta _{4j-3}(x),\beta _{4j-1}(x),\beta (x)}$ as in \ref{4.23}, following the earlier notation of \ref{3.8} and \ref{4.6}. If $\alpha (0)=a$, we also write
 $$E_{4i-3,4i-1,\alpha }'=E_{4i-3,4i-1,x }'=E_{\zeta _{4i-3}(x),\zeta _{4i-1}(x)}$$
Since $\zeta _i(x)$ depends only on $\alpha (j)$ for $j\leq t$,  we see that $E_{4i-3,4i-1,x}'=E_{4i-3,4i-1,\alpha }'$ depends only on $\alpha (j)$ for $j\leq t$, so long as  $\vert w_{k(2i,x)-1}'(x)\vert \leq 2\cdot \vert v_t\vert  -7$. In particular, $E_{1,3,x}'$ is depends only on $\alpha (0)$. We also define
 $$\psi _{4i-3,4i-1,x}'=\psi _{\zeta _{4i-3}(x),\zeta _{4i-1}(x)}.$$
  We know from \ref{4.13}  that $\zeta _i(x)=\beta _i(x)$ for $i\leq 4$ if $\alpha (0)=a$.  Thus, $E_{1,3,x}=E_{1,3,x}'$ if $\alpha (0)=a$. We also know that $(\beta _t(x)):1\leq t\leq 4)$ is of type C if $\alpha (0)=b$, with $\beta _t(x)=\zeta _t(x)$ for $t\leq 2$ and $\beta _t(x)=\zeta (x_{\alpha _1,r})$ for $t=3$, $4$,where $\alpha _1(0)=a$ and $\alpha _1(k)=\alpha (k)$ for $k\geq 1$.

 \begin{lemma}\label{5.7}
Let $x=x_{\alpha ,r}$ for any $\alpha $ with $\alpha (0)=a$, and any $r\geq 1$.

For all $i$,
 \begin{equation}\label{5.7.1}k(2i,x)=k(2i+1,x)\leq k(2i+2,x),\end{equation}
 \begin{equation}\label{5.7.2}\zeta _{4i}(x)<\zeta _{4i+s}(x){\rm{\ for\ all\ }}s\geq 1.\end{equation}
 \begin{equation}\label{5.7.3}k(2i+2,x)>k(2i+1,x)\Rightarrow k(2i+2,x)=k(2i,x)+2{\rm{\ for\ all\ }}s\geq 4\end{equation}
 \begin{equation}\label{5.7.4}k(2i+2,x)=k(2i+1,x)\Leftrightarrow \zeta _{4i+2}(x)<\zeta _{4i+s}(x){\rm{\ for\ all\ }}s\geq 3.\end{equation}
 
\begin{equation}\label{5.7.5}k(2i+2,x)>k(2i,x)\Leftrightarrow E_{4j+1,4j+3,x}'\subset  E_{4i+1,4i+3,x}'{\rm{\ for\ all\ }}j>i\end{equation}
\begin{equation}\label{5.7.6}k(2i+2,x)=k(2i,x)\Leftrightarrow E_{4j+1,4j+3,x}'\cap  E_{4i+1,4i+3,x}'=\emptyset {\rm{\ for\ all\ }}j>i\end{equation}
\end{lemma}

\begin{proof} It is easily checked that (\ref{5.7.1}) and (\ref{5.7.2}) hold for the paths $\zeta _{4i+s}(x)$ which were found in \ref{5.5}. Also we have $k=k(2i+2,x)=k(2i,x)+2$ for exactly all the pairs $(\zeta _{4i+3}(x),\zeta _{4i+1}(x))$ for which $w_{k-1}'(\zeta _{4i+3}(x)=w_{k-1}'(x)$, as we saw in \ref{5.5}. We saw in \ref{5.5} that (\ref{5.7.3}) holds for these pairs and (\ref{5.7.4}) holds for the others. Then (\ref{5.7.5}) and (\ref{5.7.6}) follow immediately from (\ref{5.7.1}) to (\ref{5.7.4}).

\end{proof}

\begin{lemma}\label{5.8}

Once again, let $x=x_{r,\alpha }$ with $\alpha (0)=a$. The finite sequence  of sets $\{ E_{4i_j-3,4i_j-1,x}':j\leq n\} $ is strictly decreasing, that is, satisfies 
$$E_{4i_{j+1}-3,4i_{j+1}-1,x}'\subset E_{4i_j-3,4i_j-1,x}',\ \ E_{4i_{j+1}-3,4i_{j+1}-1,x}'\neq E_{4i_j-3,4i_j-1,x}'$$
 for all $1\leq j<n$, if and only if $k(2i_j,x)<k(2i_j+2,x)$ for all $j<n$, that is, if and only if $w_{k-1}'(\zeta _{4i_j-3}(x))=w_{k-1}'(x)$ for $k=(2i_j,x)$ for all $1\leq j<n$.

Now let  $\ell \geq 1$, and let $n_1(\ell )=n_1(\ell ,\alpha )$ be the largest even integer such that 
$$|w_{k(j,x)-1}'(\zeta _{2j-1}(x))|=|w_{k(j,x)-1}'(\zeta _{2j}(x))|<|v_\ell |$$
for $2j\leq n_1-1$.  (By \ref{5.5}, $n_1(1)=4$ and $n_1(\ell )=2^{\ell }+2^{\ell -1}$ for $\ell \geq 2$ is independent of $\alpha $). Then the following hold.
\begin{itemize}
\item[1.] Let $u=dd$ or $dv_0t_0$ or $av_0t_0$ or $bv_0t_0$. Then
$$D(t_{\ell ,\alpha }u)\subset E_{2n_1-3,2n_1-1,x}'\ \Leftrightarrow  \alpha (\ell )=\begin{cases}a{\rm{\ if\ }}\ell =1,\\ b{\rm{\ if\ }}\ell >1.\end{cases}$$

Exactly similar results hold for $u_{\ell ,\alpha }$.
\item[2.] Let  $\ell >1$ and let $\alpha _1$ and $\alpha _2$ be such that $\alpha _1(\ell ')=\alpha _2(\ell ')$ for $\ell '<\ell $, but $\alpha _1(\ell )=b$ and $\alpha _2(\ell )=a$, so that  $t_{\ell ,\alpha _1}=t_{\ell ,\alpha _2}$. Write $x_1=x_{r,\alpha _1}$ and $x_2=x_{r,\alpha _2}$.  Write $n_1(\ell )=n_{1}$. Then 
$$D(v_{\ell ,\alpha _1})\subset E_{2n_{1}-3,2n_{1}-1,x_1}',$$
$$D(v_{\ell ,\alpha _2})\subset E_{2n_{1}-3,2n_{1}-1,x_2}'\subset E_{2n_{1}-3,2n_{1}-1,x_1}',$$
but
$$D(v_{\ell ,\alpha _1})\cap E_{2n_{1}-3,2n_{1}-1,x_2}'=\emptyset .$$

\end{itemize}

Hence, the number of sets $E_{4i-3,4i-1,x_1}'$ containing  $D=D(t_{\ell ,\alpha _1}L_3L_2)=D(t_{\ell ,\alpha _2}L_3L_2)$ is one more than the number of sets $E_{4i-3,4i-1,x_2}'$ containing $D$. These numbers are $2^{\ell -1}$ and $2^{\ell -1}-1$ respectively.  In contrast, the number of sets $E_{4i-3,4i-1,x_1}'$ containing $D'=D(v_{\ell ,\alpha _1})$ is the same as the number of sets $E_{4i-3,4i-1,x_2}'$ containing $D(v_{\ell ,\alpha _2})$. These numbers are both $2^{\ell -1}$.

  \end{lemma} 
 
 \begin{proof} The first statement follows immediately from \ref{5.7}, in particular, from (\ref{5.7.5}) and (\ref{5.7.6}).
 
 Now we consider the proof of 1. First consider the case $\ell =1$. Since we are taking $\alpha (0)=a$, the set $E_{1,3}=E_{1,3}'$ is the same whether $\alpha (1)=a$ or $b$. In both cases, we have $w_2'(\zeta _i)=v_0L_3L_2R_3L_3^3$ for $3\leq i\leq 6$ and $k(\zeta _{2i-1},\zeta _{2i})=3$ for $2\leq i\leq 4$. If $\alpha (1)=a$ then 
 $$w_2'(\zeta _7)=w_2'(\zeta _8)=v_0t_0aL_3=v_0L_3L_2R_3L_3^2a^2L_3.$$
 In this case, for any path $\zeta \in R_{m,0}$ such that $\zeta $ intersects $D(t_1u)$ for any $u$ as in 1, all unit-disc crossings from the second onwards  (determined by $w_2'(\zeta )$) are in $D(t_1)=D(t_{1,\alpha })$, and this second unit-disc crossing is between those of $\zeta _6$ and $\zeta _7$ (determined by $w_2'(\zeta _6)$ and $w_2'(\zeta _7)$). We therefore have $\zeta _6<\zeta <\zeta _7$, and all of $\zeta $, from the second unit-disc crossing onwards, is in $E_{5,7}'$ and we therefore have $D(t_1u)\subset E_{5,7}'$ for all such $u$.
 
   If $\alpha (1)=b$, then
  $$w_2'(\zeta _7)=w_2'(\zeta _8)=v_0t_0L_2R_3L_3=v_0L_3L_2R_3L_3^2aL_2R_3L_3.$$
 We still have $\zeta _6<\zeta _7$. Now for any other path $\zeta $ such that  $\zeta $ intersects  $D(t_1)=D(t_{1,\alpha })$,we have $w_2'(\zeta )=w_2'(\zeta _7)=w_2'(\zeta _8)$ and, since the third unit-disc crossing $\zeta _8$ is the furthest right of all such paths, the third unit disc crossing of $\zeta $ is to the left of the third unit-disc-crossing of $\zeta _8$, but to the right of the second unit-disc-crossing of $\zeta $. We therefore have $\zeta _8<\zeta $ and hence $\zeta $ does not intersect $E_{5,7}'$ in this case, and we deduce that $D(t_1u)\cap E_{5,7}'=\emptyset $ in this case. 
 
 Now let $\ell \geq 2$. Then, writing $n_1(\ell )=n_1$,
 $$k(n_1-1,x)=k(n_1,x)-2=k(\zeta _{2n_1,-3},\zeta _{2n_1-1}),$$
  and we have 
 $$\zeta _{2n_1-1}<\zeta _{2n_1}<\zeta _{2n_1-3}<\zeta _{2n_1-2},$$
which is different from the case $\ell =1$. Write $k=k(n_1,x)$. Write $t_{\ell ,\alpha }=zdd$. For all the choices considered of $u$, we have 
 $$w_{k-1}'(t_{\ell ,\alpha }u)=zdL_2R_3L_3u_1,$$
  for a word $u_1$ depending on $u$. We also have
  $$w_{k-1}'(\zeta _{2n_1-1})=w_{k-1}'(\zeta _{2n_1})=\begin{cases}zdaL_3{\rm{\ if\ }}\alpha (\ell )=a,\\ zdL_2R_3L_3{\rm{\ if\ }}\alpha (\ell )=b.\end{cases}$$
  For any path $\zeta \in R_{m,0})$ which intersects $D(t_{\ell ,\alpha }u)$, we have
  $$w_{k-2}'(\zeta )=w_{k-2}'(\zeta _{2n_1-1}),$$
  and hence, for the same reasons as when $\ell =1$ we have $\zeta _{2n_1 }<\zeta $ if $\alpha (\ell )=b$ and $\zeta <\zeta _{2n_1-1}$ if $\alpha (\ell )=a$. But because of the difference between the ordering of $\zeta _{2n_1-3}$ and $\zeta _{2n_1-1}$ from the case of $\ell =1$, we have
  \begin{equation}\label{5.8.1}D(t_{\ell ,\alpha }u)\subset E_{2n_1-3,2n_1-1,\alpha }'{\rm{\ if\ }} \alpha (\ell )=b,\end{equation}
   \begin{equation}\label{5.8.2}D(t_{\ell ,\alpha }u)\cap E_{2n_1-3,2n_1,\alpha }'=\emptyset{\rm{\ if\ }}\alpha (\ell)=a.\end{equation}
 Now we consider the proof of 2.  If $x_1=x_{r,\alpha _1}$ and $x_2=x_{r,\alpha _2}$ where $\alpha _1(\ell ')=\alpha _2(\ell ')$ for $\ell '<\ell $ but $\alpha _1(\ell )=b$ and $\alpha _2(\ell )=a$, then for the same reasons as before, we have
  $$E_{2n_1-3,2n_1-1,x_2}'\subset E_{2n_1-3,2n_1-1,x_1}'.$$
  and
  $$D(v_{\ell ,\alpha _j})\subset E_{2n_{1}-3,2n_{1}-1,x_j}',$$
  for $j=1$, $2$, and
  $$D(v_{\ell ,\alpha _1})\cap E_{2n_1-3,2n_1-1,\alpha _2}'=\emptyset .$$
  
From (\ref{5.8.1}) and (\ref{5.8.2}),  the number of sets $E_{4i-3,4i-1,x_1}'$ which contain $D(t_{\ell ,\alpha _1})$, which is a decreasing sequence of sets ending in $E_{2n_1-3,2n_1-1,x_1}'$, is one more than the number of sets $E_{4i-3,4i-1,x_2}'$ which contain $D(t_{\ell ,\alpha _2})$ (which is again a decreasing sequence of sets ending in $E_{2n_1-3,2n_1-1,x_2}'$). We have seen that the  first number is the number of $i$ with $2i\leq n_1$ such that $\zeta _{4i-3}<\zeta _{4i-5}$, that is the number of $i$ such that $E_{4i-3,4i-1,x_1}'$ is of type A, that is, by \ref{5.5}, $2^{\ell -1}$ sets if $\ell \geq 2$ and just one set if $\ell -1$. The number of sets $E_{4i-3,4i-1,x_2}'$ which contain $D(t_{\ell ,\alpha _2})$ is one less if $\ell \geq 2$ and one more if $\ell =1$, that is, $2^{\ell -1}-1$ sets if $\ell \geq 2$ and two sets if $\ell =1$.
   
   \end{proof}
 
The following  lemma follows from \ref{4.10}. 
\begin{lemma}\label{5.9}Let $x=x_{r,\alpha }$ for any $\alpha $ and $r$. If $uE_{4j-3,4j-1,x}'$ intersects $I(\zeta _{4i-1-2s}(x),\zeta _{4i-2s}(x))$ for at least one of $s=0$ or $1$ then
\begin{equation}\label{5.9.1}I(\zeta _{4i-1}(x),\zeta _{4i}(x))\cup I(\zeta _{4i-3}(x),\zeta _{4i-2}(x))\cup E_{4i-3,4i-1,x}'\subset uE_{4j-3,4j-1,x}'\end{equation}
\end{lemma}

\section{The $\beta $ sequence for $x_{r,\alpha}$}\label{5.10}

We now interpret \ref{4.19} for $\gamma =\gamma ^{r, \alpha }$, for any $\alpha \in \{ a,b\} ^{\mathbb N}$, and obtain more precise information about the paths $(\beta _{2i-1}(\gamma ),\beta _{2i}(\gamma ))$. We write $x=x_{r,\alpha }$ when no confusion can arise. In particular, we give some information about how the sequences $\{ (\beta _{2i-1}(\gamma ),\beta _{2i}(\gamma ))\} $ and  
$$\{ (\psi _{1,4\lfloor (i-1)/2\rfloor -1,\gamma }(\beta _{2i-1}(\gamma )),\psi _{1,4\lfloor (i-1)/2\rfloor -1,\gamma }(\beta _{2i}(\gamma )))\} $$
 differ from the sequence $\{ (\zeta _{2j-1}(x),\zeta _{2j}(x))\} $. We sometimes write $\beta _i(x_{r,\alpha })$ or $\beta _i(x)$ for $\beta _i(\gamma )=\beta _i(\gamma ^{r,\alpha })$, or even $\beta _i$, if no confusion can arise.

\begin{subsectiontheorem}\label{5.10.1} Let $x=x_{r,\alpha }$ for any $r$ and $\alpha $ and et $\gamma =\gamma ^{r,\alpha }$. Then for any $i\ge 1$ with $2i\le n(\gamma )$, paired  satellite arcs for $\psi _{1,4i-1,\gamma }$ of type 3b) as in \ref{4.19} do not occur. Consequently every quadruple $(\beta _{4i+t}(\gamma ):-3\le t\le 0)$ is either an upper or lower $\zeta $ quadruple, or a type C quadruple. Type C quadruples only occur if $\alpha (0)=b$. Moreover, the following hold. Cases 2 and 3 only occur for satellite quadruples over $x$. 
\begin{itemize}
\item[1.] $\beta _{4i+t}(x)=\zeta _{4j+t}(x')$ for $-3\leq t\leq 0$ and for $x'$ and $j$ depending on $i$, where $(\zeta _{4j+t}(x'):-3\le t\le 0)$ is an upper $\zeta $ quadruple for $x'$.  If $\alpha (0)=a$, then every subword $v_0t_0$ in $x$ is a subword in the same position in $x'$.  If $\alpha (0)=b$, then every subword $w_0$ in $x$ is replaced by $v_0$ in the same position in $x'$. Once that is done, some occurrences of $a$ or $b$ preceding occurrences of $v_0t_0$ may be replaced by $b$ or $a$
\item[2.] $I(\beta _{4i+t-1}(x),\beta _{4i+t}(x))=uI(\zeta _{4j+t-1}(x'),\zeta _{4j+t}(x'))$ for $t=-2$ and $t=0$, and for $x'$ and $j$ depending on $i$, where all letters of $u$ are in $\{ L_3,L_2,R_3\} $, and an odd number of letters are in $\{ L_3,L_2\} $, and with the same restrictions on $ux'$ as in 1.
\item[3.] For some $k\leq k(2i-1,x)$, the value of $w_{k}'(\beta _{4i+t}(x))$ (for $-3\leq t\leq 0$) is maximal or minimal subject to the value of $w_{k-1}'(\beta _{4i+t}(x))$

\end{itemize}

For any $1\le \ell \le i$, the support of $\psi _{4\ell -3,4i-1,x}$ can only intersect $D(x')$ for any $x'$as above, if $(\beta _{4i+t}(x):-3\le t\le 0)$ is of form 1 or 2.  Moreover, in each of these cases, the pair of principal arcs of $\psi _{1,4i-5,x}\cdot (\beta _{4i+t-1}(x),\beta _{4i+t}(x))$ (for $t=-2$, $0$) is of one of the following forms, up to isotopy preserving $S^1\cup Z_\infty $:
\begin{enumerate}[a)]
\item $I(\zeta _{4\ell +t-1}(x),\zeta _{4\ell +t}(x))$ for some $\ell $, and $t=-2$, $0$, or is in the image of  $I(\zeta _{4\ell +t-1}(x'),\zeta _{4\ell +t}(x'))$ under a single disc exchange corresponding to a basic exchange $a\leftrightarrow b$, where $x'$ is the image of $x$ under that exchange;
\item both arcs are  parallel to $I(\zeta _{4\ell -3}(x),\zeta _{4\ell -2}(x))$ for some $\ell $;
\item both arcs are parallel to a common  arc of $\zeta _{4\ell +s}(x)$ for $-3\leq s\leq 0$ for some $\ell $ where the arc is between $S^1$ crossings, up to or before the $k(4\ell -3,x)-1$'th crossing.
\end{enumerate}

Similar properties hold if $\psi _{4\ell -3,4i-1}$ is replaced by $\psi '$ for any decomposition $\psi _{4\ell -3,4i-1}=\psi ''\circ \psi '$, where each of $\psi '$ and $\psi ''$ is a composition of disc exchanges. 
\end{subsectiontheorem}

{\textbf{Remark}} The difference between 1 and 2 is the difference between the supports of disc exchanges spanning the upper and lower circle. The extra possibility at the end of a) arises because of the differences in the sequences $\zeta _j(x_{r,\alpha })$ and $\zeta _j(x_{r,\alpha '})$ if $\alpha (k)\neq \alpha '(k)$ for some $k$.

\begin{proof} Naturally, this is proved by induction. We can start the induction because the results are true for $i=1$.  Now suppose the results are true  for $(\beta _{4j+t}(\gamma ):-3\le t\le 0)$, and for $\psi _{1,4j-1,\gamma }$, for all $j< i$. As in \ref{4.19}, we aim to  obtain the results for $(\beta _{4i+t}(\gamma ):1\le t\le 4)$ and for $\psi _{1,4i+3,\gamma }$, by using the information already obtained about $\psi _{1,4i-5,\gamma }^{-1}$, and about the disc exchanges of which it is a composition. We need to consider the action of $\psi _{1,4i-5,\gamma }^{-1}$ on $D(x)$. 

We define $x'$ to be {\em{of permissible form}} if $x'$ is obtained from $x_r$ by replacing some occurrences of $a$ or $b$ preceding occurrences of $v_0t_0$ by $b$ or $a$. This is equivalent to the description of $x'$ in terms of $x_{r,\alpha }$ in statement 1 of the theorem. By the inductive hypothesis, the support of each disc exchange in the composition which can have an effect on $D(x)$ is of the form $C(ya\leftrightarrow yb,E_{4j-3,4j-1,z}')$ where $E_{4j-3,4j-1,z}'=E_{\zeta _{4j-3}(z),\zeta _{4j-1}(z)}$ and $y$ is some admissible word. Then we see  that in order to have an effect on $D(x)$, the first disc exchange $\xi _1$ in the composition whose support intersects $D(x)$  must have $ya$, $yb$ and $z$ of permissible  form, and then $\xi _1(D(x))=D(x^1)$ where $x^1$ is also of permissible form. The exchanges for $\beta _j(\gamma )$ are all basic exchanges, because $w_1'(\beta _j)=w_1'(\beta _3)=L_3$ for all $j\ge 3$.  Write $\psi _{1,4i-5}^{-1}\xi _n\circ \cdots \xi _1$ on $D(x)$.  and write $\xi _{\ell ,1}=\xi _\ell \circ \cdots \xi _1$.  Let the support of $\xi _\ell $ be $C(y_\ell a\leftrightarrow y_\ell b,E_j)$, where $E_j=E_{4j_1(\ell )-3,4j_1(\ell )-1,z_\ell }$ for $e_\ell \subset \{ a,b\} $, so that $\xi _{\ell ,1}(D(x))=D(x^{\ell })$ and $x^\ell =y_{\ell +1}e_{\ell +1}'z_{\ell +1}$ for $e_{\ell +1}'\in \{ a,b\} $, as it is in the support of the next disc exchange $\xi _{\ell +1}$. We do not have $E_j$ of the form $E_{4n -3,4n-1,\gamma }$ for $(\beta _{4n +t}:-3\le t\le 0)$ as in 3 in the statement of the lemma, because in that case, by the form of $y_\ell e_\ell z_\ell$, $uE_{4n-3,4n-1,x}$ does not intersect $D(y_\ell e_\ell z_\ell )$ for any $\ell $. Because of the form of the words $y_\ell e_\ell b_\ell $, the set $y_\ell e_\ell E_\ell $ contains any set $E_{4j_2-3,4j_2-1,x^{\ell +1}}'$ that it intersects, and if the two sets intersect and are on the same level, so that the topmost principal pairs of the two coincide, then they coincide. It follows that principal arcs of type 3b) of \ref{4.19} do not occur for $\xi _{\ell ,1}$ for any $\ell $, and hence they  do not occur for $\psi _{1,4i-5,\gamma }$.

The types 2 and 3 of quadruples $(\beta _{4i+t}(\gamma ):-3\le t\le 0)$  only occur because they are in the boundaries of sets $C(v\leftrightarrow v',E)$, and hence give satellite arcs. 
\end{proof}

\subsection{Satellites without further effect}\label{5.10.2}
From Theorem \ref{5.10.1} there are two kinds of satellite quadruples for $\gamma =\gamma ^{r,\alpha }$: those as in 2 or 3 of the lemma. If $(\beta _{4i+t}(\gamma ):-3\le t\le 0)$ is as in 2 of the theorem, then some of the disc exchanges in the composition for $\psi _{4i-3,4i-1}^{-1}$ might be in the composition for $\psi _{1,4j-1}^{-1}$ which is not the identity on $D(x_{r,\alpha })$. But for quadruples $(\beta _{4i+t}(\gamma ):-3\le t\le 0)$ as in 3 of the theorem this is impossible. We therefore refer to these as {\em{satellite (quadruples) without further effect.}}

The result of \ref{5.10.1} can be improved. We have the following. In particular, this implies that $x'$ in the statement of \ref{5.10.1} can almost always be taken of the form $x_{r,\alpha '}$ for some $\alpha '\in \{ a,b\} ^{\mathbb N}$. For an $\alpha \in \{ a,b\} ^{\mathbb N}$, we write $n(\ell ,\alpha )$  for the largest  integer such that 
$$D(v_{\ell -1 ,\alpha }t_{\ell -1,\alpha })\subset \psi _{1,4\lfloor n(\ell )-1)/2\rfloor -1}(D(\beta _{2n(\ell )-1},\beta _{2n(\ell )})).$$ 

\begin{subsectiontheorem}\label{5.11}  Fix $\alpha \in \{ a,b\} ^{\mathbb N}$. Then,  for $2p\le n(\ell  )$, 
$$\psi _{1,4p-1,\alpha }\mid \psi _{1,2n(\ell )-1,\alpha }^{-1}(D(v_{\ell -1,\alpha }t_{\ell -1,\alpha }))$$
 can be written in the form 
$$\xi _{q_{p,\ell } ,p\ell }'\circ \cdots \xi _{1,p,\ell}',$$ 
where each $\xi _{j,p,\ell }'$ is a composition of disc exchanges $\xi _{i,j,p,\ell }$ with supports $C(z_{j,\ell }a\leftrightarrow z_{j,\ell }b,E_{i,j,p,\ell }')$ where the following hold.
\begin{itemize}
\item[1.]   $|z_{j,\ell}|$ is increasing with $j$;
\item[2.]  $z_{j,\ell }e$ is a prefix  of $x$ for one of $e=a$ or $b$;
\item[3.] Each $E_{i,j,p,\ell }$ is of the form $E_{4j_1-3,4j_1-1,x_{i,j,p,\ell }}'$ for some $j_1$ depending on $i$, where  $x_{i,j,p,\ell }=x_{\theta _{i,j,p,\ell }}$ for $\theta _{i,j,p,\ell }\in \{ a,b\} ^{\mathbb N}$ which can be taken to be $\alpha $ except when $|z_{j,\ell }a|=|v_k|$ for some $k$, when we can take $\theta _{i,j,p,\ell }(s)=\alpha (s)$ for $s\ne k$.
\item[4.]  
$$\psi _{1,4p-1,\alpha }^{-1}(D(v_{\ell -1,\alpha }t_{\ell -1 ,\alpha }a)\cup D(v_{\ell -1 ,\alpha }t_{\ell -1,\alpha }b))$$
$$=D(v_{\ell -1,\theta }t_{\ell -1,\theta }a)\cup D(v_{\ell -1,\theta }t_{\ell -1,\theta }b)$$
 for some $\theta =\theta _{\ell ,p,\alpha }\in \{ a,b\} ^{\mathbb N}$ such that $\theta _{\ell _1,p_1,\alpha }(i)=\theta _{\ell _ 2,p_2,\alpha }(i)$ for $i\le \ell _1\le \ell _2$ and $n(\ell _1)\le 2{\rm{min}}(p_1,p_2)$. 
 \item 5. A similar result to 4 holds for $D(v_{\ell -1,\alpha }t_{\ell -1,\alpha }eu)$ where $e\in \{ a,b\} $ and $u$ is any prefix of $t_{\ell ,\alpha }$. This time, the image of $D(v_{\ell -1,\alpha }t_{\ell -1,\alpha }eu)$ under $\psi _{1,4s-1,\alpha }^{-1}$ is of the form $D(v_{\ell -1,\theta }t_{\ell -1,\theta }e'u'$ where $e'\in \{ a,b\} $ and $u'$ is a prefix of $t_{\ell ,\theta }$. \end{itemize}
\end{subsectiontheorem}
\begin{proof} [1.] The proof involves rearranging the disc exchanges in the composition 
$$\psi _{1,4p-1,\alpha }=\psi _{1,3,\alpha }\circ \cdots \circ \psi _{4r-3,4r-1,\alpha }.$$
Since each $\psi _{4i-3,4i-1,\alpha }$ is a composition of disc exchanges, we already have a representation of $\psi _{1,4p-1,\alpha }$ as a composition of disc exchanges $\xi _q'\circ \cdots \xi _1'$. For such a  composition we write $\xi _{i,1}'=\xi _i'\circ \cdots \circ \xi _1'$. First, we summarise what we know from \ref{5.10}, given that we are looking at the restriction of $\psi _{1,4r-1,\alpha }$ to $\psi _{1,4r-1}^{-1}(D(v_{\ell ,\alpha }t_{\ell ,\alpha })$. It follows from this that each disc exchange in the composition for $\psi _{4j-3,4j-1}$, for each $j$, has support of the form $C(za\leftrightarrow zb,E_{4j_1-3,4j_1-1,y})$ for some $j_1$, $z$ and $y$, where each of $z$ and $y$ is obtained from a prefix of $x_{r,\alpha }$ (or, equivalently, of $x_r$) by replacing occurrences of $a$ or $b$ before some subwords $v_{1,\alpha }$ or $t_{1,\alpha }$ by $b$ or $a$.  So the same is true for the exchanges in the full decomposition for $\psi _{1,4n(\ell )-1,\alpha }$. So the support of $\xi _j'$ is of the form $C(z_j'a\leftrightarrow z_j'b,E_{4j_1-3,4j_1-1,y_j'})$, where each $z_j'$ and $y_j'$ have the properties just stated for $z$ and $y$ respectively. However, we do not have from this that $|z_j'|$ is non-decreasing  with $j$. In order to get this, we need to be able to commute  a disc exchange $\xi _{j_1,i_1}$ in the composition  for  $\psi _{4i_1-3,4i_1-1,x}$ with $\xi _{j_2,i_2}$ in the composition for $\psi _{4i_2-3,4i_3-1,x}$, if $i_1<i_2$ and the supports are $C_1=C(z_{j _1,i_1}'a\leftrightarrow z_{j_1,i_1}'b,E_{4i_1-3,4i_1-1,x})$ and $C_2=C(z_{j_2,i_2}'a\leftrightarrow z_{j_2,i_2}'b,E_{4i_2-3,4i_2-1,x})$ with $|z_{j_1,i_1}'|<|z_{j_2,i_2}'|$. We only have to do this restricted to $\xi '(D(v_{\ell ,\alpha }t_{\ell ,\alpha }))$, where $\xi '$ is the composition of disc exchanges applied before $\xi _{\ell _2,i_2}$. Commutativity holds if $C_2$ is contained in $z_{j_1,i_1}'eE_{4i_1-3,4i_1-1,x}$ (for $e=a$ or $b$) whenever it intersects it. This is true if 
$$D(z_{j_2,i_2}')\cap z_{j_1,i_1}'eE_{4i_1-3,4i_1-1,x}\ne \emptyset \Rightarrow D(z_{j _2,i_2}')\subset z_{j_1,i_1}'eE_{4i_1-3,4i_1-1,x}.$$
Write $y_s''$ and $y_s'$ for the longest and shortest words respectively with
$$D(y_s'')\subset E_{4i_s-3,4i_s-1,x}\subset D(y_s').$$
Then it suffices to show that if $i_1<i_2$ and $|z_{j_1,i_1}|<|z_{j_2,i_2}'|$ and $z_{j_1,i_1}'ey_1'$ is a prefix of $z_{j _2,i_2}'$ or $z_{j_2,i_2}$ is a prefix of $z_{j_1,i_1}'ey_1''$  then $z_{j_1,i_1}'ey_1''$ is a prefix of $z_{j_2,i_2}'$. To do this we use properties of the set $A$ of words $y$ which are obtained from $x_r$ by replacing some occurrences of $a$ preceding subwords $v_1$ or $t_1$ by $b$. For each $\alpha $, $x_{r,\alpha }$ is such a word, but the set of words $y$ is considerably larger than the set of words $x_{r,\alpha }$. Nevertheless, any such word has  strings of $d$'s  in certain positions, and this means that that the set $A$ still has strong non-recurrence properties. In fact, the following holds. Suppose $z$ and $y$ are both prefixes of words in $A$. Then $zy$ can only be a prefix of a word in $A$ if $z$ has the same length as a word of $x_r$ which ends in $v_i$ for some $i\ge k$, where $k$ is such that $|v_{k-1}|< |y|\le |v_k|$. The reason is simply that, for any $y\in A$, a string of  at least $k$ $d$'s only occurs in the same position as in $x_r$, and any occurrence  in $x_r$ is at the end of a subword $v_i$ for some $i\ge k$. 

So now we apply this with $z=z_{\ell _1,i_s}'e_s$ and $y=y_s'$ for $s=1$, $2$ under the assumption that $|z_{j_1,i_1}|<|z_{j_2,i_2}'|$ and $i_1<i_2$. We let $k_s$ be such that $|v_{k_s-1}|\le |y_s'|<|v_{k_s}|$ Then the assumption that $i_1<i_2$ means that $k_1\le k_2$. In fact if $(\beta _{4i_1+t}-3\le t\le 0)$ is non-satellite then $|y_1'|<|y_2'|$ ,and if $(\beta _{4i_1+t}:-3\le t\le 0)$ is satellite then from the way that satellites are produced we see that $|v_{k-1}|<|y_1'|\le |v_k|$ whenever $|v_{k-1}|<|y_2'|\le |v_k|$, that is $k_1\le k_2$. It follows that if $|z_{j_1,i_1}|<|z_{j_2,i_2}'|$, then  
$$|z_{j_1,i_1}e_1|+|v_{k_1}|\le |z_{j_1,i_1}e_1|+|v_{k_2}|<|z_{j_2,i_2}e_2|$$
and hence 
$$|z_{j_1,i_1}e_1y_1''|<|z_{j_2,i_2}e_2|$$
as required.

The proof of 2, 3, 4 and 5 follow by induction. Once we know $\psi _{1,4n(\ell )-1,x}^{-1}$, then we know the prefixes $z_{j,\ell +1}$ for all $j$ with $|z_{j,\ell +1}a|< |v_{\ell +1}|$. We have $z_{j,\ell }=z_{j,\ell +1}$ for $|z_{j,\ell +1}a|<|v_{\ell +1}|$. Moreover the disc exchanges in the decomposition for $\psi _{1,4r-1,\alpha }$ with supports of the form $C(z_{j,\ell }a\leftrightarrow z_{j,\ell }b,E)$ are the same for any $|z_{j,\ell }|\le |v_\ell |$ if $2n(\ell )\le r$. Now the disc exchanges in the composition for $\psi _{1,4r-1}^{-1}$ which affect the choice of $\theta =\theta _{\ell ,r,\alpha }$ are those corresponding to prefixes $|z_{j,\ell }|\le |v_{\ell }|$. It follows that $\theta _{\ell ,r,\alpha }(i)=\theta _{\ell ,s,\alpha }(i)$ if $i\le \ell $ and $2n(\ell )\le {\rm{min}}(r,s)$, and also $\theta _{\ell _1,r,\alpha }(i)=\theta _{\ell _2,r,\alpha }(i)$ if $i\le \ell _1\le \ell _2$ and $2n(\ell _1)\le r$. The result follows.\end{proof}

From \ref{5.11.1}, and in particular from 4 and 5 of \ref{5.11.1},  we have the following corollary. The proof is immediate. 

\begin{corollary}\label{5.12} Let $\alpha \in \{ a,b\} ^{\mathbb N}$ and let $x=x_{r,\alpha }$.  Let $(\beta _{4i+t}(x):-3\le t\le 0)$ be any quadruple which is not a satellite quadruple without further effect, such that  with $n(\alpha ,\ell )<2i$ for some  $\ell \ge 1$. Then there is $\theta =\theta _{i,\alpha }\in \{ a,b\} ^{\mathbb N}$  and $j=j(i)$ such that $(\beta _{4i+t}(x):-3\le t\le 0)$ is an  $\zeta $ quadruple $(\zeta _{4j+t}(x_{r,\theta _i}):-3\le t\le 0)$, which is an upper quadruple if $(\beta _{4i+t}(x):-3\le t\le 0)$ is non-satellite. Moreover $\theta _{i_1,\alpha }(s)=\theta _{i_2,\alpha }(s)$ for all $s\le \ell -1$ if $n(\alpha ,\ell )\le {\rm{min}}(2i_1,2i_2)$.

Moreover 
$$(\beta _{2n(r+1,\alpha )-1}(\gamma ),\beta _{2n(\gamma )}(\gamma ))=(\zeta _{2n_1(r+1,y)-1}(y),\zeta _{2j}(y)).$$ for $y=x_{r,\theta }$  and  $\theta =\theta _{\lfloor n(\alpha ,r+1)/2\rfloor ,\alpha }$.
\end{corollary}

It follows from Theorem \ref{5.2} that $\theta _{\lfloor n(\alpha ,r+1)/2\rfloor ,\alpha }$ is independent of $\alpha $, or can be assumed so, that is, $\theta _{\lfloor n(\alpha ,r+1)/2\rfloor ,\alpha }(s)$ is independent of $\alpha $ for $0\le s\le r$.  We write $\alpha ^*$ for this common value of $\theta _{\lfloor n(\alpha ,r+1)/2\rfloor ,\alpha }$. Then combining the above with Theorem \ref{5.2}, we have the following corollary.

\begin{corollary}\label{5.13} There is $\alpha ^*\in \{ a,b\} ^{\mathbb N}$ with the following property. For all $\alpha \in \{ a,b\} ^{\mathbb N}$, we have$\theta _{i,\alpha }(s)=\alpha ^*(s)$ for $s\le \ell -1$, if $n(\ell ,\alpha )\le 2i$, and $\theta _{i,\alpha }(s)=\alpha (s)$ for all $s\le r$ if $i=\lfloor n(r+1,\alpha )/2\rfloor $.\end{corollary}

I have not been able to compute $\alpha ^*$. This is mainly because of the difficulty of detemining  accurately the satellites with effect. The only thing which seems certain is that $\alpha ^*(0)=\alpha ^*(1)=a$ and $\alpha ^*(2)=b$. As the notation suggests, $\alpha ^*$ is independent of $r$, that is, the values of $\alpha ^*(s)$ for $s\le r$ are the same for $x_{\alpha ,q}$ for all $q\ge r$, and, of course, for all $\alpha \in \{ a,b\} ^{\mathbb N}$.

\chapter{The main theorems}\label{6}

The theorems from which our main Theorem \ref{1.9} is deduced are the following.

  \begin{theorem}\label{6.1} The set of  capture paths $\gamma \in  {\cal{Z}}_{m}(3/7,+,+,0)$ for which there is no other inequivalent capture path $\zeta $ such that $\sigma _{\zeta }\circ s$ and $\sigma _{\gamma }\circ s$ are Thurston equivalent, is of positive density bounded from $0$. 
 
 More precisely, there are capture paths $\gamma \in {\cal{Z}}_{m}(3/7,+,+,0)$ such that the following hold. There is no other capture path $\zeta \in {\cal{Z}}(3/7,+,+,0)$ (up to equivalence) such that $\sigma _\gamma \circ s$ and $\sigma _\zeta \circ s$ are Thurston equivalent, and there is a decreasing sequence of neighbourhoods $U$ of a point   $P$ in the boundary of the gap of $L_{3/7}$ containing the endpoint $x$ of  $\gamma $ such that 
 $$\lim _{U\to P}\lim _{n\to \infty }\dfrac{\# (V(U)\cap Z_n)}{\#(U\cap Z_n)}=1$$
 where $V(U)\subset U$, and $y\in V(U)\cap Z_\infty  $ if and only if the following hold. The point  $y$ is the endpoint of a capture path $\gamma (y)\in \in {\cal{Z}}(3/7,+,+,0)$ such that there is no other capture path $\zeta \in {\cal{Z}}(3/7,+,0)$ (up to equivalence) such that $\sigma _{\gamma (y)}\circ s$ and $\sigma _{\zeta }\circ s$ are Thurston equivalent. 
 \end{theorem}
 
 \begin{theorem}\label{6.2} For each integer $r\geq 0$, and for all sufficiently large $m$ given $r$, the set of capture paths $\gamma =\gamma _1\in  {\cal{Z}}_{m}(3/7,+,+,0)$ such that $\sigma _{\gamma }\circ s$ is equivalent to $\sigma _{\gamma _i}\circ s$ for exactly $2^{r+1}$ inequivalent capture paths $\gamma _i$ for $1\leq i\leq 2^{r+1}$, and with $\gamma _i\in {\cal{Z}}_m(3/7,+,+,0)$ for $1\leq i\leq 2^{r+1}$, and to no path  in ${\cal{Z}}_m(3/7,+,-,0)$ with $S^1$ crossing at $e^{2\pi it}$ for $t\in [\frac{39}{56},\frac{5}{7})$,  is of positive density bounded from $0$. 
 
 More precisely there is a capture path $\gamma _1\in  {\cal{Z}}_{m}(3/7,+,+,0)$ such that $\sigma _{\gamma }\circ s$ is equivalent to $\sigma _{\gamma _i}\circ s$ for exactly $2r$ inequivalent capture paths $\gamma _i$ for $1\leq i\leq 2r$, with endpoints $x_i$, and $\gamma _i\in {\cal{Z}}(3/7,+,+,0)$ for all $i\leq 2^{r+1}$, and there are neighbourhoods $U_i$ of points  $P_i$ in the boundaries of the gaps of $L_{3/7}$ containing the points $x_i$, and a sequence of neighbourhoods of $P_i$ in $U_i$, converging to $P_i$ such that, for $U_i'$ in this sequence,
 $$\lim _{U_i'\to P_i}\lim _{n\to \infty }\dfrac{\#(V(U_i':U_1,\cdots U_{2N})\cap Z_n)}{\#(U_i'\cap Z_n)}=1$$
 where $U_i'$ is in this sequence of neighbourhoods, and $V(U_i':U_1,\cdots U_{2r})$ is the set of points $y_i\in U_i'\cap Z_\infty $ with the following property. The point $y_i$ is the endpoint of a capture path $\gamma (y_1)\in {\cal{Z}}(a_1,+,+)$ and there are exactly $2^{r+1}$ inequivalent capture paths $\gamma (y_j)$ such that $y_j\in U_j$ for $1\leq j\leq 2^{r+1}$ and $y_i\in U_i'$ and $\sigma _{\gamma (y_i)}\circ s$ are all Thurston equivalent, but there is no other inequivalent capture path $\zeta \in {\cal{Z}}(3/7,+,0)$ such that $\sigma _\zeta \circ s$ and the $\sigma _{\gamma (y_i)}\circ s$ are Thurston equivalent. Moreover there are neighbourhoods  $U_j''(U_i')$ of $y_j$ with $U_i''(U_i')=U_i'$ such that
 $$ V(U_i':U_1,\cdots U_{2r})=V(U_i',U_1''(U_i'),\cdots U_{2r}''(U_i'))$$
 and 
 $$\lim _{U_i'\to y_i}U_j''(U_i)=y_j.$$
 \end{theorem}
 
 \section{Reductions}\label{6.3}
 The more precise statements make it likely that these are really the same theorem. This is in fact the case, given Theorem \ref{5.2}, which gives us capture paths $\gamma _i$ for $1\leq i\leq 2^{r+1}$ to work with, as in the more precise statement of \ref{6.2}. We will make choice of  the path $\gamma $, its endpoint $x$ and the point $P$ of \ref{6.1}, in \ref{6.6}. Similarly, we will make a precise choice of the $\gamma _i$, the endpoints $x_i$ and the  $P_i$ of \ref{6.2}, in \ref{6.7}. We will choose $P$ of the same preperiod as $x$, and similarly for the $P_i$. 
 
 The usual way to construct sets of positive density is to  use a convergent product: for example, a product of the form
 $$\prod _{n=1}^\infty (1-\exp (-(\log n)^\alpha )$$
for any $\alpha >1$. This is the type of convergent product that we will use.

Let $\gamma \in {\cal{Z}}_m(3/7,+,+)$. A sufficient condition for $\sigma _{\gamma }\circ s$ not to be Thurston equivalent to $\sigma _{\gamma '}\circ s$ for any inequivalent capture path $\gamma '$ is that, for all $i$,
\begin{equation}\label{6.3.1}(\beta _{2i-1}(\gamma ),\beta _{2i}(\gamma ))=(\beta _{2i'-1}(\gamma '),\beta _{2i'}(\gamma '))\Rightarrow i=i'\wedge \beta _j(\gamma )=\beta _j(\gamma ')\forall j\leq 2i.\end{equation}
To see that this is sufficient, the pairs $(\beta _{2i-1}(\gamma ),\beta _{2i}(\gamma ))$ are defined so that $\beta =\beta _{2n(\gamma )-1}(\gamma )$ is the path in $R_{m,0}$ such that $\sigma _{\beta }\circ s$ is Thurston equivalent to $\sigma _\gamma \circ s$. Therefore, if $\sigma _{\gamma }\circ s$ and $\sigma _{\gamma '}\circ s$ are Thurston equivalent, we have $n(\beta )=n(\gamma )=n$ and 
$$\beta =\beta _{2n-1}(\gamma )=\beta _{2n-1}(\gamma ')$$
Since $(\beta _{2n-1}(\gamma ),\beta _{2n}(\gamma ))$ is an adjacent pair in $R_{m,0}$ with $\beta _{2n-1}(\gamma )<\beta _{2n}(\gamma )$, and similarly for $\gamma '$, we then have $\beta _{2n(\gamma )}(\gamma )=\beta _{2n(\gamma ')}(\gamma ')$, and hence 
(\ref{6.3.1}) implies that 
$$\beta _i(\gamma )=\beta _i(\gamma '){\rm{\ for\ all\ }}i\leq 2n(\beta ).$$
Now (\ref{6.3.1}) also  implies that
\begin{equation}(\beta _{2i-1}(\gamma ),\beta _{2i}(\gamma ))=(\beta _{2i'-1}(\gamma '),\beta _{2i'}(\gamma '))\Rightarrow \psi _{4j-3,4j-1,\gamma }=\psi _{4j-3,4j-1,\gamma '}\end{equation}
for all $2j\le i$. Hence if the condition holds we also have
$$\psi _{1,2i-1,\gamma }=\psi _{1,2i-1,\gamma '}.$$
In particular if (\ref{6.3.1}) holds for $\gamma $ then for $n=n(\gamma )=n(\gamma ')$ we have 
\begin{equation}\beta _{2n-1}(\gamma )=\beta _{2n-1}(\gamma ')\Rightarrow \psi _{1,2n-1,\gamma }=\psi _{1,2n-1,\gamma '}.\end{equation}
But then, denoting the endpoints of $\beta =\beta _{2n-1}(\gamma )$, and $\gamma $ and $\gamma '$ by $p(\beta )$ and so on (as usual), we have
$$p(\gamma )=\psi _{1,2n-1}^{-1}(p(\beta ))=p(\gamma ')$$ and hence $\gamma =\gamma '$. The method for proving Theorems \ref{6.1} and \ref{6.2} is essentially to prove that Condition (\ref{6.3.1}) holds with positive density for all $i$. This will be done by showing that, if $A_i$ is the set for which (\ref{6.3.1}) holds, then, for a suitable sequence $n_i$ converging to $\infty  $, the density of $A_{n_{i+1}}$ in $A_{n_i}$ approaches $1$ as $i$ tends to $\infty $, sufficiently fast to yield a set of positive density. 



Naturally enough, the positive density condition that we seek will be obtained by proving that a set of conditions which imply the sought condition, hold with positive density. The theorem which we will prove, which will imply both theorems \ref{6.1} and \ref{6.2}, is the following. Here, $D(\zeta _1,\zeta _2)$ is as in \ref{4.2}.

\begin{theorem}\label{6.4}  Let $U\subset \{ z:{\rm{Im}}(z)>0\} $ be an open set   in $Z_\infty $, and let ${\cal{U}}$ be the set of capture paths in ${\cal{Z}}(3/7,+,+,0)$ with endpoints in $U$. Let $n_0\geq 1$. Suppose that
 $\beta _i(\gamma )$ is constant for $\gamma \in {\cal {U}}$ and $i\leq 4n_0+2$, so that $\psi _{1,4i-1,\gamma }$ is also constant for $2i\le n_0$. Write
$$(\beta _{0,1},\beta _{0,2})=(\beta _{4n_0+1}(\gamma ),\beta _{4n_0+2}(\gamma ))$$
for any $\gamma \in {\cal{U}}$, and write 
$$D_0=D(\beta _{0,1},\beta _{0,2}).$$

Then 
$$\limsup _{n\to \infty }d_n=d>0,$$
where:
\begin{itemize}
\item  $d_n=\dfrac{\# (U_1\cap Z_n)}{\# (U\cap Z_n)}$;
\item  ${\cal{U}}_1$ is the set of capture paths in ${\cal{Z}}(3/7,+,+,0)$ with endpoints in $U_1$ and $\gamma \in {\cal{U}}_1$;
\item $U_1$ is the subset of $U$ on which the following hold for suitable constants.
\end{itemize}
 There are  strictly increasing sequences of integers $n_i$ and $i(r)$, depending on $\gamma $, but with $n_0$ as before and, such that the following hold. 
Write
$$(\beta _{i,1},\beta _{i,2})=(\beta _{4n_i+1}(\gamma ),\beta _{4n_i+2}(\gamma )),$$
 $$D_i=D(\beta _{i,1},\beta _{i,2}),$$
 $$E_i=\cup \{  E_{\omega _1,\omega _3}:(\omega _t:1\leq t\le 4){\rm{\ is\ a\ quadruple,\ }}(\omega _1,\omega _2)=(\beta _{i,1}(\gamma ),\beta _{i,2}(\gamma ))\} .$$
 Let $u_j$ be the longest word such that
 $$D(\beta _{2j-1},\beta _{2j})\subset D(u_j).$$
 For suitable $\alpha _1$ and $\alpha _2>1$, we  have
 $$|u_{n_{j+1}+1}|-|u_{2n_j+1}|\le (\log |u_{2n_j+1}|)^{\alpha _1}$$
 $$|u_{2n_{i(r+1)}+1}|\ge \exp( (\log |u_{2n_{i(r)}+1}|)^{\alpha _2})$$

 Also, the following hold: 
\begin{enumerate}[a)]
\item For all $i$, if $C_1$ is a component of 
$$\cup \{ C:C{\rm\ is\ a\ component\ of\ }(\sigma _{\alpha  _{1,4\ell -3}}\circ s)^{-n}(E_{4\ell -3,4\ell -1}){\rm\ with\ }\ell \le n_i,\ n>0,D_i\not\subset C\} ,$$
or if $(\beta _{2j-1},\beta _{2j})$ is an inactive satellite, and 
$C_1$ is a component of 
$$\cup \{ C:C{\rm\ is\ a\ component\ of\ }(\sigma _{\alpha _{1,4\ell -3}}\circ s)^{-n}(E_{4\ell -3,4\ell -1}){\rm\ with\ }2\ell <j,\ n>0,D_i\not\subset C\} $$
then for any principal arc $I$ for an adjacent pair in $R_{m,0}$ which is inside $D_i$ and  over $D_{i_1}$  for some $i_1>i$, we have
$$(C_1\cap I=\emptyset )\vee (C_1\subset D_i\wedge (C_1\subset E_{2n_i-3,2n_i-1}\vee C_1\cap E_{2n_i-3,2n_i-1}=\emptyset )).$$

\item  For any $r$, $\ell $ and $j$ with $2n_{i(r-1)}<2\ell <j\le 2n_{i(r+1)}$, such that $(\beta _{2j-1},\beta _{2j})$ is not an inactive satellite, and any $n>0$ and any component $C$ of $(\sigma _{\beta _{4\ell -3}}\circ s)^{-n}(E_{4\ell -3,4\ell -1})$, or $\sigma _{\beta _{i(r-1),1}}\circ s)^{-n}(D_{i(r-1)})$, we have 
$$C\cap I(\beta _{2j-1},\beta _{2j})=\emptyset .$$
\end{enumerate}

Moreover, given $n_0$ we can find $n_j$ for $j\ge 1$ so that  $d$ is a function $d(d^0)$ of $d^0=\limsup _{m\to \infty }d_{m,0}$, where 
$$d_{m,0}=\sup _{D_0\subset U}\dfrac{\# ((D_0\setminus \cup _{n>0}s^{-n}(E_0))\cap Z_m)}{\# (D_0\cap Z_m)},$$
and
$$\lim _{d^0\to 1}d(d^0)=1.$$
\end{theorem}
\section{Why \ref{6.4} implies \ref{6.1} and \ref{6.2}}\label{6.5}
Of course, it is not immediately obvious that Theorem \ref{6.4} implies \ref{6.1} and \ref{6.2}. The main purpose of the conditions on $U_1$ and ${\cal{U}}_1$ in \ref{6.4} is to bound the effect of the homeomorphisms $\psi _{1,4i-1,\gamma }$  on sets $I(\beta _{2j-1}(\gamma ),\beta _{2j}(\gamma ))$ for $\gamma \in {\cal{U}}_1$. The following lemma is the key to showing that \ref{6.4} implies \ref{6.1} and \ref{6.2}.

We say that a pair $(\beta _{2j-1},\beta _{2j})$ is an {\em{inactive satellite}} for $\gamma $ if 
$$\psi _{1,4\lfloor j/2\rfloor -1}^{-1}(p(\gamma ))\notin D(\beta _{2j-1},\beta _{2j}).$$

\begin{lemma}\label{6.5.1}The conditions on ${\cal{U}}_1$ in \ref{6.4} imply the following for $\gamma \in {\cal{U}}_1$, where $n_i$ is the sequence, which depends on $\gamma $, which is given in 2 of \ref{6.4}. 
\begin{itemize}
\item[1.]  For all $j \geq 2n_i$, 
$$D(\beta _{2j-1}(\gamma ),\beta _{2j }(\gamma ))\subset D_i.$$
Moreover  $\psi _{4\ell -3,4\ell -1,\gamma }$ is the identity on $I(\beta _{2j-1}(\gamma ),\beta _{2j}(\gamma ))$ for all $n_i\leq \ell \le n_{i+1}$ and $2n_k \leq j\leq 2n_{k +1}$, whenever $i(r)\le i<k \le i(r+2)$ for some $r$, and provided that $(\beta _{2j-1},\beta _{2j})$ is not an inactive satellite.
\item[2.] If
\begin{equation}\label{6.5.1.1}\beta _\ell (\zeta )=\beta _\ell (\gamma ){\rm{\ for\ all\ }}\ell \leq 4n_0=4n_0(\gamma ),\end{equation}
and
\begin{equation}\label{6.5.1.2}I(\beta _{2i-1}(\zeta  ),\beta _{2i}(\zeta ))\subset D_k,\end{equation}
for some $i$ and $k$, then 
\begin{equation}\label{6.5.1.3}\beta _j(\zeta )=\beta _j(\gamma ){\rm{\ for\ all\ }}j\leq 2n_{k}.\end{equation}
\item[3.] If $\zeta \in {\cal{U}}$ and  $\beta _{2n(\gamma )}(\gamma )=\beta _{2n(\zeta )}(\zeta )$ then $\zeta =\gamma $.
\end{itemize}
\end{lemma}
\begin{proof}
1. The proof of the first statement is by induction on $j$. The pair $(\beta _{2j-1},\beta _{2j})$ is the pair $(\omega _1,\omega _2)$ with the property that $I(\omega _1,\omega _2)$ is mapped by $\psi _{1,4\lfloor (j-1)/2\rfloor -1}$ to cross $\gamma $ immediately inside the principal arc of $(\psi _{1,4\lfloor (j-2)/2\rfloor -1}(\beta _{2j-3}),\psi _{1,4\lfloor (j-1)/2\rfloor -1}(\beta _{2j-2}))$. So we need to show that $I(\omega _1,\omega _2)\subset D_i$. So we consider the homeomorphism 
$$\psi _{1,4\lfloor (j-1)/2\rfloor -1}=\psi _{1,4n_i-1}\circ \psi _{4n_i+1,4\lfloor (j-1)/2\rfloor -1}.$$
By the inductive hypothesis, 
$$\psi _{4n_i+1,4\lfloor (j_1-1)/2\rfloor -1}={\rm{\ identity\ on\ }}I(\beta _{2j_1-1},\beta _{2j_1}){\rm{\ for\ all\ }}2n_i+1\le j_1<j.$$ So 
$$\psi _{1,4\lfloor (j-1)/2\rfloor -1}=\psi _{1,4n_i-1}{\rm{\ on\ }}I(\beta _{2j_1-1},\beta _{2j_1}){\rm{\ for\ all\ }}2n_i+1\le j_1<j.$$
 By condition b), this will also be true for $j_1=j$, provided $(\beta _{2j-1},\beta _{2j})$ is not an inactive satellite. So now we need to consider the action of $\psi _{1,4n_i-1}^{-1}$ on the region bounded by 
 $$I(\psi _{1,4\lfloor (j-2)/2\rfloor -1}(\beta _{2j-3}),\psi _{1,4\lfloor (j-1)/2\rfloor -1}(\beta _{2j-2})).$$
We use \ref{3.8}. Note that 
$$\psi _{1,4n_i-1}^{-1}=\psi _{4n_i-3,4n_i-1}^{-1}\circ \cdots \circ \psi _{1,3}^{-1}$$ 
  Therefore, the  homeomorphism $\psi _{1,4n_i-1}^{-1}$ is a composition of disc exchanges $\xi _1\circ \cdots \circ \xi _n$, with supports of the form $C(v_i\leftrightarrow v_i',E_i')$. For some of these disc exchanges, one of the exchanged topological discs contains $D_i$. For the others, any component of the union which intersects $D_i$, and which also intersects some $I(\beta _{2j_1-1},\beta _{2j_1})$ for some $j_1<j$, is contained in $D_i$. This means that the first disc exchanges map $\psi _{1,2n_i-1}(D_i)$ to $D_i$ and then after that the only relevant disc exchanges preserve $D_i$. By the condition on $I(\omega _1,\omega _2)$ in a), it follows that $I(\beta _{2j-1},\beta _{2j})\subset D_i$. To see that $I(\beta _{2j-1},\beta _{2j})\subset D_i$ if $(\beta _{2j-1},\beta _{2j})$ is an inactive satellite, we use the other part of condition a). This time the first disc exchanges which are applied are all in the composition for $\psi _{1,2n_i-1}^{-1}$, and the second group are in the composition for $\psi _{1,4\lfloor (j-1)/2\rfloor -1}^{-1}$, and the second part of condition a) ensures that $I(\beta _{2j-1},\beta _{2j})\subset D_i$.

2. Suppose that (\ref{6.5.1.1}) holds and (\ref{6.5.1.2}) holds for some $i$ and $k$. Then (\ref{6.5.1.2}) also holds for $k'$ replacing $k$, for any $k'\leq k$. Therefore we can proceed by induction, both on $j$ and on $k$. We can assume that (\ref{6.5.1.1}) and (\ref{6.5.1.2}) hold, and that (\ref{6.5.1.3}) holds for $k-1$ replacing $k$. Then for all $j$ with $2n_{k-1}\leq j\leq i$ such that $(\beta _{2j-1},\beta _{2j})$ is not an inactive satellite, 
$$\psi _{1,4n_{k-1}-1,\gamma }I(\beta _{2j-1}(\gamma ),\beta _{2j}(\gamma ))=\psi _{1,4\lfloor (j-1)/2\rfloor -1,\gamma }(I(\beta _{2j-1}(\gamma ),\beta _{2j}(\gamma ))$$
and these are mapped over 
$$\psi _{1,4\lfloor (i-1)/2\rfloor -1,\zeta }(I(\beta _{2i-1}(\zeta ),\beta _{2i}(\zeta ))\subset \psi _{1,4n_{k-1}-1,\gamma }(D_k).$$
It follows by induction on $k$ that, for $4n_{k-1}\leq j\leq 4n_k$,
$$\beta _{j}(\zeta )=\beta _j(\gamma )$$
and for $n_{k-1}\leq  j \leq n_k$,
$$\psi _{1,4j -1,\zeta }=\psi _{1,4j-1,\gamma },$$
which give (\ref{6.5.1.3})  for $j$. The inductive step for inactive satellites relies on the fact that an inactive  satellite $I(\beta _{2j-1}(\zeta ),\beta _{2j}(\zeta ))$ is obtained from the support of $\psi _{1,2\lfloor (j-1)/2\rfloor -1,\zeta }$ and similarly for $\gamma $ replacing $\zeta $. So once again, the proof is inductive. 

\noindent 3. By 2, (\ref{6.5.1.3}) holds for the largest $j$ with $2n_j\leq n(\gamma  )$. So 
$$\psi _{1,4n_j-1,\gamma }=\psi _{1,4n_j-1,\zeta }=\psi _{1,4i-1,\gamma }{\rm{\ on\ }}I{\rm{\ for\ all\ }}2n_j\leq 2i\leq n(\beta ).$$
For any $i>2n_1$, and any $i'>n_1$, by induction on $i$,  the set $I(\beta _{2i-1}(\zeta ),\beta _{2i}(\zeta ))$ and  the support of $\psi _{4i'-3,4i'-1,\zeta }$, are both contained in 
$$\cup _{n>0}(\sigma _{\beta ^j}\circ s)^{-n}(E_j)$$
which is disjoint from 
$$I=I(\beta _{2n(\gamma )-1}(\gamma ),\beta _{2n(\gamma )}(\gamma ))=I(\beta _{2n(\zeta )-1}(\zeta ),\beta _{2n(\zeta )}(\zeta )).$$ 
It follows that
$$\psi _{1,4\lfloor n(\gamma )/2\rfloor -1,\gamma }(I)=\psi _{4n_j-3,4n_j-1,\gamma }(I)=\psi _{1,4\lfloor n(\zeta )/2\rfloor -1,\gamma }(I)$$
and hence since one endpoint of $\psi _{1,4\lfloor n(\gamma )/2\rfloor -1}(I)$ is the endpoint of $\gamma $ by the definition of $I(\beta _{2n(\gamma )-1}(\gamma ),\beta _{2n(\gamma )}(\gamma ))$ and the same endpoint is the endpoint of $\zeta $ by the corresponding definition for $\zeta $, we have $\gamma =\zeta $.
\end{proof}

In order to use this lemma to prove Theorems \ref{6.1} and \ref{6.2}, we need to know that the conditions of \ref{6.4}  can be satisfied. This is given by the following theorem.

\begin{subsectiontheorem}\label{6.5.2} 
\begin{itemize}
\item[1.] It is possible to choose a set $$U\subset \{ z:{\rm{Im}}(z)>0\} $$
 and  an even $n_0>0$ such that the following hold for $U$ and for the set ${\cal U}$ of capture paths in ${\cal Z}(3/7,+,+,0)$ with endpoints in $U$.
\begin{enumerate}[(i)]
\item $\beta _j(\gamma )=\beta _j$ is constant for $\gamma \in {\cal{Z}}(3/7,+,+,0)$ with $j\le 4n_0+2$, so that
$$D_0=D(\beta _{4n_0+1}(\gamma ),\beta _{4n_0+2}(\gamma )$$
is independent of $\gamma \in {\cal U}$.
\item  If $\zeta \in {\cal{Z}}(3/7,+,0)$ and  $\sigma _\gamma \circ s$ and $\sigma _\zeta \circ s$ are Thurston equivalent for some $\gamma \in {\cal{U}}$ then $\zeta \in {\cal{Z}}(3/7,+,+,0)$ or $\zeta $ crosses $S^1$ at $e^{2\pi it}$ for $t\in (\frac{19}{28},\frac{39}{56})$. If $\zeta \in {\cal{Z}}(3/7,+,+,0)$ and 
$$D(\beta _{2n(\zeta )-1}(\zeta ),\beta _{2n(\zeta )}(\zeta ))\subset D_0,$$ 
then $\zeta \in {\cal{U}}$, and consequently
$$\beta _i(\zeta ) =\beta _i{\rm{\  for\ }}i\le 4n_0.$$
\end{enumerate}

\item[2.]  Fix any integer $r\ge 1$. Then we can find sets $U_i\subset \{ z:{\rm{Im}}(z)>0\} $, for $1\le i\le 2^{r+1}$, and even integers $n_{0,i}$ such that the following hold, given $d^0>1$, where ${\cal U}_i$ is the set of capture paths in ${\cal Z}(3/7,+,+,0)$ which endpoints in $U_i$.
\begin{enumerate}[(i)]
\item For each $i\le 2^{r+1}$ and each $j\le 2n_{0,i}+2$, the map $\gamma \mapsto  \beta _j(\gamma )=\beta _{j,i}$ is constant on ${\cal U}_i$. 
\item The map $(i,\gamma ) \mapsto (\beta _{2n_{0,i}+1}(\gamma ),\beta _{2n_{0,i}+2}(\gamma ))$ is constant on $\cup _{j=1}^{2^{r+1}}(\{ j\} \times {\cal U}_j$
Consequently
$$(\beta _{0,1},\beta _{0,2})=(\beta _{4n_{0,i}+1}(\gamma ),\beta _{4n_{0,i}+2}(\gamma ))$$
is independent of $i$ and of $\gamma \in {\cal U}_i$, and therefore
$$D_0=\psi _{1,4n_{0,i}-1}^{-1}(U_i)$$
is independent of $i$.

\item If $\zeta \in {\cal{Z}}(3/7,+,0)$   is such that $\sigma _\gamma \circ s$ and $\sigma _\zeta \circ s$ are Thurston equivalent for some $\gamma \in \cup _i{\cal{U}}_i$, then $\zeta \in {\cal{Z}}(3/7,+,+,0)$ or $\zeta $ crosses $S^1$ at $e^{2\pi it}$ for $t\in (\frac{19}{28},\frac{39}{56})$. If  $\zeta \in {\cal{Z}}(3/7,+,+,0)$ and $\gamma \in {\cal U}_i$ for some $i\le 2^{r+1}$, and
$$D(\beta _{2n(\zeta )-1}(\zeta ),\beta _{2n(\zeta )}(\zeta ))\subset D_0,$$ 
then $\zeta \in {\cal U}_j$ for some $j\le 2^{r+1}$, and consequently
$$\beta _i(\zeta )=\beta _{i,j}{\rm{\ for\ }}i\le 4n_{0,i}.$$
\end{enumerate}
\end{itemize}
In both cases we can choose $U$, or the $U_i$, so that 
$$d^0\le \limsup _{m\to \infty }\frac{\# ((D_0\setminus \cup _{n>0}s^{-n}(E_0))\cap Z_m}{\#(D_0\cap Z_m)},$$
where 
$$E_0=\cup \{ E_{\omega _1,\omega _3}:(\omega _1,\omega _2)=(\beta _{0,1},\beta _{0,2}).$$

\end{subsectiontheorem}

If we can prove both this theorem and Theorem \ref{6.4} then we have proved Theorems \ref{6.1} and \ref{6.2}. Of course this theorem has two cases, the first of which is needed for \ref{6.1}, and the second for \ref{6.2}. Each of these two cases has two subcases, as we shall see. It seems likely that for $\gamma \in {\cal{U}}$ or $\gamma \in {\cal{U}}_i$, there is no capture path $\zeta \in {\cal{Z}}(3/7,+,-,0)$ such that $\sigma _\zeta \circ s$ and $\sigma _\gamma \circ s$ are Thurston equivalent. Indeed it seems likely that there is also no such $\zeta \in {\cal{Z}}(3/7,+,p)$ for any $p>0$. But it has not been possible to completely prove this because of incomplete knowledge of $R_{m,p}$ for $p\ge 1$: an issue which has also affected all our earlier analysis.
 
\section{The choice of  $U$ for Theorem \ref{6.1}}\label{6.6}

We make a quite explicit example, with $n_0=1$. We choose $\beta _1$ and $\beta _2$ so that $\beta _1$ is the maximum path in $R_{m,0}$ with first unit-disc-crossing to the right of $D(u)$ ,where $u=(L_3L_2R_3)^2L_3^3L_2R_3L_3^4L_2$ and therefore $\beta _2$ is the minimum path in $R_{m,0}$ with first unit-disc-crossing to the left of $D(u)$. We then have $w_1'(\beta _3)=w_1'(\beta _4)=L_3$. Whatever the choice of $w_2'(\beta _3)$, there is a set $C(u_1\leftrightarrow u_2,E_{\beta _1,\beta _3})$ which intersects $D(\beta _1,\beta _2)$ without being contained in it, where 
$$u_1=(L_3L_2R_3)^2L_3^3L_2R_3L_3^3,\ u_2=(L_3L_2R_3)^2L_3L_2BCL_1R_2R_3L_2R_3.$$
  We choose
$$w_2'(\beta _3)=w_2'(\beta _4)=u_2L_3xL_3$$
 for a suitable nontrivial $x$. Then $E_{\beta _1,\beta _3}$ contains the subset of $D(L_3L_2R_3)$ which starts at the right-hand edge of $D(L_3)$ and ends at the first unit-disc-crossing of $\beta _1$ and $\beta _2$. It follows that the second unit-disc-crossing of $\beta _3$ and $\beta _4$ is in the interior of $u_2E$. Thus $(\beta _t:1\leq t\leq 4)$ is a type AC quadruple, in the language of section \ref{4}. In \ref{4.15}, it was shown that for any capture path $\gamma $ with $\beta _i(\gamma )=\beta _i$ for $i\le 4$, we have 
 $$I(\beta _{2i-1}(\gamma ),\beta _{2i}(\gamma ))\subset  \cup _{n\ge 0}u_2^nD(\beta _1,\beta _2).$$
  In fact we will make conditions that ensure that, for all $i\ge 3$,
  $$I(\beta _{2i-1}(\gamma ),\beta _{2i}(\gamma ))\subset D(\beta _3,\beta _4)$$
  To start the process, we choose
  $$D(\beta _5,\beta _6)\subset E_{1,3}\cap D(\beta _3,\beta _4)\setminus \cap _{n>o}s^{-n}(E_{1,3}\cup D(\beta _1,\beta _2)),$$
  so that 
  $$\psi _{1,3}\mid D(\beta _5,\beta _6)=\xi \mid D(\beta _5,\beta _6),$$
  where $\xi $ is the single disc exchange with support $C(u_1\leftrightarrow u_2,E_{1,3})$.  
 
 Since $n_0=1$, we have $D_0=D(\beta _5,\beta _6)$, and $E_0$ is the union of all quadruples $(\omega _t:1\leq t\leq 4)$ with $(\omega _1,\omega _2)=(\beta _5,\beta _6)$.
 It is clear that by choice of $(\beta _5,\beta _6)$, we can ensure that 
 $$\psi _{1,3}(D(\beta _5,\beta _6))\subset \{ z:{\rm{Im}}(z)>0\} ,$$
  for example by choosing 
  $$D(\beta _5,\beta _6)\subset (u_2u_3X),$$
  where $X=BC$ or $UC$, and $u_3$ has an even number of letters in $\{ L_3,L_2\} $ if $X=BC$, and an odd number if $X=UC$. This is because $u_1$ has an odd number of letters in $\{ L_3,L_2\} $.  Also, we can choose $(\beta _5,\beta _6)$ so that
 $$\limsup _{m\to \infty }\dfrac{\#(\cup _{n>0}s^{-n}(E_0)\cap Z_m)}{\# (D_0\cap Z_m)}$$
is arbitrarily close to $1$. Now we choose
$$U=\psi _{\beta _1,\beta _3}(D(\beta _5,\beta _6)),$$ 
and choose ${\cal{U}}$ to be the set of capture paths in ${\cal{Z}}(3/7,+,+,0)$ with endpoints in $U$. Then $\beta _i=\beta _i(\gamma )$ for all $\gamma \in {\cal{U}}$, and $\psi _{1,3,\gamma }=\psi _{\beta _1,\beta _3}$ for all $\gamma \in {\cal{U}}$. Therefore the conditions at the beginning of \ref{6.4} are satisfied.  For such a $(\beta _5,\beta _6)$, we can take $d^0$ arbitrarily close to 1.  Therefore, the condition at the end of \ref{6.4} is satisfied. 

Now we need to prove Theorem \ref{6.5.2}. So  we need to show that any $\zeta\in {\cal{Z}}(3/7,+,0)$ such that $D(\beta _{2n(\zeta )-1}(\zeta ),\beta _{2n(\zeta )}(\zeta ))\subset D_0$,  we have $\zeta \in {\cal{U}}$. We have
 $${\cal{Z}}(3/7,+,0)={\cal{Z}}(3/7,+,+,0)\cup {\cal{Z}}(3/7,+,-,0).$$
 First we deal with the case of ${\cal{Z}}(3/7,+,+,0)$. We will deal with the case of ${\cal{Z}}(3/7,+,-,0)$ in \ref{6.8}.

\begin{lemma}\label{6.6.1} If $\zeta \in {\cal{Z}}(3/7,+,+,0)$ and $D(\beta _{2n(\zeta )-1}(\zeta ),\beta _{2n(\zeta )}(\zeta ))\subset D_0$, then $\beta _i(\zeta )=\beta _i(\gamma )$ for $i\le 4$ for any $\gamma \in {\cal{U}}$. \end{lemma}
\begin{proof}  By the conditions on $\beta _i(\gamma )$ for $i\le 6$ and $\gamma \in {\cal{U}}$, we see that $(\beta _i(\gamma ):1\le i\le 4)$ is a quadruple of type AC. By  \ref{4.17}, and the choice of $D_0=D(\beta _5(\gamma ),\beta _6(\gamma )$, we have $I(\beta _{2i-1}(\gamma ),\beta _{2i}(\gamma ))\subset D(\beta _3(\gamma ),\beta _4(\gamma )$ for all $i\ge 2$ and all $\gamma \in {\cal{U}}$. Also by the choice of $\gamma $ and by \ref{4.16}, \ref{4.17} and \ref{4.18}, $D(\beta _{2i-1}(\gamma ),\eta _{2i}(\gamma )$ does not intersect any of the sets of the form
$$D(\beta _1(\zeta ),\beta _2(\zeta ))\cap E_{1,3,\zeta }\cup \cup _{n\ge 1}u^nE_{1,3,\zeta }$$
for any $\zeta \in {\cal{Z}}(3/7,+,+,0)$ such that $(\beta _i(\zeta ):1\le 4)$ is different from $(\beta _i(\gamma ):1\le i\le 4)$. So if $D(\beta _{2n(\zeta )-1}(\zeta ),\beta _{2n(\zeta )}(\zeta ))\subset D_0$, we must have $\beta _i(\gamma )=\beta _i(\zeta )$ for $i\le 4$. 
\end{proof}
 
 \section{The choice of $U$ for  Theorem \ref{6.2}}\label{6.7}
 
 We use the examples constructed in Chapter \ref{5}. So fix $r\ge 1$. We consider $\gamma ^{r,\alpha }$ for $\alpha \in \{ a,b\} ^{\mathbb N}$. The definition of $\gamma ^{r,\alpha }$ depends only on $\alpha (k)$ for $0\le k\le r$. For any $\alpha \in \{ a,b\} ^{\mathbb N}$, let  $n(k,\alpha )$ be as in \ref{5.11}, that is, $n(k,\alpha )$ is the largest (automatically even) integer $2i$ such that 
 $$I(\psi _{1,4i-1,\alpha }\cdot (\beta _{4i-1}(\alpha )),\beta _{4i}(\alpha ))\not \subset D(v_{k,\alpha })$$
if $1\le k\le r$, and $n(r+1,\alpha)$ is the largest integer $i$ such that $(\beta _{2i-1}(\alpha ),\beta _{2i}(\alpha ))$ is defined. 
 
  We have seen that the branched coverings $\sigma _{\gamma ^{r,\alpha }}\circ s$ are all Thurston equivalent for all $\alpha \in \{ a,b\} ^{\mathbb N}$, that is,
 $$\beta _{2n(r+1,\alpha  )}(\gamma ^{r,\alpha }))=\beta _{2n(r+1,\alpha ^\infty )}(\gamma ^{r,\alpha ^\infty }){\rm{\ for\ all\ }}\alpha ,$$
 where $\alpha ^\infty $ is defined by $\alpha ^\infty (k)=a$ for all $k\in \mathbb N$. 
 But otherwise the sequences $(\beta _i(\gamma ^{r,\alpha }))$ and $(\beta _i(\gamma ^{r,\alpha ^\infty }))$ are quite different. Define 
 $$\Delta _{i,\alpha }=D(\beta _{2n(i,\alpha )-1}(\gamma ^{r,\alpha }),\beta _{2n(i,\alpha )}(\gamma ^{r,\alpha }))$$
 for $i\le r$ and 
 $$\Delta _{r+1,\alpha }=D(\beta _{2n(r+1,\alpha )-3}(\gamma ^{r,\alpha }),\beta _{2n(r+1,\alpha )=2}(\gamma ^{r,\alpha }))$$
 $$\Delta _{i,\alpha }'=\psi _{1,2n(i,\alpha )-5,\gamma ^{r,\alpha }}(\Delta _{i,\alpha })$$
 for $i\le r$ and similarly for $i=r+1$. Although there is no reason why this should be true in general, we have seen that $\Delta _{i,\alpha }'=D(\omega _1,\omega _2)$, for some adjacent pair $(\omega _1,\omega _2)$ in $R_{m,0}$ (depending on $i$ and $\alpha $).
 
  Let $\zeta _i(x)$ be the sequence defined in \ref{4.4} for some words $x$. In particular, $\zeta _i(x)$ defined for all $i$ if $x=x_{r,\alpha }$ with $\alpha (0)=a$. Let $\alpha ^*$ be as in \ref{5.13}. We define
 $$n_1(i,*)=n_1(i,\alpha ^*),\ \ \ x^*=x_{r,\alpha ^*},\ \ \Delta _{i,*}=\Delta _{i,\alpha ^*},$$
where the first two of these definitions are as in \ref{5.8} and \ref{5.12} respectively. For $1\le k\le r$, $n(k,\alpha )$, as given above, is also the largest even integer $2i$ such that 
$$D(\beta _{4i-1}(\alpha ),\beta _{4i})(\alpha ))\not \subset D(v_{k-1,*}t_{k-1,*}a)\cup D(v_{k-1,*}t_{k-1,*}b).$$ Define $\alpha ^{*,i}$ by 
$$\alpha ^{*,i}(j)=\alpha ^*(j){\rm{\ if\ }}j\neq i,\ \ \alpha ^{*,i}(i)\neq \alpha ^*(i),$$
and define 
$$x^{*,i}=x_{r,\alpha ^{*,i}},\ \ \Delta _{i,i,*}=\Delta _{i,\alpha ^{*,i}}.$$
 Then 
 $$n_1(i,\alpha ^{*,i})=n_1(i,*).$$
We saw in Chapter \ref{5}, in particular in \ref{5.12} and \ref{5.13}, that for any $\alpha \in \{ a,b\} ^{\mathbb N}$, 
\begin{equation}\label{6.7.0.1}\begin{array}{l}(\beta _{2n(i,\alpha )-1}(x_{r,\alpha }),\beta _{2n(i,\alpha )}(x_{r,\alpha }))=(\zeta _{2n_1(i,*)}(x^*),\zeta _{2n_1(i,*)}(x^*))\\ 
{\rm{\ or\ }}{\rm{\ if\ }}(\zeta _{2n_1(*)-1}(x^{*,i}),\zeta _{2n_1(i,*)}(x^{*,i}))),\end{array}\end{equation}
but it is not always clear which of these is true.
 We have
 $$I(\beta _{2n(\gamma ^{r,\alpha  })-1}(\gamma ^{r,\alpha }),\beta _{2n(\gamma ^{r,\alpha  })}(\gamma ^{r,\alpha }))\subset D(v_0t_0a)=D(v_1)$$
 for all $\alpha \in \{ a,b\} ^{\mathbb N} $. Now $x_{r,\alpha }$ consists only of the letters $L_3$, $L_2$ and $R_3$, apart from the last letter $C$.  The number of letters of $x_{r,\alpha }$ in $\{ L_3,L_2\} $ is odd. The same are true for $x^*$. Therefore, writing $x_{r,\alpha }=yC$, and $x^*=y^*C$, we have 
 $$D(yBC)\cup D(y^*BC)\subset \{ z:{\rm{Im}}(z)>0\} .$$
 For any capture path $\gamma \in {\cal{Z}}(3/7,+,+,0)$ with endpoint in $D(yBC)$, it is implicit in the results of Section \ref{5} that 
 $$\beta _i(\gamma )=\beta _i(\gamma ^{r,\alpha })=\beta _i(\alpha ){\rm{\ for\ }}i\le 2n(\gamma ^{r,\alpha } )-2,$$
 and hence, for $n=n(\gamma ^{r,\alpha })-1$, 
 $$D(\beta _{2n-1}(\gamma ),\beta _{2n}(\gamma )=D{r+1,*}.$$
Write $x^*= y^*C$. Our candidate for the set $U$ is a set of the form 
$$D(\omega  _1,\omega _2)=\psi _{1,2n(r+1,\alpha )-3}(D(\eta _1,\eta _2)),$$
 where $D(\eta _1,\eta _2)\subset D(y^*BC)$ and $D(\omega _1,\omega _2)\subset D(yBC)$. We choose such a set $U$ with 
 $$\psi _{1,2n(r+1,\alpha )-1}^{-1}(U)\subset \{ z:{\rm{Im}}(z)>0\} .$$
 It is certainly possible to choose $U$ so that this is true. 
We can also choose $U$ so that  $d^0=d^0(U)$ arbitrarily close to $1$, by choosing  $U$ so that 
$$U\subset D(yBCL_1R_1^N)$$
for an arbitrarily large $N$.  We will then have
$$D_0=D(\beta _{4n_0+1},\beta _{4n_0+1})\subset D_{r+1,*}=D(\beta _{4n_0-1},\beta _{4n_0}).$$
 It remains to prove Theorem \ref{6.5.2} for this choice of $U$. So as usual we let ${\cal{U}}$ be the set of capture paths in ${\cal{Z}}(3/7,+,+,0)$ with endpoints in $U$. We write $4n_0=2n(r+1,\alpha )-2$. Suppose that $\zeta \in {\cal{Z}}(3/7,+,0)$ and
 
 \begin{equation}\label{6.7.1.1}D(\beta _{2n(\zeta )-1}(\zeta ),\beta _{2n(\zeta )}(\zeta ))\subset D_0=\Delta _{r+1,*}.\end{equation}
  We want to show that  there is some $\theta  \in \{ a,b\} ^{\mathbb N}$ such that  
 $$\beta _{i}(\zeta )=\beta _i(\gamma ^{r,\theta }){\rm{\ for\ }}i\le 2n(r+1,\theta )-2.$$
   Once again, we consider the case of $\zeta \in {\cal{Z}}(3/7,+,+,0)$ first, and we will deal with $\zeta \in {\cal{Z}}(3/7,+,-,0)$ in \ref{6.8}. 

\begin{subsectiontheorem}\label{6.7.1} Let ${\cal{U}}$ be as above. Let $\gamma \in {\cal{U}}$ and $\zeta \in {\cal{Z}}(3/7,+,+,0)$ such that (\ref{6.7.1.1}) holds. Then there is $\alpha $ and $2n_0=n(r+1,\theta )$ such that
$$\beta _i(\zeta )=\beta _i(\gamma )=\beta _i(\gamma ^{r,\theta }){\rm{\ for\ }}i\le 4n_0.$$\end{subsectiontheorem}

\begin{proof}
The proof, as usual, is an induction.We shall show the following.
\begin{itemize}
\item[A-1]
 If $r\ge 1$, then $p(\zeta )\in D(v_0t_0)\cup D(w_0t_0)$  where, as usual, $p(\zeta )$ denotes the endpoint of $\zeta $ and $\beta _i(\zeta )=\beta _i(\gamma ^{r,\theta })$ for $1\le i\le 4$, for some $\theta \in \{ a,b\} ^{\mathbb N}$. Moreover 
$$D(\beta _{2j-1}(\zeta ),\beta _{2j}(\zeta ))\subset D(\beta _3(\zeta ),\beta _4(\zeta ))$$
for all $j\ge 2$.
\item[A-2] For each $i$ with $1\le i\le r$ there exists $\theta \in \{ a,b\} ^{\mathbb N}$ such that $p(\zeta )\in D(v_{i,\theta })$ and $\beta _j(\zeta )=\beta _j(\theta)$ for $j\le n(i,\theta )$.
Moreover
$$D(\beta _{2j-1}(\zeta ),\beta _{2j}(\zeta ))\subset D(v_{i-1,*}t_{i-1,*}a)\cup D(v_{i-1,*}t_{i-1,*}b)$$
for all $j> n(i,\theta )$. 
\item[A-3] $p(\zeta )\in D(v_{r,\theta }v_{r-1,\theta }t_{r-1,\theta })$ and $\beta _j(\zeta )=\beta _j(\gamma ^{r,\theta})$ for all $j\le n(\gamma ^{r,\theta })-2$. \end{itemize}

To prove A-1, we use \ref{4.15}, and \ref{4.16} to \ref{4.18}, which give information about the sequence $\beta _i(\zeta )$. By (\ref{6.7.1.1}), we need
$$\Delta _{r+1,*}\subset  (D(\beta _1,\beta _2)\cap E_{1,3})\cup \cup _{n\ge 1}u^nE_{1,3}.$$ 
where $u\leftrightarrow u'$ is a basic exchange for $\beta _1$ of the type in (\ref{3.11.4}). Since $\Delta _{r+1,*}\subset D(v_0t_0)$ if $r\ge 1$ and $\Delta _{r+1,*}\subset D(v_0L_3L_2R_3L_3^3)$ if $r=0$, this gives just two choices for  $(\beta _1,\beta _2)$: namely, $\beta _1$ is maximal with $w_1'(\beta _1)$ to the right of $D(v_0L_3L_2R_3L_3^3L_2)$ or $D(w_0L_3L_2R_3L_3^3L_2)$, that is, 
$$(\beta _1(\zeta ),\beta _2(\zeta ))=(\beta _1(\theta),\beta _2(\theta ))$$
for some $\theta \in \{ a,b\} ^{\mathbb N}$. Of course, as a result of this, $(\beta _1(\zeta ),\beta _2(\zeta ))$ only depends on $\theta (0)$, but then we can also choose $\theta $ so that  
$$(\beta _3(\zeta ),\beta _4(\zeta ))=(\beta _3(\gamma ^{r,\theta }),\beta _4(\gamma ^{r,\theta })).$$

By \ref{4.16} and \ref{4.18}, and the choice of $\beta _i(\zeta )$ for $i\le 4$, we then have
$$D(\beta _{2j-1}(\zeta ),\beta _{2j}(\zeta ))\subset D(L_3L_2L_3L_3)\setminus D((L_3L_2R_3)^2L_3)$$
for all $j\ge 1$. We now claim that we have 
$$D(\beta _{2j-1}(\zeta ),\beta _{2j}(\zeta ))\subset D(\beta _3(\zeta ),\beta _4(\zeta ))\subset D(v_0t_0)$$
for all $j\ge 2$. If this is not true for a least $j$, then $j=2j_0+1$ is odd by \ref{4.19}, and   there is an exchange $u_1\leftrightarrow u_1'$ for $\beta _3$ (and for $\beta _{j_1}$ for $3\le j_1<4j_0$) such that 
$$u_1(E_{1,3})\cup D(\beta _3,\beta _4))\subset D(\beta _3,\beta _4)$$
Then $u_1$ does not exist if $r=0$ and  $v_0t_0a$ or $v_0t_0b$ must be a proper prefix of $u_1$ if $r\ge 1$. But then the same  $v_0t_0a$ or $v_0t_0b$ is also a prefix of $u_1'$and from considering sets $C(u_2\leftrightarrow u_2',E')$ which can intersect $C(u_1\leftrightarrow u_1',E)$, for $E$ and $E'\subset D(L_3L_2L_3)L_3)\setminus D(L_3L_2R_3)^2L_3)$, we see that the component of the support of $\psi _{5,4j_0-1,\zeta }$ which intersects $D(\beta _{4j_0+1}(\zeta ),\beta _{4j_0+2}(\zeta ))$ is contained in $D(u_1)\cup D(u_1')$. The same is true for any set $C(u_2\leftrightarrow u_2',E_{1,3})$ which intersects $D(\beta _{4j_0+1}(\zeta ),\beta _{4j_0+2}(\zeta ))$, apart from $C(v_0\leftrightarrow w_0,E_{1,3})$ which contains it. It follows that 
$$D(\beta _{4j_0+1}(\zeta ),\beta _{4j_0+2}(\zeta ))\subset D(u_1)\cup D(u_1')\subset D(\beta _3,\beta _4).$$
Hence we have
$$D(\beta _{2j-1}(\zeta ),\beta _{2j}(\zeta ))\subset D(\beta _3(\zeta ),\beta _4(\zeta )){\rm{\ for\ all\ }}j\ge 2.$$
This completes the proof of [A-1]. So now we consider [A-2].
By a similar argument, if $C(u_1\leftrightarrow u_1',E)$ is the  first exchange to intersect $D(\beta _{2j-1}(\zeta ),\beta _{2j}(\zeta ))$ for a minimal $j$, where $v_0t_0a$ or $v_0t_0b$ is a proper prefix of $u_1$ then we have
$$D(\beta _{2j-1}(\zeta ),\beta _{2j}(\zeta ))\subset D(u_1)\cup D(u_1').$$
Then by induction, 
\begin{equation}\label{6.7.1.2}D(\beta _{2\ell -1}(\zeta ),\beta _{2\ell }(\zeta ))\subset D(u_1)\cup D(u_1'){\rm{\ for \ all\ }}\ell \ge j.\end{equation} 
By (\ref{6.7.1.1}) this can only be true for $u_1$ a prefix of $x_{r,*}$.  Also we are assuming that $v_0t_0a=v_1$ or $v_0t_0b$ is a proper prefix of any such $u_1$. The first such prefix of $x_{r,*}$ is $u_1=v_1t_1a$ or $v_1t_1b$ (since $v_{1,*}=v_1$ and $t_{1,*}=t_1$), that is, $u_1\leftrightarrow u_1'$ is $v_{2,*}\leftrightarrow v_1t_1b$ (interchanging $u_1$ and $u_1'$ if necessary). It follows that for all sufficiently large $\ell $,
\begin{equation}\label{6.7.1.3}D(\beta _{2\ell -1}(\zeta ),\beta _{2\ell }(\zeta ))\subset D(v_1t_1a)\cup D(v_1t_b).\end{equation}
It follows by induction on $j\ge 4$ that that there is $\theta $ such that 
$$\beta _j(\zeta )=\beta _j(\theta ){\rm{\  for\ }}j\le n(1,\theta ),$$
and (\ref{6.7.1.3}) holds for all $\ell \ge n(1,\theta )$. This completes the proof of [A-2] for $i=1$. The proof of [A-2] for a general $i$ is similar. We assume inductively that [A-2] is proved for $i-1$. Then, as before, we can prove that 
$$D(\beta _{2j-1}(\zeta ),\beta _{2j}(\zeta ))\subset D(v_{i-2,*}t_{i-2,*}a)\cup D(v_{i-2,*}t_{i-2,*}b)$$
for all $j\ge n(i-1,\theta )$. For if this is not true then as before we can find $2j_0>n(i,\theta )$ an exchange $u_1\leftrightarrow u_1'$ for $\beta _{j}$ for all $j$ with $2n(i-1,\theta )<j\le 4j_0$ such that $v_{i-2,*}t_{i-2,*}e$ is a proper prefix of $u_1$ for one of $e=a$ or $e=b$, and such that
$$D(\beta _{2j-1}(\zeta ),\beta _{2j}(\zeta ))\subset D(u_1)\cup D(u_1'){\rm{\ for\ all\ }}j>2j_0.$$
This is only possible if $u_1\leftrightarrow u_1'$ is one of the exchanges for which this happens for $(\beta _{2j-1}(\theta),\beta _{2j}(\theta))$, which means end  in $v_{k-1,*}t_{k-1,*}e$ for some $k\le i-1$, or in $t_{i-1,*}e$, for $e\in \{ a,b\} $. So then by induction on $j$ we can define $\theta (i)\in \{ a,b\} $ so that
$$\beta _j(\zeta )=\beta _j(\theta ){\rm{\ for\ all\ }}j\le 2n(i,\theta ),$$
and the rest of [A-2] for $i$ is proved as before. 

The proof of [A-3] is similar.

\end{proof}

\section{Paths in ${\cal{Z}}(3/7+,-,0)$}\label{6.8}

In order to complete the proof of \ref{6.5.1}) and \ref{6.5.2}, it suffices to show that for any capture paths $\gamma \in {\cal{U}}$ and $\zeta \in {\cal{Z}}(3/7,+,-,0)$, the maps $\sigma _\gamma \circ s$ and $\sigma _{\zeta }\circ s$ are not Thurston equivalent. In the process, we shall obtain some  restrictions on the set of  paths $\beta \in R_{m,0}$ such that $\sigma _\beta \circ s$ is equivalent to $\sigma _{\zeta }\circ s$, for $\zeta \in {\cal{Z}}(3/7,0,+,-,0)$. These restrictions appear to be of interest in their own right. More specifically, we shall prove the following. 
\begin{subsectiontheorem} Let $\zeta \in  {\cal{Z}}(3/7,+,-,0)$ which crosses $S^1$ at $e^{2\pi it}$ for $t\in [\frac{39}{56},\frac{5}{7})$ and $\beta \in R_{m,0}$ such that $\sigma _\zeta \circ s$ and $\sigma _\beta \in s$ are Thurston equivalent. Then
$$w_1'(\beta )=L_3{\rm{\ or\ }}L_3L_2R_3L_2u{\rm{\ or\ }}L_3L_2UCu{\rm{\ or\ }}L_3L_2R_3L_3^2u$$
for some $u$, and if  $w_1'(\beta )=L_3$ then
$$w_2'(\beta )=L_3^4u{\rm{\ or\ }}BCu{\rm{\ or\ }}L_3(L_2R_3)^2u$$
for some $u$.\end{subsectiontheorem}
For the sets ${\cal{U}}$ that we have chosen, for theorem \ref{6.1} (chosen in \ref{6.6}) we have 
$$w_1'(\beta _{2n(\gamma )}(\gamma ))=(L_3L_2R_3)^2L_3^3$$
 for all $\gamma \in {\cal{U}}$ and for the set chosen for Theorem \ref{6.2} (chosen in \ref{6.7}) we have 
 $$w_1'(\beta _{2n(\gamma )}(\gamma ))=L_3,\ \ w_2'(\beta _{2n(\gamma )}(\gamma ))=v_0t_0a=L_3L_2R_3L_3(L_2R_3)^2u$$
 for some $u$. It follows that $\sigma _{\gamma }\circ s$ is not equivalent to $\sigma _\zeta \circ s$ for any $\zeta \in {\cal{Z}}(3/7+,-,0)$ which crosses $S^1$at $e^{2\pi it}$ for $t\in [\frac{39}{56},\frac{5}{7})$ .

\begin{proof}

We shall choose a suitable $\psi $ and $\alpha $ for which 
$$(s,Y_m)\simeq _\psi (\sigma _\alpha \circ s,Y_m),$$
where $\simeq $ denotes Thurston equivalence. We shall then consider
$$\zeta ^1=\alpha *\psi (\zeta ).$$
 The basic reason for considering such a $\zeta ^1$ was explained in \ref{2.4}: $\sigma _\zeta \circ s$ and $\sigma _{\zeta ^1}\circ s$ are Thurston equivalent, and this is realised by composing the above equivalence on the left by $\sigma _\zeta $. The aim is to choose $\psi $ and $\alpha $ so that the first two unit-disc crossings of $\zeta ^1$ coincide with those of a path in $R_{m,0}$. If that is so, then we shall show that the path $\beta \in R_{m,0}$ such that $\sigma _{\beta }\circ s$ is equivalent to $\sigma _\zeta \circ s$ has the same first unit disc  crossing as $\zeta ^1$, and the same second unit -disc crossing if $w_1'(\beta )=w_1'(\zeta ^1)=L_3$.  
 
Let $\zeta $ cross the unit circle at $e^{2\pi it}$ for $t\in (\frac{19}{28},\frac{5}{7})$.  Then $t\in (\frac{39}{56},\frac{5}{7})$ or $t\in (\frac{19}{28},\frac{39}{56})$. We can discount $\zeta $ with $t=\frac{39}{56}$ because this point is in the boundary of the gap containing the point of preperiod $2$ labelled by $L_3L_2C$, and it is then clear that $\sigma _\zeta \circ s$ is not equivalent to $\sigma _\gamma \circ s$ for  any $\gamma \in {\cal{U}}$, because the endpoints of paths in ${\cal{U}}$ are in $U$, and all points in $U$ have preperiod higher than $2$, in fact, at least $11$. 

First suppose that  $t\in (\frac{39}{56},\frac{5}{7})$. 
Let $\psi _{m,2/7}$ and $\alpha _{m,2/7}$ be as in \ref{3.15} (and as in 3.3 and 7.2 of  \cite{R5}), that is, these homeomorphisms and closed loops arise from considering the equivalence between the type II captures $\sigma _{\zeta _{5/7}}^{-1}\circ \sigma _{\beta _{5/7}}\circ s$ and $\sigma _{\zeta _{2/7}}^{-1}\circ \sigma _{\beta _{2/7}}\circ s$. Our initial choices of $\psi $ and $\alpha $ are:  
$$\psi = \psi _{m,2/7}$$ 
and 
$$\alpha =\alpha _{m,2/7}.$$
This is just a first choice, which will need to be modified later.

Now let $t\in (\frac{19}{28},\frac{39}{56})$. We use the notation of \ref{3.16} (and of  3.3 of \cite{R6} --  see the last paragraph of that section). Let $q=\frac{17}{56}$. This time our initial choices of $\psi $ and $\alpha $ are:
$$\psi =\psi _{m,q},$$
$$\alpha =\alpha _{m,q}.$$
Once again, this is just a first choice, which will need to be modified later.

First, we consider the case $t\in (\frac{39}{56},\frac{5}{7})$ and $\psi =\psi _{m,2/7}$. Referring to the definitions in \ref{3.15}, we see that
$\psi _{2,2/7}$ is the identity on $\beta _{5/7}$, and hence $\alpha _{2/2/7}$ is an arbitrarily small perturbation of $\beta _{2/7}*\overline{\beta _{5/7}}$. For any path $\zeta $ based at $\infty $ we also have 
$$\alpha _{m,2/7}*\psi _{m,2/7}(\zeta )=\beta _{2/7}*\psi _{m,2/7}(\overline{\beta _{5/7}}*\zeta ).$$
For $\zeta \in {\cal{Z}}(3/7,+,0)$, the homeomorphism $\psi _{2,2/7}$ is the identity on $\overline{\beta _{5/7}}*\zeta $ and hence for all $m\ge 3$,
$$\alpha _{m,2/7}*\psi _{m,2/7}(\zeta )=\beta _{2/7}*\xi _{m-1,2/7}\circ \psi _{m-1,2/7}(\overline{\beta _{5/7}}*\zeta ).$$
Referring again to \ref{3.15}, we see that the support of $\xi _{2,2/7}$ intersects $\psi _{2,2/7}(\overline{\beta _{5/7}}*\zeta )$ with intersection contained in $D(BC)\cap C_{2,2/7}$. Since $\xi _{2,2/7}$ is an anticlockwise twist in the annulus $C_{2,2/7}$, this implies that  $\psi _{3,2/7}(\overline{\beta _{5/7}}*\zeta )$ intersects the unit circle in the boundary of the central gap  at $e^{2\pi i(11/14)}$ rather than $e^{2\pi i(5/7)}$, and the intersection with the unit disc nearer to the endpoint of $\zeta $ is unchanged. The support of $\xi _{3,2/7}$ does not intersect $\psi _{3,2/7}(\overline{\beta _{5/7}}*\zeta )$. Therefore
$$\psi _{3,2/7}(\overline{\beta _{5/7}}*\zeta )=\psi _{4,2/7}(\overline{\beta _{5/7}}*\zeta )$$
Referring to the description of $\xi _{4,2/7}$ restricted to $C_{4,2/7}$ in \ref{3.15}, we see that there are no intersections between $C_{4,/7}$ and $\psi _{4,2/7}(\overline{\beta _{5/7}}*\zeta )$ because of the restriction that $t>\frac{39}{56}$.  However, there is an intersection with $A_{4,2/7}$, if  and only if $\zeta $ crosses the unit circle between $e^{2\pi i(79/112)}$ and $e^{2\pi i(5/7)}$, that is, if and only if the endpoint of $\zeta $ is in $D(L_3L_2R_3L_3)$. For the moment we assume that the  endpoint of $\zeta $ is in $D(L_3L_2R_3L_3)$. Then, because $\xi _{4,2/7}(D(L_3L_2R_3L_3))=D(L_3^4)$, with the movement being along an arc with endpoints in the lower unit circle, we see that   $\psi _{5,2/7}(\overline{\beta _{5/7}}*\zeta )$ intersects the unit disc in $D(L_3^4)$ rather than in $D(L_3L_2R_3L_3)$, and crosses the unit circle between $e^{2\pi i(37/56)}$ and $e^{2\pi i(75/112)}$ instead of between $e^{2\pi i(79/112)}$ and $e^{2\pi i(5/7)}$. There are actually two components of intersection with the unit disc in $D(L_3^4)$, because the action of $\xi _{4,2/7}$ pushes $\psi _{4,2/7}(\overline{\beta _{5/7}}*\zeta )$ up through $D(L_3^4L_2)$, so that the last component of intersection of $\psi _{5,2/7}(\overline{\beta _{5/7}}*\zeta )$ enters $D(L_3^4L_2)$ from the top of the unit circle. There is  no intersection between $\psi _{5,2/7}(\overline{\beta _{5/7}}*\zeta )$ and $C_{5,2/7}$. In fact there is never any intersection between $A_{i,2/7}$ and $C_{i+1,2/7}$, because there is no intersection between $A_{1,2/7}$ and $C_{2,2/7}$. There is, however  an intersection between $\psi _{5,2/7}(\overline{\beta _{5/7}}*\zeta )$  and $A_{5,2/7}$ if and only if the endpoint of $\zeta $ is in $D(L_3L_2R_3L_3^2)$, that is, if and only if the endpoint of $\psi _{5,2/7}(\zeta )$ is in $D(L_3^5)$. In this case the action of $\xi _{5,2/7}$, with $\xi _{5,2/7}(D(L_3^5))=D(BCL_1R_2R_3L_3)$,  gives that the endpoint of $\psi _{6,2/7}(\zeta )$ is in $D(BCL_1R_2R_3L_3)$. There is then no intersection between $\psi _{6,2/7}(\overline{\beta _{5/7}}*\zeta )$  and $C_{6,2/7}$, and there is an intersection between $\psi _{6,2/7}(\overline{\beta _{5/7}}*\zeta )$  and $A_{6,2/7}$ if and only if the endpoint of $\zeta $ is in $D(L_3L_2R_3L_3^3)$, that is, if and only if the endpoint of $\psi _{6,2/7}(\zeta )$ is in $D(BCL_1R_2R_3L_3^2)$. Using the move $\xi _{6,2/7}(D(BCL_1R_2R_3L_3^2))=D(L_3(L_2R_3)^2L_3)$, the endpoint of $\psi _{7,2/7}(\zeta )$ is in $D(L_3(L_2R_3)^2L_3)$. Continuing in this way, the endpoint of $\psi _{m,2/7}(\zeta )$ is in $D(u)$, where
\begin{equation}\label{6.8.1.1}\begin{array}{l}u=u_1L_2=L_3^4L_2(R_3L_3L_2)^r,{\rm{\ or\ }}BCL_1R_2R_3L_3L_2(R_3L_3L_2)^r,\\
{\rm{\ or\ }}L_3(L_2R_3)^2L_3L_2(R_3L_3L_2)^r,\end{array}\end{equation}
if $m=5+3r$ or $6+3r$ or $7+3r$ respectively, for any $r\ge 0$.  The action of $\xi _{i,2/7}$ restricted to components of $C_{i,2/7}$ is relevant for $i\ge 8$. Note that, for the first two choices of $u$, for any extension $uv$ of $u$,  we have $w_1'(uv)=L_3$, and $u_1$ is a prefix of $w_2'(uv)$. This is not true for the third choice of $u$.

The path $\zeta ^1=\alpha _{m,2/7}*\psi _{m,2/7}(\zeta )$ obtained in this way has first unit-disc crossing along $D((L_3L_2R_3)^\infty )$, but then departs from all paths $\beta \in R_{m,0}$ with $w_1'(\beta )=L_3$, because of the action of $\xi _{2,2/7}$ restricted to $C_{2,2/7}$ on $\overline{\beta _{5/7}}*\zeta $. In order to determine the next iteration $\zeta ^2$ towards a path in $R_{m,0}$, we consider the word $w=w(\zeta ^1)=w(\psi _{m,2/7}(\zeta ))$ representing the endpoint of $\psi _{m,2/7}(\zeta )$. If $L_3^4$ or $L_3^2BC$ or $BC$ is a prefix of $w$, when we define
$$\zeta ^2=\alpha ^1*\psi ^1(\zeta ^1),$$
where 
$$(s,Y_m)\simeq _{\psi ^1}(\sigma _{\alpha ^1}\circ s,Y_m),$$
and  $\beta ^1$ is the unique path in $R_{m,0}$ which shares an endpoint with $\zeta ^1$ and with $w_1'(\beta ^1)=L_3$, and $\alpha ^1$ is an arbitrarily small perturbation of $\beta ^1*\overline{\zeta ^1}$ such that the disc (or union of discs) bounded by $\alpha ^1$ does not contain the common endpoint of $\beta ^1$ and $\zeta ^1$. 
Therefore, 
$$\zeta ^2=\beta ^1*\overline{\zeta ^1}*\psi ^1(\zeta ^1).$$
 Since $u_1$ is a prefix of $w_2'(w(\zeta ^1))$, all of $\beta ^1$ beyond the first unit disc crossing is in $D(u)$. Therefore all of $\alpha ^1$, apart from the first and last unit-disc crossing along $D(L_3L_2R_3)^\infty )$, is in
 $$D(u^4)\cup D(u^5)\cup D(u^6),$$
 where the $u^i$ are the different possibilities for $u$ in (\ref{6.8.1.1}) with $4\le m\le 6$, that is, 
$$u^4=L_3^4L_2,\ \ u^5=BCL_1R_2R_3L_3L_2,\ \ u^3=L_3(L_2R_3)^2L_3L_2.$$
From this, we can estimate the size of the components of the support of $\psi ^1$ which intersect $\zeta ^1$.  The difference between $\zeta ^1$ and $\beta ^1$ is contained in a small part of the supports of the homeomorphisms $\xi _{m,2/7}$ for $m\ge 4$, in fact, for each one, a single component of intersection with the unit disc  together with a single adjacent component of intersection with the complement of the unit disc. The path $\alpha ^1$ has a first and last crossing of the unit disc along the leaf $D((L_3L_2R_3)^\infty )$ between $e^{\pm 2\pi i(2/7)}$. It follows that the only component of ${\rm{supp}}(\psi ^1)$ which is not contained in $D(u^m)$ (for $4\le m\le 6$) has one intersection with $D(u^m(R_3L_3L_2)^\infty )$. and then the  later intersections are with $D(v^mv^nL_3)$ for $4\le n\le 6$, where $v^i$ is obtained from $u^i$ by deleting the last two letters $L_3L_2$. So this component of the support is disjoint from $\zeta ^1$.  It follows that the path $\zeta ^2$ is in $D(u)$ from the second unit-disc crossing onwards. 

The case when $u=L_3(L_2R_3)^2L_3L_2(R_3L_3L_2)^r$ is a bit more involved. In this case, because of the definition of $w_1'$ (in \ref{3.6}), it is possible to have $w_1'(w)\ne L_3$. In order to see what happens, we need to look at the longest prefix $u_1u_2$ of $w$ such that $u_2$ has all letters in $\{ L_3,L_2,R_3\} $ and is of the form 
$$u_2=(L_2R_3)^{r_1}L_2{\rm{\ for\ }}r_1\ge 0,$$ or 
$$u_2=(L_2R_3)^{r_1}L_3^{r_2}{\rm{\  for \ }}r_1\ne 1{\rm{\ and\ }}r_2\ge 1.$$
 In the first case, $u_2$ is followed in $w$ by $BC$ or $UC$, and in the second case it is followed in $w$ by $L_2$. In the first case we do have $w_1'(w)=L_3$, we define $\alpha ^1$, $\psi ^1$ and $\zeta ^2$ as above, and it does follow, much as in the previous cases, that $\zeta ^2$ is in $D(L_3(L_2R_3)^2L_3)\cup D(L_3L_2R_3L_2UC)$ from the second unit-disc crossing onwards. In the second case, we can also define $\alpha ^1$, $\psi ^1$ and $\zeta ^2$ as above, in the cases when $w_1'(u_1u_2L_2)=L_3$, and we have the same bounds on $\zeta _2$ from the second unit-disc crossing onwards.  It remains to consider the cases when $w_1'(u_1u_2L_2)\ne L_3$. This happens precisely when $r_1+r_2$ is odd and $\le 3$. We have:
 $$w_1'(u_1u_2L_2)=\begin{cases}L_3(L_2R_3)^2(L_3L_2R_3)^{r-1}L_3{\rm{\ if\ }}r_1=0{\rm{\ and\ }}r_2=1,\\ u_1L_2R_3L_3^{r_2-2}{\rm{\ if\ }}r_1=0{\rm{\ and\ }}r_2>1,\\ u_1(L_2R_3)^{r_1-2}{\rm{\ if\ }}r_1\ge 2{\rm{\ and\ }}r_2=1,\\ u_1(L_2R_3)^{r_1-1}{\rm{\ if\ }}r_1\ge 2{\rm{\ and\ }}r_2=3,\\ u_1(L_2R_3)^{r_1}L_3^{r_2-3}{\rm{\ if\ }}r_1\ge 2{\rm{\ and\ }}r_2\ge 4.\end{cases}$$
 In all these cases, we have 
 $$w_1'(w)=w_1'(u_1u_2L_2)u_3$$
 for a word $u_3$ which is often empty. In all these cases we consider the leaf $\ell =D(w_1'(u_1u_2L_2)(L_2R_3L_3)^\infty )$ of $L_{3/7}$. Exactly one endpoint $e^{2\pi ip}$ is the first  intersection point on the unit circle of at least one path in $R_{m,0}$. We then define $\alpha ^1$ to be the anticlockwise closed loop based at $\infty =v_2$ which goes along  $\ell $ and then back through the disc along the  leaf from $e^{2\pi i(5/7)}$ to $e^{2\pi i(2/7)}$. We then define $\psi ^1$ and $\zeta ^2$ as above. It is possible that $w_1'(\zeta ^2)\ne w_1'(u_1u_2L_2)$, and if the first unit disc crossing of $\zeta ^1$ is in the support of $\psi ^1$, in the support of a disc exchange in the composition for $\psi ^1$ which has connecting arc along the top of the unit circle. But if that happens then both components of intersections with the unit disc are in $D(w_1'(u_1u_2L_2))$, and there is no movement out of this set by other disc exchanges in the composition. So the first unit disc-crossing of $\zeta ^2$ is in $D(w_1'(u_1u_2L_2))$,  and the process can then be repeated.

Now we need to consider what happens for $t\in (\frac{39}{56},\frac{79}{112})$. The case of $t=\frac{79}{112}$ is not a problem because if $\zeta $ crosses the circle at $e^{2\pi i(79/112)}$ then the end of $\zeta $ is in the gap in $D(L_3L_2R_3L_2)$ which also intersects the unit circle at $e^{2\pi i(33/112)}$. There is a capture path $\beta \in R_{m,0}$  with $w(\beta )=L_3L_2R_3L_2C$. The intersection if $\beta $ with the unit circle is of course with the upper unit circle. The intersection point is not $e^{2\pi i(33/112)}$ (as was noted in passing in Chapter 7 of \cite {R5}) but at the leftmost point $e^{2\pi i(67/224)}$. For this $\beta $ we have $w_1'(\beta )=L_3(L_2R_3)^3L_3$, which satisfies the conditions of the theorem. 

Otherwise for $t\in (\frac{39}{56},\frac{79}{112})$, in some contrast to what was done before, we define $w=w(\zeta )$ and $x=w_1'(w)$. Then $\ell =x(L_2R_3L_3)^\infty $ is a vertical leaf. Let $e^{2\pi ip}$ be the  point of $\ell $ which is the starting point of at least one path in $R_{m,0}$.  Our new choice of $\psi $ and $\alpha $ --- with $\zeta ^1=\alpha *\psi (\zeta )$ --- are
$$\psi =\psi _{m,p},$$
$$\alpha =\alpha _{m,p}.$$

The structure of $\psi _{m,p}$ is similar to that of $\psi _{m,2/7}$, as noted in \ref{3.16} (and earlier, in 3.3 of \cite{R5}). In particular, $\psi _{i+1,p}=\xi _{i,p}\circ \psi _{i,p}$ for all $i\ge 1$, where the support of $\xi _{i,p}$ is $A_{i,p}\cup C_{i,p}$ and 
$$A_{i,p}=(\sigma _{\beta _p}\circ s)^{1-i}(A_{1,p}),\ \ C_{i,p}=(\sigma _{\beta _p}\circ s)^{1-i}(C_{1,p}),$$
using the notation of \ref{3.16}. The annuli of $A_{i,p}$ are thinner than the corresponding annuli of $A_{i,2/7}$ but the annuli of $C_{i,p}$ are thicker than the corresponding annuli of $C_{i,2/7}$. It is not accurate to say that the annuli ``correspond,  because it is not true (for example) that each annulus component of $A_{i,2/7}$ contains a unique annulus component of $A_{i,p}$. But this is true  for components of $A_{i,2/7}\cap \{ z:|z|\le 1\} $ and components of $A_{i,p}\cap \{ z:|z|\le 1\} $ and even on some (or most) unions of components of intersection with, alternately, $\{z:|z|\le 1\} $ and $\{ z:|z|\ge 1\} $. 

In many cases, $w(\zeta ^1)=w(\psi _{m,p}(\zeta ))\in U^x$. We consider this case first. Then as when $x=L_3$, $\zeta ^1$ is likely to diverge from all paths of $R_{m,0}$ after the first unit-disc-crossing. In this case, we define $\beta ^1$ to be the unique path in $R_{m,0}$with $w_1'(\beta ^1)=x$ and such that $\beta ^1$ shares an end point of $\zeta ^1$ $\alpha ^1$ to be a perturbation of $\beta ^1*\zeta ^1$ which bounds a disc (or union of discs) disjoint from the common endpoint of $\beta ^1$ and $\zeta ^1$. We then define
$$\zeta ^2=\alpha ^1*\psi ^1(\zeta ^1).$$
 As in the case when $x=L_3$, the action of $\psi ^1$ on $\zeta ^1$ does not cancel the first  unit-disc crossing of $\alpha ^1$, and keeps the second unit-disc crossing in $U^x$. Therefore the first two unit-disc crossings of $\zeta ^2$ coincide with those of a path of $R_{m,0}$ with $w_1'=x$.

Now suppose that $w(\zeta ^1)\notin U^x$.  There are more possibilities than in the case $x=L_3$, when $w(\zeta ^1)\notin U^{L_3}$ occurred in essentially just one way. The new possibilities come from the difference between  the homeomorphisms $\psi _{i,2/7}$ and $\psi _{i,p}$, more precisely, in the difference between the annuli in $A_{i,2/7}$ and $C_{i,2/7}$, and the annuli in $A_{i,p}$ and $C_{i,p}$. For any $p>\frac{2}{7}$, the component of $C_{2,p}$ which intersects $D(C)$ also has intersection with interior with both $D(L_3)$ and $D(R_3)$.  Write 
$$x=x_1L_3^2(L_2R_3L_3)^n{\rm{\ for\ some\ }}n\ge 0,\ \  x_2=x_1L_3,\ \ x_2'=x_1L_2,$$
 or 
 $$x=x_1R_3L_2R_3(L_3L_2R_3)^nL_3{\rm{\ for\ some\ }}n\ge 0,$$
 $$ x_2=x_1(R_3L_2)^2R_3,\ \ x_2'=x_1R_3L_2,$$
 or
 $$x=x_1R_3,\ x_2=xL_2R_3,\ \ x_2'=x_1.$$
Note that 
$$C_{2,2/7}\cap D(C)=C_{2,p}\cap D(C),$$
and 
$$A_{4,2/7}\cap D(L_3)\subset A_{4,p}\cap D(L_3)\cup C_{2,p}\cap D(L_3))\cup L_3L_2(C_{2,p}\cap D(R_3)).$$
Hence the same is true if we prefix by $x_2$. There is an important difference between $C_{5,2/7}$ and $C_{5,p}$, which arises because $\beta _{2/7}$ crosses $C_{2,2/7}$, but $\beta _p$ does not intersect $C_{2,p}$. This means that, if $x_2'$ has length $i$, there is a component of $C_{i,p}$ which wraps round the unique point of $Z_i\cap D(x_2'C)$, and $\xi _{i,p}$ restricted to this component is simply a  half twist round this annulus. Furthermore, for each $n>0$, there is a component of $C_{i+3n,p}$ within this component of $C_{i,p}$, such that $\xi _{i+3n,p}$ restricted to this component is a fractional $1/2^n$ twist.  
 The sets $w(\zeta )$ for which $D(w(\zeta ^1))\not\subset  U^x$ are of three types:
 \begin{itemize}
 \item[1.]  $D(w(\zeta ^1))\subset  D(x_2'X)\subset x_2'C_{2,p}$, where $X=UC$ or $BC$, depending on whether $x_2'$ has an even or odd number of letters in $\{ L_3,L_2\} $. If this happens, we have $w(\zeta )\in D(x_2'X')$ where $X'=BC$ or $UC$ respectively, so that $X'\ne X$. 
 \item[2.] $D(w(\zeta ^1))\subset  D(y_1u_2)$, where $y_1=x_2L_3(L_2R_3)^2L_3(L_2R_3L_3)^r$ for some maximal  $r\ge 0$,  and $u_2$ with all letters in $\{ L_3,L_2,R_3\} $ ending in $L_2$ or $L_3$, of one of the two forms as in the case $x=L_3$, and $w_1'(y_1u_2X)\ne x$, where $X\in \{ BC,UC\} $ or $X=L_2$, whichever is appropriate. In this case, it is always also true that $w_1'(y_1u_2X,x)\ne x$.  Also, we have $D(w(\zeta ))\subset  D(x_2L_3L_2R_3L_3^{3r+3})$. But in that case, $x_2L_3L_2R_3L_3^3$ is a proper prefix of $x=w_1'(w(\zeta ))$. This contradicts the forms given above for $x_2$ in terms of $x$. So, in fact, this case does not occur.
 \item[3.] $D(w(\zeta ^1))\subset x_2L_3^2(C_{2,p}\cap D(L_3))$. In this case, we  have $D(w(\zeta ))\in x_2L_3L_2(C_{2,p}\cap D(R_3))$, and $w_1'(w(\zeta ^1))\ne w_1'(w(\zeta ))$. 
  \end{itemize}

 In cases 1 and 3, noting that case 2 does not occur, we define $v_1=w_1'(\zeta ^1)=w_1'(\psi _{m,p}(\zeta ))$ and let $\ell _1=D(v_1(L_2R_3L_3)^\infty )$. Let $e^{2\pi ip_1}$ be the unique endpoint of $\ell _1$ which is the first $S^1$ intersection point of at least one path in $R_{m,0}$,  and define
$$\psi ^1=\psi _{m,p_1},\ \ \alpha ^1=\alpha _{m,p_1},\ \ \ \zeta ^2=\alpha ^1*\psi ^1(\zeta ^1 ).$$
In case 1 above we have $w(\psi _{m,p_1}(\zeta ^2))\subset D(x_2'X)$, because $D(x_2'X)\subset C_{2,p_1}$,  and in case 3 we have $D(w(\zeta ^2))\subset x_2L_3^2(C_{2,p_1}\cap D(L_3))$, because $e^{2\pi ip_1}$ is to the left of $e^{2\pi ip}$, and consequently $C_{2,p_1}\supset C_{2,p}$. 
If $w_1'(\zeta ^2)=x^2$ then all of $\zeta ^2$, from the first intersection with $S^1$ onwards, is in $U^{x^2}$, but, as in the case $x^2=L_3$, the path $\zeta ^2$ diverges from all paths in $R_{m,0}$ after the first unit disc crossing, and we need to adjust by defining 
$$\zeta ^3=\alpha ^2*\psi ^2(\zeta ^2),$$
 where, similarly to before, $\alpha ^2$ is a perturbation of $\beta ^2*\overline{\zeta ^2}$, where $\beta ^2\in R_{m,0}$ is a path with $w_1'(\beta ^2)=x^2$ and the second unit disc crossing coincides with the second unit-disc-crossing of $\zeta ^2$. As before,  all of $\zeta ^3$ also, from the first intersection with the unit disc onwards, is contained in $U^{x^2}$. If $w_1'(\zeta ^2)\ne x^2$, then we define $p_3$ similarly to $p_2$, as the unique endpoint of $D(w_1'(\zeta ^2)(L_2R_3L_3)^\infty )$ which is the first $S^1$ intersection point of at least one path in $R_{m,0}$, and then define $\alpha ^3$, $\psi ^3$ and $\zeta ^3$ similarly to $\alpha ^2$, $\psi ^2$ and $\zeta ^2$. Then we iterate this process.

We finish with a discussion of the case of $t\in (\frac{19}{28},\frac{39}{56})$, where the proof cannot be completed simply because we have no comprehensive description of the paths in $R_{m,1}$. Recall that this time we are taking $\psi =\psi _{m,q}$ for $q=\frac{17}{56}$. As before, the support of $\xi _{m,q}$ is a union of annuli $A_{m,q}\cup C_{m,q}$.  The annulus components of $C_{1,q}$ are larger than the annulus components of $C_{1,2/7}$ while the annulus $A_{1,q}$ is smaller than the annulus $A_{1,2/7}$. The first intersection between $\psi _{i,q}(\zeta )$ and $A_{i,q}\cup C_{i,q}$ is for $i=4$, as $\zeta $ intersects the component of $C_{4,q}$ which surrounds the point in $s^{-2}(0)$ with code $L_3L_2C$. As the twist in components of $C_{4,q}$ is anticlockwise,  and is a half-twist, we see that $\alpha _q*\psi _{5,q}(\zeta )$ has first intersection point with $S^1$ at a point $e^{2\pi iu}$ for $u\in (\frac{17}{56},\frac{9}{28})$. As for $\zeta ^1=\alpha _{m,q}*\psi _{m,q}(\zeta )$, one of two things happens. Either 
\begin{equation}\label{6.8.5}\zeta ^1\cap \{ z:|z|\le 1\} \subset D(L_3L_2UC),\end{equation}
or, in a few cases, the endpoint of  $\psi _{5,q}(\zeta )$ is in $C_{7,q}$, and is subject to a further quarter twist, and possibly even further fractional twists. If this happens then the path $\zeta ^1=\alpha _{m,q}*\psi _{m,q}(\zeta )$ has first unit disc crossing either along the leaf joining $e^{\pm 2\pi i(9/28)}$, with the remainder of the path in $D(L_3L_2BCL_1R_2BC)$, or with all of the path beyond the first unit disc crossing in 
$$D(L_3^2(L_3L_2R_3)^3).$$
 In both cases, we have 
\begin{equation}\label{6.8.6}\zeta ^1\cap \{ z:|z|\le 1\} \subset U^1\end{equation}
We have not defined $U^1$ in this work, but it is mentioned in \ref{3.4}. In  both cases, all further iterations give paths with the same properties. If (\ref{6.8.5}) holds, the path $\beta \in R_{m,0}$ such that $\sigma _\beta \circ s$ and $\sigma _\zeta \circ s$ are Thurston equivalent satisfies $w_1'(\beta )\in D(L_3L_2UC)$. If (\ref{6.8.6}) holds, it seems likely that there is no path $\beta \in R_{m,0}$ such that $\sigma _\beta \circ s$ and $\sigma _\gamma \circ s$ are Thurston equivalent, but that there is such a path in $R_{m,1}$, by continuing the iteration. But it is not possible to be certain of this, because the paths in $R_{m,1}$ have never been completely described.\end{proof}

\section{Idea of the proof of \ref{6.4}}\label{6.9}

The proof of the main theorems \ref{6.1} and \ref{6.2} is now reduced to proving \ref{6.4}.
We start with a set ${\cal{U}}$ of capture paths $\gamma $ such that $\beta _i(\gamma )$ is constant for $\gamma \in {\cal{U}}$ and $i\le 4n_0+2$. We have
$$D_0=D(\beta _{4n_0+1}(\gamma ),\beta _{4n_0+2}(\gamma ))=D(\beta _{0,1}(\gamma ),\beta _{0,2}(\gamma ))$$
for any $\gamma \in {\cal{U}}$. 

 Our set $U_1$ will be of the form
 $$U_1=\cap _{i=0}^\infty U_{1,i}$$
 where $U_{1,0}=U=D_0$ and $U_{1,i}$, where $U_{1,i}$ is a union of sets $\psi _{1,4n_i-1}^{-1}(D_i)$, where $\beta _j(\gamma )$ is constant on $D_i$ for $j\le 4n_i$, and hence $\psi _{1,4n_i-1}^{-1}$ is also constant on $D_i$ and $\psi _{4n_{i'}+1,4n_i-1}^{-1}$ is constant on $D_i$ for all $0\le i'<i$, and hence $D_{i'}$ is constant on $D_i$ for all $0\le i'<i$, the numbers $n_{i'}$ are constant on $D_i$ for all $i'\le i$, and the sets $E_{4j-3,4j-1}$ are constant on $D_i$ for all $j\le 4n_i$. Then $U_{1,i}$ will be the union of sets $\psi _{1,4n_i-1}(D_i)$ such that the conditions of 2a), 2b), 2c) of \ref{6.4} which concern $\beta _j$ for $j\le 4n_i$, $D_{i'}$ for $i'\le i$, and $E_{4j-3,4j-1}$ for $j\le n_i$, all hold. We then need to show that
 $$\prod _{i=0}^N\dfrac{\# (U_{1,i+1}\cap Z_n)}{\# (U_{1,i}\cap Z_n)}\cdot \#(U\cap Z_n)$$
 is bounded from $0$ for all sufficiently large $n$. So we need to show that this product is convergent. 
 
 Obviously we are going to use induction. Given $U_{1,i(r)}$,  such that conditions 1 to 3 hold for sets $D_i$ for $i\le i(r)$ such that $U_{1,i}$ is a union of sets $\psi _{1,4n_i-1}(D_i)$, we need to choose $U_{1,i(r+1)}$ so that the conditions are satisfied.  So for each fixed $D_i(r)$ in $\psi _{1,4n_{i(r)-1}}^{-1}(U_i)$, we need to determine the possible $D_{i}$ for $i<i(r+1)$ with $D_{i+1}\subset D_i\subset D_{i(r)}$, and for this set of $D_{i(r+1)}$, and $i(r)\le i<i(r+1)$, it suffices to have
 $$\# (D_{i+1})\cap Z_n)\ge (1-
 \delta_i) \cdot \# (D_{i}\cap Z_n),$$
 where $\sum _i\delta _i<\infty $, because 
 $$\#(\psi _{1,4n_{i(r)}-1}(D_{i(r)}\cap Z_n))=\#(D_{i(r)}\cap Z_n),$$
 and similarly for $i(r+1)$.
 
To give the idea of the proof, it is easiest to describe what could be done if the situation were somewhat simpler. 
\begin{enumerate}[(i)]
\item Suppose it were true that $\psi _{\omega _1,\omega _3}={\rm identity}$ on all principal arcs $I(\zeta _1,\zeta _2)\subset D_{i(r)}$ for which $I(\zeta _1,\zeta _2)$ is not contained in $D(u)$ for some word $u$ with $|u|>\exp (|u|)^{\alpha _1}$  and for all $(\omega _t:1\le t\le 4)$ with $E_{\omega _1,\omega _3}\subset D_{i(r)}$, again with $E_{\omega _1,\omega _3}$ not contained in any such $D(u)$. 
\item Suppose also that it were true that if $C$ were a component of $(\sigma _{\beta _{4\ell -3}}\circ s)^{-n}(E_{4\ell -3,4\ell -1})$ for $\ell <n_i$ then either $C\cap D_{i(r)}=\emptyset $ or $D_{i(r)}\subset C$. \end{enumerate}

Then it would be easy to choose $D_{i+1}$ and $D_j$ for $i(r-1)\le i<j\le i(r+1)$. But none of these is true. There is always at least one component $C$ of $(\sigma _{\beta _{1}}\circ s)^{-n}(E_{1,3})$ for some $n$, which is not contained in $D_{i(r)}$, and corresponding components of $(\sigma _{\beta _{4\ell -3}}\circ s)^{-n}(E_{4\ell -3,4\ell -1})$,  for each $\ell \le n_{i}$. Indeed $C\cap Z_n$ is a positive proportion of $D_i\cap Z_n$, bounded from $0$. So we need to choose $u_{n_i}$ and  $D_i$ carefully to ensure that $I(\beta _{2j -1},\beta _{2j })$ is disjoint from this $C$, for $j\le 2n_i$. The size of $|u_{2n_{i+1}+1}|-|u_{2n_i+1}|$ is chosen so that such a choice will be possible, with a set of positive density.

First we consider (i).  Every type 2 arc in $D_i$ is crossed by some component $C$ of $(\sigma _{\omega _1}\circ s)^{-n}(E_{\omega _1,\omega _3})$ for any quadruple $(\omega _i)$. But the crossing is usually on the central arc $\ell (C)$ of $C$and if $E_{\omega _1,\omega _3}$ is already small then whenever $(\omega _t:\le t\le 4)=(\beta _{-4+4l+t}:1\le t\le 4)$ and some $\ell $, it is possible to bound the size of the component of 
$$\cup _{0<m,n_i<j}(\sigma _{\beta _{4j-3}}\circ s)^{-m}(E_{4j-3,4j-1})$$
 containing $C$. In fact, we shall do that as follows. We define $k(C)=k(\omega _1,\omega _2)$, where $\omega _1$ and $\omega _2$ are the two paths in $R_{m,0}$ which bound $C$, as in \ref{4.13}. Then $k(C)-1$ is constant on each of the two components $E_1$ and $E_2$ of ${\rm int }(C)$, and if $\{ D(v):w_{k(C)+2}'(v)=u\} $ intersects $C$, it is contained in $C$. We shall ensure that if $k(C)\ge 3$ then any component of $(\sigma _{\beta _{4\ell -3}}\circ s)^{-n}(E_{4\ell -3,4\ell -1})$ which intersects $C$ is contained in
 $$C_1\cup \bigcup \{ D(v):k(v)\ge k(C)+2,\ w_{k(C)+2}'(v)=u{\rm\ for\ some\ }D(u)\subset C_1\} $$
 where $C_1$ is the corresponding component of $(\sigma _{\beta _{1}}\circ s)^{-n}(E_{1,3})$, or 
 $$C_1\cup C_2\cup \bigcup \{ D(v):k(v)\ge k(C)+2,\ w_{k(C)+2}'(v)=u{\rm\ for\ some\ }D(u)\subset C_1\cup C_2\} $$
 where $C_2$ is a component of $(\sigma _{\beta _{1}}\circ s)^{-p}(E_{1,3})$ which intersects $E_{1,3}$. This helps us to obtain the required bound on satellites. We aim to show that, for $n_{i(r-1)}\le \ell $ any component of $(\sigma _{\beta _{4\ell -3}}\circ s)^{-n}(E_{4\ell -3,4\ell -1})$ is disjoint from $I(\beta _{2j-1},\beta _{2j})$ for $n(i(r-1))\le 2j\le n(i(r+1))$, except when $(\beta _{2j-1},\beta _{2j})$ is purely satellite.
 
Now we consider (ii). There is always at least one such $C$ such that $C\cap D_i\ne \emptyset $ and $C\setminus D_i\ne \emptyset $. In fact, if $C$ is a component of $(\sigma _{\beta _1}\circ s)^{-n}(E_{1,3})$, then the proportion of points  in $C\cap D_i\cap Z_k$  is a  positive proportion of points in $D_i\cap Z_k$ for all sufficiently large $k$, independent of $i$ and $D_i$. But the proportion of sets $D_{i}$ for which it will be impossible to choose $n_{i+1}$ and  $D_{i+1}$ such that $C\cap D(\beta _{2j-1},\beta _{2j})\cap C=\emptyset $ for $j\ge 2n_{i+1}$ is of proportion $<\lambda ^{(\log n_i)^\alpha }$ for some fixed $\lambda <1$, and we can choose $\alpha >1$. In this case, this is $<n_i^{d\log \lambda }$ for any $d>0$ for sufficiently large $i$, and since $n_i>i$, we have
$$\sum _i\lambda ^{(\log n_i)^\alpha }<+\infty .$$
In all cases, we will make our initial choice so that every component $C'$ of $(\sigma _{\beta _{4\ell -3}}\circ s)^{-n'}(E_{4\ell -3,4\ell -1})$ which intersects this $C$ is either contained in it or contains $D_i$. So then if $C'\subset C$ we can ensure that $D_{i+1}\cap C=\emptyset $, and hence $D_{i+1}\cap C'=\emptyset $. 

The method which is used is loosely, but quite closely, related to standard techniques of proving non-uniform hyperbolicity properties --- as in \cite{R8}, for example, but the idea probably occurs in every single paper on the subject. We need to exclude certain behaviours of the critical orbit. Early and long recurrence needs to be excluded, which, in general terms, is what is needed to prove non-uniform hyperbolicity -- in the cases which can be handled by methods currently available. Something more needs to be excluded because of the group action, but everything that has to be excluded is related to early long recurrence --- what were known as ``followers'' in \cite{R8}.

\section{Proof of \ref{6.4} }\label{6.10}

One basic idea of the proof is to translate everything into symbolic dynamics. This means that we use the admissible words $w$  in our standard alphabet 
$$\{ BC,UC,L_i,R_j:1\le i,j\le 3\} $$
and the sets $D(w)$. The number of admissible words of length $n$ is boundedly proportional to $2^n$, and if $w_0$ and $z_0$ are admissible and of lengths $m$ and $p$ respectively, then the number of admissible words of length $n$ and starting with $w_0$ and ending with $z_0$ is boundedly proportional to $2^{n-m-p}$. 

Fix $\gamma \in {\cal{U}}$ and write $\beta _j=\beta _j(\gamma )$. As in \ref{6.4}, we define $u_j$ to be the longest word such that
$$D(\beta _{2j-1},\beta _{2j})\subset D(u_j).$$
One thing which we need to determine is how  $u_j$ differs from a prefix of $w(\gamma)$. 
Putting $j=2n_i+1$, we have
$$D_i\subset D(u_{2n_i+1}).$$
The conditions of \ref{6.4} ensure that  $u_{2n_i+1}$ is a prefix of $u_j$ for all $j\ge 2n_i+1$. In particular, $u_{2n_i+1}$ is a prefix of $u_{n(\gamma )}$ for all $i\ge 0$ and $u_{2n_j+1}$ is a prefix of $u_{2n_i+1}$ for all $j\le i$. 
It is not the case, however, that $u_{j_1}$ is a prefix of $u_{i_1}$ for all $j_1\le i_1$, and not even if $i_1=2n_i+1$ for some $i$. Still, if $j$ and $i$ are the largest integers with $2n_j+1\le j_1$ and $2n_i+1\le i_1$, then $u_{2n_j+1}$ is a prefix of $u_{2n_i+1}$, and hence is also a prefix of $u_{i_1}$. Since $|u_{j_1}|-|u_{2n_j}|=o(|u_{j_1}|)$, this is useful information. 

We also define $u_{j,1}$ to be the longest word such that
$$\psi _{1,4\lfloor j/2\rfloor -1}(D(\beta _{2j-1},\beta _{2j}))\subset D(u_{j,1}),$$
and $u_{2j-1,-}$ is the longest word with
$$E_{4j-3,4j-1}\subset D(u_{2j-1,-}).$$
For the choice of $(\beta _i:1\le i\le 4)$ for proving Theorem \ref{6.1}, there is no single nonempty word $u_{1,-}$ with $E_{1,3}\subset D(u_{1,-})$. So in this case, we define $u_{1,-}$ to be  a pair of words $(u_{1,1,-},u_{1,2,-})$ with $D(u_{1,1,-})\subset D(L_3)$ and $D(u_{1,2,-})\subset D(BC)$, where $u_{1,1,-}$ and $u_{1,2-}$ are the longest possible with 
$$E_{1,3}\subset D(u_{1,1,-})\cup D(u_{1,2,-}).$$
We also define 
$$D(u_{1,-})=D(u_{1,1,-})\cup D(u_{1,2,-}).$$

Our set $U_1$, which we want to construct to prove Theorem \ref{6.4}, will be a union of sets $u_{2n(\gamma ),1}$. By the inductive construction there is a one-to-one correspondence between the words $u_{j,1}$ and $u_j$. We will be choosing the set $U_1$ to be a union of sets $D(u_{2n(\gamma ),1})\cap Z_\infty $ by choosing the union of sets $D(u_{2n(\gamma )})\cap Z_\infty )$ and we want the union of sets $D(u_{2n(\gamma )}\cap Z_m)$ to be of positive density for all sufficiently large $m$. So,  of course, we need to control the difference between $u_j$ and $u_{j,1}$. 

We want to ensure that any changes between $u_j$ and $u_{j,1}$ can be controlled by looking either at the different prefixes $u_{2n_i+1}$ in $u_j$, or beyond the largest such prefix. In order to do this, we shall ensure the following.
 
\begin{enumerate}[(i)]
\item  For $j>2n_i+1$, for any suffix $a_1$ of $a$ such that $a\leftrightarrow b$ is a basic exchange with corresponding letters of $a$ and $b$ all different, any  occurrence of $a_1u_{1,-}$ which is a subword of $u_j$, and  which starts before the end of  the prefix  $u_{2n_i+1}$ in $u_j$, is a subword of $u_{2n_i+1}$.
\item No non-prefix subword of $u_{2n_i+1}$ is equal to $u_{2n_{i(r-1)}+1}$, for $i\le i(r+1)$.\end{enumerate}

We claim that these two conditions are enough to ensure the conditions a) and b) of \ref{6.4}. The reason is as follows.  The only way that a) or b) can fail is if $D(v)$ is moved outside $D(u_{2n_i+1})$ by $\psi _{4n_i+1,4\ell -1}$ for some $\ell >n_i$, and for some word $v$ which is an extension of $u_{2n_i+1}$, but a prefix of $u_{2n_{i+1}+1}$. Now $\psi _{4\ell +1,4\lfloor (j-1)/2\rfloor -1}$ is a composition of disc exchanges for each $2\ell <j-1$. If we look just at the disc exchanges in the composition which have an effect on $D(u_j)$, we see that the support of each one is either contained in $D(u_{2n_i+1}$ or contains it, for each $i$ with $2n_i+1\le j$. So the disc exchanges which have an effect on $D(u_j)$ can be permuted, and we can group together those whose supports contain $D(u_{2n_i+1})$ together, for each $i$. Then we can permute the order in which these occur. We can apply this corresponding to the largest $i$ first, and then apply them in decreasing order of $i$. This makes sense, because then all the disc exchanges, in the composition of $\psi _{4\ell +1,4\lfloor (j-1)/2\rfloor -1}$ which have an effect on $D(u_j)$, have support contained in $D(u_{2n_i+1})$ for any $n_i\le 2\ell $. So $\psi _{4\ell +1,4\lfloor (j-1)/2\rfloor -1}(D(u_j))$ is of the form $D(u_{j,\ell +1}$ for some word $u_{j,\ell +1}$ where $u_{j,\ell +1}$ has $u_{2n_i+1}$ as a prefix for any $2n_i\le \ell +1$. For those disc exchanges which have support between $D(u_{2n_i+1})$ and $D(u_{2n_{i+1}+1})$, it is not possible to commute in any order, but the first exchange which is applied can only give a change before an occurrence of $u_{1,-}$, or before a word arrived at by a rather special sequence of disc exchanges preceding some $u_{2\ell -1,-}$. 

Clearly, in order to satisfy these conditions, we need to choose the $n_i$ and $u_{2n_i+1}$ suitably, given $u_j$. In fact, we will choose $n_i$ and  $u_{2n_i+1}$, for $i\le i(r)$, given an extension $u$ of $u_{2n_{i(r)}+1}$. From the statement of \ref{6.4}, we want the conditions
$$|u_{2n_{j+1}+1}|-|u_{2n_j+1}|\le (\log |u_{2n_j+1}|)^{\alpha _1}=o(|u_{2n_j+1}|)$$
and 
$$|u_{2n_{2i(r+1)}+1}|\ge \exp (\log |u_{2n_{2i(r)}+1}|)^{\alpha _2})$$
From this, and the corresponding lower bound on $|u_{2n_{i(r)}+1}|$ in terms of $|u_{2n_{i(r-1)}+1}|$, we have
$$|u_{2n_{i(r+1)}+1}|\ge \exp (\log |u_{2n_{i(r-1)}+1}|^{\alpha _2^2})$$

Recall that we are aiming to construct $U_1$ as an intersection of sets $U_{1,r}$. We shall construct $U_{1,r}$ as a union of sets $D(u)$ for prefixes $u_{2n_i+1,1}$ have been defined for $i\le i(r)$ with 
$$|u|-|u_{2n_{i(r)}+1}|\le  (\log |u_{2n_{i(r)}+1}|)^{\alpha _1/4})$$
We shall do this by constructing $U_{1,r}'$ with
$$U_{1,r}=\psi _{4n_{i(r)}+1}(U_{1,r}').$$
The method is inductive. For a sequence of integers $m_r$, we will choose $U_{1,r}'$ to be a union of sets $D(u)$ with $|u|=m_r$. We start by taking $m_0=|u_{2n_0+1}|$ but for $r>0$, we will have $m_r>|u_{2n_{i(r)}+1}|$ for each choice of $u_{2n_{i(r)}+1}$. Since we can take $|u_{2n_0+1}|$ arbitrarily large, we can also take $m_0$ arbitrarily large. We shall construct $U_{1,r+1}'$ as a union of sets $D(v)$ where each $D(v)$ is a subset of one of the sets $D(u)$ in the union $U_{1,r}'$, that is, $v$ has some $u$ as a prefix, where $D(u)$ is part of the union $U_{1,r}'$. We shall ensure that
$$\# (U_{1,r+1}'\cap Z_{m_{r+1}})\ge (1-  \delta _r')\#(U_{1,r}'\cap Z_{m_{r+1}}).$$
 where $\sum \delta _r'$ converges, and, in fact, can be taken arbitrarily small.

 So it remains to show that, given $u$ such that $D(u)$ is one of the sets in the union   $U_{1,r}'$, we can find a high density set of extensions $v$ for which the properties are satisfied. Note that $i(r)$, $n_{i(r)}$ and $u_{2n_{i(r)}+1}$ all depend on $u$. The numbers $i(r+1)$, $n_i$ for $i(r)<i\le i(r+1)$, and the prefixes $u_j$ for $2n_i(r)+1<j\le 2n_{i(r+1)}+1$, will all depend on the extension $v$ of $u$. But we need to find a good proportion of extensions $v$ of $u$ for which these numbers and prefixes can be defined with the appropriate properties. We will choose $m_{r+1}$ so that
 $$\exp ((\log m_{r})^{\alpha _2/2}\le m_{r+1}\le \exp \log (m_{r}^{2\alpha _2}).$$
 Note that this also implies that
 $$\exp ((\log m_{r-1})^{\alpha _2^2/4}\le m_{r+1}\le \exp \log (m_{r-1}^{\alpha _2^2}),$$
 or otherwise put,
 $$(\log m_{r+1})^{\alpha _2^{-2}}\le m_{r-1}\le (\log m_{r+1})^{4\alpha _2^{-2}}.$$

 So fix $u$. Then  $i(r)$ is defined depending on $u$ and $n_i$ and $u_{2n_i+1}$ are defined for each $i\le i(r)$, depending on $u$, such that $u_{2n_i+1}$ is a prefix of $u$. To satisfy conditions (i) and (ii), we need to obtain a lower bound on the number of $v$ extending $u$ such that
 \begin{itemize}
 \item every subword of length $(1/10)(\log m_r)^{\alpha _1}$ in the extension of $u$ in $v$ contains a word $L_3x_0$ such that $D(L_3x_0)\cap D(u_{1,-})=\emptyset $ and $u_{1,-}$ is not a suffix of $L_3x_0$;
 \item $u_{2n_{i(r-1)}+1,-}$ is not a subword of the extension of $u$ in $v$.
 \end{itemize}
 If $r=0$ then in the second condition we simply replace $u_{2n_{i(r-1)}+1,-}$ by $u_{2n_{i(r)}+1,-}$.
 
These conditions are sufficient to ensure (i) and (ii) respectively. In order to ensure (i) we can take a particular choice of $L_3x_0$. A choice of $L_3x_0$ which will work for all our different choices of $u_{1,-}$ is $(L_3(L_2R_3)^3)^2$. So now we need to estimate the proportion of extensions $v$ of $u$ of length $m_{r+1}$ for which these conditions are not satisfied. The proportion of words of length $N$ which do not contain an occurrence of $L_3x_0$ is $\le e^{-N\lambda }$ for some $\lambda >0$ which depends only on $L_3x_0$ -- which we have fixed. We want to avoid such words, with $N\ge (1/10)(\log m_r)^{\alpha _1}$, over length $m_{r+1}$.  So the proportion of words that we want to avoid for (i) is at most
$$m_{r+1}\exp (-(\lambda /10)(\log m_r)^{\alpha _1}\le \exp ((\log m_r)^{2\alpha _2}-(\lambda /10)(\log m_r)^{\alpha _1}).$$
This proportion is at most
$$ \exp (-(\log m_r)^{\alpha _1/2}$$
provided that $\alpha _2<\alpha _1/3$, that is, provided that $\alpha _1$ is sufficiently large given $\alpha _2$. The only requirement on $\alpha _2$ is that $\alpha _2>1$.
Similarly for  (ii), for which it suffices to avoid $u_{2n_{i(r-1)}+1,-}$, the proportion of words which we want to avoid is, for a suitable $\lambda >0$,
$$\le m_{r+1}\exp (-\lambda m_{r-1})\le \exp ((\log m_{r-1})^{4\alpha _2^2}-\lambda m_{r-1}).$$
Assuming, as we may do, that $m_{r-1}$ is sufficiently large given $\lambda $ and $\alpha _2$, this is at most
$$ \exp (-\lambda m_{r-1}/2).$$
So altogether the proportion that has to be avoided is
$$\sum _{r\ge 0}\exp (-(\log m_r)^{\alpha _1/2})+\sum _{r\ge 0}\exp (-\lambda m_{r-1}/2)$$
where in the second sum we define $m_{-1}=m_0$. This sum can be taken arbitrarily small by taking $m_0$ arbitrarily large for fixed $\lambda >0$ and $\alpha _1>2$. This completes the construction of $U_{1,r+1}$ from $U_{1,r}$, and hence completes the proof.

\bibliographystyle{amsplain}

\bibliography{capturemapdriver}
}
\end{document}